\documentclass[english]{amsart}
\usepackage{amsmath}

\usepackage{graphicx}

\usepackage{amssymb}
\usepackage{svn}

\usepackage{stackrel}

\usepackage{amscd}
\usepackage[all,cmtip,line]{xy}

\SVN $LastChangedDate: 2012-10-30 15:27:28 +0000 (Tue, 30 Oct 2012) $
\SVN $LastChangedRevision: 28 $

\def\AA{\mathbb{A}}

\def\CC{\mathbb{C}}

\def\FF{\mathbb{F}}

\def\PP{\mathbb{P}}
\def\QQ{\mathbb{Q}}
\def\RR{\mathbb{R}}

\def\ZZ{\mathbb{Z}}

\newcommand{\calC}{{\mathcal{C}}}

\newcommand{\onto}{\twoheadrightarrow}

\newcommand{\lra}{\longrightarrow}

\newcommand{\dR}{\mathrm{dR}}

\newcommand{\mat}[4]{\left( \begin{array}{cc} {#1} & {#2} \\ {#3} & {#4}
\end{array} \right)}

\newcommand{\smat}[4]{{\mbox{\scriptsize $\mat{{#1}}{{#2}}{{#3}}{{#4}}$}}}
\newlength{\ownl}

\newcommand{\Aut}{{\operatorname{Aut}\,}}

\newcommand{\coker}{{\operatorname{coker}\,}}

\newcommand{\End}{{\operatorname{End}\,}}

\newcommand{\Fil}{{\operatorname{Fil}\,}}

\newcommand{\Frob}{{\operatorname{Frob}}}

\newcommand{\gr}{{\operatorname{gr}\,}}

\newcommand{\Hom}{{\operatorname{Hom}\,}}

\newcommand{\Res}{{\operatorname{Res}}}

\newcommand{\Spec}{{\operatorname{Spec}\,}}
\newcommand{\SPEC}{\mathbf{Spec}}

\newcommand{\Tan}{{\operatorname{Tan}}}
\newcommand{\tr}{{\operatorname{tr}\,}}
\newcommand{\Tr}{{\operatorname{Tr}\,}}

\newcommand{\uno}{{\operatorname{uno}\,}}

\newcommand{\GL}{\operatorname{GL}}

\newcommand{\SL}{\operatorname{SL}}

\newcommand{\cris}{{\operatorname{crys}}}

\newcommand{\red}{{\operatorname{red}}}

\newcommand{\st}{{\operatorname{st}}}

\newcommand{\A}{{\mathbb{A}}}

\newcommand{\C}{{\mathbb{C}}}
\newcommand{\D}{{\mathbb{D}}}

\newcommand{\F}{{\mathbb{F}}}

\newcommand{\Q}{{\mathbb{Q}}}

\newcommand{\Z}{{\mathbb{Z}}}

\newcommand{\CA}{{\mathcal{A}}}
\newcommand{\CB}{{\mathcal{B}}}
\newcommand{\CE}{{\mathcal{E}}}
\newcommand{\CF}{{\mathcal{F}}}
\newcommand{\CG}{{\mathcal{G}}}
\newcommand{\CH}{{\mathcal{H}}}
\newcommand{\CI}{{\mathcal{I}}}
\newcommand{\CJ}{{\mathcal{J}}}
\newcommand{\CK}{{\mathcal{K}}}
\newcommand{\CL}{{\mathcal{L}}}
\newcommand{\CM}{{\mathcal{M}}}
\newcommand{\CN}{{\mathcal{N}}}
\newcommand{\CO}{{\mathcal{O}}}
\newcommand{\CP}{{\mathcal{P}}}

\newcommand{\CR}{{\mathcal{R}}}

\newcommand{\CT}{{\mathcal{T}}}

\newcommand{\CV}{{\mathcal{V}}}

\newcommand{\cD}{\mathcal{D}}

\newcommand{\gP}{{\mathfrak{P}}}

\newcommand{\gc}{{\mathfrak{c}}}
\newcommand{\gd}{{\mathfrak{d}}}

\newcommand{\gm}{{\mathfrak{m}}}
\newcommand{\gn}{{\mathfrak{n}}}

\newcommand{\gp}{{\mathfrak{p}}}

\newcommand{\Qbar}{{\overline{\Q}}}

 \newcommand{\p}{\mathfrak{p}}

\newcommand{\Qpbar}{\overline{\Q}_p}

\newcommand{\Fpbar}{\overline{\F}_p}

\newcommand{\Ver}{\operatorname{Ver}}
\newcommand{\Sym}{\operatorname{Sym}}

\usepackage{hyperref}

\newcommand{\f}{\mathbf{f}}

\newcommand{\bk}{\mathbf{k}}

\newcommand{\bm}{\mathbf{m}}

\newcommand{\bone}{\mathbf{1}}
\newcommand{\bO}{\mathbf{0}}

\newcommand{\Nm}{\mathrm{Nm}}
\newcommand{\uhp}{\mathfrak{H}}
\newcommand{\Shom}{{\mathcal{H}om}}
\newcommand{\Sext}{{\mathcal{E}xt}}

\newcommand{\Lie}{{\mathcal{L}ie}}
\newcommand{\crys}{\mathrm{crys}}

\newcommand{\dr}{{\operatorname{dR}}}

\newcommand{\im}{{\operatorname{im}}}

\newcommand{\ol}{\overline}

\newcommand{\pp}{\prime\prime}

\theoremstyle{plain}
\newtheorem{theorem}{Theorem}[subsection]
\newtheorem{proposition}[theorem]{Proposition}
\newtheorem{corollary}[theorem]{Corollary}
\newtheorem{lemma}[theorem]{Lemma}

\newtheorem{ithm}{Theorem}

\newtheorem{icor}[ithm]{Corollary}

\theoremstyle{definition}
\newtheorem{remark}[theorem]{Remark}

\newcommand{\begpf}{\noindent{\bf Proof.}\enspace}
\newcommand{\epf}{{\ifhmode\unskip\nobreak\hfil\penalty50 \hskip1em
\else\nobreak\fi \nobreak\mbox{}\hfil\mbox{$\square$} \parfillskip=0pt
\finalhyphendemerits=0 \par\vskip5pt}}

\begin{document}

\title[Kodaira--Spencer isomorphisms and degeneracy maps]{Kodaira--Spencer isomorphisms
and degeneracy maps on Iwahori-level Hilbert modular varieties:\\ the saving trace}
\author{Fred Diamond}
\email{fred.diamond@kcl.ac.uk}
\address{Department of Mathematics,
King's College London, WC2R 2LS, UK}

\subjclass[2010]{11G18 (primary), 11F32, 11F41, 14G35 (secondary)}

\begin{abstract}
We consider integral models of Hilbert modular varieties with Iwahori level structure at primes over $p$,
first proving a Kodaira--Spencer isomorphism that gives a concise description of their dualizing sheaves.
We then analyze fibres of the degeneracy maps to Hilbert modular varieties of level prime to $p$ and
deduce the vanishing of higher direct images of structure and dualizing sheaves, generalizing prior work
with Kassaei and Sasaki (for $p$ unramified in the totally real field $F$).  We apply the vanishing results to
prove flatness of the finite morphisms in the resulting Stein factorizations, and combine them with the Kodaira--Spencer isomorphism to simplify and
generalize the construction of Hecke operators at primes over $p$ on Hilbert modular forms (integrally and mod $p$).
\end{abstract}



\maketitle


\section{Introduction}

The motivation for this paper is two-fold: the construction of Galois representations and 
of Hecke operators at primes over $p$ in the setting of mod $p$ Hilbert modular forms
for a totally real field $F$.  By mod $p$ Hilbert modular forms, we mean sections
of automorphic line bundles of arbitrary weight over mod $p$ Hilbert modular varieties
of level prime to $p$.  We will explain this below more precisely, discuss what was done previously,
and describe how the results here 
complete the picture for Hecke operators and 
feed into an argument that does so for Galois representations.

First let us recall the nature of the difficulty in constructing Hecke operators at $p$ in
characteristic $p$.  One issue already arises in the setting of classical modular forms.
For simplicity assume $N > 4$ and let $X_1(N)$ denote the modular curve $\Gamma_1(N)\backslash\uhp^*$,
so that $X_1(N)$ may be viewed as parametrizing pairs $(E,P)$, where $E$ is a generalized elliptic
curve over $\CC$ and $P$ is a point of $E$ of order $N$ (see for example \cite{DR}).  
Similarly for $p$ not dividing $N$ and $\Gamma_1(N;p) = \Gamma_1(N) \cap \Gamma_0(p)$, the modular curve
$X_1(N;p) = \Gamma_1(N;p)\backslash\uhp^*$ parametrizes data of the form $(E,P,C)$, where $(E,P)$
is as above and $C$ is a subgroup of $E$ of order $p$, and there are two natural degeneracy
maps $\pi_1,\pi_2: X_1(N;p) \to X_1(N)$, with $\pi_1(E,P,C) = (E,P)$ and $\pi_2(E,P,C) = (E/C, P\bmod C)$.
The space of modular forms of weight $k$ with respect to $\Gamma_1(N)$ may be identified with
$H^0(X_1(N),\omega^{\otimes k})$, where the line bundle 
$\omega = 0^*\Omega^1_{E/X_1(N)}$ is the pull-back
of $\Omega^1_{E/X_1(N)}$ along the zero section
$0: X_1(N) \to E$.   The Hecke
operator $T_p$ can then be defined as the composite
$$\begin{array}{rl}  H^0(X_1(N),\omega^{\otimes k}) & \longrightarrow 
    H^0(X_1(N;p), \pi_2^*\omega^{\otimes k})  \\  &\longrightarrow
    H^0(X_1(N;p), \pi_1^*\omega^{\otimes k}) \longrightarrow
    H^0(X_1(N), \omega^{\otimes k})\end{array}$$
{\em divided by $p$}, where the first map is pull-back, the second is induced by the universal isogeny
over $X_1(N;p)$, and the third is the trace.
To define $T_p$ integrally or in characteristic $p$, one 
can work instead with integral models of the
curves and show (as in \cite[\S4.5]{con2}, for example) that the resulting composite morphism 
\begin{equation} \label{eqn:classical} \pi_{1,*} \left(\pi_2^* \omega^{\otimes k}\right) \longrightarrow \pi_{1,*} \left(\pi_1^* \omega^{\otimes k}\right)
    \longrightarrow \omega^{\otimes k} \end{equation}  
of sheaves is divisible by $p$ (assuming $k \ge 1$).

This paper offers an alternative to the standard approach just described to the construction of Hecke operators at $p$. 
In particular, it allows for a more general and direct definition of the morphism which produces (\ref{eqn:classical}) after multiplication by~$p$.  Our perspective thus breaks the mindset reflected in \cite[\S3.4]{cal}, which described all prior constructions by saying that 
{\em the ``correct'' definition of $T$ involves first defining a map coming from a correspondence and then showing that it is ``divisible'' by the correct power of $p$.}

To construct integral and mod $p$ Hecke operators at primes over $p$ over more general Shimura
varieties, one encounters the further difficulty that the degeneracy maps no longer necessarily extend to
finite flat morphisms on the usual integral models, so the definition of a trace morphism requires a more
sophisticated application of Grothendieck--Serre duality.  Emerton, Reduzzi and Xiao take this approach
in \cite{ERX} to defining Hecke operators at primes over $p$ for
Hilbert modular forms,  integrally and mod $p$, but for a restricted set of weights.   On the other hand,
Fakhruddin and Pilloni set up a general framework in \cite{FP} and prove results for Hilbert modular forms
that are optimal in terms of the weights\footnote{The results in \cite{FP} are only formulated integrally,
necessitating a parity condition on the weight, but the hypothesis is not essential to the methods there.}
considered, but they assume $p$ is unramified in $F$.

This paper contains both a conceptual innovation and a technical improvement on previous
work, yielding an optimal result.   The innovation takes the form of a Kodaira--Spencer
isomorphism describing the dualizing (or canonical) sheaf on integral models of Iwahori-level
Hilbert modular varieties.  In addition to being strikingly natural and simple to state, it can also
be viewed as encoding integrality properties of Hecke operators at $p$. The technical 
advance is a cohomological vanishing theorem for the dualizing sheaf, showing that its higher
direct images relative to certain degeneracy maps are trivial.  This generalizes such a result
in \cite{DKS}, where it is proved under the assumption that $p$ is unramified in $F$, and
we apply it here to define Hecke operators integrally and in characteristic $p$, obtaining them
directly from morphisms of coherent sheaves rather than (rescaled) morphisms of complexes 
in a derived category, as in \cite{ERX} and~\cite{FP}.

The cohomological vanishing is also useful in proving the existence of Galois representations
associated to mod $p$ Hilbert modular eigenforms of arbitrary weight.  The construction of such Galois
representations is proved independently by Emerton, Reduzzi and Xiao in \cite{ERX0}, and by Goldring
and Koskivirta in \cite{GK}, under parity conditions on the weight.  These conditions are inherited from
the obvious parity obstruction to algebraicity for automorphic forms on $\Res_{F/\Q}\GL_2$; however they are unnecessary,
and unnatural, in the consideration of automorphic forms in finite characteristic.   The construction of
Galois representations for mod $p$ Hilbert modular forms of arbitrary weight is carried out in \cite{DS1}
under the assumption that $p$ is unramified in $F$, with the cohomological vanishing result in \cite{DKS}
playing a crucial role in the argument.
The generalization proved in this paper can similarly be used,
in conjunction with the methods of \cite{DS1} and \cite{ERX0},
to prove the existence of Galois representations associated to mod $p$ Hilbert modular eigenforms in full generality.

The cohomological vanishing results also have applications to the study of integral models for
Iwahori-level Hilbert modular varieties.  Indeed similar arguments to those for the dualizing sheaves
also show that the higher direct images of their structure sheaves vanish, and together with
Grothendieck--Serre duality, this implies the flatness of the finite morphisms in the Stein factorizations of the degeneracy
maps.  As a result we obtain Cohen--Macaulay models for the Iwahori-level varieties which are
finite flat over the smooth models for the varieties of prime-to-$p$ level.

\medskip

We now describe our results in more detail.  

We fix a prime $p$ and a totally real field $F$, and let $\CO_F$ denote the ring of integers of $F$.
Fix also embeddings $\Qbar \to \Qpbar$ and $\Qbar \to \CC$, and let $\CO$ denote the ring of integers
of a sufficiently large finite extension $K$ of $\QQ_p$ in $\Qpbar$.  

Let $U$ be a sufficiently small open compact subgroup of $\GL_2(\AA_{F,\f})$ containing $\GL_2(\CO_{F,p})$.
A construction of Pappas and Rapoport~\cite{PR} yields a smooth model $Y$ over $\CO$ for the Hilbert
modular variety with complex points
$$\GL_2(F)_+\backslash (\uhp^\Sigma \times \GL_2(\A_{F,\f})/U),$$
where $\uhp$ is the complex upper-half plane and $\Sigma$ is the set of embeddings $F \to \Qbar$.
The scheme $Y$ is equipped with line bundles\footnote{We remark that $\delta$ is trivializable, but only
non-canonically, and we systematically incorporate it into constructions to render them Hecke-equivariant.}
$\omega$ and $\delta$ arising from its interpretation as a coarse moduli space parametrizing abelian
schemes with additional data including an $\CO_F$-action.  Since $Y$ is smooth over $\CO$, its
relative dualizing sheaf $\CK_{Y/\CO}$ is identified with $\wedge^{[F:\QQ]}\Omega^1_{Y/\CO}$, and in~\cite{RX}
(see also \cite[\S3.3]{theta}),  Reduzzi and Xiao establish an integral version of the Kodaira--Spencer isomorphism,
taking the form
$$\CK_{Y/\CO}   \cong \delta^{-1}\omega^{\otimes2} .$$

Let $\gp$ be a prime of $\CO_F$ containing $p$.  Denote its residue degree by $f_\gp$ and let
 $$U_0(\gp) = \left\{\,\left. \left(\begin{array}{cc} a&b \\ c& d \end{array}\right) \in U \,\right|\, c_\gp \in \gp \CO_{F,\gp} \,\right\}.$$
Augmenting the moduli problem for $Y$ with suitable isogeny data then produces 
a model $Y_0(\gp)$ for the Hilbert modular variety of level $U_0(\gp)$,
equipped with a pair of degeneracy maps $\pi_1,\pi_2:Y_0(\gp) \to Y$.  The morphisms $\pi_1$ and
$\pi_2$ are projective, but not finite or flat unless $F_\gp = \QQ_p$.

The scheme $Y_0(\gp)$ is syntomic over $\CO$, hence has an invertible dualizing sheaf $\CK_{Y_0(\gp)/\CO}$.
Our Iwahori-level version of the Kodaira--Spencer isomorphism is the following (Theorem~\ref{thm:KSIwa}; see there
for the precise meaning of Hecke-equivariant).
\begin{ithm} \label{ithm:KS} There is a Hecke-equivariant isomorphism
$$\CK_{Y_0(\gp)/\CO} \cong  \pi_2^*(\delta^{-1}\omega) \otimes \pi_1^*\omega.$$
\end{ithm}
 The idea of the proof is that, for sufficiently small $U$, the Hecke correspondence $(\pi_1,\pi_2): Y_0(\gp) \to Y \times Y$
 is a closed immersion whose conormal bundle can be described using deformation-theoretic considerations similar to
 those applied to the diagonal embedding $Y \to Y \times Y$
in the proof of the Kodaira--Spencer isomorphism for level prime to $p$.   We mention also that 
 Theorem~\ref{ithm:KS} can be expressed in terms of upper-shriek functors as an isomorphism
 $$\pi_2^!\omega  \cong \pi_1^*\omega,$$
 and that  $\pi_1^*\delta$ and $\pi_2^*\delta$ are canonically isomorphic, so that $\pi_1$ and $\pi_2$ can 
 be interchanged in these statements.

We now describe our cohomological vanishing results.  For simplicity we focus on the vanishing of 
$R^i\pi_{1,*}\CK_{Y_0(\gp)/\CO}$ for $i > 0$.  Since $\pi_{1,K}: Y_0(\gp)_K \to Y_K$ is finite,  the problem
reduces to proving the vanishing of $R^i\overline{\pi}_{1,*}\CK_{\overline{Y}_0(\gp)/\Fpbar}$, where
$\overline{\pi}_1: \overline{Y}_0(\gp)  \to \overline{Y}$ is the reduction\footnote{We work here over
$\Fpbar$, slightly deviating from the notation in the paper.} of $\pi_1$.  Then $\overline{Y}_0(\gp)$
is a local complete intersection which may be written as a union of smooth subschemes $\overline{Y}_0(\gp)_J$
indexed by subsets $J$ of $\Sigma_\gp$, where $\Sigma_\gp$ is the set of embeddings $F_\gp \to \Qpbar$.
The technical heart of this paper is a complete description of the fibres of the restriction of $\overline{\pi}_1$
to $\overline{Y}_0(\gp)_J$, proving the following (Corollary~\ref{cor:fibres}; see there for the definition of
$m$ and $\delta$, and Theorem~\ref{thm:fibres} for an even more precise version).
\begin{ithm}  \label{ithm:fibres}  Every non-empty fibre of $\overline{Y}_0(\gp)_J \to \overline{Y}$ is
isomorphic to $(\PP^1)^m \times S^\delta$, where $S = \Spec(\Fpbar[T]/T^p)$ and $m$ and $\delta$
are determined by $J$.
\end{ithm}
This generalizes Theorem~D of \cite{DKS}, where it is proved under the assumption that $p$ is unramified in $F$.
Work in this direction, for arbitrary behavior of $p$, was also carried out\footnote{There is however a serious gap in the argument in \cite{ERX}; see Remark~\ref{rmk:ERX} below.} in \cite[\S4]{ERX}.  The approach taken
here to the general case relies on a brute force analysis of the local deformation theory of the fibres, undertaken
in \S\S4.3--4.5.  With this in hand, the proof is similar to the one in \cite{DKS}, as is the deduction of cohomological
vanishing via the ``dicing'' argument introduced there.   Furthermore we observe here that similar arguments give the
vanishing of $R^i\pi_{1,*}\CO_{Y_0(\gp)}$ for $i > 0$, and combining these results with Grothendieck--Serre
duality gives the following (see Corollaries~\ref{cor:Iwa} and~\ref{cor:flat0}):
\begin{ithm}  \label{ithm:main} The sheaves $R^i\pi_{1,*}\CK_{Y_0(\gp)/\CO}$ and $R^i\pi_{1,*}\CO_{Y_0(\gp)}$
vanish for $i > 0$, and are locally free of rank $1 + p^{f_\gp}$ if $i=0$.  Furthermore there
is a Hecke-equivariant isomorphism
$$\pi_{1,*} \CO_{Y_0(\gp)}  \stackrel{\sim}{\longrightarrow}  \Shom_{\CO_Y}(\pi_{1,*}\CK_{Y_0(\gp)/\CO},\,\CK_{Y/\CO}).$$
\end{ithm}
As in \cite{DKS}, we in fact treat the degeneracy map $\varphi: Y_1(\gP) \to Y$ where $\gP$ is the radical of $p\CO_F$
and $Y_1(\gP)$ is a model for the Hilbert modular variety of level $U_1(\gP)$,  this being what is needed for the
construction of Galois representations.  The case where $\gP$ contains the radical, for example $\gP = \gp$, is
a mild generalization, and this implies the result for $\pi_1$. It follows also that the models for the Hilbert modular
varieties of level $U_1(\gp)$ and $U_0(\gp)$ defined by
$$\SPEC(\varphi_*\CO_{Y_1(\gp)}) \quad \mbox{and} \quad \SPEC(\pi_{1,*}\CO_{Y_0(\gp)})$$
are finite and flat over $Y$, and hence Cohen--Macaulay.  
We remark that if $p$ is unramified in $F$, such a model for level $U_1(p)$ (and indeed $U_1(p^r)$) was constructed by Kottwitz and Wake in~\cite{KW}; it is natural to ask how these models
may be related.  Since our (finite flat over $Y$) models are defined via Stein factorization and $Y_0(\gp)$ is normal, we obtain the following, perhaps surprising, further consequence (Corollary~\ref{cor:normal}), pointed out
to us by G.~Pappas:
\begin{icor}  \label{icor:normal}  The normalization of $Y$ in $Y_0(\gp)_K$ is flat over $Y$.
\end{icor}

Note also that combining Theorems~\ref{ithm:KS} and~\ref{ithm:main} with the Kodaira--Spencer isomorphism over $Y$
and the isomorphism $\pi_1^*\delta \cong \pi_2^*\delta$ gives an isomorphism
$$\pi_{1,*}\CO_{Y_0(\gp)}  \cong  \Shom_{\CO_Y} (\pi_{1,*}\pi_2^*\omega,\,\omega),$$
and hence a canonical morphism 
$$\pi_{1,*}\pi_2^*\omega \longrightarrow \omega$$
which we call the {\em saving trace}.    Over $Y_K$ it coincides with the composite 
$$(\pi_{1,*}\pi_{2}^*\omega)_K  \longrightarrow  (\pi_{1,*}\pi_1^*\omega)_K  \longrightarrow \omega_K$$
{\em divided by $p^{f_\gp}$}, where the first morphism is induced by the universal isogeny and
the second is the trace relative to $\pi_{1,K}$.  (See (\ref{eqn:compatibility}); note that for $F = \QQ$
this recovers the divisibility of (\ref{eqn:classical}) by $p$ in the case $k=1$.)

The saving trace can be used to define the Hecke operator $T_\gp$ 
on Hilbert modular forms with coefficients in an arbitrary $\CO$-algebra $R$;
for simplicity we restrict our attention in this introduction to $R = \Fpbar$.
First recall that for any $\bk = (k_\theta)_{\theta\in \Sigma} \in  \ZZ^\Sigma$ (and sufficiently small $U$),
there is an associated automorphic bundle $\overline{\CA}_{\bk}$ over $\overline{Y}$
(denoted $\CA_{\bk,\bO,\Fpbar}$ in \S\ref{sec:bundles}; in particular $\overline{\CA}_{\bone}
 = \overline{\omega} = \omega_{\Fpbar}$).   If $F \neq \QQ$, so that the Koecher Principle holds, then 
$M_{\bk}(U;\Fpbar) = H^0(\overline{Y},\overline{\CA}_{\bk})$
is the space of Hilbert modular forms over $\Fpbar$ of weight $\bk$ and level $U$;
if $F = \QQ$, one needs to compactify $\overline{Y}$ and extend $\overline{\CA}_k = 
\overline{\omega}^{\otimes k}$ to define $M_k(U;\Fpbar)$.  Taking the direct limit over all
(sufficiently small) $U$ containing $\GL_2(\CO_{F,p})$, we obtain a $\GL_2(\AA_{F,\f}^{(p)})$-module
$M_{\bk}(\Fpbar)$.   We remark that even forms of paritious weight $\bk$ 
in characteristic $p$ do not necessarily lift to characteristic zero;
indeed this is already the case for $F = \QQ$ and $k=1$.  

If $k_\theta \ge 1$ for all $\theta\in \Sigma_\gp$ (a mild hypothesis in view of
the main result of \cite{DK2} and 
\cite[Prop.~1.13]{DDW}), then the universal isogeny induces a morphism 
$\overline{\pi}_{2}^*\overline{\CA}_{\bk - \bone} \to \overline{\pi}_{1}^*\overline{\CA}_{\bk-\bone}$
over $\overline{Y}_0(\gp)$.  Twisting its direct image by the saving trace 
$\overline{\pi}_{1,*}\overline{\pi}_2^*\overline{\omega} \cong (\pi_{1,*}\pi_2^*\omega)_{\Fpbar} \to \overline{\omega}$
over $Y_{\Fpbar}$  then yields a morphism
$$\overline{\pi}_{1,*} \overline{\pi}_{2}^*\overline{\CA}_{\bk}  \longrightarrow \overline{\CA}_{\bk},$$
and hence an endomorphism $T_\gp$ of $M_{\bk}(U;\Fpbar)$ (see (\ref{eqn:Tp})).
Taking the limit over $U$ thus defines a $\GL_2(\AA_{F,\f}^{(p)})$-equivariant endomorphism
$T_\gp$ of $M_{\bk}(\Fpbar)$.  See Theorem~\ref{thm:hecke} for a statement applicable to
more general weights, coefficients and cohomology degrees.

 \bigskip
 
\noindent {\bf Acknowledgements:}  The author would like to thank Payman Kassaei and Shu Sasaki
for helpful discussions related to this work in the context of our collaborations, and George Pappas
for correspondence pointing out Corollary~\ref{icor:normal}.
The author is also grateful to Liang Xiao for communication clarifying several points in relation to an error in \cite{ERX}, and to David Loeffler, Vincent Pilloni and Preston Wake for illuminating comments on an earlier version of the paper.
Finally the we thank the referee for their careful reading, and for numerous comments that helped improve the exposition.

\section{Hilbert modular varieties} \label{sec:HMV}

\subsection{Notation} \label{ssec:notation}
Fix a prime $p$ and a totally real field $F$ of degree $d = [F:\QQ]$.
We let $\CO_F$ denote the ring of integers of $F$, $\gd$ its different, $S_p$ the set of primes over $p$.
For each $\gp \in S_p$, we write $F_\gp$ for the
completion of $F$ at $\gp$, and we let $F_{\gp,0}$ denote its maximal unramified subextension, $\FF_\gp$ the
residue field $\CO_F/\gp$, $f_\gp$ the residue degree $[F_{\gp,0}:\Q_p] = [\FF_\gp:\FF_p]$,  and 
$e_\gp$ the ramification index $[F_\gp:F_{\gp,0}]$.  We also fix a choice of totally positive $\varpi_\gp \in F$
so that $v_\gp(\varpi_\gp) = 1$ and $v_{\gp'}(\varpi_\gp) = 0$ for all other $\gp' \in S_p$.
In particular $\varpi_\gp$ is a uniformizer in $F_\gp$, and
for each $\gp \in S_p$, and let $E_\gp(u) \in W(\FF_\gp)[u]$ denote its minimal polynomial over $F_{\gp,0}$,
whose ring of integers we identify with $W(\FF_\gp)$.

We adopt much of the notation and conventions of \cite{theta} for 
notions and constructions associated to embeddings of $F$.  In particular,
for each $\gp \in S_p$, we let $\Sigma_{\gp}$ denote the set of embeddings
$F_{\gp} \to \Qpbar$, and we identify $\Sigma := \coprod_{\gp \in S_p} \Sigma_{\gp}$
with the set of embeddings $F \hookrightarrow \Qpbar$ via the canonical bijection.
We also let $\Sigma_{\gp,0}$ denote the set of embeddings $F_{\gp,0} \to \Qpbar$,
which we identify with the set of embeddings $W(\FF_\gp) \to W(\Fpbar)$, or equivalently
$\FF_\gp \to \Fpbar$, and we let $\Sigma_0 = \coprod_{\gp \in S_p} \Sigma_{\gp,0}$. 

For each $\gp \in S_p$, we fix a choice of embedding $\tau_{\gp,0}  \in \Sigma_{\gp,0}$,
and for $i \in \Z/f_\gp\Z$, we let $\tau_{\gp,i} = \phi^i\circ\tau_{\gp,0}$ where $\phi$
is the Frobenius automorphism of $\Fpbar$.   We also fix an ordering 
$\theta_{\gp,i,1},\theta_{\gp,i,2},\ldots,\theta_{\gp,i,e_\gp}$ of the embeddings 
$\theta \in \Sigma_\gp$ restricting to $\tau_{\gp,i}$, so that
$$\Sigma = \{\,\theta_{\gp,i,j}\,|\,\gp \in S_p, i \in \ZZ/f_\gp\ZZ, 1 \le j \le e_\gp\,\}.$$
Finally we let $\sigma$ denote the permutation of $\Sigma$ defined by 
$\sigma(\theta_{\gp,i,j}) = \theta_{\gp,i,j+1}$ if $j \neq e_\gp$ and
$\sigma(\theta_{\gp,i,e_\gp}) = \theta_{\gp,i+1,1}$.

Choose a finite extension $K$ of $\Q_p$ sufficiently large to contain the images
of all $\theta \in \Sigma$; let $\CO$ denote
its ring of integers, $\varpi$ a uniformizer, and $k$ its residue field.  

For $\tau \in \Sigma_{\gp,0}$, we define $E_\tau(u) \in \CO[u]$ to be
the image of $E_\gp$ under the homomorphism induced by $\tau$, so that
$$\CO_F \otimes \CO  =  
\bigoplus_{\gp \in S_p} \bigoplus_{\tau \in \Sigma_{\gp,0}} \CO_{F,\gp} \otimes_{W(\FF_\gp),\tau} \CO
\cong \bigoplus_{\tau \in \Sigma_0}  \CO[u]/(E_\tau).$$
For any $\CO_F \otimes \CO$-module $M$, we obtain a corresponding decomposition
$M = \bigoplus_{\tau \in \Sigma_{0}} M_\tau$, and we also write $M_{\gp,i}$ for
$M_\tau$ if $\tau = \tau_{\gp,i}$.  Similarly if $S$ is a scheme over $\CO$ and
$\CM$ is a quasi-coherent sheaf of $\CO_F \otimes \CO_S$-modules on $S$,
then we write $\CM_\tau = \CM_{\gp,i}$ for the corresponding summand of $\CM$.

For each $\theta = \theta_{\gp,i,j} \in \Sigma$, we factor $E_{\tau} = s_\theta t_\theta$
where $\tau = \tau_{\gp,i}$,
\begin{equation} \label{eqn:st} \begin{array}{rcccc}
s_{\theta} &=&s_{\tau,j}& =& (u-\theta_{\gp,i,1}(\varpi_\gp))\cdots(u-\theta_{\gp,i,j}(\varpi_\gp))\\
\mbox{and}\quad t_{\theta} &=& t_{\tau,j}& = & (u-\theta_{\gp,i,j+1}(\varpi_\gp))\cdots(u-\theta_{\gp,i,e_\gp}(\varpi_\gp));\end{array}\end{equation}
note that the ideals $(s_\theta)$ and $(t_\theta)$ in $\CO[u]/(E_\tau)$ are each other's annihilators, and that the
corresponding ideals in $\CO_{F,\gp} \otimes_{W(\FF_\gp),\tau} \CO$ are independent of the choice of the
uniformizer $\varpi_\gp$.


\subsection{The Pappas-Rapoport model} \label{subsection:PRmodel} 

We now recall the definition, due to Pappas and Rapoport, for smooth integral models
of Hilbert modular varieties of level prime to $p$,

We let $\A_{F,\f} = F \otimes \widehat{\ZZ}$ denote the finite adeles of $F$, and we let
$\A_{F,\f}^{(p)} = \Q \otimes \widehat{\CO}_F^{(p)} = F \otimes \widehat{\ZZ}^{(p)}$, where 
$\widehat{\CO}_F^{(p)}$ (resp.~$\widehat{\ZZ}^{(p)}$)  is the prime-to-$p$
completion of $\CO_F$ (resp.~$\ZZ$), so $\widehat{\ZZ}^{(p)} = \prod_{\ell\neq p} \ZZ_\ell$ and
$\widehat{\CO}_F^{(p)} = \CO_F \otimes \widehat{\ZZ}^{(p)} = \prod_{v\nmid p} \CO_{F,v}$. 

Let $U$ be an open compact subgroup of 
$\GL_2(\widehat{\CO}_F) \subset \GL_2(\A_{F,\f})$ 
of the form $U_pU^p$, where $U_p = \GL_2(\CO_{F,p})$.
Consider the functor which associates, to a locally Noetherian $\CO$-scheme $S$, the
set of isomorphism classes of data $(A,\iota,\lambda,\eta,\CF^\bullet)$, where:
\begin{itemize}
\item $s:A \to S$ is an abelian scheme of relative dimension $d$;
\item $\iota: \CO_F\to \End_S(A)$ is a homomorphism;
\item $\lambda$ is an $\CO_F$-linear quasi-polarization of $A$ such that for each connected component
$S_i$ of $S$, $\lambda$ induces an isomorphism $\gc_i\gd \otimes_{\CO_F}  A_{S_i} \to A_{S_i}^\vee$
for some fractional ideal $\gc_i$ of $F$ prime to $p$;
\item $\eta$ is a level $U^p$ structure on $A$, 
i.e., a $\pi_1(S_i,\overline{s}_i)$-invariant $U^p$-orbit of
$\widehat{\CO}_F^{(p)}$-linear isomorphisms
$$\eta_i :(\widehat{\CO}_F^{(p)})^2 \to \gd \otimes_{\CO_F} T^{(p)}(A_{\overline{s}_i})$$
for a choice of geometric point $\overline{s}_i$ on each connected component $S_i$ of $S$,
where $T^{(p)}$ denotes the product over $\ell \neq p$ of the $\ell$-adic Tate modules,
and $g \in U^p$ acts on $\eta_i$ by pre-composing with right multiplication by $g^{-1}$;
\item $\CF^\bullet$ is a collection of Pappas--Rapoport filtrations, i.e., for
each $\tau = \tau_{\gp,i} \in \Sigma_0$, an increasing filtration of 
$\CO_{F,\gp} \otimes_{W(\FF_\gp),\tau} \CO_S$-modules 
$$0 = \CF_\tau^{(0)} \subset \CF_\tau^{(1)} \subset \cdots 
   \subset \CF_\tau^{(e_\gp - 1)} \subset \CF_\tau^{(e_\gp)} 
   = (s_*\Omega_{A/S}^1)_\tau$$
such that for $j=1,\ldots,e_{\gp}$, the quotient
$ \CF_\tau^{(j)}/\CF_\tau^{(j-1)}$
is a line bundle on $S$ on which $\CO_F$ acts via $\theta_{\gp,i,j}$.
\end{itemize}

If $U^p$ is sufficiently small, then the functor is representable by an infinite disjoint union $\widetilde{Y}_U$
of smooth, quasi-projective schemes of relative dimension $d$ over $\CO$.  Furthermore we have an action of
$\nu \in \CO_{F,(p),+}^\times$ on $\widetilde{Y}_U$ defined by composing the quasi-polarization with $\iota(\nu)$,
and the action factors through a free action of $\CO_{F,(p),+}^\times/(U\cap\CO_F^\times)^2$ by which
the quotient is representable by a smooth quasi-projective scheme of relative dimension $d$ over $\CO$,
which we denote by $Y_U$.   We also have a natural right action of $g \in \GL_2(\A_{F,\f}^{(p)})$ on the inverse
system of schemes $\widetilde{Y}_U$ defined by pre-composing the level structure $\eta$ with right-multiplication
by $g^{-1}$, and the action descends to one on the inverse system of schemes $Y_U$.  Furthermore we have
a compatible system of isomorphisms of the $Y_U(\CC)$ (for any choice of $\CO \to \CC$) with the Hilbert modular varieties
$$\GL_2(F)_+\backslash (\uhp^\Sigma \times \GL_2(\A_{F,\f})/U)$$
under which the action of $\GL_2(\A_{F,\f}^{(p)})$ corresponds to right multiplication.

\subsection{Automorphic bundles} \label{sec:bundles}
Let $\underline{A} = (A,\iota,\lambda,\eta,\CF^\bullet)$ denote the universal object over $\widetilde{Y}_U$.
Recall that $\CH^1_\dr(A/\widetilde{Y}_U)$ is a (Zariski-)locally free sheaf of rank two 
$\CO_F \otimes \CO_{\widetilde{Y}_U}$-modules, and hence
decomposes as $\oplus_{\tau \in \Sigma} \CH^1_\dr(A/\widetilde{Y}_U)_\tau$, where each 
$\CH^1_\dr(A/\widetilde{Y}_U)_\tau$ is locally free of rank two over $\CO_{\widetilde{Y}_U}[u]/(E_\tau)$.

For $\tau = \tau_{\gp,i}$ and $\theta = \theta_{\gp,i,j}$, we let $\CL_\theta$ denote the line bundle
$\CF_\tau^{(j)}/\CF_\tau^{(j-1)}$ on $\widetilde{Y}_U$.  We let $\CG_\tau^{(j)}$ denote the pre-image of $\CF_\tau^{(j-1)}$ in
$\CH^1_\dr(A/\widetilde{Y}_U)_\tau$ under $u - \theta(\varpi_\gp)$, so that $\CP_\theta :=  \CG_\tau^{(j)}/\CF_\tau^{(j-1)}$
is a rank two vector bundle on $\widetilde{Y}_U$, containing $\CL_\theta$ as a sub-bundle.
We let $\CM_\theta = \CP_\theta/\CL_\theta = \CG_\tau^{(j)}/\CF_\tau^{(j)}$ and 
$\CN_\theta = \CL_\theta \otimes_{\CO_{\widetilde{Y}_U}} \CM_\theta = \wedge^2_{\CO_{\widetilde{Y}_U} }\CP_\theta$.
For any $\bk$, $\bm \in \ZZ^\Sigma$, we let $\widetilde{\CA}_{\bk,\bm}$ denote the line bundle
$$\bigotimes_{\theta\in \Sigma}\left( \CL^{\otimes k_\theta}_\theta \otimes_{\CO_{\widetilde{Y}_U}} \CN_\theta^{\otimes m_\theta} \right)$$
on $\widetilde{Y}_U$.  We also let $\widetilde{\omega} = \widetilde{\CA}_{{\bf{1}},\bf{0}}$ and 
$\widetilde{\delta} = \widetilde{\CA}_{\bf{0},\bf{1}}$.  Note that $\widetilde{\omega}$ may be identified
with $\wedge^d_{\CO_{\widetilde{Y}_U}} (s_*\Omega^1_{A/\widetilde{Y}_U})$; we remark that,
less obviously, $\widetilde{\delta}$ may be identified with $\wedge^{2d}_{\CO_{\widetilde{Y}_U}} (\CH^1_{\dr}(A/\widetilde{Y}_U))$
(and more naturally with $\mathrm{Disc}_{F/\QQ} \otimes \wedge^{2d}_{\CO_{\widetilde{Y}_U}} (\CH^1_{\dr}(A/\widetilde{Y}_U))$).

There is a natural action  of $\CO_{F,(p),+}^\times$ on the vector bundles $\CF_\tau^{(j)}$ and $\CG_\tau^{(j)}$ over
the one on $\widetilde{Y}_U$, inducing actions on $\CL_\theta$, $\CP_\theta$, $\CM_\theta$ and $\CN_\theta$;
however the restriction to $(U\cap \CO_F^\times)^2$ is non-trivial, so it fails to define descent data on these
sheaves to $Y_U$ (see \cite[\S3.2]{theta}).  More precisely, if $\nu = \mu^2$ for some $\mu \in U \cap \CO_F^\times$,
then $\nu$ acts on $\CL_\theta$, $\CP_\theta$ and $\CM_\theta$ (resp.~$\CN_\theta$) by $\theta(\mu)$ (resp.~$\theta(\nu)$).
Thus if $\bk + 2\bm$ is parallel, in the sense that $k_\theta + 2 m_\theta$ is independent of $\theta$, then
we obtain descent data on $\widetilde{\CA}_{\bk,\bm}$ (provided $\Nm_{F/\Q}(U\cap \CO_F^\times) = \{1\}$),
and we let $\CA_{\bk,\bm}$ denote the resulting line
bundle on $Y_U$; note in particular this applies to $\widetilde{\omega}$ and $\widetilde{\delta}$, yielding the
line bundles on $Y_U$ which we denote $\omega$ and $\delta$.  More generally if $R$ is an $\CO$-algebra in which the
image of $\prod_\theta \theta(\mu)^{k_\theta + 2m_\theta}$ is $1$ for all $\mu \in U \cap \CO_F^\times$, then
$\widetilde{\CA}_{\bk,\bm,R}:= \widetilde{\CA}_{\bk,\bm} \otimes_{\CO} R$ descends to a line bundle on 
$Y_{U,R}: = Y_U \times_{\CO} R$ which we denote by $\CA_{\bk,\bm,R}$;
note in particular that if $p^N R = 0$ for some $N$, then this applies to all $\bk,\bm$ whenever $U$ is sufficiently small
(depending on $N$).

We also have natural left action of $\GL_2(\A_{F,\f}^{(p)})$ on the various vector bundles
over its right action on the inverse system $\widetilde{Y}_U$.  More precisely suppose that $g \in \GL_2(\A_{F,\f}^{(p)})$
is such that $g^{-1}Ug \subset U'$, where
$U$ and $U'$ are sufficiently small open compact subgroups of $\GL_2(\A_{F,\f})$ containing $\GL_2(\CO_{F,p})$.
Let $\underline{A}' = (A',\iota',\lambda',\eta',\CF^{\prime\bullet})$ denote the universal object over $\widetilde{Y}_{U'}$,
and similarly let $\CL_\theta'$, etc., denote the associated vector bundles.
As described in \cite[\S3.2]{theta}, the morphism $\widetilde{\rho}_g: \widetilde{Y}_U \to \widetilde{Y}_{U'}$
is associated with  a (prime-to-$p$) quasi-isogeny $A \to \widetilde{\rho}_g^*A'$ inducing isomorphisms
$\widetilde{\rho}_g^*\CF_\tau^{\prime(j)} \stackrel{\sim}{\to} \CF_\tau^{(j)}$, which in turn give rise to isomorphisms 
$\widetilde{\rho}_g^*\CL'_\theta \stackrel{\sim}{\to} \CL_\theta$, etc., satisfying the usual compatibilities.
Furthermore, if $\bk + 2\bm$ is parallel, then the resulting isomorphisms
$\widetilde{\rho}_g^*\widetilde{\CA}'_{\bk,\bm} \stackrel{\sim}{\to} \widetilde{\CA}_{\bk,\bm}$ descend to 
isomorphisms $\rho_g^*{\CA}'_{\bk,\bm} \stackrel{\sim}{\to} {\CA}_{\bk,\bm}$, where $\rho_g: {Y}_U \to {Y}_{U'}$
is the morphism obtained by descent from $\widetilde{\rho}_g$, and more generally we obtain 
isomorphisms $\rho_g^*{\CA}'_{\bk,\bm,R} \stackrel{\sim}{\to} {\CA}_{\bk,\bm,R}$ whenever the
image in $R$ of $\prod_\theta \theta(\mu)^{k_\theta + 2m_\theta}$ is $1$ for all $\mu \in U' \cap \CO_F^\times$.

If $R = \CC$ and $\bk + 2\bm$ is parallel, then we recover the usual automorphic line bundles on $Y_U(\C)$ whose
global sections are Hilbert modular forms of weight $(\bk,\bm)$ and level $U$, along with the 
usual\footnote{Up to a factor of $||\det||$, depending on normalizations.}
Hecke action of $\GL_2(\A_{F,\f}^{(p)})$ on their direct limit over $U$.

Finally recall from \cite[\S3.2]{theta} that the quasi-polarization and Poincar\'e duality induce
a perfect alternating pairing on $\CP_\theta$ (depending on the choice of $\varpi_\gp$),
and hence a trivialization of $\CN_\theta$, but their products do not descend to
trivializations of the bundles $\CA_{{\bf{0}},\bm,R}$ over $Y_U$ (if $\bf{m} \neq \bf{0}$).
These are however torsion bundles, which are furthermore non-canonically trivializable
for sufficiently small $U$.  The Hecke action on their global sections is given by
\cite[Prop.~3.2.2]{theta}.

\subsection{Iwahori--level structure} \label{sec:iwahori}

For $\gp \in S_p$, we let $I_0(\gp)$ denote the Iwahori subgroup
$$\left\{\,\left. \left(\begin{array}{cc} a&b \\ c& d \end{array}\right) \in \GL_2(\CO_{F,\gp}) \,\right|\, c \in \gp \CO_{F,\gp} \,\right\}.$$
In this section we recall the definition of suitable integral models of Hilbert modular varieties with such level structure
at a  set of primes over $p$.

Fix an ideal $\gP$ of $\CO_F$ containing the radical of $p\CO_F$, so that $\gP = \prod_{\gp \in P} \gp$
for some subset $P$ of $S_p$.  We are mainly interested in the cases $P =  \{\gp\}$ and $P = S_p$, but
the additional generality introduces no difficulties, and it may be instructive (and amusing) to note how
some of our results specialize to well-known ones in the case $P = \emptyset$.
Let $U$ be an open compact subgroup of $\GL_2(\A_{F,\f})$ as above, so that $U = U_p U^p$ where $U_p = \GL_2(\CO_{F,p})$
and $U^p$ is a sufficiently small open subgroup of $\GL_2(\widehat{\CO}_F^{(p)})$, and we let 
$$U_0(\gP) = \{ \, g \in U \,|\, \mbox{$g_\gp \in I_0(\gp)$ for all $\gp|\gP$}\,\}.$$d
We let $\varpi_\gP = \prod \varpi_\gp$, $f_\gP = \sum f_\gp$, $d_\gP = \sum e_\gp f_\gp$, $\Sigma_\gP = \coprod \Sigma_{\gp}$
and $\Sigma_{\gP,0} = \coprod \Sigma_{\gp,0}$, with the product, sums and unions taken over the set of $\gp$ dividing $\gP$.

For a locally Noetherian $\CO$-scheme $S$, we consider the functor which associates to $S$ the set of isomorphism classes
of triples  $(\underline{A}_1,\underline{A}_2,\psi)$, where $\underline{A}_i = (A_i,\iota_i,\lambda_i,\eta_i,\CF_i^\bullet)$
define elements of $Y_U(S)$ for $i = 1,2$ and $\psi:A_1 \to A_2$ is an isogeny of degree $p^{f_\gP}$ such that
\begin{itemize}
\item $\ker(\psi) \subset A_1[\gP]$;
\item $\psi$ is $\CO_F$-linear, i.e., $\psi\circ \iota_1(\alpha) = \iota_2(\alpha)\circ\psi$ for all $\alpha \in \CO_F$;
\item $\lambda_1\circ \iota_1(\varpi_\gP) = \psi^\vee \circ \lambda_2 \circ \psi$;
\item $\psi\circ \eta_1 = \eta_2$ (as $U^p$-orbits on each connected component of $S$);
\item $\psi^*\CF_2^\bullet \subset \CF_1^\bullet$, i.e., $\psi^*\CF_{\tau,2}^{(j)} \subset \CF_{\tau,1}^{(j)}$ for all
$\gp \in S_p$, $\tau \in \Sigma_{0,\gp}$ and $j = 1,\ldots,e_{\gp}$.
\end{itemize}

For such an isogeny $\psi$, consider also the isogeny $\xi:A_2 \to \gP^{-1}\otimes_{\CO_F} A_1$ such that
$\xi\circ\psi: A_1 \to \gP^{-1} \otimes_{\CO_F} A_1$ is the canonical isogeny with kernel $A_1[\gP]$.  The compatibility
with the quasi-polarizations $\lambda_1$ and $\lambda_2$ then implies the commutativity of the resulting diagram
of $p$-integral quasi-isogenies (over each connected component of the base $S$):
$$\xymatrix{ A_2 \ar[r]^-{\xi}\ar[d]_{\lambda_2} & \gP^{-1}\otimes_{\CO_F} A_1 \ar[d]^{\varpi_\gP \otimes \lambda_2} \\
(\gd \otimes_{\CO_F} A_2)^\vee \ar[r]^{(1\otimes\psi)^\vee} & (\gd \otimes_{\CO_F} A_1)^\vee.}$$
This in turn implies the commutativity of the diagram
$$\xymatrix{ \gP\otimes_{\CO_F} \CH_{\tau,1} \ar[r]^-{\xi_\tau^*}  \ar[d]_{\varpi_\gP^{-1} \otimes \mu_{\tau,1}}  &  \CH_{\tau,2} \ar[d]^{\mu_{\tau,2}} \\
\CH_{\tau,1}^\vee \ar[r]^{(\psi_\tau^*)^\vee} & \CH_{\tau,2}^\vee}$$
for each $\tau \in \Sigma_0$, where $\CH_{\tau,i}$ is the $\tau$-component of the locally free $\CO_F\otimes \CO_S$-module
$\CH^1_\dr(A_i/S)$ for $i=1,2$, the superscript ${}^\vee$ denotes its $\CO_S[u]/(E_\tau)$-dual, and the $\mu_{\tau,i}$ are the
$\tau$-components of the isomorphisms obtained from the polarizations and Poincar\'e duality (see \cite[(4)]{theta}).
The condition that $\psi^*\CF_{\tau,2}^{(j)} \subset \CF_{\tau,1}^{(j)}$ for all $j$ therefore implies that 
$\xi_\tau^*(\gP\otimes_{\CO_F}(\CF_{\tau,1}^{(j)})^\perp) \subset (\CF_{\tau,2}^{(j)})^\perp$, so it follows from
\cite[Lemma~3.1.1]{theta} that
\begin{equation}\label{eqn:dual}  \xi_\tau^*(\gP\otimes_{\CO_F}\CF_{\tau,1}^{(j)}) \subset (\CF_{\tau,2}^{(j)}).\end{equation}

The functor defined above is representable by a scheme $\widetilde{Y}_{U_0(\gP)}$, projective over $\widetilde{Y}_U$ relative
to either of the forgetful morphisms.\footnote{The compatibility in (\ref{eqn:dual}) is included in other references as a further requirement
on the isogeny $\psi$ in the moduli problem defining the model.  As we have shown that it is automatically satisfied, it follows that our
definition agrees with the one used in \cite[\S7]{DK2}
for example, except that there it is assumed that $\gP = \gp$, $U = U(\gn)$ for some $\gn$ and the polarization ideals are fixed (and normalized differently).
Our $\widetilde{Y}_{U_0(\gp)}$ is therefore isomorphic to an infinite disjoint union of the schemes denoted $Y^{\rm PR}$ in {\em loc.~cit.}
in the case $U = U(\gn)$, and its quotient by a finite \'etale cover if $U(\gn) \subset U$.}  
In the next section, we prove that the schemes $\widetilde{Y}_{U_0(\gP)}$ are syntomic of relative 
dimension $d$ over $\CO$ (see also \cite[Prop.~3.3]{ERX}).
Note that we again have a free action of $\CO_{F,(p),+}^\times/(U\cap\CO_F^\times)^2$ for which the quotient is representable by a 
quasi-projective scheme of relative dimension $d$ over $\CO$, which we denote by $Y_{U_0(\gP)}$, and a natural right action of
$\GL_2(\A_{F,\f}^{(p)})$ on the inverse system of the $\widetilde{Y}_{U_0(\gP)}$ (over $U^p$) descending to one on the inverse
system of $Y_{U_0(\gP)}$.  Furthermore the schemes $Y_{U_0(\gP)}$ are independent of the choices of $\varpi_\gp$, and we again
have isomorphisms of their sets of complex points with Hilbert modular varieties
$$Y_{U_0(\gP)}(\CC) \cong \GL_2(F)_+\backslash (\uhp^\Sigma \times \GL_2(\A_{F,\f})/U_0(\gP))$$
under which the action of $\GL_2(\A_{F,\f}^{(p)})$ corresponds to right multiplication.

Consider also the forgetful morphisms $\widetilde{\pi}_i: \widetilde{Y}_{U_0(\gP)} \to \widetilde{Y}_U$ sending $(\underline{A}_1,\underline{A}_2,\psi)$ to
 $\underline{A}_i$ for $i =1,2$, and their product 
 $$\widetilde{h}: = (\widetilde{\pi}_1,\widetilde{\pi}_2): \widetilde{Y}_{U_0(\gP)} \to \widetilde{Y}_U \times_{\CO} \widetilde{Y}_U,$$
 descending to morphisms $\pi_i:Y_{U_0(\gP)} \to Y_U$ and $h: Y_{U_0(\gP)} \to Y_U \times_{\CO} Y_U$, which are again
 independent of the choices of $\varpi_\gp$.
 
 \begin{proposition} \label{prop:immersion}  For sufficiently small $U$, the morphisms $\widetilde{h}$ and $h$
are closed immersions.
 \end{proposition}
 \begpf  Since $\widetilde{h}$ is projective,
 it suffices to prove that it is injective on geometric points and their tangent spaces.
 Furthermore the assertion for $h$ follows from the one for $\widetilde{h}$.
 
 Suppose then that $(\underline{A}_1,\underline{A}_2,\psi)$ and $(\underline{A}_1',\underline{A}_2',\psi')$ correspond
 to geometric points of $\widetilde{Y}_{U_0(\gP)}$ with the same image under $\widetilde{h}$.  Thus there are
 isomorphisms $f_i:A_i \to A_i'$ for $i=1,2$ which are compatible with all auxiliary data.
 Let $\xi:A_2' \to A_1$ be the unique isogeny such that $\xi\circ\psi'\circ f_1 = p$,
 and consider the $\CO_F$-linear endomorphism $\alpha:=  \xi\circ f_2 \circ \psi$ of $A_1$ (where we identify $\CO_F$
 with a subalgebra of $\End(A_1)$ via $\iota_1$).  We wish to prove that $\alpha = p$, as this implies that
 $f_2 \circ \psi = \psi'\circ f_1$, giving the desired injectivity on geometric points. 
 
 First note that the compatibility of $\psi'$ and $f_1$ with quasi-polarizations implies that 
 $\xi^\vee \circ \lambda_1 \circ \varpi_\gP \circ \xi = p^2\lambda_2'$.
 Combining this with the compatibility of $\psi$ and $f_2$ with quasi-polarizations it follows that 
 $\alpha^\vee \circ \lambda_1 \circ \alpha = p^2\lambda_1$.  In particular $F(\alpha)$ is stable
 under the $\lambda_1$-Rosati involution of $\End^0(A_1)$, which sends $\alpha$ to $p^2\alpha^{-1}$.
 By the classification of endomorphism algebras of abelian varieties, it follows that
 either $\alpha \in F$, in which case $\alpha = \pm p$, or $F(\alpha)$ is a quadratic CM-extension of
 $F$, in which case $\alpha$ is in its ring of integers and $\alpha\overline{\alpha} = p^2$.
 
 Next note that the compatibility of $\psi$, $\psi'$, $f_1$ and $f_2$ with level structures implies
 that if $U \subset U(\gn)$, then $\alpha\circ\overline{\eta}_1 = p\overline{\eta}_1$, 
 where $\overline{\eta}_1: (\CO_F/\gn)^2 \cong A_1[\gn]$ is the isomorphism
 induced by $\eta_1$.  It follows that $\alpha - p$ annihilates $A_1[\gn]$, and therefore
 so does $\overline{\alpha} - p$ (writing $\overline{\alpha} = \alpha$ if $\alpha \in F$).
 Letting $\beta$ denote the element $\alpha + \overline{\alpha} - 2p \in \CO_F$, it follows
 that $\beta$ annihilates $A_1[\gn]$, and hence $\beta \in \gn$.  Furthermore since
 $\alpha$ is a root of 
 $$X^2 - (\beta + 2p)X + p^2$$
 and either $\alpha = \pm p$ or $F(\alpha)$ is a CM-extension of $F$, it follows
 that $|\theta(\beta)| \le 4p$ for all embeddings $\theta:F \hookrightarrow \RR$.
 If $\gn$ is such that $N_{F/\QQ}(\gn) > (4p)^{[F:\QQ]}$, this implies that
 $\beta = 0$, and hence $\alpha = p$.
 
 We have now proved that $\widetilde{h}$ is injective on geometric points.
 The injectivity on tangent spaces is immediate from the Grothendieck--Messing Theorem,
a version of which we recall below for convenience and for later use.
 \epf
 
Let $S$ be a scheme and $i:S \hookrightarrow T$ a nilpotent divided power thickening.
If $t: B \to T$ is an abelian scheme of dimension $d$, then the restriction over $S$ of the
crystal $R^1t_{\crys,*}\CO_{B,\crys}$
is canonically identified with $R^1{\overline{t}}_{\crys,*}\CO_{\overline{B},\crys}$, where $\overline{t}: \overline{B} \to S$ is the base-change of
$t:B \to T$.  The image of $t_*\Omega^1_{B/T}$ under the resulting isomorphism
\begin{equation} \label{eqn:crysDR}
\CH^1_\dr(B/T) \stackrel{\sim}{\longrightarrow} (R^1t_{\crys,*}\CO_{B,\crys})_{T}
\stackrel{\sim}{\longrightarrow} (R^1\overline{t}_{\crys,*}\CO_{\overline{B},\crys})_{T}\end{equation}
is thus an $\CO_T$-subbundle $\CV_B$ of $(R^1\overline{t}_{\crys,*}\CO_{\overline{B},\crys})_{T}$
whose restriction to $S$ corresponds to $\overline{t}_*\Omega^1_{\overline{B}/S}$ under the canonical isomorphism
$$\CH^1_\dr(\overline{B}/S) \stackrel{\sim}{\longrightarrow} (R^1\overline{t}_{\crys,*}\CO_{\overline{B},\crys})_{S}  
\longrightarrow i^*(R^1\overline{t}_{\crys,*}\CO_{\overline{B},\crys})_{T}.$$
The Grothendieck--Messing Theorem (as in \cite[Ch.V, \S4]{gro}) states that the functor sending $B$ to
$(\overline{B},\CV_B)$ defines an equivalence between the categories of abelian schemes over $T$
and that of pairs $(A,\CV)$, where $s: A \to S$ is an abelian scheme and $\CV$ is an $\CO_T$-subbundle
of $(R^1s_{\crys,*}\CO_{A,\crys})_{T}$ such that $i^*\CV$ corresponds to $s_*\Omega^1_{A/S}$ under
the canonical isomorphism of $i^*(R^1s_{\crys,*}\CO_{A,\crys})_{T}$ with $\CH^1_\dr(A/S)$.
A morphism of pairs $(A,\CV_1) \to (A',\CV')$ being a pair of morphisms $(A \to A',
\CV' \to \CV)$ satisfying the evident compatibility, it follows that if $t:B \to T$ and $t':B' \to T$ are abelian schemes and $\psi:\overline{B} \to \overline{B}'$ is a morphism of their 
base-changes to $S$, then $\psi$ extends (necessarily uniquely) to a morphism $B \to B'$ if and only if the morphism
$$\psi_{\crys}^*: (R^1\overline{t}'_{\crys,*}\CO_{\overline{B}',\crys})_{T}
\longrightarrow (R^1\overline{t}_{\crys,*}\CO_{\overline{B},\crys})_{T}$$
sends $\CV'$ to $\CV$.
In particular, the functor
sending $B$ to $\overline{B}$ is faithful, which is all that is needed in the proof
of Proposition~\ref{prop:immersion}.

For later reference, we note if the data $\underline{A} = (A,\iota,\lambda,\eta,\CF^\bullet)$
corresponds to an element of $\widetilde{Y}_U(S)$, then to give a lift to an element of
$\widetilde{Y}_U(T)$ is equivalent to giving a lift $\CE^\bullet$ of the Pappas--Rapoport filtrations to
$(R^1s_{\crys,*}\CO_{A,\crys})_{T}$, by which we mean a collection of 
$\CO_{F,\gp} \otimes_{W(\FF_{\gp}),\tau} \CO_T$-submodules 
$$0 = \CE_\tau^{(0)} \subset \CE_\tau^{(1)} \subset \cdots 
   \subset \CE_\tau^{(e_{\gp} - 1)} \subset \CE_\tau^{(e_{\gp})}$$
of  $(R^1s_{\crys,*}\CO_{A,\crys})_{T,\tau}$ for each $\tau = \tau_{\gp,i}$
such that
\begin{itemize}
\item $\CE_\tau^{(j)}/\CE_\tau^{(j-1)}$ is a line bundle on $T$ on which $\CO_F$ acts via $\theta_{\gp,i,j}$.
\item $\iota^*\CE_\tau^{(j)}$ corresponds to $\CF_\tau^{(j)}$ under the canonical isomorphism
$$\iota^*(R^1s_{\crys,*}\CO_{A,\crys})_{T} \cong (R^1s_{\crys,*}\CO_{A,\crys})_{S} \cong \CH^1_\dr(A/S)$$
\end{itemize}
for $j=1,\ldots,e_{\gp}$.   The bijection is defined by sending the data of a lift $\underline{\widetilde{A}} = 
(\widetilde{A},\widetilde{\iota},\widetilde{\lambda},\widetilde{\eta},\widetilde{\CF}^\bullet)$
to the lift of filtrations corresponding to $\widetilde{\CF}^\bullet$ under the canonical isomorphism 
$$\CH^1_\dr(\widetilde{A}/T) \cong 
(R^1\widetilde{s}_{\crys,*}\CO_{\widetilde{A},\crys})_{T} \cong (R^1s_{\crys,*}\CO_{A,\crys})_{T}.$$
The injectivity of the map is a straightforward consequence of the Grothendieck--Messing Theorem;
for the surjectivity one needs also to know that the quasi-polarization $\lambda$
extends to the lift of $A$ associated to $(A,\CV)$ with $\CV = \oplus_\tau \CE_\tau^{(e_{\gp})}$,
which is ensured by \cite[Prop.~2.10]{Vol}.  Since the theorem provides an equivalence of categories, it follows also that if a triple $(\underline{A}_1,\underline{A}_2,\psi)$
corresponds to an element of $\widetilde{Y}_{U_0(\gP)}(S)$, then to give a lift to an element of
$\widetilde{Y}_{U_0(\gP)}(T)$ is equivalent to giving lifts $\CE_i^\bullet$ 
of the Pappas--Rapoport filtrations $\CF_i^\bullet$ to $(R^1s_{i,\crys,*}\CO_{A_i,\crys})_{T}$
for $i=1,2$, such that $\CE_1^\bullet$ and $\CE_2^\bullet$ are compatible with 
$$\psi_{\crys,T}^*: (R^1s_{2,\crys,*}\CO_{A_2,\crys})_{T}  \longrightarrow
 (R^1s_{1,\crys,*}\CO_{A_1,\crys})_{T}.$$

\subsection{Local structure: an example}  \label{sec:lcieg}
In the next section, we will recall the analysis of the local structure of $\widetilde{Y}_{U_0(\gP)}$. First however,
at the suggestion of the referee, we consider the case where $F$ is a quadratic extension ramified at $p$.  This will already illustrate the key ideas and techniques; the general case is then mainly a matter of transforming them into an inductive argument.

To further fix ideas and make the analysis more concrete, let 
$F = \QQ(\sqrt{p})$ and $\gP = \gp = \varpi_{\gp}\CO_F$, where 
$\varpi_{\gp} = \sqrt{p}$.  We also let $\CO = \ZZ_p[\varpi]$
and $\Sigma = \{\theta_1,\theta_2\}$, where $\varpi^2 = p$,
$\theta_1(\varpi_\gp) = \varpi$ and $\theta_2(\varpi_\gp) = -\varpi$,
so that $\CO_F \otimes \CO \cong \CO[u]/(u^2-p)$, in which
$s_1 = u - \varpi$ and $t_1 = u + \varpi$ in the notation of
\S\ref{ssec:notation}.

Suppose now that $y \in \widetilde{Y}_{U_0(\gp)}(\F_p)$, and
let $S$ denote the local ring $\CO_{\widetilde{Y}_{U_0(\gp)},y}$
and $(\underline{A},\underline{A}',\psi)$ denote the corresponding triple over $S$.  Thus the $S[u]/(u^2-p)$-module $H^0(A,\Omega^1_{A/S})$ is free of rank two over $S$, equipped
with a filtration of $S$-modules
$$ 0 = F^{(0)} \subset F^{(1)} \subset F^{(2)} = H^0(A,\Omega^1_{A/S})$$
such that $L_1:=F^{(1)}$ and $L_2:= F^{(2)}/F^{(1)}$ are each free of rank one, $(u-\varpi)F^{(1)} = 0$ and $(u+\varpi)F^{(2)} \subset F^{(1)}$.  Viewing $H^0(A,\Omega^1_{A/S})$ as a submodule of the free rank two $S[u]/(u^2-p)$-module $H^1_{\dR}(A/S)$, and letting $G^{(1)} = (u+\varpi)H^1_{\dR}(A/S)$ (or equivalently, the kernel of $u-\varpi$ on $H^1_\dr(A/S)$)
and $G^{(2)} = (u+\varpi)^{-1}F^{(1)}$ (i.e., the
preimage of $F^{(1)}$ in $H^1_\dr(A/S)$ under $u+\varpi$),
we obtain inclusions of free $S$-modules
$$\begin{array}{ccccccccc}
&&&&F^{(2)}&&&&\\
&&&\rotatebox[origin=c]{45}{$\subset$}&& \rotatebox[origin=c]{-45}{$\subset$}&&&\\
0&\subset&F^{(1)}&&&&G^{(2)}&\subset&H^1_\dr(A/S)\\
&&&\rotatebox[origin=c]{-45}{$\subset$}&& \rotatebox[origin=c]{45}{$\subset$}&&&\\
&&&&G^{(1)}&&&&
\end{array}$$
such that each successive quotient is free of rank one over $S$.

Furthermore the free rank two $S$-modules $P_1:= G^{(1)}$ and
$P_2:= G^{(2)}/F^{(1)}$ are equipped with perfect alternating pairings $\langle\cdot,\cdot\rangle_i$ whose construction we briefly recall (see \cite[\S3.1]{theta} for more details).  Firstly, Poincar\'e duality and the polarization on $A$ yield an isomorphism
$$H^1_\dr(A/S)  \stackrel{\sim}{\lra}
 \Hom_S(\gd^{-1}\otimes_{\CO_F}H^1_\dr(A/S),S) 
 \stackrel{\sim}{\longleftarrow} \Hom_{\CO_F\otimes S}
 (H^1_\dr(A/S),\CO_F\otimes S),$$
and hence a perfect alternating pairing $\langle\cdot,\cdot\rangle_{\dr}$ on $H^1_\dr(A/S)$ over 
$\CO_F\otimes S = S[u]/(u^2 - p)$.
Furthermore one finds that this induces an isomorphism
$$G^{(1)}  \stackrel{\sim}{\lra}
\Hom_S
 (H^1_\dr(A/S)/(u-\varpi)H^1_\dr(A/S),
 (u+\varpi)(\CO_F\otimes S)),$$
and our perfect pairing on $P_1 = G^{(1)}$ (over $S$) is then
obtained from the isomorphisms
$H^1_\dr(A/S)/(u-\varpi)H^1_\dr(A/S)\stackrel{\sim}{\lra}
 G^{(1)}$ and $S\stackrel{\sim}{\lra} (u+\varpi)(\CO_F\otimes S)$
defined by multiplication by $u+\varpi$.
On the other hand,
one finds that $F^{(1)}$ and $G^{(2)}$ are orthogonal complements, and the perfect pairing on $P_2 = G^{(2)}/F^{(1)}$ is then
obtained from the resulting isomorphism
$$G^{(2)}/F^{(1)}  \stackrel{\sim}{\lra}
\Hom_S (G^{(2)}/F^{(1)},
 (u-\varpi)(\CO_F\otimes S)).$$

Similarly we have the filtration
$$ 0 = F'^{(0)} \subset F'^{(1)} \subset F'^{(2)} = H^0(A',\Omega^1_{A'/S}),$$
which we use to define $S$-modules $L_i'\subset P_i'$
for $i=1,2$ such that $L_i'$ and $P_i'/L_i'$,
free of rank one and $P_i'$ is equipped with a perfect alternating
pairing.  Furthermore, the $S[u]/(u^2-\varpi)$-linear map
$\psi_\dr^*:H^1_\dr(A'/S) \to H^1_\dr(A/S)$ induces
morphisms $\psi_i^*:P'_i \to P_i$ restricting to $L_i' \to L_i$,
and the equation 
$\psi\circ\lambda'\circ\psi^\vee = \varpi_{\gp}\lambda$
(where $\lambda$ and $\lambda'$ are the polarizations on
$A$ and $A'$) implies that 
$$\langle \psi_\dr^*(x),\psi_\dr^*(z) \rangle_{\dr}
  = u \langle x,z \rangle'_{\dr},$$
and hence $\langle \psi_i^*(x),\psi_i^*(z)\rangle_i
 = \theta_i(\varpi_\gp)\langle x,z\rangle'_i
 = \pm \varpi\langle x,z\rangle'_i$
for $i=1,2$ and $x,z\in P_i'$.

Consider the behavior of the maps $\psi_i^*$ at the closed point, i.e., $\psi_{i,0}^*:{P}'_{i,0} \to {P}_{i,0}$, where we
use ${\cdot}_0$ to denote $\cdot\otimes_S\F_p$.  
We find that
${\psi}_{i,0}^*$ has rank one (see the argument in the general 
case), and since ${\psi}_{i,0}^*({L}_{i,0}') \subset {L}_{i,0}$, it follows that
${\psi}_{i,0}^*({L}_{i,0}') = 0$ or ${\psi}_{i,0}^*({P}_{i,0}') = {L}_{i,0}$.
To fix ideas even further, let us suppose that both these hold for
$i = 1$, and that only the first holds for $i=2$.
We will prove that in this case the completion of $S$
at its maximal ideal $\gn$ is isomorphic to the $\CO$-algebra
$$\widehat{R} := \CO[[X_1,X_1',X_2]]/(X_1X_1'+\varpi).$$
In order to do so, we will construct a homomorphism 
$\widehat{R} \to \widehat{S}$ using a parametrization of the
$S$-lines $L_i \subset P_i$ and $L_i' \subset P_i'$ with respect to suitable choices of bases for the ambient $S$-planes.
We will then use the Grothendieck--Messing Theorem to prove inductively that the resulting homomorphism
$R_n \to S_n$ is an isomorphism for all $n\ge 1$, where
$S_1 = S/(\varpi,\gn)$, $S_n = S/\gn^n$ for $n\ge 2$, and
$R_n$ is defined similarly.

First note that our assumptions on the $\psi_{i,0}^*$ imply that $F_{0}'^{(2)} = uH^1_\dr(A_0'/\F_p)$. Indeed, we have
$${F}_0'^{(1)} \subset \ker(\psi_{\dr,0}^*)
\quad\mbox{and}\quad
{\psi}_{\dr,0}^*({F}_0'^{(2)}) \subset {F}_0^{(1)}
  = {\psi}_{\dr,0}^*({G}_0'^{(1)}),$$
so that ${F}_0'^{(2)} \subset {G}_0'^{(1)}$, and
comparing dimensions gives equality.  On the other hand,
$${F}_0^{(1)} \subset {\psi}_{\dr,0}^*({G}_0'^{(2)})
  = \psi_{\dr,0}^* (u^{-1}{F}_0'^{(1)})
   \subset {G}_0^{(1)},$$
and comparing dimensions implies 
${\psi}_{\dr,0}^*({G}_0'^{(2)})= {G}_0^{(1)}$,
but ${\psi}_{\dr,0}^*({G}_0'^{(2)})\not\subset 
{F}_0^{(2)}$, so ${F}_0^{(2)} \neq {G}_0^{(1)}$.  It follows that ${F}_0^{(2)}$ is a cyclic $\F_p[u]/(u^2)$-module and ${F}_0^{(1)} = u{F}_0^{(2)}$.  We may therefore choose bases $(x,z)$ for $H^1_\dr(A_0/\F_p)$ and
$(x',z')$ for $H^1_\dr({A}_0'/\F_p)$ as $\F_p[u]/(u^2)$-modules
so that
$${F}_0^{(1)} = \langle ux \rangle,\,
{F}_0^{(2)} = \langle x \rangle,\,
{F}_0'^{(1)} = \langle ux' \rangle \,\,\,\mbox{and}\,\,\,
{F}_0'^{(2)} = \langle ux',uz' \rangle$$
(where $\langle\cdot\rangle$ denotes generation as $\F_p[u]/(u^2)$-modules).
Elementary manipulations show that we can furthermore choose 
these bases so that the matrix of $\psi_{\dr,0}^*$ is
$\smat{0}{1}{-u}{0}$ and
 $\langle x,z \rangle_{\dr,0} = \langle x',z' \rangle'_{\dr,0} = 1$.

We then have bases for $P_{1,0}$ and $P_{1,0}'$
defined by
$$
e_{1,0} = ux,\,\, f_{1,0} = uz, \quad\mbox{and}\quad
e'_{1,0} = ux',\,\, f'_{1,0} = uz',$$
so that
$L_{1,0} = \F_p e_{1,0}$, $L_{1,0}' = \F_p e_{1,0}'$,
$\langle e_{1,0},f_{1,0} \rangle_{1,0}
 = \langle e'_{1,0},f'_{1,0} \rangle'_{1,0}=1$,
and the matrix of 
$\psi_{1,0}^*$ takes the form $\smat{0}{1}{0}{0}$.  
Consider the isomorphisms $P_{1,0} \otimes_{\F_p} S_1 \cong P_{1,1}$ and $P'_{1,0} \otimes_{\F_p} S_1 \cong P'_{1,1}$ 
obtained from the canonical isomorphisms
$$H^1_{\dr}(A_0/\F_p)\otimes_{\F_p} S_1  \stackrel{\alpha}{\longrightarrow} H^1_{\dr}(A_1/S_1)
\quad\mbox{and}\quad
H^1_{\dr}(A'_0/\F_p)\otimes_{\F_p} S_1 \stackrel{\alpha'}{\longrightarrow} H^1_{\dr}(A'_1/S_1)$$
(systematically using $\cdot_n$ for $\otimes_S S_n$).
Their functoriality properties further ensure that
the matrix of $\psi_{1,1}^*$ has the same form as above with respect
to the corresponding bases for $P_{1,1}$ and $P_{1,1}'$,
which we denote $(e_{1,1},f_{1,1})$ and $(e_{1,1}',f_{1,1}')$,
and that $\langle e_{1,1},f_{1,1} \rangle_{1,1}
 = \langle e'_{1,1},f'_{1,1} \rangle'_{1,1}=1$.
However $e_{1,1}$ and $e'_{1,1}$ are no longer bases
for $L_{1,1}$ and $L_{1,1}'$; instead we have 
$L_{1,1} = S_1(e_{1,1} - s_{1,1}f_{1,1})$ and $L_{1,1}'  =
 S_1(e_{1,1}'-s_{1,1}'f_{1,1})$ for some (unique) 
$s_{1,1},s_{1,1}' \in \gn S_1$.  An elementary matrix calculation
then shows that we may lift the chosen bases $(e_1,f_1)$ for $P_1$ and $(e_1',f_1')$ for $P_1'$ over $S$ so that $\langle e_1,f_1 \rangle_1 = \langle e_1',f_1' \rangle_1' = 1$ and the resulting
matrix of $\psi_1^*$ is $\smat{0}{1}{-\varpi}{0}$.  We then
have $L_1 = S(e_1 - s_1 f_1)$ and $L_1' = S(e_1' - s_1'f_1)$
for some uniquely determined $s_1,s_1' \in \gn$, and the fact that
$\psi_1^*(L_1') \subset L_1$ means that
$$\mat{0}{1}{-\varpi}{0}\left(\begin{array}{c}1\\-s_1'\end{array}\right) = \left(\begin{array}{c}-s_1'\\-\varpi\end{array}\right)
 = -s_1'\left(\begin{array}{c}1\\-s_1\end{array}\right),$$
and therefore $s_1s_1' = -\varpi$.

Similarly, using the bases for $P_{2,0}$ and $P_{2,0}'$ defined by
$$ e_{2,0} = x + F_0^{(1)},\,\,f_{2,0} = uz + F_0^{(1)},\quad \mbox{and}\quad
e'_{2,0} = uz' + F_0'^{(1)},\,\,f'_{2,0} = -x' + F_0'^{(1)}$$
gives $L_{2,0} = \F_p e_{2,0}$, $L_{2,0}' = \F_p e_{2,0}'$,
$\langle e_{2,0},f_{2,0} \rangle_{2,0}
 = \langle e'_{2,0},f'_{2,0} \rangle'_{2,0}=1$,
and the matrix of 
$\psi_{2,0}^*$ takes the form $\smat{0}{0}{0}{1}$.  
Note however that the canonical isomorphism $\alpha$
does not (necessarily) send $L_{1,0} = F^{(1)}_0$
to $L_{1,1} = F^{(1)}_1$, and therefore does not
yield an isomorphism between $P_{2,0}\otimes_{\F_p}S_1$
and $P_{2,1}$.  However letting $S_{1/2} = S/I$, where
$I = (s_1,s_1',\gn^2)\subset (\varpi,\gn^2)$, we find
that $\alpha$ and $\alpha'$ still induce isomorphisms $P_{2,0}\otimes_{\F_p}S_{1/2}\cong P_{2,1/2}$ and
$P'_{2,0}\otimes_{\F_p}S_{1/2}\cong P'_{2,1/2}$.
We may then lift the resulting bases for $P_{2,1/2}$ and
$P'_{2,1/2}$ to ones, say $(e_2,f_2)$ for $P_2$ and 
$(e_2',f_2')$ for $P_2'$, such that $\langle e_2,f_2 \rangle_2
 = \langle e_2',f_2' \rangle_2' = 1$ and the matrix for $\psi_2^*$
has the form $\smat{-\varpi}{0}{0}{1}$.  Defining $s_2,s_2' \in \gn$ by $L_2 = S(e_2-s_2 f_2)$ and $L_2' = (e_2' - s_2'f_2')$,
the fact that $\psi_2^*(L_2') \subset L_2$ now translates into the equation $s_2' = -\varpi s_2$.

We now define the homomorphism $\rho: \widehat{R} \to \widehat{S}$ by $X_1 \mapsto s_1$, $X_1' \mapsto s_1'$ and $X_2 \mapsto s_2$, and we will sketch the proof that it is an isomorphism.

We first prove that $\rho$ is surjective, or equivalently, 
that the $\F_p$-vector space $\gn/(\varpi,\gn^2)$ is 
spanned by $s_{1,1}$, $s_{1,1}'$ and $s_{2,1}$, or equivalently, if $\delta:S_1 \to T:=  \F_p[\epsilon]/(\epsilon^2)$
is an $\F_p$-algebra homomorphism such that 
$\delta(s_{1,1}) = \delta(s_{1,1}') = \delta(s_{2,1}) = 0$, then 
$\delta$ factors through $\F_p$.  This in turn is equivalent to the assertion that for such a $\delta$, the triple
$(\underline{\widetilde{A}},\underline{\widetilde{A}}',\widetilde{\psi}) := 
{\delta}^*(\underline{A}_1,\underline{A}_1',\psi_1)$, 
is isomorphic to the base-change (via $\F_p \hookrightarrow T$)
of $(\underline{A}_{0},\underline{A}_{0}',\psi_{0})$. 
By the Grothendieck--Messing Theorem\footnote{As described at the end of \S\ref{sec:iwahori}, but in the simplest case, namely for the thickening $\Spec(\F_p) \hookrightarrow \Spec(T)$, and using the canonical isomorphisms $H^1_\crys(A_0/T)\cong H^1_\dr(A_0/\F_p)\otimes_{\F_p} T$ and $H^1_\crys(A'_0/T)\cong H^1_\dr(A'_0/\F_p)\otimes_{\F_p} T$ to reinterpret (\ref{eqn:crysDR}).}, this in turn is equivalent to the assertion that for $i = 1,2$, 
the $T$-modules
$ F_1^{(i)}\otimes_{S_1} T$
and  $F_{1}'^{(i)}\otimes_{S_{1}} T$
correspond (respectively) to $F_0^{(i)}\otimes_{\F_p} T$
and $F_0'^{(i)}\otimes_{\F_p} T$ under the canonical isomorphisms
$$H^1_{\dr}(A_0/\F_p)\otimes_{\F_p} T 
  \cong H^1_{\dr}(\widetilde{A}/T)\quad
 \mbox{and} \quad H^1_{\dr}(A'_0/\F_p)\otimes_{\F_p} T 
  \cong H^1_{\dr}(\widetilde{A}'/T)$$
(induced by $\alpha$ and $\alpha'$).
Since $F^{(1)}_1 = L_{1,1} = S_1(e_{1,1}-s_{1,1}f_{1,1})$ and
$\delta(s_{1,1}) = 0$, we have that 
$ F_1^{(1)}\otimes_{S_1} T = T(e_{1,1}\otimes 1)$, and
by construction, $e_{1,1}$ corresponds
to $e_{1,0} \otimes 1$ under $\alpha$.
Similarly, we see that $F_{1}'^{(1)}\otimes_{S_{1}} T$
corresponds to $F_0'^{(1)}\otimes_{\F_p} T$.
Note also that $\delta$ factors through $S_{1/2}$, so the
same argument shows that 
$$(F_1^{(2)}/F_1^{(1)})\otimes_{S_1}T = L_{2,1}\otimes_{S_1} T
 = L_{2,1/2}\otimes_{S_{1/2}}T$$
corresponds to $L_{2,0}\otimes_{\F_p} T = (F_0^{(2)}/F_0^{(1)})\otimes_{\F_p}T $ under the isomorphism
induced by $\alpha$, and hence that $ F_1^{(2)}\otimes_{S_1} T$ corresponds to $ F_0^{(1)}\otimes_{\F_p} T$.
Finally the fact that
$L'_{2,1/2}\otimes_{S} T = T(e'_{2,0}\otimes 1)$ follows automatically from the description of $\psi_2^*$, so we similarly obtain the desired conclusion for $F_1'^{(2)}\otimes_{S_1} T$.

In order to prove that $\rho_1:R_1 \to S_1$ is an isomorphism, we could again consider the map on tangent spaces, but in order to be more indicative of the inductive step treating $\rho_n$, we will interpret the argument as the construction of a surjective homomorphism $S_1 \to R_1$.  Note that if 
$\widetilde{y} \in \widetilde{Y}_{U_0(\gp)}(R_1)$ 
is a lift of $y$, then the induced map 
$S = \CO_{\widetilde{Y}_{U_0(\gp)},y} \to R_1$ factors through $S_1$.  Furthermore, by the Grothendieck--Messing Theorem (with $i:\Spec(\F_p) \hookrightarrow \Spec(R_1)$), to give such a lift is equivalent to giving lifts
$$\widetilde{F}^{(1)} \subset \widetilde{F}^{(2)} \subset H^1_{\dr}(A_0/\F_p)\otimes_{\F_p}R_1\quad\mbox{and}\quad
\widetilde{F}'^{(1)} \subset \widetilde{F}'^{(2)} \subset H^1_{\dr}(A'_0/\F_p)\otimes_{\F_p}R_1
$$
of the Pappas--Rapoport filtrations compatible with 
$\psi_{\dr,0}^*\otimes 1$.  More precisely, we require that the $\widetilde{F}^{(i)}$ and $\widetilde{F}'^{(i)}$ (for $i=1,2$) be free $R_1$-modules\footnote{The first bullet renders this equivalent to the successive quotients being free over $R_1$, and the second implies that they are in fact $R_1[u]/(u^2)$-modules.} such that
\begin{itemize}
\item $\widetilde{F}^{(i)}\otimes_{R_1}\F_p = F_0^{(i)}$ and
$\widetilde{F}'^{(i)}\otimes_{R_1}\F_p = F_0'^{(i)}$;
\item $u\widetilde{F}^{(i)} \subset \widetilde{F}^{(i-1)}$
and $u\widetilde{F}'^{(i)} \subset \widetilde{F}'^{(i-1)}$
(where $\widetilde{F}^{(0)} = \widetilde{F}'^{(0)} = 0$);
\item $(\psi_{\dr,0}^*\otimes 1)(\widetilde{F}'^{(i)}) \subset \widetilde{F}^{(i)}$.
\end{itemize}
We will write down such lifts in terms of the bases
$(\widetilde{x},\widetilde{z})$ and $(\widetilde{x}',\widetilde{z}')$, where $(x,z)$ and $(x',z')$ are the bases
already chosen for $H^1_\dr(A_0/\F_p)$ and $H^1_\dr(A_0'/\F_p)$ and
$\widetilde{x} = x \otimes 1$, etc.  Letting 
$\widetilde{P}_1 = uH^1_{\dr}(A_0/\F_p)\otimes_{\F_p} R_1$ and
$\widetilde{P}'_1 = uH^1_{\dr}(A'_0/\F_p)\otimes_{\F_p} R_1$,
note that the matrix of the $R_1$-linear map 
$\widetilde{\psi}_1:\widetilde{P}_1' \to \widetilde{P}_1$
induced by $\psi_{\dr,0}^*$ with respact to the bases
$(u\widetilde{x},u\widetilde{z})$ and $(u\widetilde{x}',u\widetilde{z}')$ is $\smat{0}{1}{0}{0}$. It follows that
$$\widetilde{F}^{(1)} = 
R_1 u(\widetilde{x} - \overline{X}_1\widetilde{z})
 \quad\mbox{and}\quad
\widetilde{F}'^{(1)} = 
R_1 u(\widetilde{x}' - 
\overline{X}'_1\widetilde{z}')$$
are free $R_1$-modules such that 
$\widetilde{\psi}_1^*(\widetilde{F}'^{(1)}) \subset
 \widetilde{F}^{(1)}$, where $\overline{X}_1$ denotes the
image of $X_1$ in $R_1$, etc.  Letting $\widetilde{P}_2$
denote the free rank two $R_1$-module
$(u^{-1}\widetilde{F}^{(1)})/\widetilde{F}^{(1)}$ and similarly
defining $\widetilde{P}'_2$, one finds that the matrix of 
the $R_1$-linear
map $\widetilde{\psi}_2:\widetilde{P}_2' \to \widetilde{P}_2$
induced by $\psi_{\dr,0}^*$ is $\smat{0}{0}{0}{1}$ with respect
to the bases 
$$(\widetilde{x} - \overline{X}_1\widetilde{z} + 
\widetilde{F}^{(1)}, u\widetilde{z} + \widetilde{F}^{(1)}) 
\quad\mbox{and}\quad
(\overline{X}_1\widetilde{x}'+u\widetilde{z}' + \widetilde{F}'^{(1)},
-\widetilde{x}'+\overline{X}_1'\widetilde{z}'
 + \widetilde{F}'^{(1)}).$$
(Note that $\overline{X}_1\widetilde{x}'+u\widetilde{z}'
 \in u^{-1}\widetilde{F}'^{(1)}$ since
$u(\overline{X}_1\widetilde{x}'+u\widetilde{z}') = 
  u\overline{X}_1(\widetilde{x}'-\overline{X}'_1\widetilde{z}')$.)
We thus find that 
$$\widetilde{F}^{(2)} = \langle \widetilde{x}-
 (\overline{X}_1+u\overline{X}_2)\widetilde{z} \rangle,
\quad\mbox{and}\quad
\widetilde{F}'^{(2)} = \langle u(\widetilde{x}'-\overline{X}'_1\widetilde{z}'),
\overline{X}_1\widetilde{x}'+u\widetilde{z}' \rangle$$
yields lifts of the Pappas--Rapoport filtrations with the
desired properties (where $\langle\cdot\rangle$ here denotes generation as $R_1[u]/(u^2)$-modules).
As already explained, this yields
a homomorphism $S_1 \to R_1$.  Furthermore our definition of the parameters $s_{1,1}$, $s_{1,1}'$ and $s_{2,1}$ ensures that 
$s_{1,1} \mapsto \overline{X}_1$, $s_{1,1}' \mapsto 
\overline{X}_1'$ and $s_{2,1} \mapsto \overline{X}_2 \bmod (\overline{X}_1,\overline{X}_1')$, so the
map is surjective.  Since $R_1$ and $S_1$ are Artinian, 
and $\rho_1:R_1 \to S_1$ is also surjective, it follows that
the lengths of $R_1$ and $S_1$ are the same, and
therefore that $\rho_1$ is an isomorphism.

Suppose now that $n \ge 1$ and that $\rho_n:R_n \to S_n$ is an isomorphism.  In order to prove that $\rho_{n+1}$ is an isomorphism, it suffices (arguing as above) to show\footnote{In this case, we have that $R_2 = R_1$ and $S_2 = S_1$, so we could assume $n \ge 2$ if we wanted.} that the composite $S \onto S_n \stackrel{\rho_n^{-1}}{\longrightarrow}R_n$ lifts to a surjective homomorphism $S \to R_{n+1}$. Using the thickening 
$\Spec(S_n) \hookrightarrow \Spec(R_{n+1})$, with
$H^1_\crys(A_n/R_{n+1})$ instead of $H^1_\dr(A_0/\F_p)\otimes_{\F_p} R_1$, the construction of lifts of the
Pappas--Rapoport filtrations on $H^1_\dr(A_n/S_n)$ is similar 
to the one above, once one has
abstracted the argument establishing the existence of suitable
bases for $\widetilde{P}_i$ and $\widetilde{P}'_i$ (see the proof in the general case).  While the isomorphism $H^1_\crys(A_0/R_1) \cong H^1_\dr(A_0/\F_p)\otimes_{\F_p} R_1$ was implicitly used in
the proof of the surjectivity of the resulting homomorphism
$S_1 \to R_1$, no such isomorphism is available relative
to the thickening $\Spec(S_n) \hookrightarrow \Spec(R_{n+1})$,
nor is it needed, since the surjectivity of the resulting
homomorphism $S_{n+1} \to R_{n+1}$ is immediate from that of its
composite with $R_{n+1} \to R_n$.

\subsection{Local structure: in general}  \label{sec:lci}
We now proceed to analyze the local structure of $\widetilde{Y}_0(\gP)$ for arbitrary $F$ and $\gP$,
proving\footnote{The result is essentially
Proposition~3.3 of \cite{ERX}, but we will give a complete proof for several reasons: 1) our different
definition of the moduli problem, 2) our greater degree of generality, and 
3) the omission of the proof in {\em loc.~cit.} that the maps
to Grassmannians induce isomorphisms on tangent spaces, where our perspective
provides the basis for the construction of an Iwahori-level Kodaira--Spencer isomorphism.}
in particular that the schemes are reduced and syntomic, i.e., flat local complete intersections, over $\CO$.

Let $y$ be a closed point of $\widetilde{Y}_{U_0(\gP)}$ in characteristic $p$.  We let $S$ denote
the local ring $\CO_{\widetilde{Y}_{U_0(\gP)},y}$, $\gn$ its maximal ideal, $S_0$ the residue field,
$S_1 = S/(\gn^2,\varpi)$ and $S_n = S/\gn^n$ for $n > 1$.  
We let $(\underline{A},\underline{A}',\psi)$ denote the corresponding triple over $S$, and
similarly write $(\underline{A}_n,\underline{A}_n',\psi_n)$ for the triple over $S_n$ for $n \ge 0$.
For each $\tau = \tau_{\gp,i} \in \Sigma_0$, 
we view the Pappas--Rapoport filtrations as being defined by $S[u]/(E_\tau)$-submodules
$F_\tau^{(j)} \subset H^0(A,\Omega^1_{A/S})_\tau$ and their quotients 
$F_{\tau,n}^{(j)} \subset H^0(A_n,\Omega^1_{A_n/S_n})_\tau$
over $S_n[u]/(E_\tau)$ for $n \ge 0$, and similarly with $A$ replaced by $A'$.
For each $\theta = \theta_{\gp,i,j} \in \Sigma$, we let 
$L_{\theta} = F_\tau^{(j)}/F_\tau^{(j-1)}$ and $P_{\theta} = G_\tau^{(j)}/F_\tau^{(j-1)}$,
where $G_\tau^{(j)}$ denotes the preimage of $F_\tau^{(j-1)}$ under multiplication by
$u - \theta(\varpi_{\gp})$ on $H^1_{\dr}(A/S)_\tau$, so that $L_\theta$ is 
a rank one summand of the free rank two $S$-module $P_\theta$.
We similarly define $L_\theta'$ and $P_\theta'$, and systematically use
$\cdot_n$ to denote $\cdot \otimes_S S_n$.

By definition, the isogeny $\psi:A \to A'$ induces morphisms 
$\psi_\theta^*:P_\theta' \to P_\theta$ for $\theta \in \Sigma$
such that $\psi_\theta^*(L_\theta') \subset L_\theta$.   
Note that $\psi_\theta^*$ is an isomorphism for $\theta \not\in \Sigma_\gP$. 
Recall also from \cite[\S3.1]{theta} that Poincar\'e duality and the quasi-polarizations
give rise to perfect alternating pairings on $P_\theta$ and $P_\theta'$, which we denote
by $\langle\cdot,\cdot\rangle_\theta$ and $\langle\cdot,\cdot\rangle_\theta'$.
Furthermore the compatibility of $\psi$ with the quasi-polarizations
implies the commutativity of the resulting diagram
$$\xymatrix{ \wedge^2_S P'_\theta  \ar[r]^-{\wedge^2\psi_\theta^*} \ar[d]^-{\wr} & \wedge^2_S P_\theta \ar[d]^-{\wr} \\
S \ar[r]^-{\theta(\varpi_\gP)}  & S.}$$

Note in particular if $\theta \in \Sigma_\gP$, then the $S_0$-linear map $\psi_{\theta,0}^*:P_{\theta,0}' \to P_{\theta,0}$
is not invertible (since $\wedge^2_{S_0}\psi_{\theta,0}^* = 0$); we will show that it has rank one.
To that end, let $W = W(S_0)$
and consider the free rank two $\CO_F\otimes W$-module
$D:=H^1_\cris(A_0/W)$.  Thus $D$ decomposes as $\oplus_\tau D_\tau$ where each $D_\tau$ is free
of rank two over $W[u]/(E_\tau)$, and we let $\widetilde{G}_\tau^{(j)}$ denote the preimage of $G_{\tau,0}^{(j)} = u^{-1}F_{\tau,0}^{(j-1)}$
in $D_\tau$ under the canonical projection to  $D_\tau/pD_\tau \cong H^1_\dr(A_0/S_0)_\tau$.
Note that if $\tau = \tau_{\gp,i}$, then we have a sequence of inclusions
$$u^{e_\gp - 1}D_\tau = \widetilde{G}_{\tau}^{(1)} \subset \widetilde{G}_\tau^{(2)} \subset \cdots
\subset \widetilde{G}_\tau^{(e_\gp-1)}\subset \widetilde{G}_\tau^{(e_\gp)}$$
such that $\dim_{S_0}(\widetilde{G}_\tau^{(j)}/pD_\tau) = 
 \dim_{S_0}(G_{\tau,0}^{(j)}) = j+1$
and the projection to $H^1_\dr(A_0/S_0)_\tau$ identifies $\widetilde{G}_\tau^{(j)}/u\widetilde{G}_\tau^{(j)}$
with $P_{\theta,0}$ for $j=1,\ldots,e_\gp$ and $\theta = \theta_{\gp,i,j}$.
Similarly decomposing $D':=H^1_\cris(A'_0/W)$ as  $\oplus_\tau D'_\tau$, and
defining $\widetilde{G}_\tau^{\prime(j)}$ as the preimage of $G_{\tau,0}^{\prime(j)}$ in $D'_\tau$,
the isogeny $\psi_0:A_0 \to A_0'$ induces an injective $\CO_F\otimes W$-linear homomorphism
$D' \to D$ whose cokernel is killed by $u$ and has dimension $f_\gP$ over $S_0$.  Furthermore the restrictions
$\widetilde{\psi}_\tau^* :D'_\tau \to D_\tau$ are isomorphisms for $\tau \not\in \Sigma_{\gP,0}$,
and restrict to homomorphisms $\widetilde{G}_\tau^{\prime(j)} \to \widetilde{G}_\tau^{(j)}$ whose
cokernel projects onto that of $\psi_{\theta,0}^*$ for each $\tau = \tau_{\gp,i} \in \Sigma_{\gP,0}$, $\theta = \theta_{\gp,i,j}$.
Therefore it suffices to prove that the cokernel of each $\widetilde{G}_\tau^{\prime(j)} \to \widetilde{G}_\tau^{(j)}$,
necessarily non-trivial if $\tau \in \Sigma_{\gP,0}$, has length (at most) one over $W$.  For $j=1$, this follows from the
identifications $u^{e_\gp - 1}D_\tau = \widetilde{G}_{\tau}^{(1)}$ and $u^{e_\gp - 1}D'_\tau = \widetilde{G}_{\tau}^{\prime(1)}$,
which implies that the sum over $\tau$ of the lengths of the cokernels is $f_\gP$.  For $j > 1$, it then follows by induction
from the injectivity of $\widetilde{\psi}_\tau^*$ and the fact that $\widetilde{G}_\tau^{(j)}/\widetilde{G}_\tau^{(j-1)}$
and $\widetilde{G}_\tau^{\prime(j)}/\widetilde{G}_\tau^{\prime(j-1)}$ each have length one.

We then define
$$\Sigma_y = \{\,\theta \in \Sigma_\gP\,|\, \im(\psi^*_{\theta,0}) = L_{\theta,0}\,\}\quad\mbox{and}\quad
\Sigma'_y = \{\,\theta \in \Sigma_\gP\,|\, \ker(\psi^*_{\theta,0}) = L'_{\theta,0}\,\}.$$
Note that $\Sigma_\gP = \Sigma_y \cup \Sigma'_y$, and define
$$\begin{array}{c} R =  \CO [X_\theta,X_\theta']_{\theta \in \Sigma}/ (g_\theta)_{\theta\in \Sigma}, \\ \ \\
\mbox{where}\quad
    g_\theta = \left\{\begin{array}{ll} 
       X_\theta - \theta(\varpi_\gP) X_\theta', & \mbox{if $\theta\in \Sigma \setminus \Sigma_y'$,}\\
       X'_\theta - \theta(\varpi_\gP) X_\theta, & \mbox{if $\theta\in \Sigma'_y \setminus \Sigma_y$,}\\
       X_\theta X'_\theta + \theta(\varpi_\gP), & \mbox{if $\theta\in \Sigma_y \cap \Sigma'_y$.} \end{array} \right. \end{array}$$
We will show that a suitable parametrization of the lines $L_\theta$ and $L_\theta'$ by the variables
$X_\theta$ and $X_\theta'$ defines an $\CO$-algebra homomorphism $R \to S$ inducing
an isomorphism $W \otimes_{W(k)} \widehat{R}_\gm \stackrel{\sim}{\to} \widehat{S}_\gn$,
where $\gm$ is the maximal ideal of $R$ generated by $\varpi$ and the variables
$X_\theta$ and $X_\theta'$ for $\theta \in \Sigma$.  (Recall that $k$ is the residue field of $\CO$, so that $R/\gm = k$, and that
$W(k)$ denotes the ring of Witt vectors of $k$.)

Let $\alpha$ denote the composite of the canonical isomorphisms
$$H^1_\dr(A_0/S_0) \otimes_{S_0} S_1
  \stackrel{\sim}{\longrightarrow}
  H^1_\cris(A_0/S_1) 
  \stackrel{\sim}{\longrightarrow}
   H^1_\dr(A_1/S_1)$$
(where we write $H^1_\cris(A_0/S_1)$ for 
$(R^1s_{0,\crys,*}\CO_{A_0,\crys})(S_1)$), 
and similarly define 
$\alpha': H^1_\dr(A_0'/S_0) \otimes_{S_0} S_1 
   \stackrel{\sim}{\longrightarrow} H^1_\dr(A_1'/S_1)$.
   Thus $\alpha$
   and $\alpha'$ are $\CO_F\otimes S_1$-linear and compatible
   with $\psi_0^*$ and $\psi_1^*$.  
It follows that for each $\tau = \tau_{\gp,i} \in \Sigma_0$,
$\alpha$ and $\alpha'$ decompose as direct sums of $S_1[u]/u^{e_{\gp}}$-linear
isomorphisms $\alpha_\tau$ and $\alpha'_{\tau}$ such that the diagrams
$$\xymatrix{  H^1_{\dr}(A_0'/S_0)_\tau \otimes_{S_0} S_1 \ar[r]^-{\psi_{\tau,0}^* \otimes 1} \ar[d]_-{\alpha_\tau'} &
H^1_{\dr}(A_0/S_0)_\tau \otimes_{S_0} S_1  \ar[d]_-{\alpha_\tau} \\
H^1_{\dr}(A_1'/S_1)_\tau \ar[r]^-{\psi_{\tau,1}^*} & H^1_\dr(A_1/S_1)_\tau}$$
commute.  We then consider the chain of ideals
$$ (\gn^2,\varpi) = I_\tau^{(0)} \subset I_\tau^{(1)} \subset \cdots \subset I_\tau^{(e_{\gp} -1)} \subset I_\tau^{(e_{\gp})} \subset \gn$$
where $I_\tau^{(j)}$ is defined by the vanishing of the maps
$$F_{\tau,0}^{(\ell)}\otimes_{S_0} S_1 \longrightarrow H^1_\dr(A_1/S_1)_\tau / F_{\tau,1}^{(\ell)}
  \quad\mbox{and}\quad
F_{\tau,0}^{\prime(\ell)}\otimes_{S_0} S_1 \longrightarrow H^1_\dr(A'_1/S_1)_\tau / F_{\tau,1}^{\prime(\ell)}$$
induced by $\alpha_\tau$ and $\alpha_\tau'$ for $\ell = 1,\ldots,j$.  For $\theta = \theta_{\gp,i,j}$, we
let $S_\theta = S/I_\tau^{(j-1)}$, so that $\alpha_\tau$ and $\alpha_\tau'$ restrict in particular to
$S_\theta[u]$-linear isomorphisms
$$F_{\tau,0}^{(j-1)} \otimes_{S_0}  S_\theta \stackrel{\sim}{\longrightarrow} F_\tau^{(j-1)} \otimes_S S_\theta
  \quad\mbox{and}\quad
    F_{\tau,0}^{\prime(j-1)} \otimes_{S_0}  S_\theta \stackrel{\sim}{\longrightarrow} F_\tau^{\prime(j-1)} \otimes_S S_\theta,$$
and hence to $S_\theta[u]$-linear isomorphisms
$$G_{\tau,0}^{(j)} \otimes_{S_0}  S_\theta \stackrel{\sim}{\longrightarrow} G_\tau^{(j)} \otimes_S S_\theta
  \quad\mbox{and}\quad
    G_{\tau,0}^{\prime(j)} \otimes_{S_0}  S_\theta \stackrel{\sim}{\longrightarrow} G_\tau^{\prime(j)} \otimes_S S_\theta.$$
We thus obtain $S_\theta$-linear isomorphisms $\alpha_\theta$ and $\alpha_\theta'$ such that the
diagram
$$\xymatrix{  P_{\theta,0}' \otimes_{S_0} S_\theta \ar[r]^-{\psi_{\theta,0}^* \otimes 1} \ar[d]_-{\alpha_\theta'} &
                  P_{\theta,0} \otimes_{S_0} S_\theta  \ar[d]_-{\alpha_\theta} \\
                  P_\theta' \otimes_{S} S_\theta \ar[r]^-{\psi_\theta^*\otimes 1} & P_\theta \otimes_S S_\theta}$$
commutes.  Furthermore $\alpha_\theta$ is compatible with the pairings $\langle \cdot,\cdot \rangle_{\theta,0}$
and $\langle \cdot,\cdot \rangle_\theta$, and similarly $\alpha_\theta'$ is compatible with
$\langle \cdot,\cdot \rangle'_{\theta,0}$ and $\langle \cdot,\cdot \rangle'_\theta$.

We claim that bases $\CB_\theta = (e_\theta,f_\theta)$ for $P_\theta$ and 
$\CB_\theta' = (e_\theta',f_\theta')$ for $P_\theta'$  may be chosen so that
\begin{itemize}
\item $e_\theta \otimes_S S_\theta = \alpha_\theta(L_{\theta,0} \otimes_{S_0} S_\theta)$ and
        $e'_\theta \otimes_S S_\theta = \alpha'_\theta(L'_{\theta,0} \otimes_{S_0} S_\theta)$;
\item $\langle e_\theta,f_\theta \rangle_\theta = \langle e'_\theta, f'_\theta \rangle'_\theta = 1$;
\item the matrix of $\psi_\theta^*$ with respect to $\CB_\theta$ and $\CB_\theta'$ is
\begin{equation}\label{eqn:matrices}
\mat{1}{0}{0}{\theta(\varpi_\gP)},\quad
\mat{\theta(\varpi_\gP)}{0}{0}{1}\quad\mbox{or}\quad\mat{0}{1}{-\theta(\varpi_\gP)}{0}
\end{equation}
according to whether\footnote{Note that if $\theta\not\in \Sigma_\gP$, then $\theta(\varpi_\gP)\in \CO^\times$
and $\psi_\theta^*$ is an isomorphism.  This is the simplest case, but we toss it in
with the case $\theta \in \Sigma_y \setminus \Sigma'_y$ for convenience; note that we could just as well
have included it with the case $\theta \in \Sigma'_y \setminus \Sigma_y$.}
 $\theta \in \Sigma \setminus \Sigma_y'$,
$\Sigma_y' \setminus \Sigma_y$ or $\Sigma_y \cap \Sigma_y'$.
\end{itemize}
Indeed first choose bases $\CB_{\theta,0} = (e_{\theta,0},f_{\theta,0})$ for $P_{\theta,0}$ and 
$\CB_{\theta,0}' = (e_{\theta,0}',f_{\theta,0}')$ for $P_{\theta,0}'$ so that $L_{\theta,0} = S_0e_{\theta,0}$,
$L_{\theta,0}' = S_0e_{\theta,0}'$, $\langle e_{\theta,0},f_{\theta,0} \rangle_{\theta,0} =
\langle e'_{\theta,0},f'_{\theta,0} \rangle'_{\theta,0} = 1$
 and $\psi_{\theta,0}^*$ has the required form ($\!\bmod\,\gn$), 
and lift the bases 
$(\alpha_\theta(e_{\theta,0} \otimes 1),\alpha_\theta(f_{\theta,0}\otimes 1))$ and
$(\alpha'_\theta(e'_{\theta,0} \otimes 1),\alpha'_\theta(f'_{\theta,0}\otimes 1))$
arbitrarily to bases $\CB_\theta$ and $\CB_\theta'$ over $S$ satisfying the condition on the
pairings.  The matrix  $T$ of $\psi_\theta^*$
with respect to $\CB_\theta$ and $\CB_\theta'$ then has the required form mod $I_\tau^{(j-1)}$ and
satisfies $\det(T) =  \theta(\varpi_\gP)$.  We may then replace $\CB_\theta$ and
$\CB'_\theta$ by $\CB_\theta U$ and $\CB_\theta' U'$ for some
$ U,\,U'\,\in \ker(\SL_2(S) \to \SL_2(S_\theta))$
so that the resulting matrix $U^{-1}TU'$ has the required form.

We now have $L_\theta = S(e_\theta - s_\theta f_\theta)$ and $L'_\theta = S(e'_\theta - s'_\theta f'_\theta)$
for unique $s_\theta, s'_\theta \in \gn$, and the fact that $\psi_\theta^*(L_\theta') \subset L_\theta$ means that
\begin{itemize}
\item $s_\theta = \theta(\varpi_\gP)s'_\theta$ if $\theta \in \Sigma  \setminus \Sigma_y'$;
\item $s'_\theta = \theta(\varpi_\gP)s_\theta$ if $\theta \in \Sigma_y' \setminus \Sigma_y$;
\item $s_\theta s'_\theta = - \theta(\varpi_\gP)$ if $\theta \in \Sigma_y \cap \Sigma'_y$.
\end{itemize}
Thus the $\CO$-algebra homomorphism $\CO [X_\theta,X_\theta']_{\theta \in \Sigma} \longrightarrow S$
defined by $X_\theta \mapsto s_\theta$, $X'_\theta \mapsto s'_\theta$ factors through $R$.
We will prove that the resulting morphism $\rho: W\otimes_{W(k)} \widehat{R}_{\gm} \to \widehat{S}_\gn$
is an isomorphism.

We let $R_1 = W \otimes_{W(k)} R/(\gm^2,\varpi)$ and $R_n = W \otimes_{W(k)} R/\gm^n$, so we obtain
morphisms $\rho_n:  R_n \to S_n$.  It suffices to prove that each of these is an isomorphism, which we
achieve by induction on $n$.

We start by showing that $\rho_1$ is surjective, i.e., that the homomorphism induced by $\rho$ on tangent spaces
is injective.  First note that $I_\tau^{(j)} = (I_\tau^{(j-1)},s_\theta,s'_\theta)$ for all $\tau = \tau_{\gp,i}$, $\theta = \theta_{\gp,i,j}$.
Indeed by construction, the ideal $I_\tau^{(j)}/I_\tau^{(j-1)}$ of $S_\theta$ is defined by the vanishing of the maps
$$L_{\theta,0}\otimes_{S_0} S_\theta \longrightarrow (P_\theta/L_\theta) \otimes_S S_\theta \quad\mbox{and}\quad
L'_{\theta,0}\otimes_{S_0} S_\theta \longrightarrow (P'_\theta/L'_\theta) \otimes_S S_\theta$$
induced by $\alpha_\theta$ and $\alpha'_\theta$, and these 
are precisely $s_\theta$ and $s'_\theta$ ($\bmod\, I_{\tau}^{(j-1)}$) relative to the bases
$\alpha_\theta^{-1}(e_\theta \otimes 1)$, $(f_\theta + L_\theta) \otimes 1$, $\alpha_\theta^{\prime-1}(e_\theta'\otimes 1)$
and $(f_\theta'+L'_\theta) \otimes 1$.  
It follows that the image of $\rho_1$ contains each $I_\tau^{(e_{\gp})}$, so letting $T$
denote the quotient of $S$ by the ideal generated by $I_\tau^{(e_{\gp})}$ for all $\tau \in \Sigma_0$, it suffices to prove
the projection $T \to S_0$ is an isomorphism, i.e., that the triple $(\underline{A},\underline{A}', \psi) \otimes_{S} T$ 
 is isomorphic to $(\underline{A}_0,\underline{A}'_0, \psi_0)\otimes_{S_0} T$.
By the definition of $T$, we have that
$F_{\tau,0}^{(j)} \otimes_{S_0} T$ corresponds to $F_{\tau,1}^{(j)} \otimes_{S_1} T$
under $\alpha_\tau \otimes_{S_1} T$ for all $\tau$ and $j$, and similarly for $F_{\tau,0}^{\prime(j)}$, 
$F_{\tau,1}^{\prime(j)}$ and $\alpha_\tau'$, so the conclusion is immediate from the equivalence
established at the end of \S\ref{sec:iwahori}.

To deduce that $\rho_1$ is an isomorphism, it suffices to prove that $\lg(R_1) \le \lg(S_1)$, which will
follow from the existence of a surjection $S \to R_1$.  We construct such a surjection by defining suitable lifts
of the Pappas--Rapoport filtrations $F_{0}^\bullet$ and $F_{0}^{\prime\bullet}$
to $H^1_\crys(A_0/R_1)$ and $H^1_\crys(A_0'/R_1)$.  For $\tau = \tau_{\gp,i}$ and $\theta = \theta_{\gp,i,j}$,
we let $R_\theta = R_1/J_\tau^{(j-1)}$ where the ideals $J_\tau^{(j)}$ are defined 
inductively by $J_\tau^{(0)} = (0)$ and $J_\tau^{(j)} = (J_\tau^{(j-1)}, X_\theta, X'_\theta)$.
We will inductively define chains of $R_1[u]/u^{e_{\gp}}$-submodules
$$\begin{array}{ll} & 0 = \widetilde{F}_\tau^{(0)} \subset \widetilde{F}_\tau^{(1)} \subset \cdots \subset \widetilde{F}_{\tau}^{(e_{\gp})} \subset H^1_\crys(A_0/R_1)_\tau \\
 \mbox{and}\quad &
0 = \widetilde{F}_\tau^{\prime(0)} \subset \widetilde{F}_\tau^{\prime(1)} \subset \cdots \subset \widetilde{F}_{\tau}^{\prime(e_{\gp})} \subset H^1_\crys(A_0'/R_1)_\tau\end{array}$$
such that the following hold for $j = 1,\ldots,e_{\gp}$ and $\theta = \theta_{\gp,i,j}$:
\begin{itemize}
\item $\widetilde{F}_\tau^{(j)}$ and $\widetilde{F}_\tau^{\prime(j)}$ are free (of rank $j$) over $R_1$;
\item $u\widetilde{F}_\tau^{(j)} \subset \widetilde{F}_\tau^{(j-1)}$, $u\widetilde{F}_\tau^{\prime(j)} \subset \widetilde{F}_\tau^{\prime(j-1)}$
    and $\psi_{0,\crys}^*(\widetilde{F}_\tau^{\prime(j)}) \subset \widetilde{F}_\tau^{(j)}$;
\item $\widetilde{F}_\tau^{(j-1)} \otimes_{R_1} R_\theta$ corresponds to $F_{\tau,0}^{(j-1)} \otimes_{S_0} R_\theta$ and
    $\widetilde{F}_\tau^{\prime(j-1)} \otimes_{R_1} R_\theta$ corresponds to $F_{\tau,0}^{\prime(j-1)} \otimes_{S_0} R_\theta$ under the canonical
     isomorphisms  $H^1_\dr(A_0/S_0)\otimes_{S_0} R_\theta \cong H^1_\crys(A_0/R_1)\otimes_{R_1} R_\theta$ and
     $H^1_\dr(A'_0/S_0)\otimes_{S_0} R_\theta \cong H^1_\crys(A'_0/R_1)\otimes_{R_1} R_\theta$, which therefore induce 
     injective $R_\theta$-linear homomorphisms
     $$\begin{array}{ll} & \beta_\theta: P_{\theta,0} \otimes_{S_0} R_\theta \longrightarrow
        (H^1_\crys(A_0/R_1)_\tau/\widetilde{F}_\tau^{(j-1)}) \otimes_{R_1} R_\theta\\ \mbox{and}\quad &
        \beta_\theta' : P_{\theta,0}' \otimes_{S_0} R_\theta \longrightarrow
         (H^1_\crys(A'_0/R_1)_\tau/\widetilde{F}_\tau^{\prime(j-1)}) \otimes_{R_1} R_\theta; \end{array}$$
\item $(\widetilde{F}_\tau^{(j)}/\widetilde{F}_\tau^{(j-1)}) \otimes_{R_1} R_\theta$ is generated over $R_\theta$ by 
  $\beta_\theta(e_{\theta,0} \otimes 1 - f_{\theta,0} \otimes X_\theta)$ and
  $(\widetilde{F}_\tau^{\prime(j)}/\widetilde{F}_\tau^{\prime(j-1)}) \otimes_{R_1} R_\theta$ is generated over $R_\theta$ 
  by  $\beta'_\theta(e'_{\theta,0} \otimes 1 - f'_{\theta,0} \otimes X'_\theta)$.
\end{itemize}

Suppose then that $1 \le j \le e_{\gp}$ and that 
$$0 = \widetilde{F}_\tau^{(0)} \subset \widetilde{F}_\tau^{(1)} \subset \cdots \subset \widetilde{F}_{\tau}^{(j-1)} \quad \mbox{and}\quad 
0 = \widetilde{F}_\tau^{\prime(0)} \subset \widetilde{F}_\tau^{\prime(1)} \subset \cdots \subset \widetilde{F}_{\tau}^{\prime(j-1)}$$
have been constructed as above.  We let $\widetilde{G}_\tau^{(j)} = u^{-1}\widetilde{F}_{\tau}^{(j-1)}$, i.e., the preimage of $\widetilde{F}_{\tau}^{(j-1)}$
in $H^1_\crys(A_0/R_1)_\tau$ under $u$, and $\widetilde{P}_\theta = \widetilde{G}_\tau^{(j)}/\widetilde{F}_{\tau}^{(j-1)}$.
Since $H^1_\crys(A_0/R_1)_\tau$ is free of rank two over $R_1[u]/(u^{e_\gp})$, we see that 
$\widetilde{F}_{\tau}^{(j-1)} \subset u H^1_\crys(A_0/R_1)_\tau$, so that $\widetilde{G}_\tau^{(j)}$ and $\widetilde{P}_\theta$
are free (of ranks $j+1$ and $2$, respectively) over $R_1$.  Similarly we let $\widetilde{G}_\tau^{\prime(j)} = u^{-1}\widetilde{F}_\tau^{\prime(j-1)}$
and $\widetilde{P}_\theta' = \widetilde{G}_\tau^{\prime(j)}/\widetilde{F}_\tau^{\prime(j-1)}$.
The conditions on $\widetilde{F}_\tau^{(j-1)}$ and $\widetilde{F}_\tau^{\prime(j-1)}$ imply that the third bullet holds, 
that $\psi_{0,\crys}^*$ induces an $R_1$-linear homomorphism $\widetilde{\psi}^*_\theta: \widetilde{P}_\theta' \to \widetilde{P}_\theta$,
and that $\beta_\theta$ and $\beta_\theta'$ define isomorphisms
$P_{\theta,0} \otimes_{S_0} R_\theta \stackrel{\sim}{\longrightarrow} \widetilde{P}_\theta \otimes_{R_1} R_\theta$ and
$P'_{\theta,0} \otimes_{S_0} R_\theta \stackrel{\sim}{\longrightarrow} \widetilde{P}'_\theta \otimes_{R_1} R_\theta$ 
compatible with $\psi^*_{\theta,0}$ and $\widetilde{\psi}^*_\theta$.
Furthermore the construction of $\widetilde{F}_\tau^{(j)}$ and $\widetilde{F}_\tau^{\prime(j)}$ satisfying the remaining
conditions is equivalent to that of invertible $R_1$-submodules $\widetilde{L}_\theta \subset \widetilde{P}_\theta$
and $\widetilde{L}'_\theta \subset \widetilde{P}'_\theta$ such that
\begin{itemize}
\item $\widetilde{\psi}_\theta^*(\widetilde{L}_\theta') \subset \widetilde{L}_\theta$;
\item $\widetilde{L}_\theta \otimes_{R_1} R_\theta = \beta_\theta(e_{\theta,0} \otimes 1 - f_{\theta,0} \otimes X_\theta)R_\theta$;
\item $\widetilde{L}'_\theta \otimes_{R_1} R_\theta = \beta'_\theta(e'_{\theta,0} \otimes 1 - f'_{\theta,0} \otimes X'_\theta)R_\theta$.
\end{itemize}

Exactly as for de Rham cohomology (see \cite[\S3.1]{theta}), we obtain a perfect pairing
$$H^1_\crys(A_0/R_1)_\tau \stackrel{\sim}{\longrightarrow} \Hom_{R_1[u]/u^{e_{\gp}}}(H^1_\crys(A_0/R_1)_\tau, R_1[u]/(u^{e_\gp}))$$
from the homomorphisms induced by the quasi-polarization on $A_0$ and Poincar\'e duality on crystalline cohomology.
Furthermore the pairing is compatible with the corresponding one on $H^1_\dr(B/R_1)_\tau$, where $B$
and the resulting quasi-polarization are associated to an arbitrarily chosen extension of 
$0  \subset \widetilde{F}_\tau^{(1)} \subset \cdots \subset \widetilde{F}_{\tau}^{(j-1)}$
to  a lift $\widetilde{F}^\bullet$ of the Pappas--Rapoport filtrations to $H^1_\crys(A_0/R_1)$.
We may therefore apply \cite[Lemma~3.1.1]{theta} to conclude that $(\widetilde{F}_\tau^{(j-1)})^\perp = u^{j+1-e_{\gp}} \widetilde{F}_\tau^{(j-1)} = u^{j-e_{\gp}}\widetilde{G}_\tau^{(j)}$,
and hence obtain a perfect alternating $R_1$-valued pairing on $\widetilde{P}_\theta$,
compatible via $\beta_\theta$ with $\langle \cdot, \cdot \rangle_{\theta,0}$.  We similarly obtain a perfect alternating
$R_1$-valued pairing on $\widetilde{P}_\theta'$ compatible via $\beta_\theta'$
with $\langle \cdot, \cdot \rangle'_{\theta,0}$.  Furthermore the compatibility of $\psi_0$ with the quasi-polarizations again implies
that $\det \widetilde{\psi}^*_\theta = \theta(\varpi_\gP)\, (=0)$ with respect to the pairings.  

The same argument as in the construction
of the bases $\CB_\theta$ and $\CB_\theta'$ then yields lifts of $\beta_\theta(e_{\theta,0} \otimes 1,f_{\theta,0} \otimes 1)$ and
$\beta_\theta'(e'_{\theta,0} \otimes 1,f'_{\theta,0} \otimes 1)$
to bases $(\widetilde{e}_\theta,\widetilde{f}_\theta)$ for $\widetilde{P}_\theta$
and $(\widetilde{e}'_\theta,\widetilde{f}'_\theta)$ for $\widetilde{P}'_\theta$ 
with respect to which the matrix of $\widetilde{\psi}_\theta^*$
is as in (\ref{eqn:matrices}).  Defining  $\widetilde{L}_\theta = (\widetilde{e}_\theta -  X_\theta \widetilde{f}_\theta)R_1$
and  $\widetilde{L}'_\theta = (\widetilde{e}'_\theta -  X'_\theta \widetilde{f}'_\theta)R_1$ then completes the construction
of $\widetilde{F}_\tau^{(j)}$ and $\widetilde{F}_\tau^{\prime(j)}$ satisfying the desired properties.

We may now apply the Grothendieck--Messing Theorem (as described at the end of \S\ref{sec:iwahori})
to obtain a triple $(\widetilde{\underline{A}},\widetilde{\underline{A}}',\widetilde{\psi})$ over $R_1$ lifting 
$(\underline{A}_0,\underline{A}_0',\psi_0)$.  Furthermore the properties of the filtrations $\widetilde{F}^\bullet$
and $\widetilde{F}^{\prime\bullet}$
ensure, by induction on $j$ for each $\tau$, that the resulting morphism $S \to R_1 \to R_\theta$ factors
through $S_\theta$ and induces isomorphisms
$$\begin{array}{ccccccc}
P_\theta \otimes_S R_\theta & \stackrel{\sim}{\longrightarrow}& \widetilde{P}_\theta \otimes_{R_1} R_\theta&&
P_\theta' \otimes_S R_\theta & \stackrel{\sim}{\longrightarrow}& \widetilde{P}'_\theta \otimes_{R_1} R_\theta\\
\cup&&\cup& \mbox{and} & \cup &&\cup \\
L_\theta \otimes_S R_\theta & \stackrel{\sim}{\longrightarrow}& \widetilde{L}_\theta \otimes_{R_1} R_\theta&&
L_\theta' \otimes_S R_\theta & \stackrel{\sim}{\longrightarrow}& \widetilde{L}'_\theta \otimes_{R_1} R_\theta\end{array}$$
sending $(e_\theta \otimes 1,f_\theta \otimes 1)$ to $(\widetilde{e}_\theta \otimes 1,\widetilde{f}_\theta \otimes 1)$ and
$(e'_\theta \otimes 1,f'_\theta \otimes 1)$ to $(\widetilde{e}'_\theta \otimes 1,\widetilde{f}'_\theta \otimes 1)$,
so that $s_\theta \mapsto X_\theta$ and $s'_\theta \mapsto X'_\theta$ ($\!\bmod\, J_\tau^{(j-1)}$).
We therefore conclude that the morphism $S \to R_1$ is surjective, and hence that $\rho_1:R_1 \to S_1$
is an isomorphism.

Suppose then that $n \ge 1$ and that $\rho_n:R_n \to S_n$ is an isomorphism.  The surjectivity of $\rho_{n+1}$ is already
immediate from that of $\rho_n$, so it suffices to prove that $\lg(R_{n+1}) \le \lg(S_{n+1})$, which will follow from the
existence of a lift of $\rho_n^{-1}$ to morphism $S \to R_{n+1}$, which necessarily factors through a surjective morphism
from $S_{n+1}$.  

As in the case $n=0$, we proceed by defining suitable lifts of the Pappas--Rapoport filtrations $F^\bullet_{n}$
and $F^{\prime\bullet}_{n}$ to the crystalline cohomology of $A_n$ and $A_n'$ evaluated on the thickening
$R_{n+1} \to R_n \stackrel{\sim}{\to} S_n$ (with trivial divided power structure), 
but the argument is now simpler since the isomorphism on tangent spaces has already been established.
More precisely, we inductively define chains of $R_{n+1}[u]/(E_\tau)$-submodules
$$\begin{array}{ll} & 0 = \widetilde{F}_\tau^{(0)} \subset \widetilde{F}_\tau^{(1)} \subset \cdots \subset \widetilde{F}_{\tau}^{(e_{\gp})} \subset H^1_\crys(A_n/R_{n+1})_\tau \\
 \mbox{and}\quad &
0 = \widetilde{F}_\tau^{\prime(0)} \subset \widetilde{F}_\tau^{\prime(1)} \subset \cdots \subset \widetilde{F}_{\tau}^{\prime(e_{\gp})} \subset H^1_\crys(A_n'/R_{n+1})_\tau\end{array}$$
such that the following hold for $j = 1,\ldots,e_{\gp}$ and $\theta = \theta_{\gp,i,j}$:
\begin{itemize}
\item $\widetilde{F}_\tau^{(j)}$ and $\widetilde{F}_\tau^{\prime(j)}$ are free (of rank $j$) over $R_{n+1}$;
\item $(u - \theta(\varpi_\gP))\widetilde{F}_\tau^{(j)} \subset \widetilde{F}_\tau^{(j-1)}$, $(u- \theta(\varpi_\gP))\widetilde{F}_\tau^{\prime(j)} \subset \widetilde{F}_\tau^{\prime(j-1)}$
    and $\psi_{n,\crys}^*(\widetilde{F}_\tau^{\prime(j)}) \subset \widetilde{F}_\tau^{(j)}$;
\item $\widetilde{F}_\tau^{(j)} \otimes_{R_{n+1}} S_n$ corresponds to $F_{\tau,n}^{(j)}$ and
    $\widetilde{F}_\tau^{\prime(j)} \otimes_{R_{n+1}} S_n$ to $F_{\tau,n}^{\prime(j)}$ under the canonical
     isomorphisms  $H^1_\dr(A_n/S_n)  \cong H^1_\crys(A_n/R_{n+1})\otimes_{R_{n+1}} S_n$ and
     $H^1_\dr(A'_n/S_n) \cong H^1_\crys(A'_n/R_{n+1})\otimes_{R_{n+1}} S_n$.
\end{itemize}

Suppose then that $1 \le j \le e_{\gp}$ and that 
$$0 = \widetilde{F}_\tau^{(0)} \subset \widetilde{F}_\tau^{(1)} \subset \cdots \subset \widetilde{F}_{\tau}^{(j-1)} \quad \mbox{and}\quad 
0 = \widetilde{F}_\tau^{\prime(0)} \subset \widetilde{F}_\tau^{\prime(1)} \subset \cdots \subset \widetilde{F}_{\tau}^{\prime(j-1)}$$
have been constructed as above.  Let $\widetilde{G}_\tau^{(j)} = (u-\theta(\varpi_\gp))^{-1}\widetilde{F}_{\tau}^{(j-1)}$,
$\widetilde{P}_\theta = \widetilde{G}_\tau^{(j)}/\widetilde{F}_{\tau}^{(j-1)}$ and similarly define $\widetilde{P}'_\theta$, 
so that $\widetilde{P}_\theta$ and $\widetilde{P}_\theta'$ are free of rank two over $R_{n+1}$.  Furthermore we have
canonical isomorphisms 
$$\beta_\theta: P_{\theta,n}  \stackrel{\sim}{\longrightarrow} \widetilde{P}_\theta \otimes_{R_{n+1}} S_n
\quad\mbox{and}\quad
\beta'_\theta: P'_{\theta,n}  \stackrel{\sim}{\longrightarrow} \widetilde{P}'_\theta \otimes_{R_{n+1}} S_n,$$ 
and $\psi_{n,\crys}^*$ induces an $R_{n+1}$-linear homomorphism $\widetilde{\psi}^*_\theta: \widetilde{P}_\theta' \to \widetilde{P}_\theta$
compatible with $\psi_{\theta,n}^*$.
The construction of the required $\widetilde{F}_\tau^{(j)}$ and $\widetilde{F}_\tau^{\prime(j)}$ is then
equivalent to that of invertible $R_{n+1}$-submodules $\widetilde{L}_\theta \subset \widetilde{P}_\theta$
and $\widetilde{L}'_\theta \subset \widetilde{P}'_\theta$ such that
\begin{itemize}
\item $\widetilde{\psi}_\theta^*(\widetilde{L}_\theta') \subset \widetilde{L}_\theta$;
\item $\widetilde{L}_\theta \otimes_{R_{n+1}} S_n = \beta_\theta(e_{\theta,n}  - s_\theta f_{\theta,n} )S_n$;
\item $\widetilde{L}'_\theta \otimes_{R_{n+1}} S_n = \beta'_\theta(e'_{\theta,n}  - s'_\theta f'_{\theta,n} )S_n$.
\end{itemize}

The same argument as for $n=0$ (with $t_j^{-1}$ generalizing $u^{j-e_{\gp}}$) yields
perfect alternating $R_{n+1}$-valued pairings on $\widetilde{P}_\theta$ on $\widetilde{P}'_\theta$
compatible via $\beta_\theta$ and $\beta_\theta'$ with $\langle \cdot,\cdot \rangle_{\theta,n}$ and
$\langle \cdot, \cdot \rangle'_{\theta,n}$, and with respect to which
$\det \widetilde{\psi}^*_\theta = \theta(\varpi_\gP)$.  We may thus 
lift $\beta_\theta(e_{\theta,n},f_{\theta,n})$ and
$\beta_\theta'(e'_{\theta,n},f'_{\theta,n})$
to bases $(\widetilde{e}_\theta,\widetilde{f}_\theta)$ for $\widetilde{P}_\theta$
and $(\widetilde{e}'_\theta,\widetilde{f}'_\theta)$ for $\widetilde{P}'_\theta$ 
with respect to which the matrix of $\widetilde{\psi}_\theta^*$
is as in (\ref{eqn:matrices}).  Defining
$\widetilde{L}_\theta = (\widetilde{e}_\theta -  X_\theta \widetilde{f}_\theta)R_{n+1}$
and  $\widetilde{L}'_\theta = (\widetilde{e}'_\theta -  X'_\theta \widetilde{f}'_\theta)R_{n+1}$.
thus yields the desired filtrations $\widetilde{F}^\bullet$ and $\widetilde{F}^{\prime\bullet}$,
and hence a lift of $(A_n,A_n',\psi_n)$ to a triple over $R_{n+1}$.

This completes the proof that $\rho$ is isomorphism, and hence that $\CO^{\wedge}_{\widetilde{Y}_{U_0(\gP)},y}$
is isomorphic to
$$W\otimes_{W(k)}  \CO [[X_\theta,X_\theta']]_{\theta \in \Sigma}/ (g_\theta)_{\theta\in \Sigma}.$$
As this is a reduced complete intersection, flat of relative dimension $d$ over $\CO$, and $\widetilde{Y}_{U_0(\gP),K}$ is
smooth of dimension $d$ over $K$, it follows that $\widetilde{Y}_{U_0(\gP)}$ is reduced and syntomic of
dimension $d$ over $\CO$, and that the same holds for $Y_{U_0(\gP)}$ since
$\widetilde{Y}_{U_0(\gP)} \to Y_{U_0(\gP)}$ is an \'etale cover.

\section{The Kodaira--Spencer isomorphism}  \label{sec:KS}

\subsection{The layers of a thickening}  \label{sec:filtration}

Recall that the Kodaira--Spencer isomorphism describes the dualizing sheaf of the smooth scheme $Y_U$ over $\CO$ 
in terms of automorphic line bundles; more precisely we have a canonical\footnote{A precise sense in which it is canonical
is its compatibility with the natural actions of $\GL_2(\A_{F,\f}^{(p)})$.  The isomorphism however depends up to sign on a
choice of ordering of $\Sigma$, which we fix in retrospect and perpetuity.} isomorphism
$$\delta^{-1} \otimes \omega^{\otimes 2}  \stackrel{\sim}{\longrightarrow}  \CK_{Y_U/\CO},$$
where $\omega$ and $\delta$ are the line bundles on $Y_U$ defined in \S\ref{sec:bundles} and 
$\CK_{Y_U/\CO} =  \wedge^d_{\CO_{Y_U}} \Omega^1_{Y_U/\CO}$ (see \cite[\S3.3]{theta}).
We will prove an analogous result for $Y_{U_0(\gP)}$, which by the results of the preceding section is
Gorenstein over $\CO$, hence has an invertible dualizing sheaf $\CK_{Y_{U_0(\gP)}/\CO}$.

First recall that the Kodaira--Spencer isomorphism for $Y_U$ was established in \cite{RX} by constructing
a suitable filtration on the sheaf of relative differentials and relating its graded pieces to automorphic line
bundles.  In \cite[\S3.3]{theta}, we gave a version of the construction which can be viewed 
as peeling off layers of the first infinitesimal thickening of the diagonal in $\widetilde{Y}_U \times_{\CO} \widetilde{Y}_U$.
We use the same approach here to describe the layers of the first infinitesimal thickening of
$\widetilde{Y}_{U_0(\gP)}$ in $\widetilde{Y}_U \times_{\CO} \widetilde{Y}_U$ along
$$\widetilde{h}: = (\widetilde{\pi}_1,\widetilde{\pi}_2): \widetilde{Y}_{U_0(\gP)} \to \widetilde{Y}_U \times_{\CO} \widetilde{Y}_U,$$
which we recall is a closed immersion for sufficiently small $U$ by Proposition~\ref{prop:immersion}.
Indeed the Kodaira--Spencer filtration for $Y_U$ may be viewed as the particular case where $\gP = \CO_F$, 
for which the proof in \cite[\S3.3]{theta} requires only minor modifications to establish the desired generalization.

We assume throughout this section that $U$ is sufficiently small that Proposition~\ref{prop:immersion} holds.
To simplify notation and render it more consistent with that of \cite[\S3.3]{theta}, we will write 
$S$ for $\widetilde{Y}_U$, $X$ for $S \times_{\CO} S$, $Y$ for $\widetilde{Y}_{U_0(\gP)}$,
$\CJ$ for the sheaf of ideals of $\CO_X$ defining the image of $\widetilde{h}:Y \hookrightarrow X$,
$\iota:Y \hookrightarrow Z = \SPEC(\CO_X/\CJ^2)$ for the first infinitesimal thickening of $Y$ in $X$,
and $\CI = \CJ/\CJ^2$ for the sheaf of ideals of $\CO_Z$ defining the image of $\iota$.
Since $X$ is smooth over $\CO$ of dimension $2d$
and $Y$ is syntomic over $\CO$ of dimension $d$ (all dimensions being relative), 
it follows that $\widetilde{h}$ is a regular immersion (see for example  \cite[\S0638]{stacks}) and its conormal bundle 
$${\widetilde{\calC}}:= \widetilde{h}^*\CJ = \iota^*\CI$$
is locally free of rank $d$ on $Y = \widetilde{Y}_{U_0(\gP)}$.

\begin{theorem} \label{thm:KSfil}  There exists  a decomposition 
$\widetilde{\calC} = \bigoplus_{\tau \in \Sigma_0} \widetilde{\calC}_\tau$,
together with an increasing filtration
$$ 0 = \Fil^0 (\widetilde{\calC}_\tau) \subset \Fil^1 (\widetilde{\calC}_\tau) \subset \cdots
\subset  \Fil^{e_\gp - 1}( \widetilde{\calC}_\tau )
  \subset  \Fil^{e_\gp}( \widetilde{\calC}_\tau) =  \widetilde{\calC}_\tau $$
for each $\tau = \tau_{\gp,i}$,  such that 
$\gr^j(\widetilde{\calC}_\tau) \cong \widetilde{\pi}_1^*\CM_\theta^{-1}  \otimes_{\CO_Y} \widetilde{\pi}_2^*\CL_\theta$
for each $j=1,\ldots,e_\gp$, where $\theta =  \theta_{\gp,i,j}$.
\end{theorem}

\begpf  
Let $A$ denote the universal abelian scheme over $S$, let $\psi:A_1 \to A_2$ denote the
universal isogeny over $Y$, and let $B_i = q_i^*A$ for $i=1,2$, where
$q_i$ denotes the restriction to $Z$ of the $i^{\mathrm{th}}$ projection $S \times_{\CO} S \to S$.  
Consider the $\CO_F\otimes \CO_Z$-linear morphism
\begin{multline}\label{eqn:psicrys}
q_2^*\CH_\dr(A/S) \cong \CH^1_\dr(B_2/Z) \cong  (R^1 s_{2,\crys,*} \CO_{A_2,\crys})_Z \\
\stackrel{-\psi_\crys^*}{\longrightarrow}  (R^1 s_{1,\crys,*} \CO_{{A}_1,\crys})_Z
 \cong \CH^1_\dr(B_1/Z) \cong q_1^*\CH^1_\dr(A/S)\end{multline}
extending\footnote{The choice of sign is made for consistency with the conventions of the
classical Kodaira--Spencer isomorphism, i.e., the case $\gP = \CO_F$, where $\widetilde{h}$ 
is the diagonal embedding, $\psi$ is the identity, and we would ordinarily consider
the resulting morphism $q_1^*\CH^1_\dr(A/S) \to q_2^*\CH^1_\dr(A/S)$.  In general the morphism induced by
$\psi$ goes in the opposite direction, so we introduce a minus sign here to obviate factors
of $(-1)^d$ in later compatibility formulas.}
$\iota^*\CH^1_\dr(B_2/Z) \cong \CH^1_\dr(A_2/Y)
\stackrel{-\psi^*}{\longrightarrow} \CH^1_\dr(A_1/Y)  \cong \iota^*\CH^1_\dr(B_1/Z)$, where the isomorphisms flanking $-\psi_\crys^*$ are the ones in (\ref{eqn:crysDR}), for $\iota:Y \hookrightarrow Z$ with $B=B_i$.
Fixing for the moment $\tau = \tau_{\gp,i} \in \Sigma_0$ and $\theta = \theta_{\gp,i,1}$, it
follows from the definition of $\CP_{\theta} = \CG_\tau^{(1)}$
that (\ref{eqn:psicrys})  restricts to an $\CO_Z$-linear morphism
$q_2^*\CP_{\theta}  \to q_1^*\CP_{\theta}$.
Furthermore the composite 
$$q_2^*\CL_{\theta}  \hookrightarrow q_2^*\CP_{\theta} \to q_1^*\CP_{\theta}
      \onto q_1^*\CM_{\theta}$$
has trivial pull-back to $Y$, so it factors through a morphism
$$\iota_*\widetilde{\pi}_2^* \CL_{\theta}  = q_2^*\CL_{\theta} \otimes_{\CO_{Z}} (\CO_{Z}/\CI)
  \longrightarrow q_1^*\CM_{\theta} \otimes_{\CO_{Z}} \CI = \iota_*\widetilde{\pi}_1^* \CM_{\theta}  \otimes_{\CO_{Z}} \CI,$$
and hence induces a morphism
$$\Xi_{\theta}:  \iota_* (\widetilde{\pi}_1^*\CM^{-1}_{\theta}\otimes_{\CO_Y} \widetilde{\pi}_2^*\CL_{\theta})  \longrightarrow  \CI.$$

We then define the sheaf of ideals $\CI_{\tau}^{(1)}$ on $Z$ to be the image of $\Xi_\theta$,
and we let $Z_{\tau}^{(1)}$ denote the subscheme of $Z$ defined by $\CI_{\tau}^{(1)}$, and
$q^{(1)}_{\tau,i}$  (for $i=1,2$) the restrictions of the projection maps to $Z_{\tau}^{(1)} \to S$.
By construction the pull-back of $\Xi_\theta$ to $Z_{\tau}^{(1)}$ is trivial, and hence so is that of the
morphism $q_2^*\CL_{\theta} \to q_1^*\CM_{\theta}$, which implies that the pull-back of (\ref{eqn:psicrys})
maps $q_{\tau,2}^{(1)*}\CL_{\theta} = q_{\tau,2}^{(1)*}\CF_\tau^{(1)}$ to 
$q_{\tau,1}^{(1)*}\CL_{\theta} = q_{\tau,1}^{(1)*}\CF_\tau^{(1)}$.  It therefore follows from the definition of
$\CG_\tau^{(2)}$ and the $\CO_F\otimes \CO_Z$-linearity of (\ref{eqn:psicrys}) that its pull-back to $Z_\tau^{(1)}$
restricts to an $\CO_{Z_\tau^{(1)}}$-linear morphism $q_{\tau,2}^{(1)*}\CG_\tau^{(2)} \to q_{\tau,1}^{(1)*}\CG_\tau^{(2)}$
(if $e_\gp > 1$), and hence to a morphism
 $$q_{\tau,2}^{(1)*}\CP_{\theta'}  {\longrightarrow}  q_{\tau,1}^{(1)*}\CP_{\theta'},$$
 where $\theta' = \theta_{\gp,i,2}$.
The same argument as above now yields a morphism
$$\Xi_{\theta'}:  \iota_* (\widetilde{\pi}_1^*\CM^{-1}_{\theta'}\otimes_{\CO_Y} \widetilde{\pi}_2^*\CL_{\theta'})  \longrightarrow  \CI/\CI_\tau^{(1)}$$
whose image is $\CI_\tau^{(2)}/\CI_\tau^{(1)}$ for some sheaf of ideals $\CI_{\tau}^{(2)} \supset \CI_{\tau}^{(1)}$ on $Z$.

Iterating the above construction thus yields, for each $\tau  = \tau_{\gp,i} \in \Sigma_{0}$, a chain of sheaves of ideals
$$0 = \CI_{\tau}^{(0)} \subset \CI_{\tau}^{(1)} \subset \cdots \subset \CI_{\tau}^{(e_\gp)}$$
on $Z$ such that for $j=1,\ldots,e_\gp$ and $\theta = \theta_{\gp,i,j}$, the morphism (\ref{eqn:psicrys}) induces 
\begin{itemize}
\item $\CO_F \otimes \CO_{Z_\tau^{(j)}}$-linear morphisms $q_{\tau,2}^{(j)*}\CF_\tau^{(j)} \longrightarrow q_{\tau,1}^{(j)*} \CF_\tau^{(j)}$
\item and surjections $\Xi_\theta: \iota_* (\widetilde{\pi}_1^*\CM_\theta^{-1} \otimes_{\CO_Y} \CL_{\theta})  \onto  \CI_{\tau}^{(j)}/\CI_{\tau}^{(j-1)}$,
\end{itemize}
where $Z_{\tau}^{(j)}$ is the closed subscheme of $Z$ defined by $\CI_\tau^{(j)}$, and
$q_{\tau,1}^{(j)}$, $q_{\tau,2}^{(j)}$ are the projections $Z_{\tau}^{(j)} \to S$.

Furthermore we claim that the map 
$$\bigoplus_{\gp \in S_p}  \bigoplus_{\tau \in \Sigma_{\gp,0}} \CI_{\tau}^{(e_\gp)}  \to \CI$$
is surjective.  Indeed let $Z'$ denote the closed
subscheme of $Z$ defined by the image, so $Z'$ is the scheme-theoretic intersection of the
$Z_{\tau}^{(e_\gp)}$.  For $i=1,2$, let
$q'_i$ denote the projection map $Z' \to S$ and $s_i':B_i' \to Z'$ the pull-back
of $s:A\to S$.   By construction (\ref{eqn:psicrys}) pulls back to a morphism 
$q_2^{\prime*}\CH^1_{\dr}(A/S) \to q_1^{\prime*}\CH^1_\dr(A/S)$ under which 
$q_2^{\prime*}\CF_\tau^{(j)}$ maps to $q_1^{\prime*}\CF_\tau^{(j)}$ for all $\tau$ and $j$.
In particular $s_{2,*}'\Omega^1_{B_2'/T} = q_2^{\prime*}(s_*\Omega^1_{A/S})$ maps to
$q_1^{\prime*}(s_*\Omega^1_{A/S}) = s'_{1,*} \Omega^1_{B_1'/T}$, 
which the Grothendieck--Messing Theorem therefore
implies is induced by an isogeny $\widetilde{\psi}: B'_1 \to B_2'$ of abelian 
schemes lifting $\psi$, which is necessarily compatible with $\CO_F$-actions, quasi-polarizations
and level structures on $\underline{B}_i' : = q_i^{\prime*}\underline{A}$.  Since
$\widetilde{\psi}$ also respects the Pappas--Rapoport filtrations $q_i^{\prime*}\CF^\bullet$, it follows
that the triple $(\underline{B}_1',\underline{B}_2',\widetilde{\psi})$ corresponds
to a morphism $r:Z' \to Y$ such that $\widetilde{h}\circ r$ is the closed immersion $Z' \hookrightarrow X$.
Since $\widetilde{h}$ is also a closed immersion, it follows that $r$ is an isomorphism,
yielding the desired surjectivity.

Now defining $\widetilde{\calC}_\tau  = \iota^*\CI_\tau^{(e_\gp)}$ and $\Fil^j(\widetilde{\calC}_\tau)= \iota^*\CI_{\tau}^{(j)}$
for each $\tau = \tau_{\gp,i}$ and $j=1,\ldots,e_\gp$, we obtain surjective morphisms
$$ \bigoplus_{\tau \in \Sigma_0} \widetilde{\calC}_\tau \onto \widetilde{\calC}, \qquad \mbox{and}\qquad
\widetilde{\pi}_1^*\CM_\theta^{-1} \otimes_{\CO_Y} \widetilde{\pi}_2^*\CL_{\theta} \onto  \gr^j (\widetilde{\calC}_\tau )
 \quad \mbox{for each $\theta = \theta_{\gp,i,j}$.}$$
Since the $\widetilde{\pi}_1^*\CM_\theta^{-1} \otimes_{\CO_Y} \widetilde{\pi}_2^*\CL_{\theta}$ are line bundles and
$\widetilde{\calC}$ is a vector bundle of rank $d$, it follows that all the maps are isomorphisms.
\epf

Rewriting the line bundles in the statement of the theorem as 
$$\widetilde{\cD}_\theta = \widetilde{\pi}_1^*\CM_\theta^{-1} \otimes \widetilde{\pi}_2^*\CL_{\theta}
  \cong \widetilde{\pi}_1^*\CN_\theta^{-1} \otimes \widetilde{\pi}_1^*\CL_{\theta} \otimes
   \widetilde{\pi}_2^*\CL_{\theta}$$
(all tensor products being over $\CO_Y$), we consequently obtain an isomorphism
$$ \bigwedge^d \widetilde{\calC} \,\,
\cong\,\, \bigotimes_{\theta \in \Sigma} \widetilde{\cD}_\theta \,\,\cong\,\,
 \widetilde{\pi}_1^* \widetilde{\delta}^{-1} \otimes \widetilde{\pi}_1^*\widetilde{\omega} \otimes \widetilde{\pi}_2^*\widetilde{\omega}.$$
 
Recall that the action of $\CO_{F,(p),+}^\times$ on $Y$ is defined by its action on the quasi-polarizations
of both abelian schemes in the triple $(\underline{A}_1,\underline{A}_2,\psi)$.  In particular the action
factors through $\CO_{F,(p),+}^\times/(U \cap \CO_F^\times)^2$ and is compatible
under $\widetilde{h}$ with the diagonal action on the product $X = S \times_{\CO} S$.
The conormal bundle $\widetilde{\calC}$ of $\widetilde{h}$ is thus equipped with an action of
$\CO^\times_{F,(p),+}$ over its action on $Y$, coinciding with the action on $\widetilde{\calC}$ obtained
from its identification with the pull-back of the conormal bundle of the closed immersion 
$h: Y_{U_0(\gP)} \hookrightarrow Y_U \otimes_{\CO} Y_U$, which we denote by $\calC$.

Recall also that we have a natural action of $\CO^\times_{F,(p),+}$ on the line bundles $\CL_\theta$
and $\CM_\theta$ over its action on $S$, under which $\nu$ acts as $\theta(\mu)$ if
$\nu = \mu^2$ for $\mu \in U \cap \CO_F^\times$.  It follows that the action factors through
$\CO_{F,(p),+}^\times/(U \cap \CO_F^\times)^2$ and hence defines descent data on the line
bundles $\widetilde{\cD}_\theta$;  we let $\cD_\theta$ denote the resulting line bundle
on $Y_{U_0(\gP)}$. It is straightforward to check that the morphisms constructed in the proof
of Theorem~\ref{thm:KSfil} are compatible with the actions of $\CO_{F,(p),+}^\times$ on
the bundles $\widetilde{\cD}_\theta$ and $\widetilde{\calC}$, so that the decomposition
and the filtrations on resulting components of $\widetilde{\calC}$ descend to ones 
over $Y_{U_0(\gP)}$, as do the isomorphisms between the graded pieces and the line
bundles $\widetilde{\cD}_\theta$.

Similarly the constructions in the proof of Theorem~\ref{thm:KSfil} are compatible with
the natural action of $\GL_2(\A_{F,\f}^{(p)})$.  More precisely suppose that $g \in \GL_2(\A_{F,\f}^{(p)})$
is such that $g^{-1}Ug \subset U'$, where $U$, $U'$ are sufficiently small open compact subgroups of 
$\GL_2(\A_{F,\f})$ containing $\GL_2(\CO_{F,p})$ (and such that Proposition~\ref{prop:immersion} also
holds for $U'$), and let $\widetilde{\rho}_g$ denote the morphisms $Y \to Y' : = \widetilde{Y}_{U'_0(\gP)}$
and $S \to S' = \widetilde{Y}_{U'}$ defined by the action on level structures as in \S\ref{subsection:PRmodel}.
Systematically annotating with $'$ the corresponding objects for $U'$, so for example writing
$\widetilde{\calC}'$ for the conormal bundle of $\widetilde{h}':Y' \hookrightarrow X'$,
Theorem~\ref{thm:KSfil} yields a decomposition $\widetilde{\calC}' = \oplus_{\tau \in \Sigma_0} \widetilde{\calC}'_\tau$,
filtrations $\Fil^j(\widetilde{\calC}_\tau')$ and isomorphisms $\gr^j(\widetilde{\calC}_\tau) \cong \widetilde{\cD}'_\theta$,
where as usual $\theta = \theta_{\gp,i,j}$ if $\tau = \tau_{\gp,i}$.    In addition to the isomorphisms
$\widetilde{\rho}_g^*\widetilde{\cD}'_\theta \stackrel{\sim}{\longrightarrow} \widetilde{\cD}$ defined as in \S\ref{sec:bundles},
we have the isomorphism $\widetilde{\rho}_g^*{\widetilde{\calC}'} \stackrel{\sim}{\longrightarrow} \widetilde{\calC}$
obtained from the Cartesian diagram
$$\xymatrix{ Y \ar[r]^-{\widetilde{h}} \ar[d]_-{\widetilde{\rho}_g} & X \ar[d]^-{(\widetilde{\rho}_g,\widetilde{\rho}_g)} \\
Y' \ar[r]_-{\widetilde{h}'} & X'.}$$
We then find that this isomorphism restricts to give
$\widetilde{\rho}_g^*\Fil^j(\widetilde{\calC}_\tau') \stackrel{\sim}{\longrightarrow} \Fil^j(\widetilde{\calC}_\tau)$
for all $\tau$ and $j$, and that the resulting diagram
$$\xymatrix{ \widetilde{\rho}_g^*\widetilde{\cD}_\theta' \ar[r]^-{\sim} \ar[d]_-{\wr} & \widetilde{\rho}_g^*\gr^j(\widetilde{\calC}'_\tau) \ar[d]^-{\wr} \\
\widetilde{\cD}_\theta \ar[r]^-{\sim} & \gr^j(\widetilde{\calC}_\tau)}$$
commutes, where as usual $\theta = \theta_{\gp,i,j}$ if $\tau = \tau_{\gp,i}$.
Combining this with the compatibility of these isomorphisms with the descent data defined by the
$\CO_{F,(p),+}^\times$-actions, we obtain analogous results for the vector bundles on $Y_{U'_0(\gP)}$ and $Y_{U_0(\gP)}$
with respect to the morphism $\rho_g: Y_{U_0(\gP)} \to Y_{U'_0(\gP)}$.  

Summing up, we have now proved the following:
\begin{theorem}\label{thm:KSfil2} Letting $\calC$ denote the conormal bundle of the closed immersion 
$h: Y_{U_0(\gP)} \hookrightarrow Y_U \otimes_{\CO} Y_U$,
there exists a decomposition 
${\calC} = \bigoplus_{\tau \in \Sigma_0} {\calC}_\tau$,
together with an increasing filtration
$$ 0 = \Fil^0 ({\calC}_\tau) \subset \Fil^1 ({\calC}_\tau) \subset \cdots
\subset  \Fil^{e_\gp - 1}( {\calC}_\tau )
  \subset  \Fil^{e_\gp}( {\calC}_\tau) =  {\calC}_\tau $$
for each $\tau = \tau_{\gp,i}$,  such that 
$\cD_\theta \cong \gr^j({\calC}_\tau)$
for each $j=1,\ldots,e_\gp$, $\theta =  \theta_{\gp,i,j}$.
Furthermore if $g^{-1}U g \subset U'$, then the canonical isomorphism $\rho_g^*\calC' \stackrel{\sim}{\longrightarrow} \calC$ restricts to
$\rho_g^*\Fil^j(\calC'_\tau) \stackrel{\sim}{\longrightarrow} \Fil^j(\calC_\tau)$
for all $\tau$ and $j$ as above, and the resulting diagram
$$\xymatrix{{\rho}_g^*{\cD}_\theta' \ar[r]^-{\sim} \ar[d]_-{\wr} 
  & {\rho}_g^*\gr^j({\calC}'_\tau) \ar[d]^-{\wr} \\
{\cD}_\theta \ar[r]^-{\sim} & \gr^j({\calC}_\tau)}$$
commutes.
\end{theorem}

\begin{corollary}  \label{cor:det}  If $U$ is sufficiently small that Proposition~\ref{prop:immersion} holds, then
there is an isomorphism
\begin{equation}\label{eqn:KSnew}
\pi_1^*\delta^{-1}  \otimes \pi_1^*\omega \otimes \pi_2^* \omega \stackrel{\sim}{\longrightarrow}  \bigwedge^d  \calC ,\end{equation}
where $\calC$ is the conormal bundle of the closed immersion $(\pi_1,\pi_2): Y_{U_0(\gP)} \hookrightarrow Y_U \times Y_U$.
Furthermore the isomorphisms are compatible with the action of $g \in \GL_2(\A_{F,\f}^{(p)})$, in the sense that if
$U'$ is also sufficiently small and $g^{-1}Ug \subset U'$, then the resulting diagram
$$\xymatrix{
\rho_g^*(\pi_1^{\prime*}\delta^{\prime-1}  \otimes \pi_1^{\prime*}\omega' \otimes \pi_2^{\prime*} \omega')
\ar[r]^-{\sim} \ar@{=}[d]
 & \rho_g^*\bigwedge^d\calC' 
  \ar@{=}[d]
  \\
\pi_1^*\rho_g^* \delta^{\prime-1}  \otimes \pi_1^* \rho_g^* \omega' \otimes \pi_2^* \rho_g^* \omega'
   \ar[d]_-{\wr}  
  & \bigwedge^d\rho_g^*\calC' 
   \ar[d]^-{\wr} 
 \\
\pi_1^*\delta^{-1}  \otimes \pi_1^*\omega \otimes \pi_2^* \omega 
\ar[r]^-{\sim} &
 \bigwedge^d\calC }$$
 commutes.
\end{corollary}

\begin{remark} \label{rmk:symmetric}
We comment briefly on the lack of symmetry between the two degeneracy maps $\pi_1$ and $\pi_2$ in
the statements.  Exchanging $\widetilde{\pi}_1$ and $\widetilde{\pi}_2$ in the definition of $\widetilde{\cD}_\theta$
gives the line bundle
 $$\widetilde{\CE}_\theta := \widetilde{\pi}_2^*\CM_\theta^{-1} \otimes \widetilde{\pi}_1^*\CL_{\theta}
  = \widetilde{\pi}_2^*\CN_\theta^{-1} \otimes \widetilde{\pi}_2^*\CL_{\theta} \otimes \widetilde{\pi}_1^* \CL_\theta$$
on $Y$, descending to a line bundle $\CE_\theta$ on $Y_{U_0(\gP)}$.  We may then use the isomorphism
  $\theta(\varpi_\gP)^{-1}(\wedge^2\psi_\theta^*):\widetilde{\pi}_2^*\CN_\theta \stackrel{\sim}{\longrightarrow} \widetilde{\pi}_1^*\CN_\theta$
to define an isomorphism $\widetilde{\cD}_\theta \stackrel{\sim}{\longrightarrow} \widetilde{\CE}_\theta$ compatible with descent data,
and hence inducing an isomorphism $\cD_\theta \cong \CE_\theta$.  The isomorphisms are furthermore compatible with the 
action of $\GL_2(\A_{F,\f}^{(p)})$, but are dependent (up to an element of $\CO^\times$) on the choices of $\varpi_\gp$.  There is
however a canonical isomorphism
$$\pi_2^*\delta\,\, \cong \,\,\Nm_{F/\Q}(\gP) \otimes_{\Z} \pi_1^*\delta\,\, \cong\,\, \pi_1^* \delta$$
making $\pi_1$ and $\pi_2$ interchangeable in the statement of Corollary~\ref{cor:det}.
\end{remark}

\subsection{Dualizing sheaves}  \label{sec:dualizing}

We now combine the Kodaira--Spencer isomorphism over $Y_U$ with the 
description of the conormal bundle obtained in the preceding section in order
to establish a Kodaira--Spencer isomorphism over $Y_{U_0(\gP)}$.

Recall that the Kodaira--Spencer filtration on $\Omega^1_{Y_U/\CO}$ yields an isomorphism
$$\wedge^d \Omega^1_{Y_U/\CO} \cong \delta^{-1} \otimes \omega^{\otimes 2}$$
(see \cite[\S2.8]{RX} and \cite[Thm.~3.3.1]{theta}), which is furthermore compatible
with the action of $\GL_2(\A_{F,\f}^{(p)})$ in the same sense as the isomorphism of
Corollary~\ref{cor:det}.  The identification
$\Omega^1_{(Y_U \times_{\CO} Y_U)/\CO} = p_1^*\Omega^1_{Y_U/\CO} \otimes_{\CO_{Y_U\times_\CO Y_U}}p_2^*\Omega^1_{Y_U/\CO}$,
where $p_1,p_2:Y_U \times_\CO Y_U \to Y_U$ are the projection maps, therefore yields an isomorphism
\begin{equation} \label{eqn:KSx2}
\wedge^{2d} \Omega^1_{(Y_U \times_\CO Y_U) /\CO} \cong p_1^*(\delta^{-1} \otimes \omega^{\otimes 2})
 \otimes p_2^*(\delta^{-1}  \otimes \omega^{\otimes 2}).\end{equation}

We recall a few general facts about dualizing sheaves (see for example \cite[Ch.~0DWE]{stacks}).
If a scheme $Y$ is Cohen--Macaulay (of constant dimension) over a base $R$, then it admits a relative
dualizing sheaf $\CK_{Y/R}$, which is invertible if and only if $Y$ is Gorenstein over $R$.   If $\rho:Y \to Y'$ is \'etale
and $Y'$ is Cohen--Macaulay over $R$, then $\CK_{Y/R}$ is canonically identified with $\rho^*\CK_{Y'/R}$, and 
with $\wedge_{\CO_Y}^n\Omega^1_{Y/R}$ if $Y$ is smooth over $R$ of dimension $n$.  More generally if 
$X$ is smooth over $R$ of dimension $n$, $Y$ is syntomic over $R$ of dimension $n-d$, and $i: Y \hookrightarrow X$
is a closed immersion with conormal bundle $\calC_{Y/X}$, then
\begin{equation} \label{eqn:cotcx}
\CK_{Y/R}\,\, \cong \,\, (\wedge_{\CO_Y}^d \calC_{Y/X})^{-1} \otimes_{\CO_Y} i^*\CK_{X/R}  \,\,\cong\,\,
                   (\wedge_{\CO_Y}^d \calC_{Y/X})^{-1} \otimes_{\CO_Y} i^*(\wedge^n_{\CO_X} \Omega^1_{X/R}).\end{equation}
Furthermore the isomorphism is compatible in the obvious sense with \'etale base-change $X \to X'$, and if $Y$ itself is
smooth over $R$, then it may be identified with the isomorphism arising from the canonical exact sequence
\begin{equation} \label{eqn:diffs} 0 \longrightarrow \calC_{Y/X} \longrightarrow i^* \Omega^1_{X/R}  
\longrightarrow \Omega^1_{Y/R} \longrightarrow 0.\end{equation}
 
 Combining the isomorphisms (\ref{eqn:KSnew}), (\ref{eqn:KSx2}) and \ref{eqn:cotcx}), we conclude that
 if $U$ is sufficiently small that Proposition~\ref{prop:immersion} holds, then we have an isomorphism
\begin{equation}\label{eqn:KStotal} \CK_{Y_{U_0(\gP)}/\CO} 
  \,\, \cong\,\, \pi_2^*\delta^{-1} \otimes \pi_1^* \omega \otimes \pi_2^* \omega. 
\end{equation}
Furthermore the compatibility of (\ref{eqn:KSnew}) and (\ref{eqn:KSx2}) with the action of $\GL_2(\A_{F,\f}^{(p)})$ and that
of (\ref{eqn:cotcx}) with \'etale base-change implies that if $U'$ is also sufficiently small for Proposition~\ref{prop:immersion}
to hold and $g \in \GL_2(\A_{F,\f}^{(p)})$ is such that $g^{-1}Ug \subset U'$, then the diagram
\begin{equation}\label{eqn:Hecke}  \begin{array}{ccc} \rho_g^* \CK_{Y_{U'_0(\gP)}/\CO} & \stackrel{\sim}{\longrightarrow} &  
\pi_2^{*}\rho_g^*\delta^{\prime-1} \otimes \pi_1^{*}\rho_g^*\omega' \otimes \pi_2^{*}\rho_g^*\omega' \\
\wr\downarrow & & \downarrow\wr \\
 \CK_{Y_{U_0(\gP)}/\CO} & \stackrel{\sim}{\longrightarrow} &  \pi_2^*\delta^{-1} \otimes \pi_1^*\omega \otimes \pi_2^*\omega \end{array} \end{equation}
commutes, where we have written $\rho_g$ for both morphisms $Y_{U_0(\gP)} \to Y_{U_0'(\gP)}$ and $Y_U \to Y_{U'}$ as usual, the top arrow is the pull-back
of (\ref{eqn:KStotal}) for $U'$, and the right arrow is induced
by the isomorphisms $\rho_g^*\delta' \stackrel{\sim}{\longrightarrow} \delta$ and $\rho_g^*\omega' \stackrel{\sim}{\longrightarrow} \omega$ defined in
\S\ref{sec:bundles}.

Note that we may therefore remove the assumption that $h:Y_{U_0(\gP)} \to Y_U \times Y_U$ is a closed immersion.  Indeed for any $U''$ sufficiently small ,
we may choose $U$ normal in $U''$ so that Proposition~\ref{prop:immersion} holds, and the commutativity of (\ref{eqn:Hecke}) for $g \in U''$ and $U = U'$
implies that (\ref{eqn:KStotal}) descends to such an isomorphism over $Y_{U''_0(\gP)}$, which is furthermore compatible with the action of
$\GL_2(\A_{F,\f}^{(p)})$ in the usual sense.  We have now proved:

\begin{theorem} \label{thm:KSIwa}   If $U$ is a sufficiently small open compact subgroup of $\GL_2(\A_{F,\f}^{(p)})$ containing $\GL_2(\CO_{F,p})$,
then there is an isomorphism
$$ \CK_{Y_{U_0(\gP)}/\CO} 
  \,\, \cong\,\, \pi_2^*\delta^{-1} \otimes \pi_1^* \omega \otimes \pi_2^* \omega;$$
furthermore the isomorphisms for varying $U$ are compatible with the action of $\GL_2(\A_{F,\f}^{(p)})$ in the sense that
diagram (\ref{eqn:Hecke}) commutes.
\end{theorem}

\begin{remark} \label{rmk:square} Recall also from Remark~\ref{rmk:symmetric} that $\pi_1^*\delta \cong \pi_2^*\delta$, so that we may exchange $\pi_1$ and
$\pi_2$ in the statement of the theorem, and deduce also the existence of isomorphisms
$$ \pi_1^*\CK_{Y_U/\CO} \otimes \pi_2^*\CK_{Y_U/\CO} \,\, \cong \,\, \CK_{Y_{U_0(\gP)}/\CO}^{\otimes 2},$$
compatible in the obvious sense with the action $\GL_2(\A_{F,\f}^{(p)})$.
\end{remark}

We will also need to relate the Kodaira--Spencer isomorphisms at level $U$ and level $U_0(\gP)$.
Consider the morphisms
$$\psi_{\theta,\CL}^*: \widetilde{\pi}_2^* \CL_\theta \longrightarrow \widetilde{\pi}_1^* \CL_\theta \quad\mbox{and}\quad
\psi_{\theta,\CM}^{*\vee}:  \widetilde{\pi}_1^*\CM_\theta^{-1} \longrightarrow \widetilde{\pi}_2^*\CM_\theta^{-1}$$
induced by the universal isogeny $\psi$ over $\widetilde{Y}_{U_0(\gP)}$.   Tensoring over all $\theta \in \Sigma$
then yields morphisms which descend to define
$${\pi}_2^* \omega \longrightarrow {\pi}_1^* \omega \quad\mbox{and}\quad
{\pi}_1^*(\delta^{-1} \otimes \omega)  \longrightarrow \pi_2^*(\delta^{-1}\otimes \omega).$$
Combining these with (\ref{eqn:KStotal}) and the isomorphisms
$\pi_j^*\CK_{Y_U/\CO} \cong \pi_j^*(\delta^{-1}\otimes \omega)  \otimes \pi_j^*\omega$
obtained by pull-back along $\pi_j$ (for $j=1,2$) of the Kodaira--Spencer isomorphism at level $U$
(i.e. the case $\gP = \CO_F$), we thus obtain morphisms
\begin{equation} \label{eqn:integral}
\pi_1^*\CK_{Y_U/\CO}  \longrightarrow \CK_{Y_{U_0(\gP)}/\CO} \quad\mbox{and}\quad
\pi_2^*\CK_{Y_U/\CO}  \longrightarrow \CK_{Y_{U_0(\gP)}/\CO} \end{equation}
(whose tensor product is $\Nm_{F/\QQ}(\gP)$ times the isomorphism of Remark~\ref{rmk:square}).

We claim that the morphisms (\ref{eqn:integral}) extend the canonical isomorphisms 
$$\pi_{1,K}^*\CK_{Y_{U,K}/K} \stackrel{\sim}{\longrightarrow} \CK_{Y_{U_0(\gP),K}/K} \quad \mbox{and}
\quad \pi_{2,K}^*\CK_{Y_{U,K}/K} \stackrel{\sim}{\longrightarrow} \CK_{Y_{U_0(\gP),K}/K}$$
induced by the \'etale morphisms $\pi_{1,K}$ and $\pi_{2,K}$ (writing $\cdot_K$ for base-changes
from $\CO$ to $K$).

We can again reduce to the case that $U$ is sufficiently small that $h$ is a closed immersion,
and replace $Y_{U,K}$ replaced by $\widetilde{Y}_{U,K}$, $Y_{U_0(\gP),K}$ by $\widetilde{Y}_{U_0(\gP),K}$, etc.
Furthermore, it suffices to prove the desired equality holds on fibres at geometric closed points of $\widetilde{Y}_{U_0(\gP),K}$.

To that end, we first recall the description on fibres of the Kodaira--Spencer isomorphism over $\widetilde{Y}_{U,K}$.
Since we are now working in characteristic zero, the filtration on $\Omega^1_{\widetilde{Y}_{U,K}/K}$ obtained from \cite[Thm.~3.3.1]{theta}
has a canonical splitting giving an isomorphism
\begin{equation} \label{eqn:KSK}  \bigoplus_{\theta\in \Sigma} \left( \CM_{\theta,K}^{-1} \otimes \CL_{\theta,K} \right) \stackrel{\sim}{\longrightarrow}  \Omega^1_{\widetilde{Y}_{U,K}/K}.\end{equation}
It follows from its construction that it is dual to the isomorphism whose fibre at the point $y \in \widetilde{Y}_U(\overline{K})$,
corresponding to a tuple $\underline{A}$, is the map
$$d: \Tan_y(\widetilde{Y}_{U,K})   \longrightarrow \bigoplus_{\theta \in \Sigma} \Hom_{\overline{K}} (H^0(A,\Omega^1_{A/\overline{K}})_\theta, H^1(A,\CO_A)_\theta)$$
with $\theta$-component $d_\theta$  described as follows:   Let $\widetilde{y} \in \widetilde{Y}_U(T)$ be a lift of $y$ to $T:= \overline{K}[\epsilon]/(\epsilon^2)$
corresponding to
data $\widetilde{\underline{A}}$ lifting $\underline{A}$, and let $\alpha_\theta$ denote the canonical isomorphism
$$H^1_\dr(\widetilde{A}/T)_\theta   \cong H^1_\crys(A/T)_\theta \cong 
H^1_\dr(A/\overline{K})_\theta \otimes_{\overline{K}} T.$$
For $e_\theta \in H^0(A,\Omega^1_{A/\overline{K}})_\theta$, let $\widetilde{e}_\theta$ be a lift
of $e_\theta$ to $H^0(\widetilde{A},\Omega^1_{\widetilde{A}/T})_\theta$.  Any two such lifts differ by an element of 
$\epsilon H^0(\widetilde{A},\Omega^1_{\widetilde{A}/T})_\theta$, which corresponds to $H^0(A,\Omega^1_{A/\overline{K}})_\theta \otimes_{\overline{K}} \epsilon T$ under the above isomorphism.
Therefore
$$\alpha_\theta(\widetilde{e}_\theta) - e_\theta \otimes 1 = f_\theta \otimes \epsilon$$
for some
$f_\theta \in H^1_\dr(A/\overline{K})_\theta$ whose image in $H^1(A,\CO_A)_\theta$ is independent of the choice
of $\widetilde{e}_\theta$, and this image is $d_\theta(\widetilde{y})(e_\theta)$.

Similarly the filtration obtained from Theorem~\ref{thm:KSfil} on the conormal bundle $\widetilde{\calC}_K$
of $\widetilde{h}:\widetilde{Y}_{U_0(\gP),K} \hookrightarrow \widetilde{Y}_{U,K} \times \widetilde{Y}_{U,K}$
has a canonical splitting giving an isomorphism
\begin{equation} \label{eqn:KSK0}
\bigoplus_{\theta\in \Sigma} \left( \widetilde{\pi}_1^*\CM_{\theta,K}^{-1} \otimes \widetilde{\pi}_2^*\CL_{\theta,K} \right) \stackrel{\sim}{\longrightarrow}  \widetilde{\calC}_K,
\end{equation}
which combined with (\ref{eqn:cotcx}) gives an exact sequence
$$0 \longrightarrow   \bigoplus_{\theta\in \Sigma} \left( \widetilde{\pi}_1^*\CM_{\theta,K}^{-1} \otimes \widetilde{\pi}_2^*\CL_{\theta,K} \right)
       \longrightarrow \widetilde{\pi}_1^*\Omega^1_{\widetilde{Y}_{U,K}/K} \otimes \widetilde{\pi}_2^* \Omega^1_{\widetilde{Y}_{U,K}/K}
       \longrightarrow \Omega^1_{\widetilde{Y}_{U_0(\gP),K}/K} \longrightarrow 0.$$
It follows from the construction (and in particular the choice of sign) in the proof of Theorem~\ref{thm:KSfil} 
that the first morphism is dual to the one whose fibre at the point
$z \in \widetilde{Y}_{U_0(\gP)}(\overline{K})$ corresponding to a tuple $(\underline{A}_1,\underline{A}_2,\psi)$
is the map
$$\Tan_{\widetilde{\pi}_1(z)}(\widetilde{Y}_{U,K}) \times \Tan_{\widetilde{\pi}_2(z)}(\widetilde{Y}_{U,K})   \longrightarrow
     \bigoplus_{\theta\in \Sigma}  \Hom_{\overline{K}} (H^0(A_2,\Omega^1_{A_2/\overline{K}})_\theta, H^1(A_1,\CO_{A_1})_\theta)$$
with $\theta$-component  described as follows:  Letting $y_1 = \widetilde{\pi}_1(z)$ correspond to $\underline{A}_1$ and
$y_2 = \widetilde{\pi}_2(z)$ correspond to $\underline{A}_2$, the lift $(\widetilde{y}_1,\widetilde{y}_2)$ is sent to
$$ d_\theta(\widetilde{y}_1) \circ \psi_{\theta,\CL}^* -  \psi_{\theta,\CM}^* \circ d_\theta(\widetilde{y}_2),$$
where we continue to write $\psi_{\theta,\CL}^*:H^0(A_2,\Omega^1_{A_2/\overline{K}})_\theta \to H^0(A_1,\Omega^1_{A_1/\overline{K}})_\theta$ and
$\psi_{\theta,\CM}^*:H^1(A_2,\CO_{A_2})_\theta \to H^1(A_1,\CO_{A_1})_\theta$ for the maps induced by $\psi$.
It follows that the diagram
$$\xymatrix{
{\displaystyle\bigoplus_{\theta\in \Sigma} }\left( \widetilde{\pi}_1^*\CM_{\theta,K}^{-1} \otimes \widetilde{\pi}_2^*\CL_{\theta,K} \right)  \ar[r]^-{\sim}  \ar[d] &
\widetilde{\calC}_K \ar[d] \\
{\displaystyle\bigoplus_{\theta\in \Sigma}} \left( (\widetilde{\pi}_1^*\CM_{\theta,K}^{-1} \otimes \widetilde{\pi}_1^*\CL_{\theta,K} )
\oplus (\widetilde{\pi}_2^*\CM_{\theta,K}^{-1} \otimes \widetilde{\pi}_2^*\CL_{\theta,K} )  \right) \ar[r]^-{\sim}  &
\widetilde{\pi}_1^*\Omega^1_{\widetilde{Y}_{U,K}/K} \bigoplus \widetilde{\pi}_2^*\Omega^1_{\widetilde{Y}_{U,K}/K} }$$
commutes, where the  horizontal isomorphisms are those of (\ref{eqn:KSK}) and (\ref{eqn:KSK0}), the left vertical map
respects the decomposition and is defined by $(1\otimes \psi_{\theta,\CL}^*,-\psi_{\theta,\CM}^{*\vee}\otimes 1)$ in the $\theta$-component,
and the right vertical map is the inclusion in (\ref{eqn:cotcx}).
Taking top exterior powers in the resulting exact sequences (and tracking signs) then yields the desired compatibility,
which we can state as follows:
\begin{proposition}  \label{prop:compatibility}  The canonical isomorphism $\pi_1^* \CK_{Y_{U,K}/K} \to \CK_{Y_{U_0(\gP),K}/K}$
extends uniquely to the morphism $\pi_1^*\CK_{Y_U/\CO} \to \CK_{Y_{U_0(\gP)}/\CO}$ corresponding under the Kodaira--Spencer
isomorphisms to the morphism $\pi_1^*(\delta^{-1}\otimes\omega) \to \pi_2^*(\delta^{-1}\otimes\omega)$
induced by $\psi$; similarly the canonical isomorphism $\pi_2^* \CK_{Y_{U,K}/K} \to \CK_{Y_{U_0(\gP),K}/K}$
extends uniquely to the morphism $\pi_2^*\CK_{Y_U/\CO} \to \CK_{Y_{U_0(\gP)}/\CO}$ corresponding 
to the morphism $\pi_2^*\omega \to \pi_1^*\omega$ induced by $\psi$.
\end{proposition}
Note that the uniqueness in the statement follows from the flatness over $\CO$ of the schemes $Y_U$ and $Y_{U_0(\gP)}$ (and the invertibility of their dualizing sheaves).

\section{Degeneracy fibres}

\subsection{Irreducible components}  \label{subsec:components}
We now turn our attention to Hilbert modular varieties in characteristic $p$, i.e., over the residue
field $k$ of $\CO$.  We fix an open compact subgroup $U \subset \GL_2(\A_{F,\f})$, as usual sufficiently
small and containing $\GL_2(\CO_{F,p})$, and we let 
$\overline{Y} = Y_{U,k}$ and $\overline{Y}_0(\gP) = Y_{U_0(\gP),k}$.

Our aim is to describe the fibres of the degeneracy map\footnote{Analogous results hold for $\overline{\pi}_2 = \pi_{2,k}$,
as can be shown using similar arguments or by deducing them from the results for $\overline{\pi}_1$ using a $w_\gP$-operator.}
 $\overline{\pi}_1 = \pi_{1,k}: \overline{Y}_0(\gP) \to \overline{Y}$,
or more precisely its restrictions to irreducible components of $\overline{Y}_0(\gP)$.  
We accomplish this by generalizing the arguments in \cite[\S7]{DKS}, where
this was carried out under the assumption that $p$ is unramified in $F$.  See also
\cite[\S4]{ERX} for partial results in this direction without this assumption.

As usual, we first carry out various constructions in the context of their \'etale covers
by fine moduli schemes.  Recall that the scheme $T: = \widetilde{Y}_{U,k}$ is equipped with a 
stratification defined by the
vanishing of partial Hasse invariants, constructed in the generality of Pappas--Rapoport models by
Reduzzi and Xiao in \cite[\S3]{RX}.   More precisely, for
$\theta \in \Sigma$, define $n_\theta = 1$ unless $\theta = \theta_{\gp,i,1}$ for some 
$\gp$, $i$, in which case $n_\theta = p$, and recall that $\sigma$ denotes the ``right shift'' 
permutation of $\Sigma$ (see~\S\ref{ssec:notation}).  We then have a surjective morphism
\begin{equation}\label{eqn:hdef}
\widetilde{h}_\theta : \CP_\theta \longrightarrow \CL_{\sigma^{-1}\theta}^{\otimes n_\theta}
\end{equation}
induced by $u$ if $n_\theta = 1$, and by the pair of surjections
\begin{equation}\label{eqn:hdef2}
\CP_\theta \,\, \stackrel{u^{e_{\gp} -1}}{\longleftarrow}  \,\,
\CH^1_\dr(A/T)_{\tau_{\gp,i}} \,\,
\stackrel{\Ver_A^*}{\longrightarrow}
\,\, \CL_{\sigma^{-1}\theta}^{\otimes p}
\end{equation}
if $n_\theta = p$, and
for each $J \subset \Sigma$, the vanishing locus of the restriction of
$\widetilde{h}_\theta$ to $\CL_\theta$ for $\theta \in J$ defines a smooth subscheme of $T$
of codimension $|J|$, which we denote $T_J$ (see \cite[\S4.1]{theta}).  

Now consider $S := \widetilde{Y}_{U_0(\gP),k}$, and let $(\underline{A},\underline{A}',\psi)$ denote
the universal triple over $S$.  Recall that for each $\theta \in \Sigma$, the isogeny $\psi$ induces an
$\CO_S$-linear morphism $\psi_\theta^* : \CP'_\theta \to \CP_\theta$; its image is a rank one
subbundle of $\CP_\theta$ if $\theta \in \Sigma_\gP$ and it is an isomorphism if $\theta \not\in \Sigma_\gP$.
Furthermore $\psi_\theta^*$ restricts to a morphism $\CL'_\theta \to \CL_\theta$ of line bundles over $S$, which
we denote $\psi_{\theta,\CL}^*$, and hence also induces one we denote $\psi_{\theta,\CM}^*:\CM_{\theta}' \to \CM_\theta$.

For $J,J' \subset \Sigma_{\gP}$, consider the closed subschemes $S_{J,J'}$ of $S$ defined\footnote{Note the
deviation here from the notational conventions of \cite{DKS} and \cite{DK2}, where these would be written as
$S_{\sigma(J'),J }$ and $Z_{\sigma(J'),J}$.}  by the vanishing of $\psi_{\theta,\CM}^*$ for all $\theta \in J$,
and $\psi_{\theta,\CL}^*$ for all $\theta \in J'$.    Note that a closed point $y$ of $S$ is in $S_{J,J'}$ if and only
if $J \subset \Sigma_y$ and $J' \subset \Sigma_y'$, where $\Sigma_y$ and $\Sigma'_y$ were defined in \S\ref{sec:lci}.
Furthermore, in terms of the parameters and basis elements chosen in \S\ref{sec:lci} for stalks at $y$, we see that 
$$\begin{array}{rl}
\mbox{if $\theta \in \Sigma_y$, then}& 

\psi_{\theta,\CM}^*(f_\theta' + L_\theta') = \left\{\begin{array}{ll} s_\theta(f_\theta + L_\theta), & \mbox{if $\theta \in \Sigma'_y$,}\\
                                                                                                    \theta(\varpi_\gP)(f_\theta + L_\theta), & \mbox{if $\theta \not\in \Sigma'_y$,}\end{array}\right.\\

\mbox{and if $\theta \in \Sigma_y'$, then}& 

\psi_{\theta,\CL}^*(e_\theta' - s_\theta' f_\theta') = \left\{\begin{array}{ll} - s_\theta'(e_\theta - s_\theta f_\theta), & \mbox{if $\theta \in \Sigma_y$,}\\
                                                                                                    \theta(\varpi_\gP)(e_\theta - s_\theta f_\theta), & \mbox{if $\theta \not\in \Sigma_y$,}\end{array}\right.

 \end{array}$$
It follows that if $y \in S_{J,J'}$, then $S_{J,J'}$ is defined in a neighborhood of $y$ by the ideal generated by
$$\{\varpi\}\quad \cup \quad \{\,s_\theta\,|\,\theta \in J\,\}\quad \cup\quad\{\,s_\theta'\,|\,\theta\in J'\}.$$

We have the following generalization\footnote{Note however that the descriptions of the cases in the proof of 
\cite[Cor.~4.3.2]{DKS} are incorrect: the first case there should
be $\theta \in (I_Q - I)\cap (J_Q - J)$, the second $(J_Q - I_Q) \cup (J-I)$, and the third $(I_Q - J_Q) \cup (I-J)$.} of \cite[Cor.~4.3.2]{DKS}:

\begin{proposition}  The scheme $S_{J,J'}$ is a reduced local complete intersection 
of dimension $d - |J\cap J'|$, and is smooth over $k$ if $J \cup J' = \Sigma_{\gP}$.
\end{proposition}
\begpf Let $y$ be a closed point of $S_{J,J'}$.  The description of the completion 
$\CO_{\widetilde{Y}_0(\gP),y}$ at the end of \S\ref{sec:lci} gives
an isomorphism
$$\CO_{S,y}^{\wedge}   \cong  k(y) [[ X_\theta,X_\theta']]_{\theta \in \Sigma}/ \langle \overline{g}_\theta \rangle_{\theta \in \Sigma}$$
under which $s_\theta \mapsto X_\theta$ and $s_\theta' \mapsto X_\theta'$, where
$$\ol{g}_\theta = \left\{\begin{array}{ll}
X_\theta & \mbox{if $\theta\in\Sigma_y - \Sigma_y'$,}\\
X'_\theta & \mbox{if $\theta\in\Sigma'_y - \Sigma_y$,}\\
X_\theta X'_\theta & \mbox{if $\theta\in\Sigma_y \cap \Sigma_y'$,}\\
X_\theta - \overline{\theta}(\varpi_\gP)X_\theta',& \mbox{if $\theta\not\in\Sigma_\gP$.}
\end{array}\right.$$
It follows that $\CO_{S_{J,J'},y}^{\wedge}$ is isomorphic to the quotient of this ring by the ideal $\langle X_\theta \rangle_{\theta \in J} + \langle X'_\theta \rangle_{\theta \in J'}$, and therefore to 
$\widehat{\bigotimes}_{\theta\not\in J \cap J'} R_\theta$, where
$$R_\theta = \left\{\begin{array}{ll}
k(y)[[X_\theta,X_\theta']]/(X_\theta X'_\theta),&\mbox{if $\theta \in (\Sigma_y - J) \cap (\Sigma_y' - J')$,}\\
\mbox{$k(y)[[X_\theta]]$ or $k(y)[[X_\theta']]$,}&
 \mbox{otherwise.}
\end{array}\right.$$
In particular, $\CO_{S_{J,J'},y}^{\wedge}$ is a reduced complete intersection of dimension $d - |J\cap J'|$, and is regular if $J \cup J' = \Sigma_\gP$.
\epf

If $J' = \Sigma_{\gP} - J$, then $S_{J,J'}$ is a union of irreducible components of $S$, which we
write simply as $S_J$; furthermore each irreducible component of $S$ is a connected
component of $S_J$ for a unique $J \subset \Sigma_{\gP}$.

For each $J \subset \Sigma_{\gP}$, we consider the restriction of $\widetilde{\pi}_{1,k}$ to $S_J$.
We first note that the morphism factors through $T_{J'}$, where $J' = \{\,\theta\in J\,|\,\sigma^{-1}\theta\not\in J\,\}$.  (Recall that $T_{J'}$ is the intersection of the vanishing loci of the partial Hasse invariants $\CL_\theta 
\stackrel{\widetilde{h}_\theta}{\longrightarrow} \CL_{\sigma^{-1}\theta}^{\otimes n_\theta}$ for $\theta \in J'$.)
Indeed since $S_J$ is reduced, it suffices to check the assertion on geometric closed points, which reduces to
the statement that if  $(\underline{A},\underline{A}',\psi)$ corresponds to an $\Fpbar$-point of $S$ such that
$\psi_{\theta,\CM}^*:M'_\theta \to M_\theta$ and $\psi_{\sigma^{-1}\theta,\CL}^*:L'_{\sigma^{-1}\theta} \to L_{\sigma^{-1}\theta}$
both vanish, then so does $\widetilde{h}_\theta:L_\theta \to L_{\sigma^{-1}\theta}^{\otimes n_\theta}$ (with the obvious notation).
If $\tau = \tau_{\gp,i}$ and $\theta = \theta_{\gp,i,j}$ for some $j > 1$, then the vanishing of $\psi_{\theta,\CM}^*$
means that 
$$\psi_{\tau}^*(u^{-1}F_{\tau}^{\prime(j-1)}) =\psi_\tau^*(G_\tau^{\prime(j)}) \subset F_{\tau}^{(j)},$$
and hence equality holds by comparing
dimensions.  On the other hand the vanishing of $\psi_{\sigma^{-1}\theta,\CL}^*$ means that 
$\psi_{\tau}^*(F_{\tau}^{\prime(j-1)}) = F_{\tau}^{(j-2)}$, and it follows that $uF_{\tau}^{(j)} = F_{\tau}^{(j-2)}$, so
$\widetilde{h}_\theta$ vanishes. Similarly  if $\tau = \tau_{\gp,i}$ and $\theta = \theta_{\gp,i,1}$, then
we find that 
$$\psi_\tau^*(u^{e_\gp-1}H^1_\dr(A'/\Fpbar)_\tau) = \psi_\tau^*(G_\tau^{\prime(1)}) = F_\tau^{(1)}$$
and $\psi_{\phi^{-1}\tau}^*(F_{\phi^{-1}\tau}^{\prime(e_\gp)}) = F_{\phi^{-1}\tau}^{(e_\gp-1)}$.  Since
$\Ver_A^*\circ \psi_\tau^* =  \psi_{\phi^{-1}\tau}^*\circ \Ver_{A'}^*$ and 
$$\Ver_{A'}^*(H^1_\dr(A'/\Fpbar)_\tau) = (F_{\phi^{-1}\tau}^{\prime(e_\gp)})^{(p)},$$
it follows that $\Ver_A^*(F_\tau^{(1)}) = u^{e_{\gp}-1} (F_{\phi^{-1}\tau}^{(e_\gp-1)})^{(p)}$,
which again means that $\widetilde{h}_\theta$ vanishes.
Similarly one finds that the restriction of $\widetilde{\pi}_{2,k}$ to $S_J$ factors
through $T_{J''}$, where $J'' = \{\,\theta\not\in J\,|\,\sigma^{-1}\theta \in J\,\}$.

Now consider the scheme
$$P_J : = \prod_{\theta \in J''} \PP_{T_{J'}}(\CP_\theta),$$
where the product is fibred over $T_{J'}$ and $\PP_{T_{J'}}(\CP_\theta)$ denotes the projective bundle parametrizing rank 
one subbundles of $\CP_\theta$ over $T_{J'}$.  (We will
freely use $\CP_\theta$ to
denote the pull-back of the rank two vector bundle associated to the universal object
over $T$ when the base and morphism to $T$ are clear from the context.)  The image
of $\psi_\theta^*: \CP_\theta' \to \CP_\theta$ is a rank one subbundle of $\CP_\theta$
over $S_J$, so it determines a morphism $S_J \to \PP_{T_{J'}}(\CP_\theta)$, and we let
$$\widetilde{\xi}_J : S_J \longrightarrow P_J$$
denote their product over $\theta \in J''$.  

The proof of \cite[Prop.~4.5]{ERX} shows
that $\widetilde{\xi}_J$ is bijective on closed points, and hence is a Frobenius factor,
in the sense that {\em some} power of the Frobenius endomorphism on $P_J$ factors
through it.  We will instead use relative Dieudonn\'e theory, as in \cite[\S7.1]{DKS}, 
to geometrize the pointwise construction in \cite{ERX} of the inverse map, showing
in particular that a {\em single} power of Frobenius suffices.

Before we carry this out in the next section, let us remark that the closed subschemes $T_{J}$ of $T$
descend to closed subschemes of $\overline{Y}$, which we denote by $\overline{Y}_{J}$, and similarly
the subschemes $S_J$ (or more generally $S_{J,J'}$) descend to subschemes of $\overline{Y}_0(\gP)$, 
which we denote $\overline{Y}_0(\gP)_J$ (or $\overline{Y}_0(\gP)_{J,J'}$).  Furthermore since 
we are now working in characteristic $p$, we may choose $U$ sufficiently small that the natural action
$U \cap \CO_F^\times$ on the vector bundles $\CP_\theta$ is trivial.  The action of
$\CO_{F,(p),+}^\times$ therefore defines descent data on the $\CP_\theta$ to a vector
bundle over $\overline{Y}$, and denoting this too by $\CP_\theta$, we obtain a morphism
$${\xi}_J : \overline{Y}_0(\gP)_J \longrightarrow 
 \prod_{\theta \in J''} \PP_{\overline{Y}_{J'}}(\CP_\theta).$$

\subsection{Frobenius factorization}  \label{sec:Frobfact}
We will define a morphism $P_J \to S$ corresponding to an isogeny from $A^{(p)}$, where $A$ is the pull-back
to $P_J$ of the universal abelian scheme over $T_{J'}$.  In order to do this, we first construct a certain
Raynaud $(\CO_F/\gP)$-module scheme over $T_{J'}$.  To that end, we need to define a line bundle
$\CA_\tau$ over $P_J$ for each $\tau \in \Sigma_{\gP,0}$, along with morphisms 
$$s_\tau:\CA_\tau^{\otimes p} \longrightarrow \CA_{\phi\circ\tau}\quad
\mbox{and}\quad t_\tau:\CA_{\phi\circ\tau} \longrightarrow \CA_{\tau}^{\otimes p}$$
such that $s_\tau t_\tau = 0$ for all $\tau$.

Suppose then that $\tau = \tau_{\gp,i} \in \Sigma_{\gP,0}$, let $\theta = \theta_{\gp,i,1}$.
and consider the exact sequence
\begin{equation}\label{eqn:Raynaud} 0 \longrightarrow \CA_\tau \longrightarrow \CP_\theta \longrightarrow \CA'_\tau \longrightarrow 0\end{equation}
over $P_J$, where $\CA_\tau$ and $\CA'_\tau$ are the line bundles defined by
\begin{itemize}
\item $\CA_\tau = \CL_\theta$ and $\CA'_\tau = \CM_\theta$ if $\theta \in J$;
\item $\CA_\tau = \ker(\widetilde{h}_\theta) \stackrel{\sim}{\leftarrow} \CM_{\sigma^{-1}\theta}^{\otimes p}$ and 
 $\CA'_\tau = \CP_\theta/\CA_\tau \stackrel{\sim}{\rightarrow} \CL_{\sigma^{-1}\theta}^{\otimes p}$
if $\sigma^{-1}\theta = \theta_{\gp,i-1,e_\gp}\not\in J$, with
$\widetilde{h}_\theta:\CP_\theta \twoheadrightarrow \CL_{\sigma^{-1}\theta}^{\otimes p}$ given by
(\ref{eqn:hdef2}) and the isomorphism $\CM_{\sigma^{-1}\theta}^{\otimes p} \stackrel{\sim}{\rightarrow} \CA_\tau$
induced by $\Frob_A^*:  (\CG_{\phi^{-1}\circ\tau}^{(e_{\gp})})^{(p)} {\rightarrow} \CP_\theta$;
\item the tautological exact sequence over $\PP_{T_{J'}}(\CP_\theta)$ if $\theta \in J''$.
\end{itemize}
Note that the first two conditions are both satisfied if $\theta \in J'$, in which case the definitions coincide since
$\ker(\widetilde{h}_\theta) = \CL_\theta$ on $T_{J'}$.

We now show that $\Frob_A^*:\CP_{\theta_{\gp,i,1}}^{(p)} \to \CP_{\theta_{\gp,i+1,1}}$
and $\Ver_A^*:\CP_{\theta_{\gp,i+1,1}} \to \CP_{\theta_{\gp,i,1}}^{(p)}$ restrict to morphisms
$\CA_\tau^{\otimes p} \to \CA_{\phi\circ\tau}$ and $\CA_{\phi\circ\tau} \to \CA_\tau^{\otimes p}$,
and hence also induce morphisms 
$\CA_\tau^{\prime\otimes p} \to \CA'_{\phi\circ\tau}$ and $\CA'_{\phi\circ\tau} \to \CA_\tau^{\prime\otimes p}$,
\begin{itemize}
\item If $\theta = \theta_{\gp,i,1} \in J$, then
$\CA_\tau = \CL_\theta = \Ver_A^*(u^{1-e_{\gp}}\CP_{\theta_{\gp,i+1,1}})$, so 
$\Frob_A^*(\CA_\tau^{\otimes p}) = 0$ and
$$\Ver_A^*(\CA_{\phi\circ\tau}) \subset \Ver_A^*(\CP_{\theta_{\gp,i+1,1}})
 = u^{e_\gp-1}(\CF_\tau^{(e_\gp)})^{(p)} \subset \CL_\theta^{\otimes p} = \CA_\tau^{\otimes p}.$$
 \item If $\theta_{\gp,i,e_\gp} \not\in J$, then $\CA_{\phi\circ \tau} = 
 \Frob_A^*((\CG_\tau^{(e_\gp)})^{(p)})$, so $\Ver_A^*(\CA_{\phi\circ\tau}) = 0$ and
 $$\Frob_A^*(\CA_\tau^{\otimes p}) \subset \Frob_A^*(\CP_{\theta_{\gp,i,1}}^{(p)})
 \subset  \Frob_A^*((\CG_\tau^{(e_\gp)})^{(p)}) = \CA_{\phi\circ \tau}.$$
\item If $\theta = \theta_{\gp,i,1} \not\in J$ and $\theta_{\gp,i,e_\gp} \in J$, then $\theta_{\gp,i,j} \in J'$
for some $j \in \{2,\ldots,e_\gp\}$, which implies that $u\CF_\tau^{(j)} = \CF_\tau^{(j-2)}$.  We deduce that
$$\CP_\theta \subset \CG_{\tau}^{(j-1)} = \CF_\tau^{(j)} \subset \CF_\tau^{(e_\gp)} \subset u\CH^1_\dr(A/P_J)_\tau$$
(where the last inclusion holds since $u^{e_\gp-1}\CF_\tau^{(e_\gp)}\subset u^{j-1}\CF_\tau^{(j)} = u^{j-1}\CG_\tau^{(j-1)} = 0$), 
and hence that $\Frob_A^*(\CA_\tau^{\otimes p}) \subset \Frob_A^*(\CP_\theta^{(p)}) \subset \Frob_A^*((\CF_\tau^{(e_\gp)})^{(p)} = 0$
and $\Ver_A^*(\CA_{\phi\circ\tau}) \subset \Ver_A^*(\CP_{\theta_{\gp,i+1,1}}) = u^{e_\gp-1}(\CF_\tau^{(e_\gp)})^{(p)} = 0$.
  \end{itemize}

 We now define $s_\tau:\CA_\tau^{\otimes p} \to \CA_{\phi\circ\tau}$ and $t_\tau:\CA_{\phi\circ\tau} \to \CA_\tau^{\phi\circ\tau}$ 
 as the restrictions of $\Frob_A^*$ and $\Ver_A^*$, yielding Raynaud data $(\CA_\tau,s_\tau,t_\tau)_{\tau \in \Sigma_{\gP,0}}$, and we let $C$ denote the
 corresponding $(\CO_F/\gP)$-module scheme over $P_J$.  We similarly obtain Raynaud data from the morphisms 
 $s_\tau':\CA_\tau^{\prime\otimes p} \to \CA'_{\phi\circ\tau}$
 and $t_\tau:\CA'_{\phi\circ\tau} \to \CA_\tau^{\prime\otimes p}$ induced by $\Frob_A^*$ and $\Ver_A^*$, and
 we let $C'$ denote the resulting $(\CO_F/\gP)$-module scheme.

Recall from \cite[\S7.1.2]{DKS} that the Dieudonn\'e crystal $\D(C)$ is isomorphic to $\Pi^*(\CA)$, with $\Phi = \Pi^*(s)$
and $V = \Pi^*(t)$, where $\CA = \oplus \CA_\tau$, $s = \oplus s_\tau$, $t = \oplus t_\tau$ (the direct sums being over
$\tau \in \Sigma_{\gP,0}$) and $\Pi^*$ is the functor denoted $\Phi^*$ in \cite[\S4.3.4]{BBM}.  (The assertion of \cite[Prop.7.1.3]{DKS}
is for $(\CO_F/p)$-module schemes, but the analogous one for $(\CO_F/\gP)$-module schemes is an immediate consequence.)
Similarly we have that $\D(C')$ is isomorphic to $\Pi^*(\CA')$ with $\Phi = \Pi^*(s')$ and $V = \Pi^*(t')$ 
(again letting $\CA' = \oplus \CA'_\tau$, etc.).
On the other hand, we have the canonical isomorphisms 
$$\D(A^{(p)}[\gP]) \cong \Pi^*\D(A[\gP])_{P_J}  \cong \Pi^*\CH^1_\dr(A/P_J) \otimes_{\CO_F} (\CO_F/\gP)$$
provided by \cite[(4.3.7.1), (3.3.7.3)]{BBM}, under which $\Phi$ (resp.~$V$) corresponds to $\Pi^*\Frob_A$
(resp.~$\Pi^*\Ver_A$).  Combined with the isomorphism
$$\CH^1_\dr(A/P_J) \otimes_{\CO_F} (\CO_F/\gP) =  \bigoplus_{\gp|\gP} \CH^1_\dr(A/P_J) \otimes_{\CO_F} (\CO_F/\gp)
\stackrel{\sim}{\longrightarrow} \bigoplus_{\gp|\gP}  \bigoplus_{i\in \ZZ/f_\gp\ZZ} \CP_{\theta_{\gp,i,1}},$$
induced by multiplication by  $u^{e_\gp-1}$ in the $\gp$-component, it follows that (\ref{eqn:Raynaud})
yields an exact sequence of Dieudonn\'e crystals
$$ 0  \longrightarrow \D(C) \longrightarrow \D(A^{(p)}[\gP]) \longrightarrow \D(C') \longrightarrow 0$$
over $P_J$.  Since $P_J$ is smooth over $k$, we may apply \cite[Thm.~4.1.1]{BM3} to conclude that
this arises from an exact sequence of $(\CO_F/\gP)$-module schemes
$$0 \longrightarrow C' \longrightarrow A^{(p)}[\gP] \longrightarrow C \longrightarrow 0.$$
We now define $A'' = A^{(p)}/C'$, and we let $\alpha:A^{(p)} \to A''$ denote the resulting isogeny and
$s'':A'' \to P_J$ the structure morphism.

We now equip $A''$ with the required auxiliary data so that the triple
$(\underline{A}^{(p)},\underline{A}'',\alpha)$ will correspond to a morphism $P_J \to S$.
First note that $A''$ inherits an $\CO_F$-action from $A^{(p)}$, as well as a unique
$U^p$-level structure with which $\alpha$ is compatible. 

The existence of a quasi-polarization
$\lambda''$ satisfying the required compatibility amounts to $\lambda^{(p)}$ inducing an
isomorphism $\gP\gc\gd\otimes_{\CO_F}A'' \stackrel{\sim}{\longrightarrow} (A'')^\vee$
over each connected component of $P_J$.  The corresponding argument in \cite[\S6.1]{theta}
carries over {\em mutatis mutandis} (with $A$ replaced by $A^{(p)}$, $A'$ by $A''$, $H$ by $C'$, $\gp$ by $\gP$ and
$\CI_i$ by $\CA_\tau^{\otimes p}$, now running over $\tau \in \Sigma_{\gP,0}$ instead of 
$i \in \ZZ/f\ZZ$) to reduce this to the orthogonality of $\CA_\tau^{\otimes p}$ and 
$u^{1-e_\gp} \CA_\tau^{\otimes p}$ under the pairing on $\CH^1_\dr(A^{(p)}/P_J)_\tau$
induced by $\lambda^{(p)}$ (via \cite[(4)]{theta}).  This in turn is equivalent to the orthogonality of
$\CA_\tau$ and $u^{1-e_\gp}\CA_\tau$ under the pairing induced by $\lambda$, which follows from
\cite[Lemma~3.1.1]{theta} via the implication that the resulting perfect pairing on 
$\CP_{\theta_{\gp,i,1}}$ is alternating.

We now define Pappas--Rapoport filtrations for $A''$, i.e., a chain of $\CO_{P_J}[u]/(u^{e_\gp})$-modules
$$0 =: \CF_\tau^{\pp(0)} \subset \CF_\tau^{\pp(1)} \subset \CF_\tau^{\pp(2)} \subset \cdots \CF_\tau^{\pp(e_\gp-1)} \subset \CF_\tau^{\pp(e_\gp)}
       := (s''_*\Omega^1_{A''/P_J})_\tau$$
for each $\tau = \tau_{\p,i} \in \Sigma_0$, such that $\CF_\tau^{\pp(j)}/\CF_\tau^{(j-1)}$ is an invertible $\CO_{P_J}$-module annihilated by $u$
for $j=1,\ldots,e_\gp$.

For each $\tau \in \Sigma_0$, consider the
morphism
$$\alpha_\tau^*:  \CH^1_\dr(A''/P_J)_\tau \longrightarrow \CH^1_\dr(A^{(p)}/P_J)_\tau = \CH^1_\dr(A/P_J)^{(p)}_{\phi^{-1}\circ\tau}$$
Zariski-locally free rank two $\CO_{P_J}[u]/(u^{e_\gp})$-modules.   Note that if $\tau \not\in \Sigma_{\gP,0}$, then
$\alpha_\tau^*$ is an isomorphism sending $(s''_*\Omega^1_{A''/P_J})_\tau$ to $(s_*\Omega^1_{A/P_J})_{\phi^{-1}\circ\tau}^{(p)}$,
and we define $\CF^{\pp(j)}_\tau = (\alpha_\tau^*)^{-1}((\CF_{\phi^{-1}\tau}^{(j)})^{(p)})$
for $j=1,\ldots,e_{\gp}-1$.  

We may therefore assume  that $\tau = \tau_{\gp,i}$ for some $\gp|\gP$, so the image of $\alpha_\tau^*$ is the $\CO_{P_J}$-subbundle
$u^{1-e_\gp}\CA_{\phi^{-1}\circ\tau}^{\otimes p}$ of rank $2e_\gp -1$, and the kernel of $\alpha_\tau^*$ is a rank one subbundle of
$\CH^1_\dr(A''/P_J)_\tau$.  For each $j=1,\ldots,e_\gp-1$, we will define $\CF_\tau^{\pp(j)}$
as an $\CO_{P_J}$-subbundle of $\CH^1_\dr(A''/P_J)_\tau$ of rank $j$, and then verify the required inclusions on
fibres at closed points.

Suppose first that $\theta_{\gp,i-1,j}\not\in J$.  Since $u^{e_\gp-1}\CF_{\phi^{-1}\circ\tau}^{(j-1)} = 0$, it follows that
$\CF^{(j-1)}_{\phi^{-1}\circ\tau}$ is contained in $u^{1-e_\gp}\CA_{\phi^{-1}\circ\tau}$, so
$(\CF^{(j-1)}_{\phi^{-1}\circ\tau})^{(p)} \subset \im(\alpha_\tau^*)$.  Therefore 
$(\alpha_\tau^*)^{-1}(\CF^{(j-1)}_{\phi^{-1}\circ\tau})^{(p)}$ is an $\CO_{P_J}$-subbundle of $\CH^1_\dr(A''/P_J)_\tau$ of
rank $j$, which we define to be $\CF^{\pp(j)}_\tau$.

Suppose next that $\theta_{\gp,i-1,j} \in J$.  We claim that
$$\CA_{\phi^{-1}\circ\tau} \subset \CF_{\phi^{-1}\circ\tau}^{(j)}  \subset \CG^{(j+1)}_{\phi^{-1}\circ\tau} \subset u^{1-e_\gp}\CA_{\phi^{-1}\circ\tau}.$$
Indeed if $\theta_{\gp,i-1,1} \in J$, then
$\CA_{\phi^{-1}\circ\tau} = \CF^{(1)}_{\phi^{-1}\circ\tau}$ and
$$\CF_{\phi^{-1}\circ\tau}^{(1)}  \subset \CF^{(j)}_{\phi^{-1}\circ\tau} \subset  \CG^{(j+1)}_{\phi^{-1}\circ\tau} \subset 
\CG^{(e_\gp)}_{\phi^{-1}\circ\tau} \subset u^{1-e_\gp}\CF^{(1)}_{\phi^{-1}\circ\tau}.$$
On the other hand if $\theta_{\gp,i-1,1} \not\in J$, then $\theta_{\gp,i-1,\ell} \in J'$ for some $\ell \in \{2,\ldots,j\}$, so that
$\CF_{\phi^{-1}\circ\tau}^{(\ell)} = \CG_{\phi^{-1}\circ\tau}^{(\ell-1)} = u^{-1}\CF_{\phi^{-1}\circ\tau}^{(\ell-2)}$,
which implies that $u^{e_\gp-1}\CG_{\phi^{-1}\circ\tau}^{(e_\gp)} = u^{e_\gp-2}\CF_{\phi^{-1}\circ\tau}^{(e_\gp-1)}
\subset u^{\ell-1}\CF^{(\ell)} = u^{\ell-2}\CF^{(\ell-2)} = 0$, and hence
$$
\CA_{\phi^{-1}\circ\tau} \subset \CG^{(1)}_{\phi^{-1}\circ\tau} \subset \CG^{(\ell-1)}_{\phi^{-1}\circ\tau} = \CF^{(\ell)}_{\phi^{-1}\circ\tau}
 \subset \CF^{(j)}_{\phi^{-1}\circ\tau} \subset 
 \CG^{(j+1)}_{\phi^{-1}\circ\tau} \subset
 \CG^{(e_\gp)}_{\phi^{-1}\circ\tau}
\subset u^{1-e_\gp}\CA_{\phi^{-1}\circ\tau}.$$
Now define $\CT_{\phi^{-1}\circ\tau}^{(j+1)}$ to be $\CF^{(j+1)}_{\phi^{-1}\circ\tau}$ if $\theta_{\gp,i-1,j+1} \in J$, and
to be the preimage in $\CG^{(j+1)}_{\phi^{-1}\circ\tau}$ of the pull-back to $P_J$ of the tautological subbundle
of $\CP_{\theta_{\gp,i-1,j+1}}$ over $\PP_{T_{J'}}(\CP_{\theta_{\gp,i-1,j+1}})$ if $\theta_{\gp,i-1,j+1}\not\in J$
(in which case $\theta_{\gp,i-1,j+1} \in J''$).  In either case
$(\CT^{(j+1)}_{\phi^{-1}\circ\tau})^{(p)}$ is contained in the image of $\alpha_\tau^*$, so that
$(\alpha_\tau^*)^{-1}(\CT^{(j+1)}_{\phi^{-1}\circ\tau})^{(p)}$ is an $\CO_{P_J}$-subbundle of $\CH^1_\dr(A''/P_J)$ of
rank $j+2$.  Furthermore we have
$$\alpha_\tau^*(u^{e_\gp-1}\CH^1_\dr(A''/P_J)_\tau) = \CA_{\phi^{-1}\circ\tau}^{\otimes p}
\subset (\CT^{(j+1)}_{\phi^{-1}\circ\tau})^{(p)},$$ so that $u^{e_\gp-1}\CH^1_\dr(A''/P_J)_\tau$ is a rank two
$\CO_{P_J}$-subbundle of $(\alpha_\tau^*)^{-1}(\CT^{(j+1)}_{\phi^{-1}\circ\tau})^{(p)}$.  We conclude that
$$\CF^{\pp(j)}_\tau := u(\alpha_\tau^*)^{-1}(\CT^{(j+1)}_{\phi^{-1}\circ\tau})^{(p)}$$
is a rank $j$ subbundle of $u\CH^1_\dr(A''/P_J)_\tau$, and hence of $\CH^1_\dr(A''/P_J)_\tau$.

To prove that the vector bundles $\CF_\tau^{\pp(j)}$ define a Pappas--Rapoport filtration, we must show that
$$u\CF^{\pp(j)}_{\tau} \subset \CF^{\pp(j-1)}_{\tau} \subset \CF^{\pp(j)}_{\tau}$$
for all $\tau = \tau_{\gp,i} \in \Sigma_0$ and $j=1,\ldots,e_\gp$, and the desired compatibility with $\alpha$ is 
the statement that $\alpha_{\phi\circ\tau}^*(\CF^{\pp(j)}_{\phi\circ\tau}) \subset (\CF_{\tau}^{(j)})^{(p)}$ for all such
$\tau$ and $j$.  These inclusions are immediate from the definitions if $\tau \not \in \Sigma_{\gP,0}$,
so we assume $\tau \in \Sigma_{\gP,0}$.

Note that since $P_J$ is smooth, it suffices to prove
the corresponding inclusions hold on fibres at every geometric closed point $y \in P_J(\Fpbar)$.
To that end, let $F_\tau^{(j)}$ (resp.~$F_\tau^{\pp(j)}$) denote the fibre at $y$ of $\CF_\tau^{(j)}$
(resp.~$\CF_\tau^{\pp(j)}$), and for each $\tau = \tau_{\gp,i}$ and $j$ such that 
$\theta = \theta_{\gp,i,j} \in J''$, let $T_\tau^{(j)}$ denote the fibre of the 
preimage in $\CG_\tau^{(j)}$ of the tautological line bundle in $\CP_{\theta_{\gp,i,j}}$ over $P_J$.

We also let $W = W(\Fpbar)$ and consider the Dieudonn\'e module
$$\D(A_y[\gp^\infty]) = H^1_\crys(A_y/W)_\gp = \bigoplus_{\tau \in \Sigma_{\gp,0}} H^1_\crys(A_y/W)_\tau$$
for each $\gp|\gP$. Thus the modules $H^1_\crys(A_y/W)_\tau$ are free of rank two over $W[u]/(E_\tau)$, 
equipped with $W[u]/(E_{\tau})$-linear maps
$$H^1_\crys(A_y/W)^\phi_{\phi^{-1}\circ\tau} \stackrel{\Phi}{\longrightarrow} H^1_\crys(A_y/W)_{\tau} 
                 \quad \mbox{and} \quad H^1_\crys(A_y/W)_{\tau} \stackrel{V}{\longrightarrow} H^1_\crys(A_y/W)^\phi_{\phi^{-1}\circ\tau}$$
induced by $\Frob_A$ and $\Ver_A$ and satisfying $\Phi V = V \Phi = p$,
where we use $\cdot^\phi$ to denote $\cdot \otimes_{W,\phi} W$.   Consider also the Dieudonn\'e module
$\D(A''_y[\gp^\infty]) = \oplus_{\tau \in \Sigma_{\gp,0}} H^1_\crys(A''_y/W)_\tau$, similarly equipped with
morphisms $\Phi$ and $V$, as well as the injective $W[u]/(E_{\tau})$-linear maps 
$$\widetilde{\alpha}_{\tau}^* : H^1_\crys(A''_y/W)_{\tau}  \longrightarrow H^1_\crys(A_y/W)_{\phi^{-1}\circ\tau}^\phi,$$
compatible with $\Phi$ and $V$.  Furthermore the morphisms $\Phi$, $V$ and $\widetilde{\alpha}_\tau^*$ are
compatible under the canonical isomorphisms $H^1_\crys(A_y/W)_\tau \otimes_W \Fpbar \cong H^1_\dr(A_y/\Fpbar)_\tau$
with the maps induced by the corresponding isogenies on de Rham cohomology.  Thus letting
$\widetilde{F}_\tau^{(j)}$ denote the preimage of $F_\tau^{(j)}$ in $H^1_\crys(A_y/W)_\tau$, and similarly
defining $\widetilde{F}_\tau^{\pp(j)}$, we are reduced to proving that
$$ u\widetilde{\alpha}_\tau^*(\widetilde{F}^{\pp(j)}_{\tau}) \subset
      \widetilde{\alpha}_\tau^*(\widetilde{F}^{\pp(j-1)}_{\tau}) \subset
      \widetilde{\alpha}_\tau^*(\widetilde{F}^{\pp(j)}_{\tau}) \subset
      (\widetilde{F}^{(j)}_{\phi^{-1}\circ\tau})^{\phi}$$
for all $\tau = \tau_{\gp,i} \in \Sigma_{\gP,0}$ and $j = 1,\ldots,e_\gp$, or equivalently that
\begin{equation}\label{eqn:nesting1} u\widetilde{F}^{\prime(j)}_{\tau} \subset
    \widetilde{E}^{(j-1)}_{\tau} \subset
    \widetilde{E}^{(j)}_{\tau} \subset
      \widetilde{F}^{(j)}_\tau,\end{equation}
 where $\widetilde{E}_\tau^{(j)} :=  (\widetilde{\alpha}_{\phi\circ\tau}^*(\widetilde{F}_{\phi\circ\tau}^{\pp(j)}))^{\phi^{-1}}$.

We use the following description of $\widetilde{E}_\tau^{(j)}$ for $\tau = \tau_{\gp,i}$ and
$j = 0,1,\ldots,e_\gp$, where we let $\widetilde{T}_\tau^{(j)}$ denote the preimage of 
 $T_\tau^{(j)}$ in $H^1_\crys(A_y/W)_\tau$ for $\theta_{\gp,i,j} \in J''$.
 \begin{itemize}
 \item Since $\widetilde{\alpha}_\tau^*(\widetilde{F}^{\pp(0)}_\tau) 
           = p\widetilde{\alpha}_\tau^*(H^1_\crys(A_y''/W)_\tau) 
           = pu^{1-e_\gp}\widetilde{A}_{\phi^{-1}\circ\tau}^\phi = u\widetilde{A}^\phi_{\phi^{-1}\circ\tau}$, where
  $$\widetilde{A}_\tau = \left\{\begin{array}{ll} \widetilde{F}^{(1)}_{\tau},& \mbox{if $\theta_{\gp,i,1} \in J$};\\
  u^{e_\gp-1} V^{-1}(\widetilde{F}_{\phi^{-1}\circ\tau}^{(e_\gp-1)})^{\phi},& \mbox{if $\theta_{\gp,i-1,e_\gp} \in J$};\\
  \widetilde{T}^{(1)}_\tau,& \mbox{if $\theta_{\gp,i,1} \in J''$};\end{array}\right.$$
we have
  $$\widetilde{E}^{(0)}_{\tau} =
  \left\{\begin{array}{ll}  \Phi((\widetilde{F}_{\phi^{-1}\circ\tau}^{(e_\gp-1)})^{\phi}),& \mbox{if $\theta_{\gp,i-1,e_\gp} \not\in J$};\\
  u \widetilde{F}^{(1)}_{\tau},& \mbox{if $\theta_{\gp,i,1} \in J$};\\
  u \widetilde{T}^{(1)}_\tau,& \mbox{if $\theta_{\gp,i,1} \in J''$}.\end{array}\right.$$
 \item If $1 \le j \le e_{\gp} -1$, then it follows from the definition of $\CF_{\phi\circ\tau}^{\pp(j)}$ that
 $$\widetilde{E}^{(j)}_{\tau} =
  \left\{\begin{array}{ll}   \widetilde{F}_{\tau}^{(j-1)},& \mbox{if $\theta_{\gp,i,j} \not\in J$};\\
  u \widetilde{F}^{(j+1)}_{\tau},& \mbox{if $\theta_{\gp,i,j+1} \in J$};\\
  u \widetilde{T}^{(j+1)}_\tau,& \mbox{f $\theta_{\gp,i,j+1} \in J''$}.\end{array}\right.$$
  (Note that if $\theta_{\gp,i,j} \not\in J$ and $\theta_{\gp,i,j+1} \in J$, then $\theta_{\gp,i,j+1} \in J'$,
  so that $u\widetilde{F}^{(j+1)}_\tau = \widetilde{F}_\tau^{(j-1)}$.)
\item Since $\Phi((\widetilde{F}_\tau^{\pp(e_\gp)})^{\phi}) = \widetilde{F}_{\phi\circ\tau}^{\pp(0)}$, we have
 $(\widetilde{E}_\tau^{(e_\gp)}) = p^{-1} (V(\widetilde{E}_{\phi\circ\tau}^{(0)}))^{\phi^{-1}}$, so
 the formulas in the case $j=0$ imply that
  $$\widetilde{E}^{(e_\gp)}_{\tau} =
  \left\{\begin{array}{ll}   \widetilde{F}_{\tau}^{(e_\gp-1)} & \mbox{if $\theta_{\gp,i,e_\gp} \not\in J$};\\
  (V(u^{1-e_\gp} \widetilde{F}^{(1)}_{\phi\circ\tau}))^{\phi^{-1}},& \mbox{if $\theta_{\gp,i+1,1} \in J$};\\
 (V( u^{1-e_\gp} \widetilde{T}^{(1)}_{\phi\circ\tau}))^{\phi^{-1}},& \mbox{if $\theta_{\gp,i+1,1} \in J''$}.\end{array}\right.$$
   \end{itemize}

Note in particular that
$$ u \widetilde{F}_\tau^{(j)} \subset \widetilde{E}^{(j)}_\tau \subset \widetilde{F}_\tau^{(j)} $$
for all $\tau = \tau_{\gp,i}$, $j = 0,\ldots,e_\gp$ (the case of $j=e_\gp$ for $\tau = \tau_{\gp,i}$
being equivalent by an application of $\Phi$ to that of $j=0$ for $\tau = \tau_{\gp,i+1}$).  The inclusions
in (\ref{eqn:nesting1}) are then immediate from the formulas
\begin{equation}\label{eqn:Fprime}
\widetilde{E}^{(j-1)}_\tau = u\widetilde{F}_\tau^{(j)} \quad \mbox{if $\theta_{\gp,i,j} \in J$,}\quad
 \mbox{and}\quad \widetilde{E}^{(j)}_\tau = \widetilde{F}_\tau^{(j-1)} \quad \mbox{if $\theta_{\gp,i,j} \not\in J$}.\end{equation}

  We have now shown that the triple $(\underline{A}^{(p)},\underline{A}'',\alpha)$ defines a morphism $\widetilde{\zeta}_J: P_J \to S$.
  Furthermore rewriting equations (\ref{eqn:Fprime}) in the form
$$\begin{array}{ccccc}
&\widetilde{\alpha}_{\phi\circ\tau}^*(u^{-1}  \widetilde{F}^{\pp(j-1)}_{\phi\circ\tau}) &=& (\widetilde{F}_\tau^{(j)})^\phi,& \mbox{if $\theta_{\gp,i,j} \in J$,}\\
 \mbox{and}& \widetilde{\alpha}_{\phi\circ\tau}^*(\widetilde{F}^{\pp(j)}_{\phi\circ\tau}) &= &(\widetilde{F}_\tau^{(j-1)})^\phi & \mbox{if $\theta_{\gp,i,j} \not\in J$},\end{array}$$
 shows that $\widetilde{\zeta}_J(y) \in S_{\phi(J)}(\Fpbar)$ for all $y \in P_J(\Fpbar)$ (where we define 
 $\phi: \Sigma_\gP \to \Sigma_\gP$ by $\phi(\theta_{\gp,i,j}) =\theta_{\gp,i+1,j}$).  Since $P_J$ is reduced, it follows that $\widetilde{\zeta}_J$ factors
 through $S_{\phi(J)}$.  Note furthermore that if $(\underline{A},\underline{A}',\psi)$ is the universal triple over $S_J$, then its pull-back by $\phi \in \Aut(k)$
 is isomorphic to the universal triple over $S_{\phi(J)}$, yielding an identification of $S_J^{(p)}$ with $S_{\phi(J)}$ under which the relative Frobenius morphism
 $S_J \to S_J^{(p)} = S_{\phi(J)}$ corresponds to the triple $(\underline{A},\underline{A}',\psi)^{(p)}$ over $S_J$.  Note also that $\phi(J)' = \phi(J')$ and
 $\phi(J)'' = \phi(J'')$, and we may similarly identify $T_{J'}^{(p)}$ with $T_{\phi(J')}$ and $P_J^{(p)}$ with $P_{\phi(J)}$, under which $\widetilde{\xi}_J^{(p)}$
 corresponds to $\widetilde{\xi}_{\phi(J)}$ and $\widetilde{\zeta}_J^{(p)}$ to $\widetilde{\zeta}_{\phi(J)}$.
 
 Finally we verify that the composites
 $$P_J \stackrel{\widetilde{\zeta}_J}{\longrightarrow} S_{\phi(J)} \stackrel{\widetilde{\xi}_{\phi(J)}}{\longrightarrow} P_{\phi(J)}
\quad \mbox{and}\quad S_{J} \stackrel{\widetilde{\xi}_{J}}{\longrightarrow} P_J \stackrel{\widetilde{\zeta}_J}{\longrightarrow} S_{\phi(J)}$$
are the Frobenius morphisms (relative to $\phi \in \Aut(k)$).  

Indeed it suffices to check this on geometric closed points, so suppose first that
$(\underline{A},(T_\theta)_{\theta \in J''})$ corresponds to an element $y \in P_J(\Fpbar)$, where $\underline{A}$ corresponds to an element 
of $T_{J'}(\Fpbar)$ and each $T_\theta$ to a line in $P_\theta = u^{-1}\widetilde{F}_\tau^{(j-1)}/\widetilde{F}_\tau^{(j-1)}$ (using the above 
notation with $y$ suppressed, so $\tau = \tau_{\gp,i}$, $\theta = \theta_{\gp,i,j}$, and $\widetilde{F}_\tau^{\bullet}$ is the preimage in
$H^1_\crys(A/W)_\tau$ of the Pappas--Rapoport filtration on $F_\tau^{(e_\gp)} = H^0(A,\Omega^1_{A/\Fpbar})_\tau \subset H^1_\dr(A/\Fpbar)_\tau$).
Its image $\widetilde{\zeta}_J(y)$ then corresponds to the triple $(\underline{A}^{(p)},\underline{A}'',\alpha)$ produced by the construction above,
and $\widetilde{\xi}_{\phi(J)}(\widetilde{\zeta}_J(y))$ is given by $(\underline{A}^{(p)}, \alpha_{\phi(\theta)}^*(P_{\phi(\theta)}''))$, where
$P''_{\phi(\theta)} = u^{-1}\widetilde{F}_{\phi\circ\tau}^{\pp(j-1)}/\widetilde{F}_{\phi\circ\tau}^{\pp(j-1)}$ and 
$\alpha^*_{\phi(\theta)}:P_{\phi(\theta)}'' \to P_\theta^{(p)}$ is induced by $\widetilde{\alpha}_{\phi\circ\tau}^*$.
From the formulas above, we have
$$\widetilde{\alpha}_{\phi\circ\tau}^*(\widetilde{F}_{\phi\circ\tau}^{\pp(j-1)}) = u^{-1}(\widetilde{E}_\tau^{(j-1)})^\phi = (\widetilde{T}_\tau^{(j)})^\phi$$
for $\theta_{\gp,i,j} \in J''$, where $\widetilde{T}_\tau^{(j)}$ is the preimage of $T_\theta$ in $H^1_\crys(A/W)_\tau$, and it follows that
$\alpha_{\phi(\theta)}^*(P_{\phi(\theta)}'') = T_\theta^{(p)}$ as required.

Suppose now that $(\underline{A},\underline{A}',\psi)$ corresponds to an element $x \in S_J(\Fpbar)$, so that if
$\tau = \tau_{\gp,i} \in \Sigma_{\gP,0}$ and $\theta = \theta_{\gp,i,j}$, then
$\widetilde{\psi}_\tau^*(u^{-1}\widetilde{F}_\tau^{\prime(j-1)}) = \widetilde{F}_\tau^{(j)}$ if $\theta \in J$ and
$\widetilde{\psi}_\tau^*(\widetilde{F}_\tau^{\prime(j)}) = \widetilde{F}_\tau^{(j-1)}$ if $\theta \not\in J$
(with the evident notation).
Its image $\widetilde{\xi}_J(x)$ then corresponds to $(\underline{A},(T_\theta)_{\theta\in J''})$ where 
$$T_\theta = \psi_\theta^*(P_\theta') = \widetilde{\psi}_\tau^*(u^{-1}\widetilde{F}_\tau^{\prime(j-1)})/\widetilde{F}_\tau^{(j-1)},$$
and we must show that the triple $(\underline{A}^{(p)},\underline{A}'',\alpha)$ obtained from it by applying $\widetilde{\zeta}_J$ is
isomorphic to $(\underline{A}^{(p)},\underline{A}^{\prime(p)},\psi^{(p)})$.  To that end note that the formulas for
$\widetilde{E}_\tau^{(j)}$ imply that $\widetilde{\psi}_\tau^*(\widetilde{F}_\tau^{\prime(j)}) = \widetilde{E}_\tau^{(j)}$ for all
$\tau = \tau_{\gp,i} \in \Sigma_{\gP,0}$ and $j=0,1,\ldots,e_\gp$ (using the equivalence provided by $\Phi$ between the cases
of $j=e_\gp$ for $\tau$ and $j=0$ for $\phi\circ\tau$).  In particular it follows from the case of $j=0$, that
$$\begin{array}{ccccc}
&  (\widetilde{\psi}_\tau^*)^\phi:  &H^1_\crys(A'/W)_\tau^\phi &\longrightarrow &H^1_\crys(A/W)_\tau^\phi \\
 \mbox{and} & \widetilde{\alpha}_{\phi\circ\tau}^*: &H^1_\crys(A''/W)_{\phi\circ\tau} &\longrightarrow&
  H^1_\crys(A/W)_\tau^\phi\end{array}$$
have the same image for all $\tau \in \Sigma_{\gP,0}$.  Combined with the fact that $\widetilde{\psi}_\tau^*$ and
$\widetilde{\alpha}_\tau^*$ are isomorphisms for $\tau \not\in \Sigma_{\gP,0}$, we conclude that the images of
$$\D(A^{\prime(p)}[p^\infty])  \longrightarrow \D(A^{(p)}[p^\infty]) \quad \mbox{and} \quad
    \D(A''[p^\infty]) \longrightarrow \D(A^{(p)}[p^\infty])$$
under the morphisms induced by $\psi^{(p)}$ and $\alpha$ are the same.  It follows that there is an isomorphism
$A^{\prime(p)} \stackrel{\sim}{\longrightarrow} A''$ such that the diagram
$$\xymatrix{&  A^{(p)}  \ar[dl]_-{\psi^{(p)}} \ar[dr]^-{\alpha} \\  A^{\prime(p)} \ar[rr]^{\sim} && A''}$$
commutes.  Its compatibility with all auxiliary structures is automatic, with the exception of the Pappas--Rapoport filtrations
for $\tau \in \Sigma_{\gP,0}$, where the compatibility is implied by the equality
$$(\widetilde{\psi}_{\tau'}^*(\widetilde{F}_{\tau'}^{\prime(j)}))^\phi = (\widetilde{E}_{\tau'}^{(j)})^\phi 
= \widetilde{\alpha}_{\phi\circ\tau}^*(\widetilde{F}_{\phi\circ\tau}^{\pp(j)})$$
for $\tau' = \phi^{-1}\circ\tau$ and $j=1,\ldots,e_\gp-1$.

We have now proved the following generalization of \cite[Lemma 7.1.5]{DKS}:
\begin{lemma}  \label{lem:Frobfact} There is a morphism $\widetilde{\zeta}_J: P_J \to S_J^{(p)}$ such that the composites
$\widetilde{\xi}_J^{(p)}\circ \widetilde{\zeta}_J:  P_J \to P_J^{(p)}$ and $ \widetilde{\zeta}_J\circ \widetilde{\xi}_J:  S_J \to S_J^{(p)}$
are the Frobenius morphisms (relative to $\phi \in \Aut(k)$).
\end{lemma}

\subsection{First-order deformations}
In this section we compute the effect of the morphisms 
$$S_J\stackrel{\widetilde{\xi}_J}{\longrightarrow} P_J  \longrightarrow T_{J'}$$
on tangent spaces with a view to generalizing the results of 
\cite[\S7.1.4]{DKS}.

We let $\FF$ denote an algebraically closed field of characteristic $p$.
First recall that if $x\in T(\FF)$ corresponds to the data $\underline{A} = (A,\iota,\lambda,\eta,F^\bullet)$, 
then the Kodaira--Spencer filtration on $\Omega^1_{T/k}$  (as described in \cite[\S3.3]{theta}) equips
the tangent space of $T$ at $x$ with a decomposition
$$\Tan_x(T)  = \bigoplus_{\tau \in \Sigma_0} \Tan_x(T)_\tau$$
and a decreasing filtration on each component such that 
$$\bigoplus_{\tau \in \Sigma_0} \Fil^{j_\tau}  (\Tan_x(T)_\tau) $$
is identified with the set of lifts $\underline{\widetilde{A}}$ of $\underline{A}$ to $T(\FF[\epsilon])$ with the property
that for all $\tau \in \Sigma_0$ and $1 \le j \le j_\tau$,  $\widetilde{F}_\tau^{(j)}$ corresponds to 
$F_\tau^{(j)} \otimes_{\FF} \FF[\epsilon]$ under the canonical isomorphism
$$H^1_{\dr}(\widetilde{A}/\FF[\epsilon]/(\epsilon^2))_\tau 
\cong H^1_\crys(\widetilde{A}/\FF)_\tau \cong H^1_\dr(A/\FF)_\tau \otimes_{\FF} \FF[\epsilon].$$
In particular, for each $\tau = \tau_{\gp,i}$ and $j = 1,\ldots,e_\gp$, the one-dimensional space
$$\gr^{j-1}  (\Tan_x(T)_\tau) := \Fil^{j-1}  (\Tan_x(T)_\tau) / \Fil^j (\Tan_x(T)_\tau)$$
is identified with the set of lines in $P_\theta \otimes_{\FF} \FF[\epsilon]$ lifting $L_\theta$,
or equivalently with $\Hom_{\FF}(L_\theta,M_\theta)$, where $\theta = \theta_{\gp,i,j}$ and $M_\theta = P_\theta/L_\theta$.

Furthermore if  $x\in T_{J'}(\FF)$, then $\underline{\widetilde{A}}$ corresponds to an element of $T_{J'}(\FF[\epsilon])$
if and only if the Pappas--Rapoport filtrations $\widetilde{F}_\tau^\bullet$ have the property that
$$\widetilde{F}_\tau^{(j)} = \left\{\begin{array}{ll} F_\tau^{(1)} \otimes_{\FF} \FF[\epsilon],&\mbox{if $j = 1$;} \\
u^{-1}\widetilde{F}_\tau^{(j-2)},&\mbox{if $j \ge 2$;}\end{array}\right.$$
for all $j$ such that $\theta = \theta_{\gp,i,j} \in J'$ (via the canonical isomorphism if $j=1$).
It follows that $\Tan_x(T_{J'})$ inherits a decomposition into components 
$\Tan_x(T_{J'})_\tau$ equipped with filtrations such that 
$$\gr^{j-1}(\Tan_x(T_{J'})_\tau)
= \left\{\begin{array}{ll} \gr^{j-1}(\Tan_x(T)_\tau),&\mbox{if $\theta_{\gp,i,j} \not\in J'$;} \\
0,&\mbox{if $\theta_{\gp,i,j} \in J'$.}\end{array}\right.$$

Suppose now that $y \in P_J(\FF)$ corresponds to the data $(\underline{A},(L_\theta'')_{\theta\in J''})$,
where $\underline{A}$ corresponds to a point $x \in T_{J'}(\FF)$ and each $L''_\theta$ is a line in $P_\theta$.
We may then decompose
$$\Tan_y(P_J)  = \bigoplus_{\tau \in \Sigma_0} \Tan_y(P_J)_\tau$$
and filter each component so that
$$\bigoplus_{\tau} \Fil^{j_\tau}  (\Tan_y(P_J)_\tau) $$
is identified with the set of lifts $(\underline{\widetilde{A}},(\widetilde{L}''_\theta)_{\theta\in J''})$ of $(\underline{A},(L''_\theta)_{\theta\in J''})$ 
to $P_J(\FF[\epsilon])$ such that the following hold for all $\tau  = \tau_{\gp,i} \in \Sigma_0$ and $1 \le j \le j_\tau$:
\begin{itemize}
\item $\widetilde{F}_\tau^{(j)}$ corresponds to $F_\tau^{(j)} \otimes_{\FF} \FF[\epsilon]$ under the canonical isomorphism;
\item if $\theta = \theta_{\gp,i,j} \in J''$, then $\widetilde{L}''_\theta$ corresponds to $L_\theta'' \otimes_{\FF} \FF[\epsilon]$
 under the resulting isomorphism 
$$\widetilde{P}_\theta = u^{-1}\widetilde{F}_\tau^{(j-1)}/\widetilde{F}_\tau^{j-1} \cong
    P_\theta \otimes_{\FF} \FF[\epsilon].$$
\end{itemize}
We then have that 
$$\dim_{\FF}\gr^{j-1}(\Tan_y(P_J)_\tau)
= \left\{\begin{array}{ll} 0,&\mbox{if $\theta_{\gp,i,j} \in J'$;}\\
2,&\mbox{if $\theta_{\gp,i,j} \in J''$;}\\
1,&\mbox{otherwise.}\\
\end{array}\right.$$
Furthermore the projection $\Tan_y(P_J) \to \Tan_x(T_{J'})$ respects the decomposition and
restricts to surjections $\Fil^j(\Tan_y(P_J)_\tau) \to  \Fil^j(\Tan_x(T_{J'})_\tau)$ inducing an
isomorphism $\gr^{j-1}(\Tan_y(P_J)_\tau) \to  \gr^{j-1}(\Tan_x(T_{J'})_\tau)$ unless 
$\theta = \theta_{\gp,i,j} \in J''$, in which case we have a decomposition
\begin{equation} \label{eqn:2d} \gr^{j-1}(\Tan_y(P_J)_\tau) = V_\theta \oplus V''_\theta\end{equation}
where $V_\theta$ (resp.~$V''_\theta$) is identified with the set of lines in 
$P_\theta \otimes_{\FF} \FF[\epsilon]$ lifting $L_\theta$ (resp.~$L_\theta''$), so that
the map to $\gr^{j-1}(\Tan_x(T_{J'})_\tau)$ has kernel $V''_\theta$ and restricts to an
isomorphism on $V_\theta$.

Finally suppose that $z \in S_J(\FF)$ corresponds to the data $(\underline{A},\underline{A}',\psi)$
and identify the tangent space $\Tan_z(S_J)$ with the set of lifts $(\underline{\widetilde{A}},\underline{\widetilde{A}}',\widetilde{\psi})$
corresponding to elements of $S_J(\FF[\epsilon])$.  We may then decompose
$$\Tan_z(S_J)  = \bigoplus_{\tau \in \Sigma_0} \Tan_z(S_J)_\tau$$
and filter each component so that
$$\bigoplus_{\tau} \Fil^{j_\tau}  (\Tan_z(S_J)_\tau) $$
is identified with the set of lifts $(\underline{\widetilde{A}},\underline{\widetilde{A}}', \widetilde{\psi})$ of $(\underline{A},\underline{A}',\psi)$ 
to $S_J(\FF[\epsilon])$ such that the canonical isomorphisms identify
$\widetilde{F}_\tau^{(j)}$ with $F_\tau^{(j)} \otimes_{\FF} \FF[\epsilon]$ and
$\widetilde{F}_\tau^{\prime(j)}$ with $F_\tau^{\prime(j)} \otimes_{\FF} \FF[\epsilon]$
for all $\tau  = \tau_{\gp,i} \in \Sigma_0$ and $1 \le j \le j_\tau$.
Note that if $\tau \not\in\Sigma_{\gP,0}$, then the isomorphism induced by $\psi$ renders the conditions for 
$\widetilde{F}_\tau^{(j)}$ and $\widetilde{F}_\tau^{\prime(j)}$ equivalent.  Suppose on the other hand that $\tau = \tau_{\gp,i} \in \Sigma_{\gP,0}$.
If $\theta = \theta_{\gp,i,j} \in J$ and $\widetilde{F}_\tau^{\prime(j-1)}$ corresponds to $F_\tau^{\prime(j-1)} \otimes_{\FF} \FF[\epsilon]$,
then the functoriality of the isomorphisms between crystalline and de Rham cohomology implies that
$\widetilde{G}_\tau^{\prime(j)} = u^{-1}\widetilde{F}_\tau^{\prime(j)}$ corresponds to 
$G_\tau^{\prime(j)} \otimes_{\FF} \FF[\epsilon]$, and consequently that $\widetilde{F}_\tau^{(j)}$ corresponds to
$F_\tau^{(j)} \otimes_{\FF} \FF[\epsilon]$.  It follows that in this case the $\FF$-linear map
$\gr^{j-1}(\Tan_z(S_J) _\tau) \to \gr^{j-1}(\Tan_x(T)_\tau)$ is trivial (where $x = \widetilde{\pi}_1$ corresponds to $\underline{A}$).
Similarly if $\theta = \theta_{\gp,i,j}\not\in J$ and $\widetilde{F}_\tau^{(j-1)}$ corresponds to $F_\tau^{(j-1)} \otimes_{\FF} \FF[\epsilon]$,
then we find that $\widetilde{F}_\tau^{\prime(j)}$ corresponds to $F_\tau^{\prime(j)} \otimes_{\FF} \FF[\epsilon]$,
so that the map $\gr^{j-1}(\Tan_z(S_J)_\tau) \to \gr^{j-1}(\Tan_{x'}(T)_\tau)$ is trivial (where $x' = \widetilde{\pi}_2(z)$ corresponds to $\underline{A}'$).
The degeneracy maps $\widetilde{\pi}_1$ and $\widetilde{\pi}_2$ therefore induce injective maps from $\gr^{j-1}(\Tan_z(S_J)_\tau)$ to
\begin{itemize}
\item $\gr^{j-1}(\Tan_{x'}(T)_\tau)  = \gr^{j-1}(\Tan_{x'}(T_{J''})_\tau)$ if $\theta \in J$;
\item $\gr^{j-1}(\Tan_{x}(T)_\tau) = \gr^{j-1}(\Tan_{x}(T_{J'})_\tau)$ if $\theta \in \Sigma_\gP$ but $\theta\not \in J$;
\item either of the above if $\theta \not \in \Sigma_\gP$.
\end{itemize}
Since $\Tan_z(S_J)$ has dimension $d$, it follows that all these maps are isomorphisms.

We now describe the map induced by $\widetilde{\xi}_J:S_J\to P_J$ on tangent spaces
in terms of its effect on filtrations and graded pieces.
\begin{lemma}  \label{lem:graded}  Let $z \in S_J(\FF)$ be a geometric point of $S_J$, and let $y = \widetilde{\xi}_J(z)$.
Then the $\FF$-linear homomorphism
$$t: \Tan_z(S_J) \to \Tan_y(P_J)$$
respects the decompositions and filtrations.
Furthermore if $\tau = \tau_{\gp,i}$ and $\theta = \theta_{\p,i,j}$, then the following hold:
\begin{enumerate}
\item If $\theta,\sigma^{-1}\theta\not\in J$, then $t$ induces an isomorphism 
$$\gr^{j-1}(\Tan_z(S_J)_\tau) \stackrel{\sim}{\longrightarrow} \gr^{j-1}( \Tan_y(P_J)_\tau).$$
\item If $\theta \in J$, then $t(\Fil^{j-1}(\Tan_z(S_J)_\tau)) \subset \Fil^j( \Tan_y(P_J)_\tau)$.
\item If $\theta,\sigma^{-1}\theta\in J$ and $j \ge 2$, then $t$ induces an isomorphism
$$\gr^{j-2}(\Tan_z(S_J)_\tau) \stackrel{\sim}{\longrightarrow} \gr^{j-1}( \Tan_y(P_J)_\tau).$$
\item If $\theta \in J''$, then $t$ induces an injection
$$\gr^{j-1}(\Tan_z(S_J)_\tau) \,\hookrightarrow\, \gr^{j-1}( \Tan_y(P_J)_\tau)$$
with image $V_\theta$ (under the decomposition (\ref{eqn:2d})), extending to an isomorphism
$$\Fil^{j-2}(\Tan_z(S_J)_\tau)/\Fil^j(\Tan_z(S_J)_\tau) \stackrel{\sim}{\longrightarrow}
  \gr^{j-1}( \Tan_y(P_J)_\tau)$$
if $j \ge 2$.
\end{enumerate}
\end{lemma}
\begpf  Suppose that $(\underline{\widetilde{A}},\underline{\widetilde{A}}', \widetilde{\psi})$ corresponds to
an element $\widetilde{z} \in \Fil^{j-1}(\Tan_z(S_J)_\tau)$, where $\tau = \tau_{\gp,i}$ and $\theta = \theta_{\gp,i,j}$,
and let $(\underline{\widetilde{A}},(\widetilde{L}_{\theta'})_{\theta'\in J''})$ correspond to
its image $\widetilde{y} \in \Tan_y(P_J)$.  We then have that $\widetilde{F}_{\tau'}^{(j')}$ (resp.~$\widetilde{F}_{\tau'}^{\prime(j')}$)
corresponds to $F_{\tau'}^{(j')} \otimes_{\FF} \FF[\epsilon]$ (resp.~$F_{\tau'}^{\prime(j')} \otimes_{\FF} \FF[\epsilon]$)
under the canonical isomorphisms
$$\begin{array}{crcl}
&H^1_\dr(\widetilde{A}/\FF[\epsilon])_{\tau'}& \cong& H^1_\dr(A/\FF)_{\tau'}\otimes_{\FF}\FF[\epsilon] \\
\mbox{(resp.}   & H^1_\dr(\widetilde{A}'/\FF[\epsilon])_{\tau'} & \cong&  H^1_\dr(A'/\FF)_{\tau'} \otimes_{\FF}\FF[\epsilon])\end{array}$$
for $\tau' = \tau$ and $j' \le j-1$, and for all $\tau' = \tau_{\gp',i'}  \neq \tau$ and $j'\le e_{\gp'}$.
It follows also that $\widetilde{L}''_{\theta'}$ corresponds to $L''_{\theta'}\otimes_{\FF} \FF[\epsilon]$ for all
$\theta' = \theta_{\gp',i',j'} \in J''$ such that $\tau' \neq \tau$, or $j' \le j$ if $\tau' = \tau$.
This proves that $t$ respects the decompositions and filtrations, and that $\theta \in J''$, then the image of
$\gr^{j-1}(\Tan_z(S_J)_\tau) \to \gr^{j-1}( \Tan_y(P_J)_\tau)$ is contained in $V_\theta$.
Furthermore, we have already seen that $\widetilde{F}_\tau^{(j)}$ corresponds to $F_\tau^{(j)}\otimes_{\FF} \FF[\epsilon]$ if $\theta \in J$,
so that 2) holds.  

For the injectivity in 1) and the first part of 4), recall that $\widetilde{F}_\tau^{\prime(j)}$ corresponds to 
$F_\tau^{\prime(j)} \otimes_{\FF} \FF[\epsilon] $ if $\theta \not\in J$, so if $\widetilde{y} \in \Fil^j(\Tan_y(P_J)_\tau)$, then 
$\widetilde{z} \in \Fil^j(\Tan_z(S_J)_\tau)$.  The assertions concerning the image then follow from the
fact that each space is one-dimensional.

Again comparing dimensions, it suffices to prove injectivity in 3) and the second part of 4).  Applying the
injectivity in the first part of 4), and then shifting indices to maintain the assumption that
$\widetilde{z} \in \Fil^{j-1}(\Tan_z(S_J)_\tau)$, both assertions reduce to the claim that if $j < e_\gp$,
$\theta \in J$ and $\widetilde{y} \in \Fil^{j+1}(\Tan_y(P_J)_\tau)$, then $\widetilde{z} \in \Fil^{j}(\Tan_z(S_J)_\tau)$.

To prove the claim, let $\widetilde{F}_\tau^{\prime\prime(j+1)}$ denote the image of 
$\widetilde{G}_\tau^{\prime(j+1)} =  u^{-1}\widetilde{F}_\tau^{\prime(j)}$ in
$\widetilde{G}_\tau^{(j+1)} =  u^{-1}\widetilde{F}_\tau^{(j)}$ under the morphism induced by $\widetilde{\psi}$.
Note that if $\sigma\theta \in J$, then $\widetilde{F}_\tau^{\prime\prime(j+1)} = \widetilde{F}_\tau^{(j+1)}$,
whereas if $\sigma\theta\not\in J$, then $\widetilde{F}_\tau^{\prime\prime(j+1)}/\widetilde{F}_\tau^{(j)} = 
\widetilde{L}''_{\sigma\theta}$.  Thus in either case, the hypothesis that $\widetilde{y} \in \Fil^{j+1}(\Tan_y(P_J)_\tau)$
implies that $\widetilde{F}_\tau^{\prime\prime(j+1)}$ corresponds to 
$F_\tau^{\prime\prime(j+1)} \otimes_{\FF} \FF[\epsilon]$.   Functoriality of the canonical isomorphism
between crystalline and de Rham cohomology then implies that their preimages under the maps induced
by $\widetilde{\psi}$ and $\psi$ correspond, i.e., that $\widetilde{G}_\tau^{\prime(j+1)}$ corresponds to
$G_\tau^{\prime(j+1)}\otimes_{\FF} \FF[\epsilon]$, and hence that  $\widetilde{F}_\tau^{\prime(j)}$ corresponds
to $F_\tau^{\prime(j)} \otimes_{\FF} \FF[\epsilon] $, so $\widetilde{z} \in \Fil^{j}(\Tan_z(S_J)_\tau)$.
\epf

\begin{lemma} \label{lem:rank}  Let $t_\tau:\Tan_z(S_J)_\tau \longrightarrow \Tan_y(P_J)_\tau$ denote the
$\tau$-component of the $\FF$-linear map $t$, with notation as in Lemma~\ref{lem:graded}.  Then
$$\begin{array}{rcl}  \ker(t_\tau) & = &  \left\{  \begin{array}{ll} \Fil^{e_\gp-1}(\Tan_z(S_J)_\tau),&\mbox{if $\theta_{\gp,i,e_\gp} \in J$;} \\
                                                                  0,&\mbox{if $\theta_{\gp,i,e_\gp} \not\in J$;} \end{array}  \right. \\
                                   \im(t_\tau) & = &  \left\{  \begin{array}{ll} \Tan_y(P_J)_\tau,&\mbox{if $\theta_{\gp,i-1,e_\gp} \not\in J$;} \\
                                                                            \Fil^{1}(\Tan_y(P_J)_\tau),&\mbox{if $\theta_{\gp,i,1} \in J$;} \\
                                                                             \widetilde{V}_{\theta_{\gp,i,1}} ,&\mbox{if $\theta = \theta_{\gp,i,1} \in J''$,} \end{array} \right. \end{array}$$
where $\widetilde{V}_{\theta_{\gp,i,1}}$ is the preimage of 
${V}_{\theta_{\gp,i,1}} \subset \gr^0(\Tan_y(P_J)_\tau)$ under t
he projection from $\Tan_y(P_J)_\tau$.

In particular $\ker(t)$ has dimension $\#\{\,(\gp,i)\,|\,\theta_{\gp,i,e_\gp}\in J\,\}$.
\end{lemma}
\begpf  It is immediate from part 2) of Lemma~\ref{lem:graded} that if
$\widetilde{z} \in \Fil^{e_\gp-1}(\Tan_z(S_J)_\tau)$, then $t(\widetilde{z}) = 0$.

Suppose on the other hand that $\widetilde{z} \in \ker(t_\tau)$.  We prove by induction on $j$ that
$\widetilde{z} \in \Fil^{j}(\Tan_z(S_J)_\tau)$ for $j = 1,\ldots,e_\gp-1$, as well as $j=e_\gp$ if
$\theta_{\gp,i,e_\gp} \not\in J$.  Indeed if $\widetilde{z} \in \Fil^{j-1}(\Tan_z(S_J)_\tau)$ and
$\theta_{\gp,i,j} \not\in J$, then since $t(\widetilde{z}) \in \Fil^j(\Tan_y(P_J)_\tau)$,
parts 1) and 4) of Lemma~\ref{lem:graded} imply that $\widetilde{z} \in \Fil^j(\Tan_z(S_J)_\tau)$.
On the other hand if $\widetilde{z} \in \Fil^{j-1}(\Tan_z(S_J)_\tau)$, $\theta_{\gp,i,j} \in J$ and
$j < e_\gp$ then since $t(\widetilde{z}) \in \Fil^{j+1}(\Tan_y(P_J)_\tau)$,
parts 3) and 4) of Lemma~\ref{lem:graded} imply that $\widetilde{z} \in \Fil^j(\Tan_z(S_J)_\tau)$.

This completes the proof of the description of the kernel of $t_\tau$, which immediately implies the formula for
the dimension of $\ker(t)$.  

The fact that the image of $t_\tau$ is contained in the described space is also immediate from 
Lemma~\ref{lem:graded}.  (Note that the first two cases in the description are not exclusive, but
are consistent since $\gr^1(\Tan_y(P_J)_\tau) = 0$ if $\theta_{\gp,i,1} \in J'$.)  Equality then follows
on noting that their direct sum over all $\tau$ has dimension equal to the rank of $t$.
\epf

Now let $Z$ denote the fibre at $x$ of $\pi_J: S_J \to T_{J'}$, defined as the restriction of $\widetilde{\pi}_{1,k}$,
and consider the tangent space $\Tan_z(Z)$, which may be identified with the kernel of the $\FF$-linear map
$\Tan_z(S_J) \to \Tan_x(T_{J'})$.  We denote the map by $s$ and write it as the composite 
$$\Tan_z(S_J) \stackrel{t}{\longrightarrow} \Tan_y(P_J)   \longrightarrow \Tan_x(T_{J'}),$$
where the second map is the natural projection. Since each map respects the decomposition
over $\tau \in \Sigma_0$, we may similarly decompose $\Tan_z(Z) = \bigoplus_{\tau} \Tan_z(Z)_\tau$
where $\Tan_z(Z)_\tau$ is the kernel of the composite 
$$s_\tau:  \Tan_z(S_J)_\tau \longrightarrow \Tan_y(P_J)_\tau   \longrightarrow \Tan_x(T_{J'})_\tau.$$
We may furthermore consider the filtration on $\Tan_z(Z)_\tau$ induced by the one on $\Tan_z(S_J)_\tau$,
so that if $z$ correspond to $(\underline{A},\underline{A}',\psi)$, then 
$\bigoplus_{\tau}\Fil^{j_\tau}(\Tan_z(Z)_\tau)$ is identified with the set of its lifts 
$(\underline{\widetilde{A}},\underline{\widetilde{A}}', \widetilde{\psi})$ to $S_J(\FF[\epsilon])$ 
such that $\underline{\widetilde{A}} = \underline{A}\otimes_{\FF} \FF[\epsilon]$ and
$\widetilde{F}_\tau^{\prime(j)}$ corresponds to $F_\tau^{(j)}\otimes_{\FF}\FF[\epsilon]$
for all $\tau \in \Sigma_0$ and $j \le j_\tau$.  (Note that this condition is automatic for $\tau \not\in \Sigma_\gP$,
in which case $\Tan_z(Z)_\tau = 0$.)

\begin{lemma}  \label{lem:tangents}  The $\FF$-vector space $\Tan_z(Z)$ has dimension $|J'| + \delta 
= |J''| + \delta$, where $\delta = \#\{\,\tau = \tau_{\gp,i}\,|\, \theta_{\gp,i,e_\gp},\theta_{\gp,i+1,1} \in J\,\}$.
More precisely, if $\tau = \tau_{\gp,i}$ and $\theta = \theta_{\gp,i,j}$, then
$$\dim_{\FF} \gr^{j-1}(\Tan_z(Z)_\tau) = \left\{\begin{array}{ll} 1,& \mbox{if $\theta = \theta_{\gp,i,j} \in J$ and either
   $\sigma\theta \not\in J$ or $j = e_\gp$;} \\
   0,&\mbox{otherwise.}\end{array}\right.$$
\end{lemma}
\begpf  It follows from the description of the image of $t_\tau$ in Lemma~\ref{lem:rank} that if $\tau = \tau_{\gp,i}$, then
$s_\tau$ is surjective unless $\theta_{\gp,i-1,e_\gp}, \theta_{\gp,i,1} \in J$, in which case the image of $s_\tau$ is
$\Fil^1(\Tan_x(T_{J'})_\tau)$, which has codimension one in $\Tan_x(T_{J'})_\tau$.  Therefore the image of $s$ has
codimension $\delta$ in $\Tan_x(T_{J'})$, and the kernel of $s$ has dimension
$$\dim_{\FF}(\Tan_z(S_J)) - \dim_{\FF}(\Tan_x(T_{J'})) + \dim_{\FF}(\coker s) = |J'| + \delta.$$

Since $\dim_{\FF} \gr^{j-1}(\Tan_z(Z)_\tau) \le \dim_{\FF} \gr^{j-1}(\Tan_z(S_J)_\tau) = 1$ for all $\tau = \tau_{\gp,i}$ and
$j = 1,\ldots,e_\gp$, and $\delta$ is the total number of $\tau$ and $j$ such that if $\theta = \theta_{\gp,i,j} \in J$ and either
$\sigma\theta \not\in J$ or $j = e_\gp$, it suffices to prove that $\gr^{j-1}(\Tan_z(Z)_\tau) = 0$ if either $\theta \not \in J$, or
$j < e_\gp$ and $\theta_{\gp,i,j}, \theta_{\gp,i,j+1} \in J$.  This is immediate from parts 1), 3) and 4) of Lemma~\ref{lem:graded},
which show that the following maps induced by $s_\tau$ are isomorphisms:
\begin{itemize}
\item $\gr^{j-1}(\Tan_z(S_J)_\tau) \longrightarrow \gr^{j-1}(\Tan_y(P_J)_\tau) \stackrel{\sim}{\longrightarrow} \gr^{j-1}(\Tan_x(T_{J'})_\tau)$
 in the case $\sigma^{-1}\theta,\theta\not\in J$;
\item $\gr^{j-1}(\Tan_z(S_J)_\tau) \longrightarrow \gr^{j}(\Tan_y(P_J)_\tau) \stackrel{\sim}{\longrightarrow} \gr^{j}(\Tan_x(T_{J'})_\tau)$ in the case $\theta,\sigma\theta\in J$, $j < e_\gp$;
\item $\gr^{j-1}(\Tan_z(S_J)_\tau) \longrightarrow V_\theta \stackrel{\sim}{\longrightarrow} \gr^{j-1}(\Tan_x(T_{J'})_\tau)$ in the case $\theta\in J''$.
\end{itemize}

\subsection{Unobstructed deformations}  \label{sec:unobstructed}
We maintain the notation of the preceding section, so that in particular $J$ is a subset of $\Sigma_\gP$, 
$J' = \{\,\theta\in J\,|\, \sigma^{-1}\theta\not\in J\,\}$, and $Z$ is the fibre of $\pi_J: S_J \to T_{J'}$
at a geometric point $x \in T_{J'}(\FF)$.  Following on from the description of the tangent space
in Lemma~\ref{lem:tangents}, we continue the study of the local structure of $Z$ at each closed point $z$ in order to
prove a generalization of \cite[Lemma~7.1.6]{DKS}, showing in particular that the reduced subscheme is smooth.
We remark however that if $\gP$ is divisible by a prime $\gp$ ramified in $F$, then the fibres over $T_{J'}$ of $\widetilde{\xi}_J$ need not be
totally inseparable, so the desired conclusions do not immediately follow as in \cite{DKS} from the commutative algebra results of \cite{KN}
and the description of the morphisms on tangent spaces.  We therefore resort to a more detailed analysis of the deformation
theory.

Fix $m \ge 0$ and let $R = \FF[[X_1,\ldots,X_m]]$ and $R_n = R/\gm^{n+1}$ for $n \ge 0$, where 
$\gm = \gm_{R} = (X_1,\ldots,X_m)$.   Similarly let $\widetilde{R} = W_2[[X_1,\ldots,X_m]]$
and $\widetilde{R}_{n} = \widetilde{R}/I^{n+1}$ where $W_2 = W_2(\FF)$ and $I = (X_1,\ldots,X_m)$,
so that $\widetilde{R}_{n}$ is flat over $W_2$ and $\widetilde{R}_{n}/p\widetilde{R}_{n} = R_{n}$.  We let 
$\widetilde{\phi} =  \widetilde{\phi}_{n}: \widetilde{R}_{n} \to \widetilde{R}_{n}$ denote the lift of the absolute 
Frobenius  $\phi = \phi_{n}$ on $R_{n}$ defined by $\widetilde{\phi}(X_j) = X_j^p$ for $j=1,\ldots,m$.  
We view $\widetilde{R}_{n+1}$ (as well as $\widetilde{R}_n$ and $R_{n+1}$) as a divided power thickening of $R_{n}$
with divided powers defined by
$$(pf+g)^{[i]} = p^{[i]}f^i + p^{[i-1]}f^{i-1}g$$
for $f \in \widetilde{R}_{n+1}$, $g \in I^{n+1}$ and $i\ge 1$.\footnote{Note that $p^{[i]} = 0$ in $W_2$ if $i>1$ and $p > 2$.}

We now denote the data corresponding to $x \in T_{J'}(\FF)$ by $\underline{A}_0$, take $(\underline{A}_n,\underline{A}_n',\psi_n)$
corresponding to an element $z_n \in Z(R_{n})$, so that $\underline{A}_n = \underline{A}_0 \otimes_{\FF} R_n$, and let
$F_n^{\bullet} = F_{0}^{\bullet} \otimes_{\FF} R_n$ and $F_{n}^{\prime\bullet}$ denote the associated Pappas--Rapoport
filtrations. Suppose that $\tau = \tau_{\gp,i}$ is such
that $\theta = \theta_{\gp,i,e_\gp} \in J$, and consider the morphisms
\begin{equation}\label{eqn:crysn}  \xymatrix{ 
H^1_\crys(A_n'/\widetilde{R}_n)_\tau \ar[r]^{\widetilde{\psi}_{\tau,n}^*}
& H^1_\crys(A_n/\widetilde{R}_n)_\tau \\
& H^1_\crys(A_n^{(p^{-1})}/\widetilde{R}_n)_\tau 
\ar[u]_{\widetilde{V}_{\tau,n}} }
\end{equation}
of free rank two $\widetilde{R}_n[u]/(E_\tau)$-modules,
where $A_n^{(p^{-1})} : = A_0^{(p^{-1})}\otimes_{\FF} R_n$ and $\widetilde{V}_{\tau,n}$ is induced by $\Ver:A_0^{(p^{-1})} \to A_0$.
Note that $H^1_\crys(A_n/\widetilde{R}_n)_\tau = H^1_\crys(A_0/W_2)_\tau \otimes_{W_2} \widetilde{R}_n$ and
$H^1_\crys(A_n^{(p^{-1})}/\widetilde{R}_n)_\tau = H^1_\crys(A_0^{(p^{-1})}/W_2)_\tau \otimes_{W_2} \widetilde{R}_n
 = H^1_\crys(A_0/W_2)^{\phi^{-1}}_{\phi\circ\tau} \otimes_{W_2} \widetilde{R}_n$, identifying $\widetilde{V}_{\tau,n}$ with 
 $\widetilde{V}_{\tau,0} \otimes_{W_2} \widetilde{R}_n$.  Note in particular that $\widetilde{V}_{\tau,n}$ has kernel and cokernel annihilated by $p$,
 and its image is $\widetilde{F}_{\tau,0}^{(e_\gp)} \otimes_{W_2} \widetilde{R}_n$, where $\widetilde{F}_{\tau,0}^{(e_\gp)}$
 is the preimage in $H^1_\crys(A_0/W_2)_\tau$ of $F_{\tau,0}^{(e_\gp)} = H^0(A_0,\Omega^1_{A_0/\FF})$ under the
 canonical projection to $H^1_\crys(A_0/\FF)_\tau \cong H^1_\dr(A_0/\FF)_\tau$.  On the other hand,
 the kernel and cokernel of $\widetilde{\psi}_{\tau,n}^*$ are annihilated by $u$, and the assumption that $\theta_{\gp,i,e_\gp} \in J$
 implies that $\psi_n$ induces a surjection 
 $$G_{\tau,n}^{\prime(e_\gp)} = u^{-1}F_{\tau,n}^{\prime(e_\gp-1)} \longrightarrow F_{\tau,0}^{(e_\gp)} \otimes_{\FF} R_n.$$
 It follows that $\widetilde{\psi}_{\tau,n}^*$ induces a surjection
  $$\widetilde{G}_{\tau,n}^{\prime(e_\gp)} \longrightarrow \widetilde{F}_{\tau,0}^{(e_\gp)} \otimes_{W_2} \widetilde{R}_n,$$
where $\widetilde{G}_{\tau,n}^{\prime(e_\gp)}$
 is the preimage in $H^1_\crys(A'_n/\widetilde{R}_n)_\tau$ of $G_{\tau,n}^{\prime(e_\gp)}$ under the
 canonical projection to $H^1_\crys(A'_n/R_n)_\tau \cong H^1_\dr(A'_n/R_n)_\tau$.   We therefore obtain
 from (\ref{eqn:crysn}) a commutative diagram of $\widetilde{R}_n[u]/(E_\tau)$-linear isomorphisms
 $$ \xymatrix{
  \widetilde{G}_{\tau,n}^{\prime(e_\gp)}/\ker(\widetilde{\psi}_{\tau,n}^*)\ar[r]^{\sim} \ar[rd]^{\sim}
   &  \widetilde{F}_{\tau,0}^{(e_\gp)} \otimes_{W_2} \widetilde{R}_n \\
&(H^1_\crys(A_0/W_2)^{\phi^{-1}}_{\phi\circ\tau} /\ker(\widetilde{V}_{\tau,0}) ) \otimes_{W_2}  \widetilde{R}_n. \ar[u]_{\wr} }
 $$
 Furthermore since $\ker(\widetilde{\psi}_{\tau,n}^*) \subset p H^1_\crys(A'_n/R_n)_\tau \subset u \widetilde{G}_{\tau,n}^{\prime(e_\gp)}$,
 we in turn obtain $R_n$-linear isomorphisms
 $$ \xymatrix{
 \widetilde{G}_{\tau,n}^{\prime(e_\gp)}/ u \widetilde{G}_{\tau,n}^{\prime(e_\gp)}\ar[r]^{\sim} \ar[rd]^{\sim}&
    (\widetilde{F}_{\tau,0}^{(e_\gp)}/ u\widetilde{F}_{\tau,0}^{(e_\gp)}) \otimes_{\FF} R_n \\ 
&(H^1_\crys(A_0/W_2)^{\phi^{-1}}_{\phi\circ\tau}  /  u H^1_\crys(A_0/W_2)^{\phi^{-1}}_{\phi\circ\tau} ) \otimes_{\FF}{R_n}.\ar[u]_{\wr}
   }$$
Finally noting that the projection $\widetilde{G}_{\tau,n}^{\prime(e_\gp)} \to G_{\tau,n}^{\prime(e_\gp)}$
identifies $ \widetilde{G}_{\tau,n}^{\prime(e_\gp)}/ u \widetilde{G}_{\tau,n}^{\prime(e_\gp)}$ with
$P'_{\theta,n} := {G}_{\tau,n}^{\prime(e_\gp)}/ F_{\tau,n}^{\prime(e_\gp - 1)}$, and that $u^{e_\gp-1}$
induces 
$$\begin{array}{rl} H^1_\crys(A_0/W_2)^{\phi^{-1}}_{\phi\circ\tau}  /  u H^1_\crys(A_0/W_2)^{\phi^{-1}}_{\phi\circ\tau} &
 = H^1_\dr(A_0/\FF)^{(p^{-1})}_{\phi\circ\tau}  /  u H^1_\dr(A_0/\FF)^{(p^{-1})}_{\phi\circ\tau} \\
 \stackrel{\sim}{\longrightarrow} &H^1_\dr(A_0/\FF)^{(p^{-1})}_{\phi\circ\tau}  [u]  = P_{\sigma\theta,0}^{(p^{-1})}\end{array}$$
(where $P_{\sigma\theta,0} = G_{\phi\circ\tau,0}^{(1)}$), we obtain from (\ref{eqn:crysn}) an isomorphism
$$   \varepsilon_\tau(z_n) :  P'_{\theta,n}  \stackrel{\sim}{\longrightarrow} P_{\sigma\theta,0}^{(p^{-1})} \otimes_{\FF} R_n $$
of free rank two $R_n$-modules.  

We say that $z_n$ is {\em unobstructed} if $\varepsilon_\tau(z_n)$ sends 
$L'_{\theta,n}$ to $L_{\sigma\theta,0}^{(p^{-1})} \otimes_{\FF} R_n$ for all $\tau = \tau_{\gp,i}$ such that
$\theta = \theta_{\gp,i,e_\gp}$ and $\sigma\theta = \theta_{\gp,i+1,1}$ are both in $J$, where
as usual $L'_{\theta,n} = F_{\tau,n}^{\prime(e_\gp)}/F_{\tau,n}^{\prime(e_\gp-1)}$ and
$L_{\sigma\theta,0} = F_{\phi\circ\tau,0}^{(1)}$.  Note that the condition on $\varepsilon_\tau(z_n)$ is equivalent to the
vanishing of the induced homomorphism
$$ L'_{\theta,n}  \longrightarrow (P_{\sigma\theta,0}/L_{\sigma\theta,0})^{(p^{-1})} \otimes_{\FF} R_n$$
of free rank one $R_n$-modules.    We let $Z(R_n)^\uno$ denote the set of unobstructed elements
in $Z(R_n)$.

Note firstly that every element $z = z_0$ of $Z(\FF)$ is unobstructed.  Indeed the condition that $\sigma\theta\in J$
implies that the image of 
$$\widetilde{\psi}^*_{\phi\circ\tau}:  H^1_\crys(A_0'/W_2) \longrightarrow H^1(A_0/W_2)$$
is $u^{1-e_\gp} \widetilde{F}_{\phi\circ\tau,0}^{(1)}$, so the claim follows from the commutativity of the
diagram
$$ \xymatrix{ H^1_\crys(A_0'/W_2)_{\phi\circ\tau}^{\phi^{-1}} \ar[rr]^{(\widetilde{\psi}^*_{\phi\circ\tau})^{\phi^{-1}}}
     \ar[d]_{\widetilde{V}_{\tau,0}'} && H^1_{\crys}(A_0/W_2)_{\phi\circ\tau}^{\phi^{-1}} \ar[d]^{\widetilde{V}_{\tau,0}} \\
     H^1_\crys(A_0'/W_2)_\tau \ar[rr]^{\widetilde{\psi}^*_\tau}  && H^1_\crys(A_0/W_2)_\tau.}$$
Furthermore the construction of $\varepsilon_\tau(z_n)$ is functorial in $R_n$ for varying $m$ and $n$,
in the sense that if
$$\begin{array}{cccc}
\alpha:  & \FF[[X_1,\ldots,X_m]]/(X_1,\ldots,X_m)^{n+1} & \longrightarrow &
\FF[[X_1,\ldots,X_{m'}]]/(X_1,\ldots,X_{m'})^{n'+1}\\
& \parallel && \parallel \\
& R_n & & R'_{n'} \end{array} $$
is an $\FF$-algebra morphism admitting a lift $\widetilde{\alpha}:\widetilde{R}_n \to \widetilde{R}'_{n'}$,
then $\varepsilon_\tau(\alpha^*z_n) = \alpha^*\varepsilon_\tau(z_n)$, so that
$\alpha^*$ sends $Z(R_n)^\uno$ to $Z(R'_{n'})^\uno$.\footnote{We could not directly define $Z(R_n)^\uno$ as the set of $R_n$-points of a closed subscheme
since the construction of $\varepsilon_\tau(z_n)$ relied on the existence and choice of flat $W_2$-lifts.  An alternative approach
would be to prove a suitable generalization of \cite[Lemma~7.2.1]{DKS} (the ``crystallization lemma''), and use it together with the flatness
of $\widetilde{\xi}_J$ to define $\varepsilon_\tau$ over $Z_\red$, without knowing a priori that it is smooth.  We opted instead here
for a bootstrapping argument that will allow us first to establish the smoothness of  $Z_\red$ so that we can apply a more direct
generalization of the crystallization lemma  proved in \cite{DKS}.}
In particular, if $z_1 \in \Tan_z(Z)$ (viewed as an element of $Z(R_1)$, with $R_1 = \FF[X]/(X^2)$), then 
$\varepsilon_\tau(z_1) - \varepsilon_\tau(z) \otimes 1$ induces an $\FF$-linear map
$$L_{\theta,0}' \cong
  L_{\theta,1}'/XL_{\theta,1}' \longrightarrow
  (P_{\sigma\theta,0}^{(p^{-1})}/L_{\sigma\theta,0}^{(p^{-1})})
   \otimes_{\F} XR_1.$$
Writing this map as $\partial_\tau(z_1) \otimes X$ for a unique
$\partial_\tau(z_1) \in \Hom_{\FF}  (L_{\theta,0}',  (P_{\sigma\theta,0}/L_{\sigma\theta,0})^{(p^{-1})} )$,
the condition on $\varepsilon_\tau(z_1)$ in the definition of unobstructedness is then equivalent to the vanishing of $\partial_\tau(z_1)$.

More generally, if $n=1$ and $m$ is arbitrary, then to give an element $z_1 \in Z(R_1)$ over $z = z_0 \in Z(\FF)$ is equivalent to giving
an $m$-tuple of tangent vectors $(z^{(1)}_1,\ldots,z_1^{(m)})$ in $\Tan_z(Z)$, or more functorially an element
of $\Hom_{\FF}(\FF^m,\Tan_z(Z))$, the induced morphism $\CO_{Z,z}^{\wedge} \to \CO_{Z,z}/\gm_z^2 \to R_1$ being surjective
if and only if the $m$-tuple is linearly independent. Furthermore any $\FF$-algebra morphism $\alpha:R_1 \to R_1'$ (for arbitrary $m,m'$)
admits a lift $\widetilde{\alpha}:\widetilde{R}_1 \to \widetilde{R}'_1$, so the function
\begin{equation}  \label{eqn:obstruction}
\partial_\tau: \Tan_z(Z)   \longrightarrow  \Hom_{\FF}  (L_{\theta,0}',  (P_{\sigma\theta,0}/L_{\sigma\theta,0})^{(p^{-1})} ) \end{equation}
is $\FF$-linear.  (Let $z_1 \in Z(R_1)$ correspond to a basis $(z^{(1)}_1,\ldots,z_1^{(m)})$ of $\Tan_z(Z)$, and
realize an arbitrary tangent vector $\sum a_i z_1^{(i)}$ as the image of $z_1$ under $X_i \mapsto a_i X$.)
The set of unobstructed elements in $\Tan_z(Z)$ is thus the kernel of the resulting $\FF$-linear map 
$$\Tan_z(Z) \longrightarrow \bigoplus_{\tau \in \Delta}
 \Hom_{\FF} (L_{\theta,0}',  (P_{\sigma\theta,0}/L_{\sigma\theta,0})^{(p^{-1})} ),$$
where $\Delta = \{\,\tau = \tau_{\gp,i}\,|\, \theta_{\gp,i,e_\gp},\theta_{\gp,i+1,1} \in J\,\}$, and is therefore a vector
space of dimension at least $|J'|$ by Lemma~\ref{lem:tangents}.  The functoriality of $\varepsilon_\tau$ also implies that an element $z_1 \in Z(R_1)$
corresponding to an $m$-tuple $(z_1^{(1)},\ldots,z_1^{(m)})$ is unobstructed if and only if each of the $z_1^{(i)}$ is,
so it follows that there is a surjective $\FF$-algebra homomorphism
\begin{equation}\label{eqn:launch}   \widehat{\CO}_{Z,z} \longrightarrow R_1 = \FF[[X_1,\ldots,X_{|J'|}]]/(X_1,\ldots,X_{|J'|})^2 \end{equation}
corresponding to an element $z_1 \in Z(R_1)^\uno$.

We shall in fact need the following more precise description of the set of unobstructed tangent vectors.
\begin{lemma}  \label{lem:unofirst}  If $z \in Z(\FF)$, then $\Tan_z(Z)^\uno = \oplus_{\tau} \Tan_z(Z)^\uno_\tau$,
where $\Tan_z(Z)^\uno$ (resp.~$\Tan_z(Z)_\tau^\uno$) is the set of unobstructed
elements of $\Tan_z(Z)$ (resp.~$\Tan_z(Z)_\tau$).  Furthermore, for each $\tau = \tau_{\gp,i}$, $\theta = \theta_{\gp,i,j}$, we have
$$\dim_{\FF} \gr^{j-1}(\Tan_z(Z)^\uno_\tau) = \left\{\begin{array}{ll} 1,& \mbox{if $\sigma\theta \in J''$;} \\
   0,&\mbox{otherwise,}\end{array}\right.$$
with respect to the filtration induced by the one on $\Tan_z(Z)_\tau$, so $\dim_\FF(\Tan_z(Z)^{\uno}) = |J'| = |J''|$.
\end{lemma}
\begpf
The assertion concerning the decomposition follows from the vanishing of the restriction of $\partial_\tau$ to
$ \Tan_z(Z)_{\tau'}$ for $\tau \in \Delta$ and $\tau' \neq \tau$,
where $\partial_\tau$ is defined in (\ref{eqn:obstruction}).  
The assertion concerning the
filtration then follows from Lemma~\ref{lem:tangents} and the injectivity of $\partial_\tau$ on
$ \Fil^{e_\gp-1}(\Tan_z(Z)_{\tau})$.
The desired vanishing and injectivity are both consequences of the claim that if $\tau \in \Delta$, $R_1 = \FF[X]/(X^2)$ and 
$z_1 \in Z(R_1)$ is a lift
of $z$ such that $F_{\tau,1}^{\prime(e_\gp-1)}$ corresponds to $F_{\tau,0}^{\prime(e_\gp-1)} \otimes_{\FF} R_1$ under
the canonical isomorphism $H^1_\dr(A_1'/R_1)_\tau \cong H^1_\crys(A_0'/R_1)_\tau \cong H^1_\dr(A_0'/\FF) \otimes_{\FF} R_1$,
then $F_{\tau,1}^{\prime(e_\gp)}$ corresponds to $F_{\tau,0}^{\prime(e_\gp)} \otimes_{\FF} R_1$ if and only if $\partial_\tau(z_1) = 0$, i.e.,
$L'_{\theta,1}$ corresponds to $L_{\sigma\theta,0}^{(p^{-1})}\otimes_{\FF} R_1$ under $\varepsilon_\tau(z_1)$.
Indeed if $z_1 \in \Tan_z(Z)_{\tau'}$ for some $\tau' \neq \tau$,
then $F_{\tau,1}'^{(j)}$ corresponds to $F_{\tau,0}^{\prime(j)} \otimes_{\FF} R_1$ 
for $j=0,\ldots,e_\gp$, and if $z_1 \in \Fil^{e_\gp -1}\Tan_z(Z)_{\tau}$,
then $F_{\tau,1}'^{(e_\gp-1)}$ corresponds to $F_{\tau,0}^{\prime(e_\gp-1)} \otimes_{\FF} R_1$, 
and $z_1$ is trivial if and only if
$F_{\tau,1}'^{(e_\gp)}$ corresponds to $F_{\tau,0}^{\prime(e_\gp)} \otimes_{\FF} R_1$

The claim follows from the commutativity of the diagram:
$$\xymatrix{P_{\theta,0}'\otimes_{\FF} R_1 \ar[rr]^-{\sim} \ar[rd]^-{\sim}_-{\varepsilon_\tau(z_0) \otimes R_1} 
&& P_{\theta,1}' \ar[ld]_-{\sim}^-{\varepsilon_\tau(z_1)} \\
& P_{\sigma\theta,0}^{(p^{-1})} \otimes_{\FF} R_1}$$
(where the top arrow is the canonical isomorphism), which in turn follows from that of
$$\xymatrix{\widetilde{G}_{\tau,0}^{\prime(e_\gp)} \otimes_{W_2} \widetilde{R}_1 \ar[rr]^-{\sim} \ar[rd]_-{\widetilde{\psi}_{\tau,0}^* \otimes \widetilde{R}_1} 
&& \widetilde{G}_{\tau,1}^{\prime(e_\gp)}  \ar[ld]^-{\widetilde{\psi}_{\tau,1}^*} \\
& \widetilde{F}_{\tau,0}^{(e_\gp)} \otimes_{W_2} \widetilde{R}_1,}$$
where the top arrow is induced by the morphisms $(\FF,W_2) \to (\FF,\widetilde{R}_1) \leftarrow (R_1,\widetilde{R}_1)$ of
divided power thickenings.  This in turn is implied by the commutativity of
$$\xymatrix{H^1_\crys(A_0'/W_2)_\tau \otimes_{W_2} \widetilde{R}_1 \ar[r]^-{\sim}  \ar[d]_-{\widetilde{\psi}_{\tau,0}^* \otimes \widetilde{R}_1} &
H^1_\crys(A_1'/\widetilde{R}_1)_\tau \ar[d]^-{\widetilde{\psi}_{\tau,1}^*} \\
H^1_\crys(A_0/W_2)_\tau \otimes_{W_2} \widetilde{R}_1 \ar[r]^-{\sim}  &
H^1_\crys(A_1/\widetilde{R}_1)_\tau ,}$$
which in turn results from $\widetilde{\psi}_{\tau,0}^*$ and $\widetilde{\psi}_{\tau,1}^*$ being evaluations of the morphism of crystals over $R_1$
induced by $\psi_1$.
\epf

Suppose now that $m$ is fixed (for example as $|J'|$) and consider the thickening $R_{n+1} \to R_n$ (with trivial divided powers), and consider the natural map $Z(R_{n+1})^\uno \to  Z(R_n)^\uno$.
We will prove (Lemma~\ref{lem:unorest} below) that this map is surjective, explaining why we call such deformations {\em unobstructed}.  Before doing so, however, we illustrate the construction of the desired liftings with some examples.

Suppose first that $[F:\Q] =2$.  If $p$ splits in $F$, then 
$J' = J'' = \emptyset$ and $\Delta = J$.  Therefore Lemma~\ref{lem:unofirst} implies that all unobstructed deformations are trivial (i.e., $Z(R_n)^\uno = Z(\F)$), so
there is nothing to prove.  The same holds if 
$p$ is inert in $F$, and either $J = \Sigma$ or $\emptyset$.
More explicitly, it is well-known that in the preceding cases, $Z$ is isomorphic to $(\Spec(\FF[X]/(X^p)))^\delta$ (see for example \cite[(43)]{DKS}), so $Z(R_n)^\uno$ consists of a single element.

Suppose on the other hand that $p$ is inert in $F$ and $J = \{\tau_0\}$, where $\gP = p\CO_F$ and 
$\Sigma = \{\,\tau_0,\tau_1\,\}$, so $J' = \{\tau_0\}$, $J' = \{\tau_1\}$ and $\Delta = \emptyset$.  In particular, the condition in the definition of unobstructed is vacuous, so $Z(R_n)^\uno = Z(R_n)$ for all $n \ge 1$.  Therefore the surjectivity of 
$Z(R_{n+1})^{\uno} \to Z(R_n)^{\uno}$ in this case amounts to the smoothness of $Z$, which follows for example from 
\cite[Lemma 7.1.4.1]{DKS} (or alternatively, from Lemmas~\ref{lem:Frobfact} and~\ref{lem:tangents} above).
To illustrate the proof of Lemma~\ref{lem:unorest} more explicitly in this case, suppose that $z_n \in Z(R_n)$ and let
$(\underline{A}_n,\underline{A}'_n,\psi_n)$ be the corresponding tuple.  So in particular $\underline{A}_n = \underline{A}_0 \otimes_{\FF} R_n$ and $\psi_n:A_n \to A'_n$ is an isogeny such that the morphisms
$$P'_{0,n}/L'_{0,n} \longrightarrow P_{0,n}/L_{0,n}
\quad\mbox{and}\quad L_{1,n}' \longrightarrow L_{1,n}$$
induced by $\psi_n$ are trivial, where 
$$P'_{i,n} = H^1_\dr(A_n'/R_n')_{\tau_i}\quad\mbox{and}\quad
P_{i,n} = H^1_\dr(A_n/R_n) = H^1_\dr(A_0/\F)_{\tau_i}\otimes_{\F} R_n$$ for $i=0,1$, and similarly
$L'_{i,n} = H^0(A_n',\Omega^1_{A_n'/R_n})$ and $L_{i,n} = H^0(A_0,\Omega^1_{A_0/\F})\otimes_{\F} R_n$.
By the Grothendieck--Messing Theorem, lifts of $\underline{A}'_n$ to $R_{n+1}$ correspond to lifts of the 
$$L'_{i,n} \subset P'_{i,n} \cong H^1_\crys(A'_n/R_n)_{\tau_i}
 \cong H^1_\crys(A'_n/R_{n+1})_{\tau_i}\otimes_{R_{n+1}}R_n$$ to 
free rank one $R_{n+1}$-submodules 
$L'_{i,n+1}$ of $P'_{i,n+1}:= H^1_\crys(A_n'/R_{n+1})_{\tau_i}$ for $i=0,1$.
Furthermore $\psi_n$ lifts to an isogeny
$\psi_{n+1}:A_{n+1} = A_0\otimes R_{n+1} \to A'_{n+1}$ (necessarily unique and compatible with
the auxiliary data) if and only if the morphisms
$$\psi_{i,n+1}^*: P'_{i,n+1} \longrightarrow
P_{i,n+1}:=   H^1_\crys(A_n/R_{n+1})_{\tau_i} \cong P_{i,0}   \otimes_{\F}R_{n+1}$$
induced by $\psi_n$ send $L'_{i,n+1}$ to $L_{i,n+1} = L_{i,0}\otimes_{\F} R_{n+1}$ for $i=0,1$.  Moreover the resulting point of $Y_{U_0(\gP)}(R_{n+1})$ lies in $S_J(R_{n+1})$, and hence in $Z(R_{n+1})$, if and only if 
$$\psi_{0,n+1}^*(P'_{0,n+1}) \subset L_{0,n+1}
 \quad\mbox{and}\quad \psi_{1,n+1}^*(L'_{1,n+1}) = 0.$$
Choosing an arbitrary lift $L'_{0,n+1}$ of $L'_{0,n}$ and setting $L'_{1,n+1} = \ker(\psi_{1,n+1}^*)$ ensures these conditions are satisfied (see the first part of the proof of Lemma~\ref{lem:unorest} for the former).  Note that these imply that $\psi_{i,n+1}^*(L'_{i,n+1}) \subset L_{i,n+1}$ for $i=0,1$, and hence supply the desired lift of $z_n$ to an element $z_{n+1} \in Z(R_{n+1})$

Continue to assume that $[F:\Q] = 2$, but suppose now that $p$ is ramified in $F$.  Let $\gP$ be the unique prime over $p$, and write $\Sigma_0=\{\tau\}$ and $\Sigma = \{\theta_1,\theta_2\}$.  If $J = \Sigma$ or $\emptyset$, then we again have $J' = J'' = \emptyset$ (and $\Delta = \Sigma_0$ or $\emptyset$, according to whether $J = \Sigma$ or $\emptyset$).  Again it follows from Lemma~\ref{lem:tangents} that all unobstructed lifts are trivial, so there is nothing to prove.  

On the other hand if $J = \{\theta_1\}$ or $\{\theta_2\}$, then $\Delta = \emptyset$, so again all lifts are unobstructed and the desired surjectivity amounts to the smoothness of $Z$, which can again be deduced from Lemmas~\ref{lem:Frobfact} and~\ref{lem:tangents}.  To see the surjectivity more explicitly, we can proceed to construct a lift of $z_n \in Z(R_n)$ to $Z(R_{n+1})$ as in the inert case, the task at hand now being to
lift the Pappas--Rapoport filtration:
$$ 0 \subset F_n'^{(1)} \subset F_n'^{(2)}
  = H^0(A_n',\Omega^1_{A_n'/R_n}) \subset H^1_\dr(A_n'/R_n)
  \cong H^1_\crys(A_n'/R_n)$$
to one in $H^1_\crys(A_n'/R_{n+1})$ with the desired properties.
In particular, if $J = \{\theta_1\}$, then we require both
$G'^{(1)}_{n+1} = u H^1_\crys(A_n'/R_{n+1})$ and 
$F'^{(2)}_{n+1}$ be sent to $F_0^{(1)} \otimes_\F R_{n+1}$ under
the induced map 
$$\psi_{n+1}^*:H^1_\crys(A_n'/R_{n+1}) \to H^1_\crys(A_n/R_{n+1}) \cong H^1_\dr(A_0/\F)\otimes_{\F} R_{n+1}.$$
Note that in this case $J'' = \{\theta_2\}$, so
$F'^{(2)}_n = u H^1_\dr(A_n'/R_n)$; furthermore the 
first part of the proof of Lemma~\ref{lem:unorest} shows that 
the image of $\psi_{n+1}^*$ is $u^{-1}F_0^{(1)} \otimes_\F R_{n+1}$.  So we obtain a lift with the desired properties by setting $F'^{(2)}_{n+1} = u H^1_\dr(A_n'/R_{n+1})$, within which we choose an arbitrary lift $F'^{(1)}_{n+1}$ of $F'^{(1)}_n$.
On the other hand, if $J = \{\theta_2\}$, then we require
$\psi_{n+1}^*$ to send $F'^{(1)}_{n+1}$ to $0$ and 
$G'^{(2)}_{n+1} = u^{-1} F'^{(1)}_{n+1}$ to 
$F^{(2)}_0\otimes_{\F} R_{n+1}$.
Note that in this case
$J' = \{\theta_2\}$, so $F^{(2)}_0 = u H^1_\dr(A_0/\F)$, 
and we obtain a lift with the desired properties by setting
$F'^{(1)}_{n+1} = \ker(\psi_{n+1}^*)$ and defining 
$F'^{(2)}_{n+1}$ by arbitrarily lifting
$L_n' := F'^{(2)}_n/F'^{(1)}_n \subset P'_n = G'^{(2)}_n/F'^{(1)}_n$ to a line $L_{n+1}'$ in 
$P'_{n+1} = G'^{(2)}_{n+1}/F'^{(1)}_{n+1}$, where 
$G'^{(2)}_\ell = u^{-1}F'^{(1)}_\ell$ for $\ell = n,n+1$.

The notion of unobstructedness can only play a more serious role if $[F:\Q] \ge 3$.  Suppose for example that $[F:\Q] = 4$ and $S_p = \{\gp\}$ with $e_\gp = f_\gp = 2$.  Write $\Sigma_0 = \{\,\tau_0,\tau_1\}$ and $\Sigma = \{\,\theta_{0,1},\theta_{0,2},\theta_{1,1},\theta_{1,2}\}$,
and consider the case $J = \{\theta_{0,1}, \theta_{1,1}, \theta_{1,2}\}$, so that
$$J' = \{\theta_{1,1}\},\quad J'' = \{\theta_{0,2}\}\quad\mbox{and}\quad \Delta = \{\tau_1\}.$$
Let $(\underline{A}_n,\underline{A}_n',\psi_n)$ correspond to 
an element $z_n \in Z(R_n)^\uno$, and consider the
associated Pappas--Rapoport filtrations 
$$0 \subset F_{i,n}'^{(1)} \subset F_{i,n}'^{(2)}
 = H^0(A_n',\Omega^1_{A_n'/R_n})_{\tau_i}
\subset H^1_{\dr}(A_n'/R_n)_{\tau_i}
\cong H^1_{\crys}(A_n'/R_n)_{\tau_i}$$
for $i=0,1$.  The task now is to lift these to filtrations
in $H^1_{\crys}(A_n'/R_{n+1})_{\tau_i}$ with the desired properties relative to the $R_{n+1}[u]/(u^2)$-linear maps
$$\psi_{i,n+1}^*: H^1_\crys(A_n'/R_{n+1})_{\tau_i}
 \longrightarrow H^1_\crys(A_n/R_{n+1})_{\tau_i} \cong
 H^1_\dr(A_0/\F)_{\tau_i} \otimes_{\F} R_{n+1}.$$
More precisely, lifting $z_n$ to $Z(R_{n+1})$ is equivalent (by the Grothendieck--Messing Theorem) to defining chains of free $R_{n+1}$-submodules 
$$0 = F_{i,n+1}'^{(0)} \subset F_{i,n+1}'^{(1)} \subset F_{i,n+1}'^{(2)}$$
of $H^1_{\crys}(A_n'/R_{n+1})_{\tau_i}$ for $i=0,1$ such that
\begin{itemize}
\item $F_{i,n+1}'^{(j)} \otimes_{R_{n+1}} R_n = F_{i,n}'^{(j)}$ for $i = 0,1$, $j=1,2$;
\item $u F_{i,n+1}'^{(1)} = 0$ and $u F_{i,n+1}'^{(2)} \subset
 F_{i,n+1}'^{(1)}$ for $i=0,1$;
\item $\psi_{n+1}^*(G_{i,n+1}'^{(j)}) \subset F_{i,0}^{(j)} \otimes_\F R_{n+1}$ if $\theta_{i,j} \in J$, where 
 $G_{i,n+1}'^{(j)}: = u^{-1} F_{i,n+1}'^{(j-1)}$;
\item $\psi_{n+1}^*(F_{0,n+1}'^{(2)}) \subset F_{0,0}^{(1)} 
  \otimes_\F R_{n+1}$;
\end{itemize}
and we must furthermore ensure that the resulting lift 
$z_{n+1} \in Z(R_{n+1})$ is unobstructed. 

Note that since $G'^{(1)}_{i,n+1} = uH^1_\crys(A_n'/R_{n+1})$, the
cases with $j =1$ in the third bullet amount to requiring the image of $\psi_{i,n+1}^*$ to be $u^{-1}F_{i,0}^{(1)}\otimes_{\F} R_{n+1}$.  This is again a special case of the assertion established in the first part of the proof of Lemma~\ref{lem:unorest}; we remark that for $i=0$ (so $\sigma^{-1}\theta,\theta \in J$ and $\phi^{-1}\circ\tau \in \Delta$, where $\theta = \theta_{0,1}$ and $\tau = \tau_0$), the argument makes crucial use of the assumption that $z_n$ is unobstructed.  Furthermore, since $\theta_{0,2} \in J''$, we
have $F'^{(2)}_{0,n} = uH^1_\dr(A'_n/R_n)_{\tau_0}$.  
It follows that all the conditions with $i=0$ are satisfied if we set $F'^{(2)}_{0,n+1} = uH^1_\crys(A'_n/R_{n+1})_{\tau_0}$ and choose within it an arbitrary lift $F'^{(1)}_{0,n+1}$ of 
$F'^{(1)}_{0,n}$.

The third bullet for $i=1$, $j =2$ translates into the 
requirement that $G'^{(2)}_{1,n+1} \subset 
 (\psi_{1,n+1}^*)^{-1}(F^{(2)}_{1,0}\otimes_\F R_{n+1})$, which forces equality since each is free of rank $3$ over $R_{n+1}$.  Note also that 
$$\psi_{1,n+1}^*(uH^1_\crys(A'_n/R_{n+1})) = 
  F^{(1)}_{1,0}\otimes_\F R_{n+1} \subset F^{(2)}_{1,0}\otimes_\F 
  R_{n+1},$$
so we obtain the desired lift of $F_{1,n}'^{(1)}$ by setting
$F'^{(1)}_{1,n+1}  
  = u(\psi_{1,n+1}^*)^{-1}(F^{(2)}_{1,0}\otimes_\F R_{n+1})$.
Finally we could choose $F'^{(2)}_{1,n+1}$ so that 
$F'^{(2)}_{1,n+1}/F'^{(1)}_{1,n+1}$ is an arbitrary lift
of the line $F'^{(2)}_{1,n}/F'^{(1)}_{1,n} \subset 
G'^{(2)}_{1,n}/F'^{(1)}_{1,n}$ to one in 
$G'^{(2)}_{1,n+1}/F'^{(1)}_{1,n+1}$, thus defining a lift 
of $z_n$ to an element $z_{n+1} \in Z(R_{n+1})$.  However the requirement that it be unobstructed determines a unique such line, specifically via the isomorphism (\ref{eqn:oblift}) constructed in the proof of Lemma~\ref{lem:unorest}.

We now return to the general setting.

\begin{lemma}  \label{lem:unorest} For all $n \ge 0$, the map $Z(R_{n+1})^\uno \to Z(R_n)^\uno$ is surjective.
\end{lemma}  
\begpf Let $(\underline{A}_n,\underline{A}'_n,\psi_n)$ correspond to an element $z_n \in Z(R_n)^\uno$.
For $\ell = n,n+1$, we let 
$$\psi_{\tau,\ell}^* : H^1_\crys(A_n'/R_{\ell})_\tau \longrightarrow H^1_\crys(A_n/R_{\ell})_\tau  \cong H^1_\crys(A_0/\FF)_\tau \otimes_{\FF} R_\ell$$
denote the morphism induced by $\psi$, so that $\psi_{\tau,n+1}^* \otimes_{R_{n+1}} R_n$ is identified with 
$\psi_{\tau,n}^*:H^1_\dr(A_n'/R_{n})_\tau \to H^1_\dr(A_n/R_{n})_\tau$ under the canonical isomorphisms
$$H^1_\crys(B/R_{n+1})\otimes_{R_{n+1}} R_n  \cong H^1_\crys(B/R_n) \cong H^1_\dr(B/R_n)$$
for $B = A_n,A'_n$.   Similarly let $\widetilde{\psi}_{\tau,\ell}^*$ denote the morphism 
$$H^1_\crys(A'_n/\widetilde{R}_\ell)_\tau \longrightarrow H^1_\crys(A_n/\widetilde{R}_\ell)_\tau \cong H^1_\crys(A_0/W_2)_\tau \otimes_{W_2} \widetilde{R}_\ell$$
for $\ell = n,n+1$, so that $\widetilde{\psi}_{\tau,n+1} \otimes_{W_2} \FF = \psi_{\tau,n+1}^*$ and 
$\widetilde{\psi}_{\tau,n+1}^* \otimes_{\widetilde{R}_{n+1}} \widetilde{R}_n = \widetilde{\psi}_{\tau,n}^*$.

Recall that to give an element of $Z(R_{n+1})^\uno$ lifting $z_n$ amounts to defining suitable chains of $R_{n+1}[u]/(u^{e_\gp})$-submodules
\begin{equation} \label{eqn:PRuno}
0 = F^{\prime(0)}_{\tau,n+1} \subset F^{\prime(1)}_{\tau,n+1} \subset \cdots F^{\prime(e_\gp)}_{\tau,n+1} \subset H^1_\crys(A_n'/R_{n+1})_\tau
\end{equation}
lifting the Pappas--Rapoport filtrations on $F^{\prime(e_\gp)}_{\tau,n}$ for all $\tau = \tau_{\gp,i} \in \Sigma_0$.
Before carrying this out, we prove that if $\tau$ is such that  $\theta = \theta_{\gp,i,1} \in J$, then the image of $\psi_{\tau,n+1}^*$ is
$u^{1-e_\gp}F_{\tau,0}^{(1)} \otimes_{\FF} R_{n+1}$.  

Suppose first that $\sigma^{-1}\theta = \theta_{\gp,i-1,e_\gp} \not\in J$, so that $\theta \in J'$ and
$(F_{\phi^{-1}\circ\tau,0}^{(e_\gp-1)})^{(p)}$ has preimage $u^{1-e_\gp}F_{\tau,0}^{(1)}$ under the morphism
$$H^1_\dr(A_0/\FF)_\tau \to H^1_\dr(A_0^{(p)}/\FF)_\tau \cong H^1_\dr(A_0/\FF)_{\phi^{-1}\circ\tau} \otimes_{\FF,\phi} \FF$$
induced by $\Ver: A_0^{(p)}  \to A_0$.  Consider the commutative diagram
$$\begin{array}{ccccc}  H^1_\crys(A_n'/R_{n+1})_\tau  & \longrightarrow & H^1_\crys(A_n^{\prime(p)}/R_{n+1})_\tau 
& \cong & H^1_\dr(A_n'/R_n)_{\phi^{-1}\circ\tau} \otimes_{R_n,\phi} R_{n+1}  \\
\downarrow && \downarrow && \downarrow \\
 H^1_\crys(A_n/R_{n+1})_\tau  & \longrightarrow & H^1_\crys(A_n^{(p)}/R_{n+1})_\tau
& \cong & H^1_\dr(A_n/R_n)_{\phi^{-1}\circ\tau} \otimes_{R_n,\phi} R_{n+1}  \\
\parallel\wr && \parallel\wr && \parallel\wr \\
H^1_\dr(A_0/\FF)_\tau \otimes_{\FF} R_{n+1}  & \longrightarrow &   
H^1_\dr(A_0^{(p)}/\FF)_\tau \otimes_{\FF} R_{n+1} & \cong&  H^1_\dr(A_0/\FF)_{\phi^{-1}\circ\tau} \otimes_{\FF,\phi} R_{n+1},
\end{array}$$
where the first horizontal maps are induced by the Verschiebung morphisms, the
horizontal isomorphisms by base-change relative to absolute Frobenius morphisms (and crystalline-de Rham isomorphisms),
top vertical maps by $\psi_n$ (and $\psi_n^{(p)}$) and vertical isomorphisms by base-change 
(and crystalline-de Rham isomorphisms again).  The image across the top of the diagram is
$F^{\prime(e_\gp)}_{\phi^{-1}\circ\tau,n} \otimes_{R_n,\phi} R_{n+1}$, and since $z_n \in S_J(R_n)$
and $\sigma^{-1}\theta \not\in J$, this maps to $F^{(e_\gp-1)}_{\phi^{-1}\circ\tau,0} \otimes_{\FF,\phi} R_{n+1}$
along the right side of the diagram. Therefore the image along the left side of the diagram must be contained
in the preimage of this, which is $u^{1-e_\gp}F_{\tau,0}^{(1)} \otimes_{\FF} R_{n+1}$, and equality follows on comparing ranks.

Suppose on the other hand that  $\sigma^{-1}\theta \in J$, so that $\phi^{-1}\circ\tau \in \Delta$.  In this case we consider the commutative diagram:
$$\xymatrix{  &\!\!\!\!\! \!\!\!\!\!  H^1_\crys(A_n^{\prime(p)}/\widetilde{R}_{n+1})_\tau \ar[dd]  \!\!\!\!\! \!\!\!\!\! & \\
\!\!\!\!\! H^1_\crys(A'_n/\widetilde{R}_{n+1})_\tau \ar[ru]  \ar[dd]  &&
H^1_\crys(A_n'/\widetilde{R}_n)_{\phi^{-1}\circ\tau} \otimes_{\widetilde{R}_n,\widetilde{\phi}} \widetilde{R}_{n+1} \ar[dd] \ar[lu]_\sim
 \\
& \!\!\!\!\!\!\!\!\!\! H^1_\crys(A_n^{(p)} / \widetilde{R}_{n+1})_\tau  \!\!\!\!\! \!\!\!\!\!   &
\\
\!\!\!\!\! H^1_\crys(A_n/\widetilde{R}_{n+1})_\tau \ar[ru] &&
H^1_\crys(A_n/\widetilde{R}_n)_{\phi^{-1}\circ\tau} \otimes_{\widetilde{R}_n,\widetilde{\phi}} \widetilde{R}_{n+1}  \ar[lu]_\sim\\
& \!\!\!\!\! \!\!\!\!\!\!\!\!\!\!  H^1_\crys(A_n^{(p^{-1})}/\widetilde{R}_n)_{\phi^{-1}\circ\tau} \otimes_{\widetilde{R}_n,\widetilde{\phi}} \widetilde{R}_{n+1}
 \ar[lu]_\sim \ar[ru] \!\!\!\!\! \!\!\!\!\! \!\!\!\!\! & \\
\!\!\!\!\! H^1_\crys(A_0/W_2)_\tau \otimes_{W_2} \widetilde{R}_{n+1} \ar[uu]^{\wr} &&
H^1_\crys(A_0/W_2)_{\phi^{-1}\circ\tau} \otimes_{W_2,\widetilde{\phi}} \widetilde{R}_{n+1}  \ar[uu]^{\wr} \\
& \!\!\!\!\! \!\!\!\!\!\!\!\!\!\!
H^1_\crys(A_0^{(p^{-1})}/W_2)_{\phi^{-1}\circ\tau} \otimes_{W_2,\widetilde{\phi}} \widetilde{R}_{n+1} \ar[ru]  \ar[lu]_\sim \ar[uu]^{\wr} &
 \!\!\!\!\! \!\!\!\!\!\!\!\!\!\!  & }$$
where the downwards arrows are induced by $\psi_n$ (and its base-changes), the northeastward arrows
by Verschiebung morphisms, all the isomorphisms are crystalline transition (and base-change) maps,
and $\widetilde{\phi}$ denotes both the restriction of the map $\widetilde{\phi}_{n+1}: \widetilde{R}_{n+1} \to \widetilde{R}_{n+1}$
to $W_2$ and its factorization through $\widetilde{R}_n \to \widetilde{R}_{n+1}$.  

We have that the image of $H^1_\crys(A_n'/\widetilde{R}_{n+1})_\tau$ in $H^1_\crys(A_n^{\prime(p)}/\widetilde{R}_{n+1})$
under the morphism induced by $\Ver:A_n^{\prime(p)} \to A_n'$ corresponds to 
$$\widetilde{F}^{\prime(e_\gp)}_{\phi^{-1}\circ\tau,n} \otimes_{\widetilde{R}_n,\widetilde{\phi}} \widetilde{R}_{n+1}
    \subset H^1_\crys(A_n'/\widetilde{R}_n)_{\phi^{-1}\circ\tau} \otimes_{\widetilde{R}_n,\widetilde{\phi}} \widetilde{R}_{n+1}$$
(this being an inclusion since $\widetilde{R}_{n+1}$ is flat over $W_2$
 and $H^1_\crys(A_n'/\widetilde{R}_n)_{\phi^{-1}\circ\tau}/\widetilde{F}^{\prime(e_\gp)}_{\phi^{-1}\circ\tau,n}$ is free over $\widetilde{R}_n/p\widetilde{R}_n = R_n$). 
Since $z_n$ is unobstructed, this submodule has the same the image along the right-hand side of the diagram as that of
$$u^{1-e_\gp}(\widetilde{F}^{(1)}_{\tau,0})^{(p^{-1})} \otimes_{W_2,\widetilde{\phi}} \widetilde{R}_{n+1}$$
under the morphism induced by $\Ver:A_0 \to A_0^{(p^{-1})}$.    It follows that the image of $H^1_\crys(A_n'/\widetilde{R}_{n+1})_\tau$
along the sides of the top left parallelogram is the same as that of 
$u^{1-e_\gp}\widetilde{F}^{(1)}_{\tau,0} \otimes_{W_2}  \widetilde{R}_{n+1} \cong u^{1-e_\gp}\widetilde{F}^{(1)}_{\tau,n+1}$ 
under the morphism 
$$H^1_\crys(A_0/W_2)_\tau \otimes_{W_2} \widetilde{R}_{n+1} \cong H^1_\crys(A_n/\widetilde{R}_{n+1})_\tau 
   \to H^1_\crys(A_n^{(p)}/\widetilde{R}_{n+1})_\tau$$
induced by Verschiebung.  Furthermore the kernel of this morphism is contained in
$p H^1_\crys(A_n/\widetilde{R}_{n+1})_\tau \subset u^{1-e_\gp}\widetilde{F}^{(1)}_{\tau,n+1}$, so we conclude that
the image of $\widetilde{\psi}_{\tau,n+1}^*$ is contained in $u^{1-e_\gp}\widetilde{F}^{(1)}_{\tau,n+1}$.  Therefore
the image of $\psi_{\tau,n+1}^*$ is contained in $u^{1-e_\gp}F^{(1)}_{\tau,n+1}$, and equality follows on comparing ranks.

We now proceed to define the lifts of the Pappas-Rapoport filtration as in (\ref{eqn:PRuno}).   If $\tau \not\in \Sigma_{\gP,0}$,
then $\psi_{\tau,n+1}^*$ is an isomorphism, under which we require that $F^{\prime(j)}_{\tau,n+1}$ correspond to 
$F^{(j)}_{\tau,n+1} = F^{(j)}_{\tau,0} \otimes_{\FF} R_{n+1}$ for all $j$.  If $\tau = \tau_{\gp,i} \in \Sigma_{\gP,0}$,
then the definition of $F^{\prime(j)}_{\tau,n+1}$ will depend on whether $\theta = \theta_{\gp,i,j}$ and $\sigma\theta$ are in $J$
(and whether $j=e_\gp$):
\begin{itemize}
\item Let $F^{\prime(0)}_{\tau,n+1} = 0$ and $G^{\prime(1)}_{\tau,n+1} = H^1_\crys(A_n'/R_{n+1})_\tau[u] = u^{e_\gp-1}H^1_\crys(A_n'/R_{n+1})_\tau$.
\item  If $\theta\not\in J$, then we let $F^{\prime(j)}_{\tau,n+1}$ denote the preimage of $F_{\tau,n+1}^{(j-1)}$ under $\psi_{\tau,n+1}^*$,
and if $j < e_\gp$, then let $G^{\prime(j+1)}_{\tau,n+1} = u^{-1}F^{\prime(j)}_{\tau,n+1}$.  Note that since 
$$F_{\tau,n+1}^{(j-1)} \subset u^{e_\gp - j + 1} H^1_\crys(A_n/R_{n+1})_\tau \subset u H^1_\crys(A_n/R_{n+1})_\tau \subset \im(\psi_{\tau,n+1}^*)$$
and $\ker(\psi_{\tau,n+1}^*)$ is free of rank one over $R_{n+1}$, it follows that $F^{\prime(j)}_{\tau,n+1}$ is free of rank $j$ over $R_{n+1}$.
In particular $F^{\prime(j)}_{\tau,n+1}$ is annihilated by $u^j$, so if $j < e_\gp$, then $G^{\prime(j+1)}_{\tau,n+1}$ is free of rank $j+2$
over $R_{n+1}$.  Furthermore $F^{\prime(j)}_{\tau,n+1} \otimes_{R_{n+1}} R_n$ is identified with the preimage of $F_{\tau,n}^{(j-1)}$
under $\psi_{\tau,n}^*$, which is $F^{\prime(j)}_{\tau,n}$ since $\theta\not\in J$, and hence also $G^{\prime(j+1)}_{\tau,n+1} \otimes_{R_{n+1}} R_n
 = G^{\prime(j+1)}_{\tau,n}$ if $j < e_\gp$.
\item  If $\theta \in J$, $\sigma\theta\in J$ and $j < e_\gp$, then let  $G^{\prime(j+1)}_{\tau,n+1}$ denote the preimage of 
$F_{\tau,n+1}^{(j+1)}$ under $\psi_{\tau,n+1}^*$ and let $F_{\tau,n+1}^{\prime(j)} = uG_{\tau,n+1}^{\prime(j+1)}$.  Note that
if $\theta_{\gp,i,1} \not\in J$, then $\theta_{\gp,i,\ell} \in J'$ for some $\ell$ such that $2 \le \ell \le j$, so
$uF_{\tau,n+1}^{(\ell)} = F_{\tau,n+1}^{(\ell -2)}$, which implies that
$$u^{e_\gp - 1} H^1_\crys(A_n/R_{n+1})_\tau \subset F_{\tau,n+1}^{(j+1)} 
 \subset u H^1_\crys(A_n/R_{n+1})_\tau \subset \im(\psi_{\tau,n+1}^*).$$
On the other hand if $\theta_{\gp,i,1} \in J$, then we have shown that 
$\im(\psi_{\tau,n+1}^*) = u^{1-e_\gp} F_{\tau,n+1}^{(1)}$, which again contains $F_{\tau,n+1}^{(j+1)}$, and
$$\psi_{\tau,n+1}^*(u^{e_\gp-1} H^1_\crys(A_n'/R_{n+1} )_\tau)
   \subset F_{\tau,n+1}^{(1)} \subset F_{\tau,n+1}^{(j+1)}.$$
Therefore in either case, $G_{\tau,n+1}^{\prime(j+1)}$ is free over $R_{n+1}$ of rank $j+2$ and contains
$u^{e_\gp-1} H^1_\crys(A_n'/R_{n+1})_\tau$, so that $F^{\prime(j)}_{\tau,n+1}$ is free of rank $j$
over $R_{n+1}$.   Furthermore $G^{\prime(j+1)}_{\tau,n+1} \otimes_{R_{n+1}} R_n$ is identified with the preimage of
$F_{\tau,n}^{(j+1)}$ under $\psi_{\tau,n}^*$, which is $G^{\prime(j+1)}_{\tau,n}$ since $\sigma\theta\in J$, and hence
also $F^{\prime(j)}_{\tau,n+1} \otimes_{R_{n+1}} R_n = F^{\prime(j)}_{\tau,n}$.
\item If $\theta \in J$ and $\sigma\theta \not\in J$, then we have defined $F^{\prime(j-1)}_{\tau,n+1}$
and $G^{\prime(j)}_{\tau,n+1}$ above, with the property that 
$P'_{\theta,n+1} := G^{\prime(j)}_{\tau,n+1}/F^{\prime(j-1)}_{\tau,n+1}$
is free of rank two over $R_{n+1}$ and annihilated by $u$.  Furthermore $P'_{\theta,n+1} \otimes_{R_{n+1}} R_n$
is identified with $P'_{\theta,n} = G^{\prime(j)}_{\tau,n}/F^{\prime(j-1)}_{\tau,n}$, and we define
$F^{\prime(j)}_{\tau,n+1}$ so that 
$$F^{\prime(j-1)}_{\tau,n+1} \subset F^{\prime(j)}_{\tau,n+1} \subset G^{\prime(j)}_{\tau,n+1}$$
and $L'_{\theta,n+1} := F^{\prime(j)}_{\tau,n+1}/F^{\prime(j-1)}_{\tau,n+1}$
is an arbitrary lift of $L'_{\theta,n} = F^{\prime(j)}_{\tau,n}/F^{\prime(j-1)}_{\tau,n}$ to a free rank
one submodule of $P'_{\theta,n+1}$.
\item  Finally if $\theta = \theta_{\gp,i,e_\gp} \in J$ and $\sigma\theta \in J$ (so $\tau\in \Delta$), then we have already defined
$F^{\prime(e_\gp-1)}_{\tau,n+1} \subset G^{\prime(e_\gp)}_{\tau,n+1}$
lifting $F^{\prime(e_\gp-1)}_{\tau,n} \subset G^{\prime(e_\gp)}_{\tau,n}$.  Furthermore we claim that
$\psi_{\tau,n+1}^*(G^{\prime(e_\gp)}_{\tau,n+1}) = F^{(e_\gp)}_{\tau,n+1}$.  Note that it suffices to prove
$\psi_{\tau,n+1}^*(G^{\prime(e_\gp)}_{\tau,n+1}) \subset F^{(e_\gp)}_{\tau,n+1}$ since $G^{\prime(e_\gp)}$
is free of rank $e_\gp+1$ over $R_{n+1}$ and contains $\ker(\psi_{\tau,n+1}^*)$, which is free of rank one.
If $e_\gp = 1$, then the desired equality is a special case of the fact that if $\theta_{\gp,i,1} \in J$, then
$\psi_{\tau,n+1}^*(H^1_\crys(A_n'/R_{n+1})_\tau) = u^{1-e_\gp}F^{(1)}_{\tau,n+1}$, so assume $e_\gp > 1$.
If $\theta_{\gp,i,e_\gp-1} \not\in J$, then $\theta \in J'$, so $uF^{(e_\gp)}_{\tau,n+1} = F^{(e_\gp-2)}_{\tau,n+1}$,
whose preimage under $\psi_{\tau,n+1}^*$ is defined to be $F^{\prime(e_\gp-1)}_{\tau,n+1}$, so $\psi_{\tau,n+1}^*$
maps $G^{\prime(e_\gp)}_{\tau,n+1} = u^{-1}F^{\prime(e_\gp-1)}_{\tau,n+1}$ to $u^{-1}F^{(e_\gp-2)}_{\tau,n+1}
 = F^{(e_\gp)}_{\tau,n+1}$.  Finally if $\theta_{\gp,i,e_\gp-1} \in J$, then we defined $G^{\prime(e_\gp)}_{\tau,n+1}$
 as the preimage of $F^{(e_\gp)}_{\tau,n+1}$.  This completes the proof of the claim, which implies that the
 restriction of $\widetilde{\psi}_{\tau,n+1}^*$ defines a surjection $\widetilde{G}^{\prime(e_\gp)}_{\tau,n+1}
 \to \widetilde{F}^{(e_\gp)}_{\tau,n+1} = \widetilde{F}^{(e_\gp)}_{\tau,0} \otimes_{W_2} \widetilde{R}_{n+1}$,
 where as usual $\widetilde{G}^{\prime(e_\gp)}_{\tau,n+1}$ denotes the preimage of
 $G^{\prime(e_\gp)}_{\tau,n+1}$ in $H^1_\crys(A_n'/\widetilde{R}_{n+1})_\tau$.  We thus obtain an
 isomorphism
 \begin{equation}  \label{eqn:oblift}
 P'_{\theta,n+1} := G^{\prime(e_\gp)}_{\tau,n+1}/F^{\prime(e_\gp-1)}_{\tau,n+1}   
    \longrightarrow P^{(p^{-1})}_{\sigma\theta,0} \otimes_{\FF} R_{n+1} \end{equation}
 exactly as in the construction of $\varepsilon_\tau(z_n)$, and we define 
 $F^{\prime(e_\gp)}_{\tau,n+1}$ so that
 $$F^{\prime(e_\gp-1)}_{\tau,n+1} \subset F^{\prime(e_\gp)}_{\tau,n+1} \subset G^{\prime(e_\gp)}_{\tau,n+1}$$
 and $F^{\prime(e_\gp)}_{\tau,n+1}/F^{\prime(e_\gp-1)}_{\tau,n+1}$ corresponds under (\ref{eqn:oblift}) to
 $L_{\sigma\theta,0}^{(p^{-1})} \otimes_{\FF} R_{n+1}$.  Note that $F^{\prime(e_\gp)}_{\tau,n+1}$ is
 free of rank $e_\gp$ over $R_{n+1}$, and lifts $F^{\prime(e_\gp)}_{\tau,n}$ since (\ref{eqn:oblift}) lifts
 $\varepsilon_\tau(z_n)$.
 \end{itemize}
 We have now defined an $R_{n+1}[u]/u^{e_\gp}$-submodule
 $F^{\prime(j)}_{\tau,n+1}$ of $H^1_\crys(A'_n/R_{n+1})_\tau$ lifting 
 $F^{\prime(j)}_{\tau,n} \subset H^1_\dr(A'_n/R_n)_\tau \cong H^1_\crys(A'_n/R_n)_\tau$
 for all $\tau = \tau_{\gp,i}$, $j=1,\ldots,e_\gp$.
 Furthermore each $F^{\prime(j)}_{\tau,n+1}$ is free of rank $j$ over $R_{n+1}$, and corresponds to $F^{(j)}_{\tau,n+1}$
 under $\psi_{\tau,n+1}^*$ if $\tau \not\in \Sigma_{\gP,0}$, so to complete the proof of the lemma it suffices to show
that the following hold for all $\tau= \tau_{\gp,i}  \in \Sigma_{\gP,0}$:
 \begin{itemize}
 \item The inclusions (\ref{eqn:PRuno}) hold with successive quotients annihilated by $u$, and hence (by the Grothendieck--Messing Theorem, as explained at the end of \S\ref{sec:iwahori})
determine data $\underline{A}'_{n+1}$
 corresponding to a point of $\widetilde{Y}(R_{n+1})$ lifting $\underline{A}'_n$ so that (\ref{eqn:PRuno}) corresponds to the Pappas--Rapoport
 filtration under the canonical isomorphism $H^1_\crys(A'_{n+1}/R_{n+1}) \cong H^1_\dr(A'_{n+1}/R_{n+1})$.
 \item If $1 \le j \le e_\gp$ and $\theta_{\gp,i,j} \not\in J$ (resp.~$\theta_{\gp,i,j} \in J$), then $\psi_{\tau,n+1}^*$ sends  
 $F^{\prime(j)}_{\tau,n+1}$ to $F^{(j-1)}_{\tau,n+1}$ (resp.~$G^{\prime(j)}_{\tau,n+1} = u^{-1}F^{\prime(j-1)}_{\tau,n+1}$
 to $F^{(j)}_{\tau,n+1}$), so that $\psi_n$ extends (again by the Grothendieck--Messing Theorem) to an isogeny $\psi_{n+1}$ and $(\underline{A}_{n+1},\underline{A}'_{n+1},\psi_{n+1})$
 corresponds to a point $z_{n+1} \in Z(R_{n+1})$ lifting $z_n$
 (where $\underline{A}_{n+1} = \underline{A}_0 \otimes_{\F} R_{n+1}$).
 \item If $\tau\in \Delta$ (i.e., $\theta_{\gp,i,e_\gp}$ and $\theta_{\gp,i+1,1}$ are both in $J$), then 
 $F^{\prime(e_\gp)}_{\tau,n+1}/F^{\prime(e_\gp-1)}_{\tau,n+1}$ corresponds
 to $(F^{(1)}_{\phi\circ\tau,0})^{(p^{-1})} \otimes_{\FF} R_{n+1}$ under $\varepsilon_\tau(z_{n+1})$,
 so that $z_{n+1} \in Z^\uno(R_{n+1})$.
 \end{itemize}
 
We start with the second assertion, which is immediate from the definition of $F^{\prime(j)}_{\tau,n+1}$ if $\theta_{\gp,i,j}\not\in J$,
so suppose that $\theta = \theta_{\gp,i,j} \in J$.   Recall we proved that if $\theta_{\gp,i,1} \in J$, then 
$\psi_{\tau,n+1}^*(H^1_\crys(A_n/R_{n+1})_\tau) = u^{1-e_\gp}F^{(1)}_{\tau,n+1}$, so it follows that
$$\psi_{\tau,n+1}^*(G^{\prime(1)}_{\tau,n+1}) =
\psi_{\tau,n+1}^*(u^{e_\gp-1}H^1_\crys(A_n/R_{n+1})_\tau) \subset F^{(1)}_{\tau,n+1}.$$
Therefore the assertion holds for $j = 1$.  Note also that if $j > 1$ and $\sigma^{-1}\theta \in J$, then the assertion is immediate
from the definition of $G^{\prime(j)}_{\tau,n+1}$.  On the other hand if $j > 1$ and $\sigma^{-1}\theta\not\in J$, then
$\theta \in J'$, so $F^{(j)}_{\tau,n+1} = u^{-1}F^{(j-2)}_{\tau,n+1}$, and therefore
$$\begin{array}{rcl} G_{\tau,n+1}^{\prime(j)} = u^{-1}F^{\prime(j-1)}_{\tau,n+1}  & = &  u^{-1}(\psi_{\tau,n+1}^*)^{-1}(F^{(j-2)}_{\tau,n+1})\\
&=& (\psi_{\tau,n+1}^*)^{-1}(u^{-1}F^{(j-2)}_{\tau,n+1}) = (\psi_{\tau,n+1}^*)^{-1}(F^{(j)}_{\tau,n+1}).\end{array}$$

For the first assertion, we need to show that if $1 \le j \le e_\gp$, then 
\begin{equation}  \label{eqn:nesting}
uF^{\prime(j)}_{\tau,n+1} \subset F^{\prime(j-1)}_{\tau,n+1} \subset F^{\prime(j)}_{\tau,n+1}.\end{equation}
Suppose that $\theta = \theta_{\gp,i,j} \not\in J$, so $F^{\prime(j)}_{\tau,n+1} = (\psi^*_{\tau,n+1})^{-1}(F_{\tau,n+1}^{(j-1)})$.
In particular if $j=1$, then $F^{\prime(1)}_{\tau,n+1} = \ker(\psi_{\tau,n+1}^*) \subset G^{\prime(1)}_{\tau,n+1}$, so 
we may assume $j \ge 2$.  If $\sigma^{-1}\theta \not\in J$, then (\ref{eqn:nesting}) is immediate from
the fact that $uF_{\tau,n+1}^{(j-1)} \subset F_{\tau,n+1}^{(j-2)} \subset F_{\tau,n+1}^{(j-1)}$, and if $\sigma^{-1}\theta\in J$,
then we have shown that $u^{-1}F_{\tau,n+1}^{\prime(j-2)} = G_{\tau,n+1}^{\prime(j-1)} = F_{\tau,n+1}^{\prime(j)}$, so
the inclusions in (\ref{eqn:nesting}) are precisely the ones in the definition of $F^{\prime(j-1)}_{\tau,n+1}$.
On the other hand if $\theta = \theta_{\gp,i,j} \in J$, then (\ref{eqn:nesting}) is immediate from the definition of
$F^{\prime(j)}_{\tau,n+1}$ if either $\sigma\theta \not\in J$ or $j = e_\gp$, so suppose that 
$\theta,\sigma\theta \in J$ and $j < e_\gp$.  Since
 $G_{\tau,n+1}^{\prime(\ell)} = (\psi_{\tau,n+1}^*)^{-1}(F^{(\ell)}_{\tau,n+1})$ for $\ell = j,j+1$, we have
 $$F^{\prime(j)}_{\tau,n+1} = u G^{\prime(j+1)}_{\tau,n+1} \subset G^{\prime(j)}_{\tau,n+1} \subset G^{\prime(j+1)}_{\tau,n+1},$$
 which in turn implies (\ref{eqn:nesting}).

 Finally the third assertion is immediate from the definition of $F^{\prime(e_\gp)}_{\tau,n+1}$ and the fact that (\ref{eqn:oblift})
 is identified with $\varepsilon_\tau(z_{n+1})$ via the canonical isomorphism 
 $$H^1_\crys(A_n/R_{n+1}) \cong H^1_\crys(A_{n+1}/R_{n+1}) \cong H^1_\dr(A_{n+1}/R_{n+1}).$$
 \epf

\subsection{Local structure}
We now use the results of the previous sections to prove a generalization of \cite[Lemma~7.1.6]{DKS}, describing
the local structure of the fibres of the projection $\pi_J: S_J \to T_{J'}$ at geometric points.  As before, we let $z \in S_J(\FF)$,
$x = \pi_J(z) \in T_{J'}(\FF)$, and $Z$ denote the (geometric) fibre of $\pi_J$ at $x$. 

\begin{lemma} \label{lem:completions}  There are isomorphisms
$$\widehat{\CO}_{Z,z} \cong \FF[[T_1,\ldots,T_\delta,X_1,...,X_{m}]]/(T_1^p,\ldots,T_{\delta}^p)
\quad\mbox{and}\quad 
\widehat{\CO}_{Z_{\red},z} \cong \FF[[X_1,...,X_m]],$$
where $m = |J'| = |J''|$ and $\delta = |\Delta|$.
In particular, $Z_\red$ is smooth of dimension $m$ over $\FF$.  Moreover the closed immersion $Z_\red \hookrightarrow Z$
identifies $\Tan_z(Z_\red)$ with $\Tan_z(Z)^\uno$.
\end{lemma}
\begpf   First consider the morphism $\chi_J = \widetilde{\zeta}_J^{(p^{-1})}:P_J^{(p^{-1})} \to S_J$, where $\zeta_J$ is the morphism
of Lemma~\ref{lem:Frobfact}; thus $\chi_J$ is finite, flat, and bijective on closed points.  Let $W$ denote
the fibre of $\pi_J\circ \chi_J$ at $x$, so that the fibre of $\chi_J$ at $x$ is a finite, flat morphism $W \to Z$ which is
bijective on closed points, and let $w = \widetilde{\xi}_J(z)^{(p^{-1})}$ denote the unique element of $W(\FF)$ mapping to $z$.
Since the completion of $W$ at $w$ is the same as its completion at the fibre of $\chi_J$ at $z$, the resulting morphism
$\widehat{\CO}_{Z,z}  \longrightarrow \widehat{\CO}_{W,w}$
is finite, flat, and in particular injective.  Furthermore choosing parameters on $P_{J,\FF}$ at $y:= \widetilde{\xi}_J(z) = w^{(p)}$
so that the morphism on completions induced by the projection $P_J \to T_{J'}$ takes the form
$$\widehat{\CO}_{T_{J',\FF},x} \cong \FF[[Y_{m+1},\ldots,Y_{d}]]
   \hookrightarrow \FF[[Y_1,\ldots,Y_d]] \cong \widehat{\CO}_{P_{J,\FF},y},$$
we see that $\widehat{\CO}_{W,w} \cong \FF[[Y_1^{1/p},\ldots,Y_d^{1/p}]]/(Y_{m+1},\ldots,Y_{d})$.
Note in particular that $\widehat{\CO}^{\red}_{W,w}$ is a domain of dimension $m$ and that if $r$ is in the nilradical of
$\widehat{\CO}_{W,w}$, then $r^p = 0$.  It follows that the same assertions hold with 
$\widehat{\CO}_{W,w}$ replaced by $\widehat{\CO}_{Z,z}$.

Now we apply the results of \S\ref{sec:unobstructed} to obtain a surjective $\FF$-algebra
morphism
\begin{equation}  \label{eqn:escape}   \widehat{\CO}_{Z,z}  \longrightarrow R = \FF[[X_1,\ldots,X_{m}]].  \end{equation}
More precisely, by Lemma~\ref{lem:unorest} there exists an element of
$$\varprojlim Z^{\uno}(R_n)  \subset \varprojlim Z(R_n) = Z(R)$$
lifting the morphism (\ref{eqn:launch}), whose surjectivity implies that of the morphism (\ref{eqn:escape})
induced by $\Spec(R) \to Z$.  Since the morphism (\ref{eqn:escape}) factors through $\widehat{\CO}_{Z,z}^\red$,
which is a domain of dimension $m$, the resulting morphism $\widehat{\CO}_{Z,z}^\red  \to R$ must be injective, and
hence an isomorphism.  Furthermore since $Z$ is Nagata, it follows that 
$$\widehat{\CO}_{Z_\red,z} = \widehat{\CO}_{Z,z}^\red \cong \FF[[X_1,\ldots, X_m]]$$
is formally smooth of dimension $m$, and hence $Z^\red$ is smooth of dimension $m$.
Note also that the image of $\Tan_z(Z_\red)$ in $\Tan_z(Z)$ is the same as that of the tangent space of $\Spec(R)$
at its closed point under the morphism induced by (\ref{eqn:escape}), which is $\Tan_z(Z)^\uno$ by construction
of (\ref{eqn:launch}).

Turning to the task of describing $\widehat{\CO}_{Z,z}$, let $\gm$ denote its maximal ideal and $\gn$ its nilradical,
and consider the exact sequence
$$0 \longrightarrow \gn/(\gn \cap \gm^2) \longrightarrow \gm/\gm^2  \longrightarrow \gm/(\gn,\gm^2) \longrightarrow 0.$$
Recall that $\gm/\gm^2$ has dimension $m+\delta$ (over $\FF$) by Lemma~\ref{lem:tangents}, and we have just seen that
$\gm/(\gn,\gm^2)$ has dimension $m$, so $\gn/(\gn\cap\gm^2)$ has dimension $\delta$.  We may therefore choose
elements $r_1,\ldots,r_\delta \in \gn$ and $s_1,\ldots,s_m \in \gm$ so that $(r_1,\ldots,r_\delta,s_1,\ldots,s_m)$
lifts a basis of $\gm/\gm^2$, and consider the surjective $\FF$-algebra morphism
$$\mu:  \FF[[T_1,\ldots,T_\delta,X_1,\ldots,X_{m}]]  \twoheadrightarrow \widehat{\CO}_{Z,z}$$
sending $T_i$ to $r_i$ for $i=1,\ldots,\delta$ and $X_i$ to $s_i$ for $i=1,\ldots,m$.  Note that the kernel of
the composite $ \FF[[T_1,\ldots,T_\delta,X_1,\ldots,X_m]] \stackrel{\mu}{\twoheadrightarrow} \widehat{\CO}_{Z,z}
\twoheadrightarrow \widehat{\CO}_{Z,z}^\red$ contains $I := (T_1,\ldots,T_\delta)$, so must in fact equal $I$,
and therefore $\ker(\mu) \subset I$.  On the other hand, from the first paragraph of the proof, we have $r_i^p = 0$
for $i=1,\ldots,\delta$, so $I^{(p)} := (T_1^p,\ldots,T_\delta^p) \subset \ker(\mu)$.  We have now shown that
$$\widehat{\CO}_{Z,z}  \cong \FF[[T_1,\ldots,T_\delta,X_1,\ldots,X_m]] / J$$
for some ideal $J$ of $\FF[[T_1,\ldots,T_\delta,X_1,\ldots,X_m]]$ such that $I^{(p)} \subset J \subset I$.

It remains to prove that $J = I^{(p)}$.  To that end consider the morphism $\widehat{\CO}_{P_{J,\FF},y} 
\to \widehat{\CO}_{S_{J,\FF},z}$ induced by $\xi_J$.  As it is a morphism of regular local $\FF$-algebras
of the same dimension, namely $d$, and its image contains that of the Frobenius endomorphism, it is (by \cite[Cor.~2]{KN}) finite flat of degree
$p^n$ where $n$ is the dimension of the kernel of the induced map on tangent spaces, which by Lemma~\ref{lem:rank}
is $\#\{\,(\gp,i)\,|\,\theta_{\gp,i,e_\gp}\in J\,\}$.  It follows that the induced morphism
$$\widehat{\CO}_{Y,y}  \longrightarrow \widehat{\CO}_{Z,z}$$
is also finite flat of degree $p^n$, where $Y \cong (\PP^1_{\FF})^m$ is the fibre of $P_J$ over $x$.
Similarly the composite 
$$\widehat{\CO}_{Y,y}  \longrightarrow \widehat{\CO}_{Z,z} \longrightarrow \widehat{\CO}_{Z_\red,z}$$
is also a morphism of regular local $\FF$-algebras of the same dimension (now $m$) whose image contains
that of the Frobenius endomorphism, so it is finite flat of degree $p^{n'}$ where $n'$ is the dimension of the 
kernel of the composite
$$\Tan_z(Z_\red) \hookrightarrow \Tan_z(Z)  \to \Tan_y(Y) \hookrightarrow \Tan_y(P_J).$$
Since the first inclusion identifies $\Tan_z(Z_\red)$ with $\Tan_z(Z)^\uno \subset \Tan_z(Z)$, it follows
from Lemmas~\ref{lem:rank} and~\ref{lem:unofirst} that $n' = n - \delta$.  Now consider the commutative
diagram:
$$\xymatrix{ & \FF[[Y_1,\ldots,Y_m]] \ar@{-->}[ld]  \ar[r]^-{\sim} \ar[d] & \widehat{\CO}_{Y,y} \ar[d] \\
Q \ar@{->>}[r] \ar@{->>}[dr] & \FF[[T_1,\ldots,T_\delta,X_1,\ldots,X_m]]/J \ar[r]^-{\sim} \ar@{->>}[d] & \widehat{\CO}_{Z,z} \ar@{->>}[d] \\
& \FF[[T_1,\ldots,T_\delta,X_1,\ldots,X_m]]/I \ar[r]^-{\sim} & \widehat{\CO}_{Z_\red,z}, } $$
where $Q = \FF[[T_1,\ldots,T_\delta,X_1,\ldots,X_m]]/I^{(p)}$ and the dashed arrow is any $\FF$-algebra
morphism lifting $\FF[[Y_1,\ldots,Y_m]] \to \FF[[T_1,\ldots,T_\delta,X_1,\ldots,X_m]]/J$.
Note that $\gr_I Q$ naturally has the structure of a finite $\widehat{\CO}_{Z_\red,z} \cong Q/IQ$-algebra, compatible
with its structure as an $\widehat{\CO}_{Y,y}$-algebra.  Furthermore $\gr_I Q$ is free of rank $p^\delta$ over
$\widehat{\CO}_{Z_\red,z}$, and hence free of rank $p^{n'}p^\delta = p^n$ over $\widehat{\CO}_{Y,y}$;
therefore $Q$ is also free of rank $p^n$ over $\widehat{\CO}_{Y,y}$.  Since the rank is the same as that of $\widehat{\CO}_{Z,z}$,
it follows that the surjection $Q\to \widehat{\CO}_{Z,z}$ is an isomorphism.
\epf

\begin{remark}  It would in fact suffice for our purposes to know that $\widehat{\CO}_{Z,z}$ is isomorphic to 
$\FF[[T_1,\ldots,T_\delta,X_1,\ldots,X_m]]/J$ for some $J$ such that $I^{(p)} \subset J \subset I$, and the
equality $J = I^{(p)}$ would follow from later considerations, but we find it interesting and satisfying to recognize it also
as a consequence of the tangent space calculations via the commutative algebra argument in the proof of
the lemma.
\end{remark}

\subsection{Reduced fibres}
In this section, we give a complete description of $Z_\red$, the reduced subscheme of the fibre $Z$ of
$\pi_J:S_J \to T_{J'}$ at any geometric point $x$ of $T_{J'}$.  This generalizes Theorem~7.2.2 of \cite{DKS}, where
it is assumed $p$ is unramified in $F$.   With the ingredients put in place in the preceding sections, we
may now proceed by an argument similar to the one in \cite[\S7.2]{DKS}.  

Recall from \S\ref{sec:Frobfact} that we have a morphism
$$\widetilde{\xi}_J:  S_J \longrightarrow P_J = \prod_{\theta \in J''} \PP_{T_{J'}}(\CP_\theta)$$
which is bijective on closed points, where the product is a fibre product over $T_{J'}$.
Taking fibres over $x \in T_{J'}(\FF)$ therefore yields such a morphism
\begin{equation}  \label{eqn:naive}  Z_\red  \longrightarrow \prod_{\theta \in J''} \PP(P_{\theta,0}), \end{equation}
where the product is over $\FF$, $P_{\theta,0} =
G_{\tau,0}^{(j)}/F_{\tau,0}^{(j-1)} = u^{-1}F_{\tau,0}^{(j-1)}/F_{\tau,0}^{(j-1)}$,
$\theta = \theta_{\gp,i,j}$, $\tau = \tau_{\gp,i}$ and $F_0^\bullet$ is
the Pappas--Rapoport filtration in the data associated to $x$.
We saw however (in the proof of Lemma~\ref{lem:completions}, where the target is denoted
$Y$) that (\ref{eqn:naive}) may not be an isomorphism, but is in fact finite flat
of degree $p^{n'}$ where $n' = \#\{\,(\gp,i)\,|\,\theta_{p,i,1} \in J''\,\}$.
We will use the approach in \cite[\S7.2.2]{DKS} to remove a factor of Frobenius from the relevant
projections $Z_\red \to \PP^1$.  We will however need the following slight generalization of the
crystallization lemma stated in \cite[\S7.2.1]{DKS}, now allowing $p$ to be ramified
in $F$ and requiring only an inclusion relation (rather than equality) between images
of components of Dieudonn\'e modules.  We omit the proof, which is essentially the
same as that of \cite[Lemma~7.2.1]{DKS}.

\begin{lemma}     \label{lem:crys}  Suppose that $S$ is a smooth scheme over an algebraically closed field $\FF$ of characteristic $p$, 
$A$, $B_1$ and $B_2$ are abelian schemes over $S$
with $O_F$-action, and $\tau \in \Sigma_{\gp,0}$.  Let $\alpha_i:A \to B_i$ for $i=1,2$ be $O_F$-linear isogenies such that
\begin{itemize}
\item $\ker(\alpha_i) \cap A[p^\infty] \subset A[p]$ for $i=1,2$, and
\item $\widetilde{\alpha}_{1,s}^* \D(B_{1,s}[p^\infty])_\tau \subset \widetilde{\alpha}_{2,s}^* \D(B_{2,s}[p^\infty])_\tau$ for all $s \in S(\FF)$
\end{itemize}
Then there is a unique morphism $\CH^1_{\dr}(B_1/S)_\tau \to \CH^1_{\dr}(B_2/S)_\tau$
of $\CO_S[u]/(u^{e_\gp})$-modules whose fibres are compatible with the injective maps
$$\D(B_{1,s}[p^\infty])_\tau \hookrightarrow \D(B_{2,s}[p^\infty])_\tau$$
 induced by $(\widetilde{\alpha}_{2,s}^*)^{-1}\widetilde{\alpha}_{1,s}^*$
for all $s \in S(\FF)$.  Furthermore if $j:\Spec(R_1) = S_1 \to S$ denotes the first infinitesimal neighborhood of $s$,
then the morphism is also compatible with the isomorphisms
$$H^1_\dr(B_{i,s}/\FF) \otimes_{\FF} R_1 \cong H^1_\dr(j^*B_i/R_1)$$
induced by their canonical isomorphisms with $H^1_\cris(B_{i,s}/R_1)$ for $i=1,2$.
\end{lemma}

Now let $(\underline{A},\underline{A}',\psi)$ denote the restriction to $Z_\red$ of the universal triple, and
suppose that $\tau = \tau_{\gp,i}$ is such that $\theta = \theta_{\gp,i,e_\gp} \in J$.  Then for any closed point
$z$ of $Z_\red$, we have
\begin{equation}\label{eqn:imageV}
\Ver^*(H^1_\dr(A_0^{(p^{-1})}/\FF)_\tau) 
 = F_{\tau,0}^{(e_\gp)} = \psi_{\tau,z}^*(G_{\tau,z}^{\prime(e_\gp)}) \subset \psi_{\tau,z}^*(H^1_\dr(A'_z/\FF)_\tau),\end{equation}
 where as usual we write $\underline{A}_0$ for the fibre of $\underline{A}$ at $z$ (or $x$) and $F^{(j)}_{\tau,0}$ for 
 (the sections of) $\CF_{\tau,0}^\bullet$, and similarly $\underline{A}'_z$ for the fibre of $\underline{A}'$,
 $F^{\prime(j)}_{\tau,z}$ for $\CF^{\prime(j)}_{\tau,z}$ and $G^{\prime(j)}_{\tau,z} = u^{-1}F^{\prime(j-1)}_{\tau,z}$.
 We may therefore apply Lemma~\ref{lem:crys}  with $S = Z_\red$, $B_1 = A_0^{(p^{-1})} \times_{\FF} Z_\red$,
 $B_2 = A'$, $\alpha_2 = \psi$ and $\alpha_1 = \Ver \times_{\FF} Z_\red$ (where $\Ver:A_0 \to A_0^{(p^{-1})}$)
 to obtain an $\CO_{Z_\red}[u]/u^{e_\gp}$-linear morphism
 $$\beta: (H^1_\dr(A_0/\FF)_{\phi\circ\tau})^{(p^{-1})} \otimes_{\FF} \CO_{Z_\red}
          = \CH^1_\dr(B_1/Z_\red)_\tau \longrightarrow \CH^1_\dr(A'/Z_\red)_\tau $$
 which at closed points is compatible in the evident sense with the injective morphisms induced by $\psi$ and $\Ver$
 on Dieudonn\'e modules and with the isomorphisms over first-order thickenings provided by crystalline-de Rham comparisons
 (cf.~diagrams (39) and (40) of \cite{DKS}).  In particular the fibre $\beta_z$ at each closed point $z$ is given by the reduction mod $p$ of the unique injective homomorphism 
$\widetilde{\beta}_z$ of free rank two $W(\F)[u]/(E_\tau)$-modules making the diagram
$$ \xymatrix{ \D(A_0[p^\infty])_{\phi\circ\tau}^{\phi^{-1}}
 \ar[rr]^{\widetilde{\beta}_z} \ar[rd]_{V} && \D(A'_z[p^\infty])_\tau \ar[ld]^{\widetilde{\psi}_{\tau,z}^*}\\
& \D(A_0[p^\infty])_\tau& }$$
commute.  It follows from (\ref{eqn:imageV}) that the image 
of $\widetilde{\beta}_z$ is the preimage 
 $\widetilde{G}_{\tau,z}^{\prime(e_\gp)}$ of ${G}_{\tau,z}^{\prime(e_\gp)}$ under reduction mod $p$, and therefore that 
of $\beta_z$ is ${G}_{\tau,z}^{\prime(e_\gp)}$.  Since this holds for each closed pont $z$, the image of $\beta$ is 
$\CG_\tau^{\prime(e_\gp)}$, and hence it induces an isomorphism
 \begin{equation}  \label{eqn:cryse1} 
 \xymatrix{
  (H^1_\dr(A_0/\FF)_{\phi\circ\tau} / u H^1_\dr(A_0/\FF)_{\phi\circ\tau} )^{(p^{-1})}   \otimes_{\FF} \CO_{Z_\red}
  \ar[d]^-{\wr}_-{u^{e_\gp-1}} \ar[r]^-{\sim}
&  \CG_\tau^{\prime(e_\gp)}/\CF_\tau^{\prime(e_\gp-1)} 
 =: \CP_\theta'. \\
  P_{\sigma\theta,0}^{(p^{-1})} \otimes_{\FF} \CO_{Z_\red} & }
\end{equation}

We can thus define the morphism
 $$ \mu_{\sigma\theta}:  Z_\red \longrightarrow \PP(P_{\sigma\theta,0}^{(p^{-1})})  $$
 such that the tautological line bundle pulls back to the subbundle corresponding to $\CL_\theta'$ under (\ref{eqn:cryse1}).
The commutativity of the diagram
$$ \xymatrix{(\widetilde{G}'^{(1)}_{\phi\circ\tau,z})^{\phi^{-1}}
  \ar@{^{(}->}[d]_{(\widetilde{\psi}^*_{\phi\circ\tau},z)^{\phi^{-1}}}
& \D(A_z'[p^\infty])_{\phi\circ\tau}^{\phi^{-1}} \ar[l]_{\sim}^{u^{e_{\gp}-1}} \ar[r]^-{\sim}_-V &
 \widetilde{F}'^{(e_{\gp})}_{\tau,z} \ar@{^{(}->}[r] & 
  \widetilde{G}'^{(e_{\gp})}_{\tau,z} \ar[d]_{\wr}^{\widetilde{\psi}_{\tau,z}^*}\\ 
   (\widetilde{G}^{(1)}_{\phi\circ\tau,0})^{\phi^{-1}}&
\D(A_0[p^\infty])_{\phi\circ\tau}^{\phi^{-1}} \ar[l]_{\sim}^{u^{e_\gp-1}} \ar[rr]^-{\sim}_-V &&
 \widetilde{F}^{(e_{\gp})}_{\tau,0}   }
$$
implies that on the fibre at $z$, the line $L_{\theta,z}'$ in $P_{\theta,z}'$ corresponds under (\ref{eqn:cryse1}) to the image of the morphism
$$(\psi^*_{\sigma\theta,z})^{(p^{-1})}: P_{\sigma\theta,z}'^{(p^{-1})} \lra P_{\sigma\theta,0}^{(p^{-1})}.$$
(cf.~the proof of \cite[Thm.~7.2.2.1]{DKS}).
If $\sigma\theta = \theta_{\gp,i+1,1} \in J$, then this image is 
$L_{\sigma\theta,0}^{(p^{-1})}$ for all $z$, so
$\mu_{\sigma\theta}$ is a constant morphism.
On the other hand, if $\sigma\theta \not\in J$, so $\sigma\theta\in J''$, then 
the composite of $\mu_{\sigma\theta}$ with the Frobenius morphism  
$$ \PP(P_{\sigma\theta,0}^{(p^{-1})}) \longrightarrow \PP(P_{\sigma\theta,0})$$
is precisely the factor indexed by $\sigma\theta$ in (\ref{eqn:naive}).

Now for any $\theta = \theta_{\gp,i,j} \in \Sigma$, let
$$P_{\theta,0}^{(n_\theta^{-1})} = \left\{\begin{array}{ll}
P^{(p^{-1})}_{\theta,0},&\mbox{if $j = 1$;}\\
P_{\theta,0},&\mbox{if $j > 1$.}\end{array}\right.$$
Replacing the $\theta$-factor in (\ref{eqn:naive}) with $\mu_\theta$ for each
$\theta = \theta_{\gp,i,j} \in J''$ such that $j = 1$ therefore yields a commutative diagram
$$\xymatrix{  Z_\red   \ar[rr]  \ar[rd]_-{(\ref{eqn:naive})} & & {\displaystyle\prod_{\theta\in J''}  \PP(P^{(n_\theta^{-1})}_{\theta,0})}  \ar[ld] \\
 & {\displaystyle\prod_{\theta\in J''}  \PP(P_{\theta,0}) } & } $$
 where the right diagonal map is the Frobenius morphism (resp.~identity) on the factors indexed by 
 $\theta = \theta_{\gp,i,j}$ such that $j=1$ (resp.~$j > 1$).  Since the schemes are all smooth over $\FF$
 and the diagonal morphisms are finite flat of the same degree (namely $p^{n'}$), it follows that the
 horizontal map is an isomorphism.
 
In order to obtain the above isomorphism, we factored out a power of Frobenius in each component of the form $\theta_{\gp,i,1}$ by describing the corresponding morphism to $\PP^1$ in terms of the data associated to $A'$ (i.e., pull-backs via $\pi_2$ of vector bundles on $\widetilde{Y}_{U,\F}$).  We can similarly describe the other components of the isomorphism in terms of the bundles associated to $A'$, again with a shift in the index, but the relation is easier to obtain since no Frobenius factor intervenes. To that end, suppose that $\tau = \tau_{\gp,i}$
and $\theta = \theta_{\gp,i,j} \in J$ for some $j < e_\gp$, and consider the $\CO_{Z_\red} [u]/(u^{e_\gp})$-linear
morphism
$$\xi_\tau^*:  \gp\CO_{F,\gp} \otimes_{\CO_{F,\gp}}  H^1_{\dr}(A/Z_\red)_\tau
 = \CH^1_\dr( (\gP^{-1} \otimes_{\CO_F} A)/Z_\red)_\tau \longrightarrow \CH^1_\dr(A'/Z_\red)_\tau$$
induced by the unique isogeny $\xi:A' \to \gP^{-1} \otimes_{\CO_F} A$ such that $\xi\circ\psi$ is the
canonical isogeny $A \to \gP^{-1} \otimes_{\CO_F} A$.  Since the image of $\CG_\tau^{\prime(j)}$
under $\psi_\tau^*$ is $\CF_\tau^{(j)}$, it follows that the image of 
$\gp\CO_{F,\gp} \otimes_{\CO_{F,\gp}} \CF_{\tau}^{(j)}$ under $\xi_\tau^*$ is $\CF_\tau^{\prime(j-1)}$.
We thus obtain an $\CO_{Z_\red}$-linear isomorphism
\begin{equation} \label{eqn:crysj} P_{\sigma\theta,0} \otimes_{\FF} \CO_{Z_\red}
  = \CP_{\sigma\theta}  \stackrel[u]{\sim}{\longrightarrow} \gp\CO_{F,\gp} \otimes_{\CO_{F,\gp}} \CP_{\sigma\theta}
  \stackrel[\xi_\tau^*]{\sim}{\longrightarrow} \CP'_\theta,\end{equation}
under which $\CL'_\theta$ corresponds to the image of the morphism $\CP'_{\sigma\theta} \to \CP_{\sigma\theta}$
induced by $\psi$.  Recall that if $\sigma\theta \not\in J$, i.e., $\sigma\theta\in J''$, then this is precisely the subbundle of 
$\CP_{\sigma\theta} = P_{\sigma,\theta} \otimes_{\FF} \CO_{Z_\red}$ defining the factor indexed
by $\sigma\theta$ in the morphism (\ref{eqn:naive}).  On the other hand if $\theta,\sigma\theta \in J$,
then this subbundle is $L_{\sigma\theta,0} \otimes_{\FF} \CO_{Z_\red}$.  
\begin{theorem}  \label{thm:reduced} For $\theta = \theta_{\gp,i,j} \in \Sigma$ such that $\sigma^{-1}\theta \in J$, let
$$ \delta_{\theta}:  P^{(n_\theta^{-1})}_{\theta,0} \otimes_{\FF} \CO_{Z_\red}  \stackrel{\sim}{\longrightarrow}  \CP'_{\sigma^{-1}\theta}$$
denote the isomorphism defined by (\ref{eqn:cryse1}) or (\ref{eqn:crysj}), according to whether $j=1$ or $j > 1$ (and applied with $\sigma^{-1}\theta$
in place of $\theta$), and let
$$ \mu_\theta : Z_\red \longrightarrow \PP(P_{\theta,0}^{(n_\theta^{-1})})$$
denote the morphism such that the tautological line bundle pulls back to the subbundle corresponding to $\CL'_{\sigma^{-1}\theta}$
under $\delta_\theta$.  Then 
\begin{enumerate}
\item \label{item:reduced} the resulting morphism 
$$  Z_\red  \stackrel{(\mu_\theta)_\theta}{\longrightarrow}   \prod_{\theta \in J''}  \PP(P_{\theta,0}^{(n_\theta^{-1})})$$
is an isomorphism,
\item and if $\sigma^{-1}\theta,\theta \in J$, then the morphism $\mu_\theta$ is the projection to the point corresponding to
$L_{\theta,0}^{(n_\theta^{-1})}$.
\end{enumerate}

\end{theorem}

\begin{remark} \label{rmk:ERX} The theorem above is asserted in \cite[\S4.9(1)]{ERX} for fibres of $\widetilde{\xi}_J$ at geometric {\em generic} points of $T_{J'}$.  There is however a serious error in their argument (as already noted in \cite[\S7]{DKS}), specifically in the unjustified claim in the paragraph after \cite[(4.9.3)]{ERX} that one  {\em can rearrange the choices of local parameters\dots}
Theorem~\ref{thm:reduced} thus fills the resulting gap in the proof of Proposition~3.19 of \cite{ERX}.
We note also that the latter is a cohomological vanishing result,
a stronger version of which follows from the results in this paper (see Remark~\ref{rmk:ERX3}), and its application in \cite{ERX} (as well as \cite{DW}) is to the construction of Hecke operators, which is ultimately achieved in greater generality in \S\ref{ss:hecke} of this paper.
\end{remark}

\subsection{Thickening}
We will now give a complete description of $Z$, generalizing \cite[Thm.~7.2.4]{DKS}.  Our approach is based on
that of \cite{DKS}: we will use the isomorphism of (\ref{eqn:cryse1}), for $\theta = \theta_{\gp,i,e_\gp} \in J$ such 
that $\sigma\theta\in J$, in order to extend the isomorphism Theorem~\ref{thm:reduced}(\ref{item:reduced}) to one between $Z$ and a suitable thickening of the target.

To that end, let $T$ denote the fibre product
$$\begin{array}{ccc}  T  & \longrightarrow &\displaystyle{ \prod_{\sigma^{-1}\theta,\theta \in J}  \PP(P_{\theta,0}^{(n_\theta^{-1})})})\\
\downarrow && \downarrow \\
\Spec \FF & \longrightarrow &\displaystyle{ \prod_{\sigma^{-1}\theta,\theta\in J}  \PP(P_{\theta,0}),} \end{array}$$
where the bottom arrow is defined by $L_{\theta,0}$ and the right vertical arrow is defined by the Frobenius morphism
$M_\theta \mapsto M_\theta^{(p)}$ on the factors indexed by $\theta$ such that $n_\theta = p$.
(Recall that these are the $\theta = \theta_{\gp,i,j}$ such that $j= 1$, and note that we could have omitted
the other factors from the definition, but we maintain them for the sake of uniformity in later formulas.)

We then let $\widetilde{Z}_\red = Z_\red \times_{\FF} T$, and let $i: Z_\red \hookrightarrow \widetilde{Z}_\red$
denote the resulting divided power thickening, i.e.,, the unique section of the projection 
$q:\widetilde{Z}_\red \to Z_\red$.  We thus have canonical isomorphisms
 $\CO_{\widetilde{Z}_\red}[u]/(u^{e_\gp})$-linear
morphisms
$$(R^1s'_{\cris,*}\CO_{A',\crys})_{\widetilde{Z}_\red,\tau} \cong q^*\CH^1_\dr(A'/Z_\red)_\tau$$
for all $\tau = \tau_{\gp,i} \in \Sigma_0$ (where we continue to let $(\underline{A},\underline{A}',\psi)$
denote the universal triple over $Z_\red$, and $s':A' \to Z_\red$ is the structure morphism).
Furthermore recall that to give a lift of the closed immersion $Z_\red \hookrightarrow Z$ to a morphism
$\widetilde{Z}_\red \to Z$ is equivalent to giving a lift of $\CF_\tau^{\prime(j)}$ (for each $\tau = \tau_{\gp,i}$,
$j= 1,\ldots,e_\gp$) to an $\CO_{\widetilde{Z}_\red}[u]/(u^{e_\gp})$-submodule $\widetilde{\CF}_\tau^{\prime(j)}$
of $q^*\CH^1_\dr(A'/Z_\red)_\tau$ such that 
\begin{itemize}
\item $\widetilde{\CF}_\tau^{\prime(j-1)} \subset \widetilde{\CF}_\tau^{\prime(j)}$ and 
 $\widetilde{\CF}_\tau^{\prime(j)}/\widetilde{\CF}_\tau^{\prime(j-1)}$ is a line bundle on $Z_\red$
 annihilated by $u$;
\item the morphism 
  $$q^*\psi_\tau^*: q^*\CH^1_\dr(A'/Z_\red)_\tau \longrightarrow q^*\CH^1_\dr(A/Z_\red)
           = H^1_\dr(A_0/\FF)_\tau \otimes_{\FF} \CO_{\widetilde{Z}_\red}$$
 sends $\widetilde{\CF}_\tau^{\prime(j)}$ to $F_{\tau,0}^{(j)} \otimes_{\FF} \CO_{\widetilde{Z}_\red}$
 if $\tau \not\in \Sigma_{\gP,0}$, sends $\widetilde{\CG}_\tau^{\prime(j)} := u^{-1}\widetilde{\CF}_\tau^{\prime(j-1)}$
 to $F_{\tau,0}^{(j)} \otimes_{\FF} \CO_{\widetilde{Z}_\red}$ if $\theta_{\gp,i,j} \in J$, and sends $\widetilde{\CF}_\tau^{\prime(j)}$
 to $F_{\tau,0}^{(j-1)} \otimes_{\FF} \CO_{\widetilde{Z}_\red}$ otherwise.
\end{itemize}
We construct such a lift by setting $\widetilde{\CF}_{\tau}^{\prime(j)} = q^*\CF_{\tau}^{\prime(j)}$ unless 
$j=e_\gp$ and $\sigma^{-1}\theta,\theta\in J$, where $\theta = \theta_{\gp,i+1,1}\in J$, in which case we define
$\widetilde{\CF}_{\tau}^{\prime(e_\gp)}$ so that 
$$\widetilde{\CL}'_{\sigma^{-1}\theta} := \widetilde{\CF}_\tau^{\prime(e_\gp)}/\widetilde{\CF}_\tau^{\prime(e_\gp-1)}
     \subset  \widetilde{\CG}_\tau^{\prime(e_\gp)}/\widetilde{\CF}_\tau^{\prime(e_\gp-1)}  = q^*\CP_{\sigma^{-1}\theta}$$
corresponds to the pull-back of the tautological line bundle on
$\PP(P_{\theta,0}^{(p^{-1})})$ under $q^*\delta_\theta$.  Note that $i^*\widetilde{\CL}'_{\sigma^{-1}\theta} = \CL'_{\sigma^{-1}\theta}$
since each corresponds to $L_{\theta,0}^{(p^{-1})} \otimes_{\FF} \CO_{Z_\red}$ under $\delta_\theta$.
The fact that the resulting sheaves $\widetilde{\CF}_\tau'^{(j)}$, including those with $j=e_\gp$, satisfy the required properties
is immediate from their definition and the corresponding properties of the $\CF_\tau'^{(j)}$.  The Grothendieck--Messing Theorem thus yields a triple
$(\underline{\widetilde{A}},\underline{\widetilde{A}}',\widetilde{\psi})$ over $\widetilde{Z}_\red$
such that $i^*\widetilde{A}' = A'$ and $\widetilde{\CF}'^\bullet$ corresponds to the Pappas--Rapoport filtration
under the canonical isomorphism
$$\CH^1_\dr(\widetilde{A}'/\widetilde{Z}_\red)  \cong 
(R^1s'_{\cris,*}\CO_{A',\crys})_{\widetilde{Z}_\red} \cong q^*\CH^1_\dr(A'/Z_\red).$$

We claim that the resulting morphism
$$\widetilde{Z}_\red \longrightarrow Z$$
is an isomorphism.  We first note that the induced morphism is injective on
tangent spaces at all closed points $\widetilde{z} = (z,t)$ in
$\widetilde{Z}_\red(\FF)  = Z(\FF) \times T(\FF)
\stackrel{\sim}{\to} Z(\FF)$.  Indeed let $R_1 = \FF[\epsilon]$ and suppose that 
$\widetilde{z}_1 = (z_1,t_1) \in Z_\red(R_1) \times T(R_1) =
\widetilde{Z}_\red(R_1)$ corresponds to an element of
the kernel of the induced map
$$\Tan_{\widetilde{z}}(\widetilde{Z}_\red) = \Tan_z(Z_\red) \times \Tan_t(T) \to \Tan_z(Z),$$
i.e., that the  triple $(\underline{\widetilde{A}}_1',\underline{\widetilde{A}}_1,\psi_1)$ associated to
$\widetilde{z}_1$ is (isomorphic to) the pull-back of the triple $(\underline{A}_0,\underline{A}_0',\psi_0)$
corresponding to $z$.  The fact that $Z_\red \hookrightarrow Z$ is a closed immersion then implies that $z_1$ is trivial,
and the fact that $\widetilde{F}_{\tau,1}^{\prime(e_\gp)}$ corresponds to $F_{\tau,0}^{\prime(e_\gp)} \otimes_{\FF} R_1$
under the canonical isomorphism $H^1_\dr(\widetilde{A}'_1/R_1) \cong H^1_\dr(A'_0/\FF) \otimes_{\FF} R_1$
implies that $\widetilde{L}'_{\sigma^{-1}\theta,1}$ corresponds to $L_{\sigma^{-1}\theta,0}^{(p^{-1})}\otimes_{\FF} R_1$ under
the isomorphisms $\widetilde{P}'_{\sigma^{-1}\theta,1} \cong P'_{\sigma^{-1}\theta,0} \otimes_{\FF} R_1
\cong P_{\theta,0}^{(p^{-1})} \otimes_{\FF} R_1$, and hence that $t_1$ is trivial.

Since $\widetilde{Z}_\red \longrightarrow Z$ is injective on tangent spaces, it is a closed immersion, and
the same argument as in the proof of \cite[Lemma~7.2.3]{DKS} shows that it is an isomorphism.
Alternatively, note that the diagram
$$\begin{array}{ccc}  \widetilde{Z}_\red & \longrightarrow & Z\\
{\scriptstyle{q}} \downarrow && \downarrow   \\
Z_\red & \longrightarrow & {\displaystyle\prod_{\theta \in J''}} \PP(P_{\theta,0}) \end{array}$$
commutes, where the right vertical map is the fibre of $\widetilde{\xi}_J$ over $x$ and the
bottom horizontal arrow is its restriction, i.e. (\ref{eqn:naive}).  Recall from the proof of Lemma~\ref{lem:completions}
that these morphisms are finite flat of degree $p^n$ and $p^{n'}$, respectively.  Since $q$
is finite flat of degree $p^\delta$, the composite morphism $\widetilde{Z}_\red \to Z \to \prod_{\theta\in J''} \PP(P_{\theta,0})$
is finite flat of degree $p^n$, and since $\widetilde{Z}_\red \to Z$ is a closed immersion, it follows that it
must be an isomorphism.  

Note in particular that the above construction extends the isomorphisms $\delta_\theta$ of (\ref{eqn:cryse1})
and (\ref{eqn:crysj}) to isomorphisms
$$\CP_{\theta,0}^{(n_\theta^{-1})}\otimes_{\FF}\CO_{\widetilde{Z}_\red}
   \stackrel{\sim}{\longrightarrow} q^*\CP'_{\sigma^{-1}\theta}  \cong
         \widetilde{\CP}'_{\sigma^{-1}\theta}$$
 over $\widetilde{Z}_\red \cong Z$ for all $\theta \in J$.  We thus obtain the following extension of Theorem~\ref{thm:reduced},
 maintaining the same notation except that now we write $\CL'_\theta \subset \CP_\theta'$ for the vector bundles 
 $\CG_\tau^{\prime(j)}/\CF_{\tau}^{\prime(j-1)}$ associated to the data $\underline{A}'$, where 
  $\tau = \tau_{\gp,i}$, $\theta = \theta_{\gp,i,j}$ and $(\underline{A},\underline{A}',\psi)$ is the universal triple over $Z$, 
 and we let
 \begin{equation} \label{eqn:epsilontheta}  \varepsilon_\theta : {\CP}'_{\sigma^{-1}\theta}    \stackrel{\sim}{\longrightarrow} 
  \CP_{\theta,0}^{(n_\theta^{-1})}\otimes_{\FF}\CO_{Z}
   \end{equation}
denote the inverse of the isomorphism just constructed over $Z$.

\begin{theorem}  \label{thm:fibres}   The morphism
$$  Z    \longrightarrow  \displaystyle \prod_{\sigma^{-1}\theta \in J}  \PP(P_{\theta,0}^{(n_\theta^{-1})}), $$
defined so the tautological line bundle on $\PP(P_{\theta,0}^{(n_\theta^{-1})})$ pulls back
to $\varepsilon_\theta(\CL'_{\sigma^{-1}\theta})$, is a closed immersion, identifying $Z$ with the
fibre over $(L_{\theta,0})_{\theta\in J-J'}$ of the morphism
$$\begin{array}{ccc} \displaystyle\prod_{\sigma^{-1}\theta \in J}  \PP(P_{\theta,0}^{(n_\theta^{-1})})  & \longrightarrow & 
\displaystyle\prod_{\theta \in J - J'}  \PP(P_{\theta,0}) \\
(M_\theta)_\theta & \mapsto & (M_\theta^{(n_\theta)})_\theta. \end{array}$$
\end{theorem}

We remark that the constructions in the proof of the theorem are independent of the quasi-polarizations in the data,
and hence compatible with the action of $\CO_{F,(p),+}^\times$ on the various schemes and vector bundles.
One therefore obtains an identical description of the geometric fibres of the morphism $\overline{Y}_0(\gP)_J  \to \overline{Y}_{J'}$
for sufficiently small $U$, with $\CL'_{\sigma^{-1}\theta} \subset \CP'_{\sigma^{-1}\theta}$ (resp.~$L_{\theta,0} \subset P_{\theta,0}$)
replaced by restrictions (resp.~fibres) of descents to $\overline{Y}_0(\gP)_{J,\FF}$ (resp.~$\overline{Y}_{J',\FF}$).
Furthermore the resulting descriptions are compatible with the action of $\GL_2(\A_{F,\f}^{(p)})$ in the obvious sense.
Since we will not make use of this, we leave it to the interested reader to make precise, and simply record the following
cruder, more immediate consequence of Theorem~\ref{thm:fibres}, generalizing Theorem~D of \cite{DKS}.

\begin{corollary}  \label{cor:fibres}  Every geometric fibre of the morphism $\overline{Y}_0(\gP)_J  \to \overline{Y}_{J'}$
is isomorphic to $(\PP^1_{\FF})^m \times (\Spec(\FF[T]/T^p))^\delta$ where $m$ and $\delta$ are as in
Lemma~\ref{lem:completions}, i.e., $m = |J'| = |J''|$ and $\delta = |\Delta| = \#\{\,\tau = \tau_{\gp,i}\,|\,\theta_{\gp,i,e_\gp}, \theta_{\gp,i+1,1}\in J\,\}$.
\end{corollary}

\section{Cohomological vanishing}

\subsection{Level $U_1(\gP)$}
One of our main results will be the vanishing of the higher direct images of both the structure and dualizing
sheaves of $Y_{U_0(\gP)}$ under the projection to $Y_U$.  We will prove this also holds with $U_0(\gP)$
replaced by its open compact subgroup $U_1(\gP):= \{\, g \in U \,|\, \mbox{$g_\gp \in I_1(\gp)$ for all $\gp|\gP$}\,\}$, where
$$I_1(\gp) = \left\{\,\left. \left(\begin{array}{cc} a&b \\ c& d \end{array}\right) \in \GL_2(\CO_{F,\gp}) \,\right|\, c, d-1 \in \gp \CO_{F,\gp} \,\right\}.$$

We first recall the definition of suitable integral models for Hilbert modular varieties of level $U_1(\gP)$, following \cite{pappas}.
We consider the functor which, to a locally Noetherian $\CO$-scheme $S$, associates the set of isomorphism classes of data
$(\underline{A},\underline{A}',\psi,P)$, where $(\underline{A},\underline{A}',\psi,P)$ corresponds to an element of
$\widetilde{Y}_{U_0(\gP)}(S)$, and $P \in A(S)$ is an $(\CO_F/\gP)$-generator of $\ker(\psi)$.
The functor is represented by an $\CO$-scheme which we denote $\widetilde{Y}_{U_1(\gP)}$, and the forgetful
morphism $\widetilde{f}:\widetilde{Y}_{U_1(\gP)} \to \widetilde{Y}_{U_0(\gP)}$ is finite flat.

The scheme $\widetilde{Y}_{U_1(\gP)}$ can be described explicitly as follows in terms of $H = \ker(\psi)$,
where $(\underline{A},\underline{A}',\psi)$ is the universal triple over $S: = \widetilde{Y}_{U_0(\gP)}$.  Firstly, the same argument
as when $p$ is unramified in $F$ shows that $H$ is a Raynaud $(\CO_F/\gP)$-module scheme (in the sense that
$H[\gp]$ is a Raynaud $(\CO_F/\gp)$-vector space scheme for each $\gp|\gP$).   The correspondence of
\cite[Thm.~1.4.1]{raynaud} thus associates to $H$ an invertible $\CO_S$-module $\CR_\tau$ for each $\tau = \tau_{\gp,i} 
\in \Sigma_{\gP,0}$, together with morphisms
$$s_\tau: \CR_\tau^{\otimes p} \longrightarrow \CR_{\phi\circ\tau}\quad\mbox{and}\quad 
    t_\tau: \CR_{\phi\circ\tau} \longrightarrow \CR_\tau^{\otimes p}$$
such that $s_\tau\circ t_\tau = w_{\gp}$ for a certain fixed element $w_\gp \in p\CO^\times$.
We then have
$$H = \SPEC ((\Sym_{\CO_S}\CR)/\CI),$$
where $\CR = \bigoplus_{\tau \in \Sigma_{\gP,0}} \CR_\tau$,
$\alpha \in \CO_F$ acts on $\CR_\tau$
as the Teichm\"uller lift of $\overline{\tau}(\alpha)$,
and $\CI$ is the sheaf of ideals generated by the $\CO_{S}$-submodules 
$(s_\tau - 1)\CR_\tau^{\otimes p}$ for $\tau \in \Sigma_{\gP,0}$.
The comultiplication on $H$ 
is given via duality by the scheme
structure on the Cartier dual 
$$H^\vee =  \SPEC ((\Sym_{\CO_S}\CR^\vee)/\CJ),$$
where $\CR^\vee  = \bigoplus_{\tau\in \Sigma_{\gP,0}} \CR_\tau^{-1}$,
and $\CJ$ is generated by
$(t_\tau - 1)(\CR_\tau^{-1})^{\otimes p}$ for $\tau \in \Sigma_{\gP,0}$, where
$t_\tau$ is regarded as a morphism $(\CR_\tau^{-1})^{\otimes p} \to \CR_{\phi\circ\tau}^{-1}$.
The same argument as in \cite[5.1]{pappas} then shows that $\widetilde{Y}_{U_1(\gP)}$
can be identified with the closed
subscheme of $H$ defined by the sheaf of ideals generated by the
$(s_\gp - 1)\left(\bigotimes _{\tau\in \Sigma_{\gp,0}}
 \CR_\tau^{\otimes(p-1)}\right)$ for $\gp|\gP$,
 where $s_\gp := \bigotimes_{\tau\in \Sigma_{\gp,0}} s_\tau$ is viewed as a morphism
 $\bigotimes_{\tau\in \Sigma_{\gp,0}}\CR_\tau^{\otimes(p-1)} \to \CO_S$.

Note that the scheme $\widetilde{Y}_{U_1(\gP)}$ is equipped with a natural action of 
$(\CO_F/\gP)^\times$ over $\widetilde{Y}_{U_0(\gP)}$, defined by $\alpha\cdot(\underline{A},\underline{A}',\psi,P)
= (\underline{A},\underline{A}',\psi,\alpha P)$.  Furthermore we may write
\begin{equation}  \label{eqn:charO} \widetilde{f}_* \CO_{\widetilde{Y}_{U_1(\gP)}} = \bigoplus_{\chi}  \CR_\chi \end{equation}
where $\chi$ runs over all characters $(\CO_F/\gP)^\times \to \CO^\times$, and each summand may be
written as
$$\CR_\chi = \bigotimes_{\tau \in \Sigma_{\gP,0}}  \CR_\tau^{\otimes m_{\chi,\tau}}$$
where the integers $m_{\tau,\chi}$ are uniquely determined by the conditions
\begin{itemize}
\item $0 \le m_{\chi,\tau} \le p - 1$ for all $\tau \in \Sigma_{\gP,0}$;
\item $m_{\chi,\tau} < p - 1$ for some $\tau \in \Sigma_{\gp,0}$ for each $\gp|\gP$;
\item $\overline{\chi}(\alpha) = \prod_{\tau} \overline{\tau}(\alpha)^{m_{\chi,\tau}}$ for all $\alpha \in (\CO_F/\gP)^\times$.
\end{itemize}
Since $\widetilde{f}$ is finite flat, we may similarly decompose the direct image of the dualizing sheaf as 
\begin{equation} \label{eqn:charK} \widetilde{f}_*\CK_{\widetilde{Y}_{U_1(\gP)}/\CO}
 = \Shom_{\CO_S} (\widetilde{f}_*\CO_{\widetilde{Y}_{U_1(\gP)}},\CK_{S/\CO})
 = \bigoplus_{\chi} \Shom(\CR_\chi,\CK_{S/\CO}).\end{equation}

Assuming as usual that $U$ is sufficiently small, and in particular $\alpha - 1 \in \gP$ for all
$\alpha \in U \cap F^\times$, the action of $\CO_{F,(p),+}^\times$ on $\widetilde{Y}_{U_1(\gP)}$ via
multiplication on quasi-polarizations factors through a free action of
 $\CO_{F,(p),+}^\times/(U\cap F^\times)^2$, and we let $Y_{U_1(\gP)}$ denote the quotient scheme.
 The morphism $\widetilde{f}$ descends to a finite flat morphism
 $f:Y_{U_1(\gP)} \to Y_{U_0(\gP)}$, so $Y_{U_1(\gP)}$ is
Cohen--Macaulay and quasi-projective of relative dimension $d$ over $\CO$.
Furthermore for sufficiently small $U$, $U'$ of level prime to $p$ and $g \in \GL_2(\A_{F,\f}^{(p)})$
such that $g^{-1}Ug \subset U'$, we obtain a finite \'etale
morphism $\tilde{\rho}_g:\widetilde{Y}_{U_1(\gP)} \to \widetilde{Y}_{U'_1(\gP)}$, descending
to ${\rho}_g:{Y}_{U_1(\gP)}  \to{Y}_{U'_1(\gP)}$ and satisfying the usual compatibilities.

Note that the assumption on $U$ implies that the
canonical isomorphism $\nu^*\underline{A} \cong \underline{A}$ over $\widetilde{Y}_{U_0(\gP)}$
is the identity on $H$ for all $\nu \in (U\cap F^\times)^2$, so $H$ descends to a Raynaud
$\CO_F/\gP$-module scheme on $Y_{U_0(\gP)}$, which we also denote by $H$.  Furthermore
the line bundles $\CR_\tau$ and morphisms $s_\tau$, $t_\tau$
all descend to give the same descriptions of $H$, $H^\vee$ and $Y_{U_1(\gP)}$
as spectra of finite flat $\CO_{Y_{U_0(\gP)}}$-algebras (with $(\CO_F/\gP)^\times$-actions)
over $Y_{U_0(\gP)}$, and we obtain analogous decompositions of $f_*\CO_{Y_{U_1(\gP)}}$
and $f_*\CK_{Y_{U_1(\gP)}/\CO}$.


\subsection{Dicing}
We will use the method of ``dicing,'' introduced in \cite[\S6.2]{DKS}, to reduce the proofs
of the cohomological vanishing results to consideration
of line bundles on irreducible components\footnote{We could alternatively have
used $\overline{Y}_0(\gP)$ and $\overline{Y}_0(\gP)_J$, as in \cite{DKS}, but for the purpose of
proving the desired vanishing results, we can replace them by \'etale covers.}
 of the special fibre of $\widetilde{Y}_{U_0(\gP)}$, or equivalently
of the schemes denoted $S_J$.

Recall from \S\ref{subsec:components} that for any subset $J \subset \Sigma_{\gP}$, 
the closed subscheme $S_J$ of $\overline{S} := \widetilde{Y}_{U_0(\gP),k}$ is defined by the vanishing of 
$$\{ \,\psi_{\theta,\CM}^*\,|\, \theta \in J \,\} \quad \cup \{\,\psi_{\theta,\CL}^*\,|\, \theta\not\in J\,\},$$
where $\psi:A \to A'$ is the universal isogeny over $\overline{S}$ and 
$$\psi_{\theta,\CM}: \CM_\theta' \to \CM_\theta,\qquad \psi_{\theta,\CL}^*: \CL_\theta' \to \CL_\theta$$
are the morphisms induced by $\psi$ on the indicated subquotients of de Rham cohomology sheaves of
$A$ and $A'$.

Following \cite[\S6.2]{DKS}, we consider the sheaves of ideals $\CI_j \subset \CO_{\overline{S}}$ corresponding
to the closed subscheme of $\overline{S}$ defined by the vanishing of 
$$\left\{ \left. \, \bigotimes_{\theta \in J}  \psi_{\CL,\theta}^* \,\right|\, J \subset \Sigma_{\gP},\,|J| = j \,\right\},$$
or equivalently the image of the morphism
$$\bigoplus_{|J| = j} \left( \bigotimes_{\theta\in J} \CL_\theta^{-1}\CL_\theta' \right)  \,\longrightarrow \, \CO_{\overline{S}}$$
induced by the $\psi_{\theta,\CL}^*$.   We thus have
$$ \CO_{\overline{S}} = \CI_0 \supset \CI_1 \supset \CI_2 \supset \cdots \CI_{d_\gP} \supset \CI_{1+d_\gP} = 0,$$
where $d_{\gP} = |\Sigma_{\gP}| = \sum_{\gp|\gP} e_\gp f_\gp$.
Letting $i_J$ denote the closed immersion $S_J \hookrightarrow \overline{S}$
and $\CI_J$ the sheaf of ideals on $S_J$ defining the vanishing locus of $\bigotimes_{\theta\in J} \psi_{\theta,\CL}^*$,
we see exactly as in \cite[\S6.2.1]{DKS} that $\CI_J$ is invertible and 
the proof of \cite[Lemma~6.2.1]{DKS} carries over {\em mutatis mutandis} to give the following:
\begin{lemma}  \label{lem:grO}  The natural map $\gr^j \CO_{\overline{S}} \longrightarrow \bigoplus_{|J| = j} i_{J,*} \CI_J$ is
an isomorphism for all $j=0,\ldots, d_\gP$.
\end{lemma}

Similarly,  letting $\CJ_J$ denote the (invertible) sheaf of ideals on $S_J$ defining the vanishing locus of\footnote{Since 
$\psi_{\theta,\CM}^*:\CM_\theta' \to \CM_\theta$ is an isomorphism for $\theta\not\in \Sigma_{\gP}$, we may view the tensor
product as being over either $\Sigma - J$ or $\Sigma_\gP - J$, and we write simply $\theta\not\in J$.}
$\bigotimes_{\theta\not\in J} \psi_{\theta,\CM}^*$, the same argument as in \cite[\S6.2.2]{DKS} yields
the following description of the graded pieces of the induced filtration 
$\Fil^j\CK_{\overline{S}/k} := \CI_j\CK_{\overline{S}/k}$ on the dualizing sheaf:
\begin{lemma}  \label{lem:grK} The natural map $\gr^j \CK_{\overline{S}/k} \longrightarrow \bigoplus_{|J| = j} i_{J,*} (\CJ_J^{-1}\CK_{S_J/k})$ is
an isomorphism for all $j=0,\ldots, d_\gP$.
\end{lemma}

Furthermore letting $\overline{S}_1 = \widetilde{Y}_{U_1(\gP),k}$, the decompositions (\ref{eqn:charO}) and
(\ref{eqn:charK}) restrict to give the decompositions
\begin{equation} \label{eqn:char2} \widetilde{f}_* \CO_{\overline{S}_1} = \bigoplus_\chi \CR_\chi  \quad\mbox{and} \quad
 \widetilde{f}_* \CK_{\overline{S}_1/k} = \bigoplus_\chi \CR_\chi^{-1} \CK_{\overline{S}/k},\end{equation}
where we again write $\widetilde{f}:\overline{S}_1 \to \overline{S}$ for the natural projection and
$\CR_\chi$ for the Raynaud bundles associated to $H = \ker(\psi)$ over $\overline{S}$.
Note also that we may now view the decompositions as being over characters 
$\chi:(\CO_F/\gP)^\times \to k^\times$.

Following the method of \cite[\S7.3]{DKS} will reduce the proof of the desired cohomological vanishing to a computation involving
the line bundles $\CI_J$, $\CJ_J$, $\CK_{{S}_J/k}$ and (the pull-back of) the $\CR_\tau$ over $S_J$.
We will now describe all of these in terms of the line bundles $\CL_\theta$, $\CL_\theta'$, $\CM_\theta$
and $\CM_\theta'$.

We start with the Raynaud bundles $\CR_\tau$ over $S_J$ for $\tau = \tau_{\gp,i} \in \Sigma_{\gP,0}$.
Suppose first that $\theta := \theta_{\gp,i,e_\gp} \not\in J$, so that $\psi_{\CL,\theta}^*$ vanishes on $S_J$,
and therefore 
$$\psi_\tau^* : \CH^1_\dr(A'/S_J)_\tau  \longrightarrow \CH^1_\dr(A/S_J)$$
 maps $\CF^{\prime(e_\gp)}_\tau$ onto $\CF^{e_\gp-1}_\tau$
(where we are now writing $(\underline{A},\underline{A}',\psi)$ for the universal object
over $S_J$).  It follows that the cokernel of $(s'_*\Omega^1_{A'/S_J})_\tau \to (s_*\Omega^1_{A/S_J})_\tau$ is 
isomorphic to the invertible sheaf $\CL_\theta = \CF_\tau^{(e_\gp)}/\CF_\tau^{(e_\gp-1)}$ on $S_J$, so
$$ \Lie(H/S_J)_\tau  = \ker( \Lie(A/S_J)_\tau \longrightarrow \Lie(A'/S_J)_\tau)$$
is isomorphic to $\CL_\theta^{-1}$, and therefore $i_J^* \CR_\tau \cong\CL_\theta$
(as in the proof of \cite[Lem.~5.1.1]{DKS}).

Similarly if $\theta = \theta_{\gp,i,e_\gp} \in J$, then $\CG^{\prime(e_\gp)}_\tau$ is the preimage of
$\CF^{(e_\gp)}_\tau$ under $\psi_\tau^*$, so the kernel of $R^1s'_*\CO_{A'} \to R^1s_*\CO_A$
is isomorphic to $\CM_\theta' = \CG_\tau^{\prime(e_\gp)}/\CF_\tau^{\prime(e_\gp)}$.
It follows that
$$ \Lie(H^\vee/S_J)_\tau  \cong \ker( \Lie((A')^\vee/S_J)_\tau \longrightarrow \Lie(A^\vee/S_J)_\tau)$$
is isomorphic to $\CM_\theta'$, and therefore so is $i_J^*\CR_\tau$.  Summing up, we have constructed
a canonical isomorphism
\begin{equation}  \label{eqn:R}
i_J^* \CR_{\tau}   \,\,\cong \,\, \left\{\begin{array}{ll}  \CL_\theta,&  \mbox{if $\theta \not\in J$.}\\
             \CM_\theta'& \mbox{if $\theta \in J$},\end{array}
 \right.
\end{equation}
where $\tau = \tau_{\gp,i} \in \Sigma_{\gP,0}$ and $\theta = \theta_{\gp,i,e_\gp}$.

Next considering $\CI_J$, its definition (and invertibility) imply that $\bigotimes_{\theta\in J} \psi_{\theta,\CL}^*$ induces
an isomorphism
\begin{equation}  \label{eqn:I}
\bigotimes_{\theta\in J}  \CL_\theta^{-1} \CL_\theta' \,\, \stackrel{\sim}{\longrightarrow} \,\,  \CI_J.\end{equation}
Similarly, we have that 
$\bigotimes_{\theta\not\in J} \psi_{\theta,\CM}^*$ induces
an isomorphism
\begin{equation} \label{eqn:J} \bigotimes_{\theta\not\in J}  \CM_\theta^{-1} \CM_\theta' \,\, \stackrel{\sim}{\longrightarrow} \,\,  \CJ_J,\end{equation}
but we will need to further relate the $\CM_\theta'$ (for $\theta \in \Sigma_{\gP} - J$) to other line bundles on $S_J$.

Let $\theta = \theta_{\gp,i,j}$ for some $\gp|\gP$, and let $\tau = \tau_{\gp,i}$,
and suppose first that $\theta,\sigma^{-1}\theta\not\in J$.
If $j > 1$, then consider the commutative diagram
$$\xymatrix{  \CP_\theta'  \ar[r]^{\psi_\theta^*}  \ar[d]_{u} & \CP_\theta \ar[d]_{u} \\
\CP_{\sigma^{-1}\theta}' \ar[r]^{\psi_{\sigma^{-1}\theta}^*} & \CP_{\sigma^{-1}\theta} .}$$
Note that the image of $u: \CP_{\theta}' \to \CP_{\sigma^{-1}\theta}'$ is $\CL_{\sigma^{-1}\theta}' = \ker(\psi_{\sigma^{-1}\theta}^*)$ 
(since $\sigma^{-1}\theta\not\in J$), so the image of $\psi_\theta^*$ is contained in the kernel of
$u: \CP_{\theta} \to \CP_{\sigma^{-1}\theta}$, which is $\CM_{\sigma^{-1}\theta} = \CG_\tau^{(j-1)}/\CF_\tau^{(j-1)}$.
Since each is a rank one subbundle of $\CP_{\sigma^{-1}\theta}$, it follows that
$$\CM_{\sigma^{-1}\theta} = \im(\psi_\theta^*) \cong  \CP'_\theta/\ker(\psi_\theta^*) = \CP_\theta'/\CL_\theta' =
\CM_\theta'.$$
Similarly if $j=1$, then we consider instead the commutative diagram
$$\xymatrix{  \CP_\theta'  \ar[r]^{\psi_\theta^*}  \ar[d] & \CP_\theta \ar[d] \\
\CP_{\sigma^{-1}\theta}^{\prime(p)}  \ar[r]^{\psi_{\sigma^{-1}\theta^{*(p)}} }& \CP_{\sigma^{-1}\theta}^{(p)},}$$
where the vertical maps are defined by the partial Hasse invariant (see (\ref{eqn:hdef}) and (\ref{eqn:hdef2})).
The same reasoning as above then shows that $\CM_\theta'$ is isomorphic to the kernel of
$$\widetilde{h}_\theta : \CP_\theta \longrightarrow \CL^{(p)}_{\sigma^{-1}\theta},$$
which in this case is isomorphic to $\CM_{\sigma^{-1}\theta}^{(p)}$, via the morphism 
induced by the restriction of $\Frob^*$ to 
$$ (\CG_{\phi^{-1}\circ\tau}^{(e_\gp)})^{(p)} = u^{-1} (\CF^{(e_\gp-1)}_{\phi^{-1}\circ\tau})^{(p)}  \longrightarrow \CH^1_\dr(A/S_J)_{\tau} [u] = \CP_\theta.$$

Suppose on the other hand that $\theta \not\in J$, but $\sigma^{-1}\theta \in J$, so $\theta \in J''$.  
Since $\underline{A}''$ corresponds to an element of $T_{J''}(S_J)$, we have that the kernel of
$$\widetilde{h}'_\theta : \CP'_\theta \longrightarrow \CL^{\prime(n_\theta)}_{\sigma^{-1}\theta}$$
is $\CL'_\theta$, and hence that $\CM_\theta'$ is isomorphic to $ \CL^{\prime(n_\theta)}_{\sigma^{-1}\theta}$.

Summing up, we have shown that if $\theta \in \Sigma_\gP - J$, then there is a canonical isomorphism
\begin{equation} \label{eqn:M}
\CM_\theta' \,\, \cong \,\,  \left\{\begin{array}{ll}
\CM_{\sigma^{-1}\theta}^{(n_\theta)}, & \mbox{if $\sigma^{-1}\theta\not\in J$} \\
\CL_{\sigma^{-1}\theta}^{\prime(n_\theta)}, &\mbox{if $\sigma^{-1}\theta \in J$}.\end{array}\right.
\end{equation}

Finally, although we will not need the following description of the line bundles $\CK_{S_J/k}$, we provide it
for completeness and coherence with the theme of Kodaira--Spencer isomorphisms.  Recall that the
isomorphism of Theorem~\ref{thm:KSIwa}  arises from a canonical isomorphism
$$\CK_{\widetilde{Y}_{U_0(\gP)}/\CO}  \,\, \cong \,\,   \bigotimes_{\theta\in \Sigma} 
\left(\widetilde{\CM}_\theta^{\prime -1}\widetilde{\CL}_\theta\right)$$
(using\, $\widetilde{\cdot}$\, for universal data over $\widetilde{Y}_{U_0(\gP)}$).   On the other
hand the proof of Lemma~\ref{lem:grK} (see \cite[\S6.2.2]{DKS}) also gives a canonical isomorphism
$$\CK_{S_J/k} \,\, \cong \,\,  \CI_J \CJ_J i_J^* \CK_{\overline{S}/k}.$$
Combining these with (\ref{eqn:I}) and (\ref{eqn:J}), we obtain a canonical isomorphism
\begin{equation}  \label{eqn:K}
\CK_{S_J/k} \,\, \cong \,\,  
\left( \bigotimes_{\theta \in \Sigma - J} \CM_\theta^{-1} \CL_\theta \right)  {\textstyle\bigotimes} 
\left( \bigotimes_{\theta \in J}  \CM_\theta^{\prime -1}  \CL_\theta' \right).
 \end{equation}
(Recall that for $\theta \not\in \Sigma_\gP$, we have $\CL_\theta \cong \CL_\theta'$ and $\CM_\theta \cong \CM_\theta'$,
so these factors could just as well have been included with $\theta \in J$ instead of $\theta \in \Sigma_\gP - J$.)

We remark that (\ref{eqn:K}) can also be proved more directly using a deformation-theoretic argument, which
furthermore produces a filtration on $\Omega^1_{S_J/k}$ analogous to the one on $\Omega^1_{Y_U/\CO}$ defined
by Reduzzi and Xiao in \cite[\S2.8]{RX} (or more precisely \cite[\S3.3]{theta} for the current set-up).

We remark also that the various isomorphisms of line bundles established in this section should be compatible
in the usual sense with descent data and Hecke action, but we will not need this for the purpose of proving the
desired cohomological vanishing theorems.

\subsection{Vanishing, duality and flatness results}  \label{sec:vanishing}
We are now ready to prove the main result, which is a generalization of Theorem E of \cite{DKS}.
Indeed all the ingredients are now in place to apply the same argument as in \cite[\S7.3]{DKS}.
Recall that $U$ is any sufficiently small open compact subgroup of $\GL_2(\CO_{F,p})$, $\gP$
is the product of any set of primes of $\CO_F$ dividing $p$, and
$$Y_{U_1(\gP)} \stackrel{f}{\longrightarrow}  Y_{U_0(\gP)} \stackrel{\pi_1}\longrightarrow Y_U$$
are the natural degeneracy maps induced by the forgetful morphisms
$\widetilde{f}: \widetilde{Y}_{U_1(\gP)} \to \widetilde{Y}_{U_0(\gP)}$ and $\widetilde{\pi}_1: \widetilde{Y}_{U_0(\gP)} \to \widetilde{Y}_U$.
We will write simply $Y$ (resp.~$Y_0(\gP)$, $Y_1(\gP)$) for $Y_U$ (resp.~$Y_{U_0(\gP)}$, $Y_{U_1(\gP)}$) when $Y$ is fixed,
and similarly abbreviate $\widetilde{Y}_U$, etc.
We let $\widetilde{\varphi}$ (resp.~$\varphi$) denote the composite $\widetilde{\pi}_1\circ \widetilde{f}$ (resp.~$\pi_1\circ f$).

\begin{theorem} \label{thm:main} The higher direct image sheaves
$$R^i \varphi_*\CK_{Y_1(\gP)/\CO} \quad\mbox{and}\quad R^i \varphi_*\CO_{Y_1(\gP)}$$
vanish for all $i > 0$.
\end{theorem}
\begpf  First note that since $\widetilde{Y} \to Y$ is \'etale and $\widetilde{Y}_1(\gP)
 = \widetilde{Y} \times_{Y} Y_1(\gP)$, we may replace $Y$ (resp.~$Y_1(\gP)$, $\varphi$)
 by $\widetilde{Y}$ (resp.~$\widetilde{Y}_1(\gP)$, $\widetilde{\varphi}$).  Next note that since
 $\widetilde{\varphi}$ is projective and $\widetilde{\varphi}_K$ is finite, we may replace these in turn 
  by their special fibres.
  Furthermore since $\widetilde{f}$ is finite, we have
  $$R^i\widetilde{\varphi}_*\CO_{\overline{S}_1} = R^i \widetilde{\pi}_{1,*} (\widetilde{f}_* \CO_{\overline{S}_1})\quad\mbox{and}\quad
      R^i\widetilde{\varphi}_*\CK_{\overline{S}_1/k} = R^i \widetilde{\pi}_{1,*} (\widetilde{f}_* \CK_{\overline{S}_1/k})$$
  (where we again let $\overline{S}_1 = \widetilde{Y}_1(\gP)_k$ and
suppress the subscript $k$ from the notation for the morphisms).
  Using the decompositions given in (\ref{eqn:charK}) and the filtrations whose graded pieces are
  described in Lemmas~\ref{lem:grO} and~\ref{lem:grK}, we are therefore reduced to proving the vanishing of
  $$R^i \pi_{J,*} (\CI_J i_J^*\CR_\chi)  \quad \mbox{and} \quad R^i \pi_{J,*} ((\CJ_J i_J^* \CR_\chi)^{-1}\CK_{S_J/k})$$
  for all $i > 0$, $J \subset \Sigma_\gP$ and $\chi: (\CO_F/\gP)^\times \to k^\times$,  where $i_J$ is the closed
  immersion $S_J \hookrightarrow \overline{S}$ and $\pi_J:S_J \to T_{J'}$ is the restriction of $\widetilde{\pi}_1$.
  
  Let us now fix $J$ and $\chi$, and write $\chi = \prod \overline{\tau}^{m_\tau}$, where the product is over
  $\tau \in \Sigma_{\gP,0}$, $0 \le m_\tau \le p -1$ for each$\tau$, and $m_\tau < p - 1$ for some $\tau \in \Sigma_{\gp,0}$
  for each $\gp|\gP$.  We now apply the results of the preceding section to write the line bundles 
  $\CI_J$, $\CJ_J$ and $i_J^*\CR_\chi$ on $S_J$ in terms of $\CL_\theta$, $\CM_\theta$, $\CL_\theta'$ and $\CM_\theta'$.
  
Firstly by formula (\ref{eqn:R}), we have
 \begin{equation} \label{eqn:R2} i_J^*\CR_\chi \cong \left(\bigotimes_{\theta\in \Sigma_\gP - J}   \CL_\theta^{\otimes m_\theta}  \right)  {\textstyle\bigotimes}
                          \left(\bigotimes_{\theta \in J}   \CM_\theta^{\prime\otimes m_\theta}\right),\end{equation}
where $m_\theta = m_\tau$ if $\theta = \theta_{\gp,i,e_\gp}$ and $\tau = \tau_{\gp,i}$, and $m_\theta = 0$
if $\theta = \theta_{\gp,i,j}$ for some $j < e_\gp$.  Combining (\ref{eqn:R2}) with (\ref{eqn:I}) and the isomorphisms
$\CL'_\theta \otimes \CM'_\theta = \CN'_\theta \cong \CN_\theta = \CL_\theta \otimes \CM_\theta$
(see Remark~\ref{rmk:symmetric}), we obtain an isomorphism
\begin{equation} \label{eqn:IR}   \CI_J i_J^*\CR_\chi \cong 
  \left(\bigotimes_{\theta\in \Sigma_\gP - J}   \CL_\theta^{\otimes m_\theta}  \right)  {\textstyle\bigotimes}
                          \left(\bigotimes_{\theta \in J} (  \CM_\theta'^{\otimes (m_\theta-1)}  \otimes \CM_\theta )\right).\end{equation}

On the other hand combining (\ref{eqn:R2}) with (\ref{eqn:J}) and (\ref{eqn:M}) gives
$$ \begin{array}{c} \CJ_J i_J^*\CR_\chi   \cong 
  \left(\displaystyle\bigotimes_{\theta,\sigma\theta \in \Sigma_\gP - J} 
   \CM_\theta^{\otimes( n_{\sigma\theta} - 1) }\otimes     \CL_\theta^{\otimes m_\theta}  \right) 
  \bigotimes
    \left(\displaystyle\bigotimes_{\theta\not\in J,\sigma\theta \in J}   \CM_\theta^{-1}  \CL_\theta^{\otimes m_\theta} \right) \\
       \bigotimes \left(\displaystyle\bigotimes_{\theta\in J,\sigma\theta \not\in J}  ( \CM_\theta'^{\otimes m_\theta} \otimes \CL_\theta'^{\otimes n_{\sigma\theta}}) \right)  
   \bigotimes     \left(\displaystyle\bigotimes_{\theta,\sigma\theta \in J}   \CM_\theta'^{\otimes m_\theta}  \right).\end{array} $$
Combined with the isomorphisms $\CN_\theta\cong \CN'_\theta$ and $\CN_\theta^{\otimes n_{\sigma\theta}} \cong \CN_{\sigma\theta}$
(see the last paragraph of \cite[\S4.1]{theta} for the latter), this simplifies to
\begin{equation} \label{eqn:JR} \begin{array}{c}
\CJ_J i_J^*\CR_\chi \cong
   \left(\displaystyle\bigotimes_{\theta,\sigma\theta\in \Sigma_\gP - J}   \CL_\theta^{\otimes (m_\theta-n_{\sigma\theta}+1)}  \right)  \bigotimes
      \left(\displaystyle\bigotimes_{\theta\not\in J,\sigma\theta\in J}   \CL_\theta^{\otimes (m_\theta+1)}  \right)  \\ \bigotimes
      \left(\displaystyle\bigotimes_{\theta\in J,\sigma\theta\not\in J}   \CM_\theta'^{\otimes (m_\theta-n_{\sigma\theta})}  \right)  \bigotimes
      \left(\displaystyle\bigotimes_{\theta,\sigma\theta\in J}   \CM_\theta'^{\otimes m_\theta}  \right).\end{array}\end{equation}
      
Since $\pi_J$ is flat and projective, it suffices to prove the vanishing of cohomology of fibres, i.e., that
$$H^i(Z,j^*(\CI_J i_J^*\CR_\chi)) = 0\quad\mbox{and}\quad H^i(Z,j^*((\CJ_J i_J^*\CR_\chi)^{-1}\CK_{S_J/k}))) = 0$$
for all $x \in T_{J'}(\Fpbar)$, where $j:Z \to S_J$ is the closed immersion of the fibre of $\pi_J$ over $x$.  Furthermore since $\pi_J$ is Cohen-Macaulay,
we have $j^*\CK_{S_J/k} = \CK_{Z/\Fpbar}$.  Since the line bundles $j^*\CL_\theta = L_{\theta,x}\otimes_{\Fpbar}\CO_Z$ and $j^*\CM_\theta = M_{\theta,x}\otimes_{\Fpbar}\CO_Z$ are trivial
on $Z$, the isomorphisms (\ref{eqn:IR}) and (\ref{eqn:JR}) give
\begin{equation} \label{eqn:IJR} j^*(\CI_J i_J^*\CR_\chi) \cong
                          \bigotimes_{\theta \in J} j^*\CM_\theta'^{\otimes (m_\theta-1)}
                          \quad\mbox{and}\quad j^*(\CJ_J i_J^*\CR_\chi)^{-1} \cong
                          \bigotimes_{\theta \in J} j^*\CM_\theta'^{\otimes (\ell_\theta)},\end{equation}
where $\ell_\theta = n_{\sigma\theta} - m_\theta$ if $\sigma\theta \not\in J$ and $\ell_\theta = - m_\theta$ if $\sigma\theta \in J$.
 
We now appeal to the description of $Z$ in Theorem~\ref{thm:fibres}, which provides an isomorphism
$$Z\,\, \stackrel{\sim}{\longrightarrow} \,\, \Spec C \times  \prod_{\theta \in J''}  \PP^1$$
for a finite (Gorenstein) $\Fpbar$-algebra $C$, such that if $\theta \in J$, then
$j^*\CM_\theta'$ corresponds to the pull-back of $\CO(1)_{\sigma\theta}$ 
(resp.~an invertible $C$-module) if $\sigma\theta \not\in J$ (resp.~$\sigma\theta\in J$).
It therefore follows from (\ref{eqn:IJR}) that
$$H^i(Z,\,j^*(\CI_J i_J^*\CR_\chi))   \cong 
H^i\left( \prod_{\theta \in J''}  \PP^1 , \bigotimes_{\theta \in J''} \CO(m_{\sigma^{-1}\theta} - 1)_\theta \right) \otimes_{\Fpbar} C,$$
which vanishes since $m_{\sigma^{-1}\theta} \ge 0$ for all $\theta \in J''$.  Similarly since $\CK_{Z/\Fpbar}$ corresponds to
the dualizing sheaf of 
$\Spec C \times  \prod_{\theta \in J''}  \PP^1$, which is isomorphic to $\otimes_{\theta\in J''} \CO(-2)_\theta$, we have
$$H^i(Z,\,j^*((\CJ_J i_J^*\CR_\chi)^{-1}\CK_{S_J/k})))  \cong 
H^i\left( \prod_{\theta \in J''}  \PP^1 , \bigotimes_{\theta \in J''} \CO(\ell_{\sigma^{-1}\theta} - 2)_\theta \right) \otimes_{\Fpbar} C,$$
which vanishes since $\ell_{\sigma^{-1}\theta} = n_\theta - m_{\sigma^{-1}\theta} \ge 1$ for all $\theta \in J''$.
\epf
 
\begin{corollary}  \label{cor:Iwa}
 The higher direct image sheaves
$$R^i \pi_{j,*}\CK_{Y_0(\gP)/\CO} \quad\mbox{and}\quad R^i \pi_{j,*}\CO_{Y_0(\gP)}$$
vanish for all $i > 0$ and $j=1,2$.
\end{corollary}
\begpf  As in the proof of the theorem, we may replace $\pi_j$ by 
$\widetilde{\pi}_j: \widetilde{Y}_0(\gP) \to \widetilde{Y}$.
By formula (\ref{eqn:charO}), the line bundle
$\CO_{\widetilde{Y}_0(\gP)}$ is a direct summand of $\widetilde{f}_*\CO_{\widetilde{Y}_1(\gP)}$,
so the vanishing of 
$$R^i\widetilde{\varphi}_*\CO_{\widetilde{Y}_1(\gP)} = R^i\widetilde{\pi}_{1,*}(\widetilde{f}_*\CO_{\widetilde{Y}_1(\gP)})$$
for all $i > 0$ implies that of $R^i\widetilde{\pi}_{1,*}\CO_{\widetilde{Y}_0(\gP)}$.
Similarly the vanishing of $R^i\widetilde{\pi}_{1,*}\CK_{\widetilde{Y}_0(\gP)/\CO}$ follows from
(\ref{eqn:charK}) and the theorem.

For the analogous assertions for $\widetilde{\pi}_2$, note that there is an automorphism 
$\widetilde{w}_{\gP}$ of $\widetilde{Y}_0(\gP)$ such that
$\widetilde{\pi}_2 = \widetilde{\pi}_1\circ \widetilde{w}_\gP$.  Indeed define $\widetilde{w}_{\gP}$ by the triple
$(\underline{A}',\gP^{-1}\otimes_{\CO_F}\underline{A},\xi)$, where $(\underline{A},\underline{A}',\psi)$ is the
universal triple over $\widetilde{Y}_0(\gP)$, $\gP^{-1} \otimes_{\CO_F}  \underline{A}$ denotes $\gP^{-1}\otimes_{\CO_F} A$
endowed with the evident auxiliary data, in particular quasi-polarization 
$$\varpi_\gP^2 \otimes \lambda : \gP^{-1} \otimes_{\CO_F} A \to \gP \otimes_{\CO_F} A^\vee = (\gP^{-1} \otimes_{\CO_F} A)^\vee,$$
and $\xi:A' \to \gP^{-1} \otimes_{\CO_F} A$ is such that $\xi\circ\psi$ is the
canonical isogeny $A \to \gP^{-1} \otimes_{\CO_F} A$.  The desired vanishing of higher direct images under $\widetilde{\pi}_2$
therefore follows from the corresponding assertions for $\widetilde{\pi}_1$ and the identifications
$\widetilde{w}_{\gP,*}\CO_{\widetilde{Y}_0(\gP)} = \CO_{\widetilde{Y}_0(\gP)}$ and
$\widetilde{w}_{\gP,*}\CK_{\widetilde{Y}_0(\gP)/\CO} = \CK_{\widetilde{Y}_0(\gP)/\CO}$ 
\epf

\begin{remark} \label{rmk:caution} We caution that the automorphism $w_\gP$ of $\widetilde{Y}_0(\gP)$ does {\em not} lift to
an automorphism of $\widetilde{Y}_1(\gP)$, so the argument handling the case of $\pi_2$ in
in the proof of the corollary cannot be used to deduce of vanishing of higher direct images under $\pi_2 \circ f$
from the theorem.
\end{remark}

\begin{corollary}  \label{cor:flat1} 
The direct image sheaves $\varphi_*\CK_{Y_1(\gP)/\CO}$ and $\varphi_*\CO_{Y_1(\gP)}$ are locally free
of rank $\prod_{\gp|\gP} (p^{2f_\gp} - 1)$ over $\CO_{Y}$.  Furthermore there is a perfect Hecke-equivariant\footnote{In the sense that diagram~(\ref{eqn:Hecke2}) in the proof commutes.} pairing
$$\varphi_*\CK_{Y_1(\gP)/\CO} \, \otimes_{\CO_{Y}} \, \varphi_*\CO_{Y_1(\gP)}
\,\longrightarrow \, \CK_{Y/\CO}.$$
\end{corollary}
\begpf  We apply Grothendieck--Serre duality to the proper morphism $\varphi$.  Since $R^i\varphi_*\CO_{Y_1(\gP)}$
vanishes for all $i > 0$, the duality isomorphism of \cite[\S0AU3(4c)]{stacks} (taking $K$ there to be $\CO_{Y_1(\gP)}[d]$) degenerates to
$$R\varphi_* \CK_{Y_1(\gP)/\CO}  \,\,\stackrel{\sim}{\longrightarrow} \,\,  R\Shom_{\CO_{Y}} (\varphi_*\CO_{Y_1(\gP)},  \CK_{Y/\CO}).$$
In particular we obtain an isomorphism
\begin{equation}
\label{eqn:duality}
\varphi_* \CK_{Y_1(\gP)/\CO}  \,\,\stackrel{\sim}{\longrightarrow} \,\,  \Shom_{\CO_{Y}} (\varphi_*\CO_{Y_1(\gP)},  \CK_{Y/\CO}).\end{equation}
Furthermore since $R^i\varphi_* \CK_{Y_1(\gP)/\CO}$ vanishes for all $i > 0$, it follows that so does
$\Sext_{\CO_{Y}}^i (\varphi_*\CO_{Y_1(\gP)},\CK_{Y/\CO})$, and hence
$$ \Sext^i_{\CO_{Y}} (\varphi_*\CO_{Y_1(\gP)},\CO_{Y})  = 0 $$
for all $i > 0$.  Since $Y$ is regular, it follows that $\varphi_*\CO_{Y_1(\gP)}$ is locally free,
with rank given by the degree of the finite morphism $\varphi_K$, namely
$$[U:U_1(\gP)] = \prod_{\gp|\gP}  [\GL_2(\CO_{F,\gp}):I_1(\gp)]  = \prod_{\gp|\gP}  (p^{2f_\gp} - 1).$$
The same argument with the roles of $\CK_{Y_1(\gP)/\CO}$ and $\CO_{Y_1(\gP)}$ exchanged
shows that $\varphi_*\CK_{Y_1(\gP)/\CO}$ is locally free of the same rank.

Furthermore it follows from the construction of the duality isomorphism that (\ref{eqn:duality}) is induced by
the trace morphism denoted $\Tr_{\varphi,\CK_{Y/\CO}}$ in \cite[\S0AWG]{stacks}, which may be viewed as
a morphism
\begin{equation}
\label{eqn:trace} \Tr:  \varphi_*\CK_{Y_1(\gP)/\CO}  \longrightarrow \CK_{Y/\CO}
\end{equation}
since $R^i\varphi_* \CK_{Y_1(\gP)/\CO}$ vanishes for all $i > 0$.

Suppose then that $U$ and $U'$ are sufficiently small and that $g \in \GL_2(\AA_{F,\f}^{(p)})$ is such that
$g^{-1} U' g \subset U'$, giving \'etale morphisms $\rho_g: Y_{U'} \to Y_{U}$ and
$\rho_g': Y_{U'_1(\gP)} \to Y_{U_1(\gP)}$ such that the diagram 
$$\xymatrix{ Y_{U'_1(\gP)} \ar[r]^{\rho_g'} \ar[d]_{\varphi'} & Y_{U_1(\gP)}\ar[d]^{\varphi}  \\
Y_{U'} \ar[r]^{\rho_g} & Y_U }$$
is Cartesian, and identifications $\rho_g^*\CK_{Y_U/\CO} = \CK_{Y_{U'}/\CO}$ and
$\rho_g'^*\CK_{Y_{U_1(\gP)}/\CO} = \CK_{Y_{U'_1(\gP)}/\CO}$.   The commutativity of the
resulting diagram
$$\xymatrix{ \rho_g^*\varphi_* \CK_{Y_{U_1(\gP)}/\CO} \ar[rr]^-{\rho_g^*(\Tr)} \ar[d]^{\wr} && \rho_g^*\CK_{Y_U/\CO} \ar@{=}[d] \\
\varphi'_*\rho'^*_g\CK_{Y_{U_1(\gP)/\CO}} \ar@{=}[r]  & \varphi'_*\CK_{Y_{U'_1(\gP)/\CO}}  \ar[r]^-{\Tr'} & \CK_{Y_{U'}/\CO} }$$
then follows from \cite[Lemma~0B6J]{stacks} and implies that of
$$\xymatrix{ \rho_g^*\varphi_*\CK_{Y_{U_1(\gP)}/\CO}\ar[r] \ar[d]^\wr & \Shom_{\CO_{Y_{U'}}} (\rho_g^*\varphi_*\CO_{Y_{U_1(\gP)}}, \rho_g^*\varphi_* \CK_{Y_{U_1(\gP)}/\CO}) \ar[r] \ar[d]^\wr &
\Shom_{\CO_{Y_{U'}}} (\rho_g^*\varphi_*\CO_{Y_{U_1(\gP)}}, \rho_g^* \CK_{Y_U/\CO}) \ar@{=}[d] 
\\ 
\varphi'_*\CK_{Y_{U'_1(\gP)}/\CO} \ar[r] \ar@{=}[d] &
\Shom_{\CO_{Y_{U'}}} (\rho_g^*\varphi_*\CO_{Y_{U_1(\gP)}}, \varphi'_* \CK_{Y_{U'_1(\gP)}/\CO}) \ar[r] &
\Shom_{\CO_{Y_{U'}}} (\rho_g^*\varphi_*\CO_{Y_{U_1(\gP)}}, \CK_{Y_{U'}/\CO})
\\
\varphi'_*\CK_{Y_{U'_1(\gP)}/\CO} \ar[r] &
\Shom_{\CO_{Y_{U'}}} (\varphi'_*\CO_{Y_{U'_1(\gP)}}, \varphi'_* \CK_{Y_{U'_1(\gP)}/\CO}) \ar[r] \ar[u]_\wr &
\Shom_{\CO_{Y_{U'}}} (\varphi'_*\CO_{Y_{U'_1(\gP)}}, \CK_{Y_{U'}/\CO}), \ar[u]_\wr}$$
where the vertical maps are base-change isomorphisms, the horizontal maps on the left are the canonical morphisms, and the ones on the right are composition with $\rho_g^*(\Tr)$ or $\Tr'$.
This in turn gives the
Hecke-equivariance of the pairing, i.e.,
the commutativity of 
\begin{equation} \label{eqn:Hecke2} \xymatrix{ \rho_g^*\varphi_*\CK_{Y_{U_1(\gP)}/\CO} \, \otimes_{\CO_{Y_{U'}}} \, \rho_g^*\varphi_*\CO_{Y_{U_1(\gP)}}
\ar[r]\ar[d]^{\wr} & \rho_g^* \CK_{Y_U/\CO} \ar@{=}[d] \\
\varphi'_*\CK_{Y_{U'_1(\gP)}/\CO} \, \otimes_{\CO_{Y_{U'}}} \, \varphi'_*\CO_{Y_{U'_1(\gP)}} \ar[r] &
\CK_{Y_{U'}/\CO}},\end{equation}
where the top (resp.~bottom) horizontal arrow is induced by $\Tr$ (resp.~$\Tr'$) and the vertical maps
are the canonical isomorphisms.
\epf

Recall that for sufficiently small $U$, the finite flat Raynaud $(\CO_F/\gP)$-module scheme $H$ over $\widetilde{Y}_0(\gP)$
descends to $Y_0(\gP)$, and hence so do the decompositions (\ref{eqn:charO}) and (\ref{eqn:charK}), giving
$$f_*\CO_{Y_1(\gP)} = \bigoplus_{\chi} \CR_\chi \quad\mbox{and}\quad 
f_*\CK_{Y_1(\gP)/\CO} = \bigoplus_{\chi}  \CR_\chi^{-1}\CK_{Y_0(\gP)/\CO},$$
where the line bundles $\CR_\chi$ (now on $Y_0(\gP)$) are indexed by the 
characters $\chi: (\CO_F/\gP)^\times \to \CO^\times$.   Theorem~\ref{thm:main} thus implies
the vanishing of $R^i\pi_{1,*}\CR_\chi$ and $R^i\pi_{1,*}(\CR_\chi^{-1}\CK_{Y_0(\gP)/\CO})$
for all $i > 0$ and characters $\chi$, and Corollary~\ref{cor:flat1} implies that
$\pi_{1,*}\CR_\chi$ and $\pi_{1,*}(\CR_\chi^{-1}\CK_{Y_0(\gP)/\CO})$
are locally free, now of rank $\prod_{\gp|\gP} (p^{f_\gp} + 1)$.  Furthermore the
same argument as in the proof of its Hecke-equivariance shows that the pairing defined
in the corollary is compatible with the natural action of $(\CO_F/\gP)^\times$ on $Y_1(\gP)$,
and hence respects the decompositions and restricts to perfect Hecke-equivariant pairings
$$\pi_{1,*}(\CR_{\chi^{-1}}^{-1}\CK_{Y_0(\gP)/\CO}) \, \otimes_{\CO_{Y}} \, \pi_{1,*} \CR_\chi
\longrightarrow \CK_{Y/\CO}$$
for all characters $\chi$. We record this in the case of the trivial character:
\begin{corollary}  \label{cor:flat0} 
The direct image sheaves $\pi_{1,*} \CK_{Y_0(\gP)/\CO}$ and $\pi_{1,*}\CO_{Y_0(\gP)}$ are locally free
of rank $\prod_{\gp|\gP} (p^{f_\gp} + 1)$ over $\CO_{Y}$, and there is a perfect Hecke-equivariant pairing
$$\pi_{1,*}\CK_{Y_0(\gP)/\CO} \, \otimes_{\CO_{Y}} \, \pi_{1,*}\CO_{Y_0(\gP)}
\,\longrightarrow \, \CK_{Y/\CO}.$$
\end{corollary}

\begin{remark}  \label{rmk:caution2} The same argument as in the proof of Corollary~\ref{cor:Iwa} shows that we can replace 
$\pi_1$ by $\pi_2$ in the statement of Corollary~\ref{cor:flat0}.  (In fact the automorphism $\widetilde{w}_\gP$
descends to an automorphism $w_\gP$ of $Y_0(\gP)$ such that $\pi_2 = \pi_1 \circ w_\gP$.) Note however that 
this argument does not allow one to replace $\pi_1$ by $\pi_2$ in the assertions for the twists by the bundles $\CR_\chi$
(cf.~Remark~\ref{rmk:caution}).
\end{remark}

Suppose now that $R$ is any (Noetherian) $\CO$-algebra and consider the base extensions
from $\CO$ to $R$ of the schemes $Y$ and $Y_i(\gP)$ for $i = 0,1$ and morphisms $\varphi$, $\pi_j$
for $j = 1,2$, which we denote $Y_{R}$, etc.  We then have:

\begin{corollary}  \label{cor:bc} For $\CE = \CO_{Y_1(\gP)}$ and $\CK_{Y_1(\gP)/\CO}$, 
the base-change morphisms
$$(R^i\varphi_*\CE)_R  \to R^i\varphi_{R,*} (\CE_R)$$
are isomorphisms for all $i \ge 0$.  In particular, $R^i\varphi_{R,*}(\CE_R) = 0$ for all $i > 0$,
$\varphi_{R,*}(\CE_R)$ is locally free over $\CO_{Y_R}$ of rank $\prod_{\gp|\gP} (p^{2f_\gp} - 1)$,
and there is a perfect Hecke-equivariant pairing
$$\varphi_{R,*}\CK_{Y_1(\gP)_R/R} \, \otimes_{\CO_{Y_{R}}} \, \varphi_{R,*}\CO_{Y_1(\gP)_R}
\,\longrightarrow \, \CK_{Y_{R}/R}.$$
Furthermore analogous assertions hold with $U_1(\gP)$ replaced by $U_0(\gP)$ and $\varphi$
replaced by $\pi_j$ for $j=1$ and $2$.
\end{corollary}
\begpf  Since $Y$ and $Y_1(\gP)$ are flat over $\CO$, the schemes $Y_{R}$ and $Y_1(\gP)$
are Tor independent over $Y$ (in the sense of \cite[Defn.~08IA]{stacks}).  The assertions in the case
of $U_1(\gP)$ are then immediate from \cite[Lem.~08IB]{stacks}, Theorem~\ref{thm:main}, 
Corollary~\ref{cor:flat1}, and the compatibility with base-change of formation of dualizing sheaves for 
Cohen--Macaulay morphisms (see \cite[Thm.~3.6.1]{conrad} or Lemmas~0E2Y 
and~0E9W of \cite{stacks}). The same argument applies for $U_0(\gP)$, but using
Corollaries~\ref{cor:Iwa} and~\ref{cor:flat0} and Remark~\ref{rmk:caution2}.
\epf

\begin{remark} \label{rmk:ERX3}
Our method also yields an improvement on the
cohomological vanishing results asserted in~\cite[Prop.~3.19]{ERX}.  More precisely, suppose that $\bk,\bm \in \ZZ^\Sigma$
and consider the automorphic line bundle $\widetilde{\CA}_{\bk,\bm}$ on $\widetilde{Y}$ (as defined in~\S\ref{ss:hecke}).
Applying the same analysis as above to the line bundle
$\widetilde{\pi}_2^* \widetilde{\CA}_{\bk,\bm}$
on $\widetilde{Y}_{U_0(\gP)}$ instead of $\CR_\chi$,
one obtains a filtration on $\widetilde{\pi}_2^* \widetilde{\CA}_{\bk,\bm,k}$ for which the graded pieces have fibres (relative to $\widetilde{\pi}_1$) described exactly as in the first formula
of (\ref{eqn:IJR}), but with the exponents\footnote{Note that the $m_\theta$ in (\ref{eqn:IJR}) are not the constituents of the
weight vector $\bm$ appearing in \S\ref{ss:hecke}.}
$m_\theta$ replaced by $n_{\sigma\theta}k_{\sigma\theta} - k_\theta$. (This refines a similar analysis carried out in 
\cite[\S4.9]{ERX} for reduced fibres at generic points, whose description
is completed by our Theorem~\ref{thm:reduced}; see Remark~\ref{rmk:ERX}.) The same argument as in the proof of
Theorem~\ref{thm:main} therefore
shows that if $n_{\sigma\theta}k_{\sigma\theta} \ge k_\theta$
for all $\theta \in \Sigma_{\gP}$, then 
$R^i\widetilde{\pi}_{1,*}(\widetilde{\pi}_2^*\widetilde{\CA}_{\bk,\bm})$ vanishes\footnote{The assertion in \cite{ERX} is that the codimension of
its support is at least $i+1$; stated in these terms, 
the codimension is in fact infinite.} for all $i > 0$. Furthermore the
same argument as in the proof of Corollary~\ref{cor:bc} shows
that, for such $\bk$, the base-change morphism 
$(\widetilde{\pi}_{1,*}\widetilde{\pi}_2^*\widetilde{\CA}_{\bk,\bm})_R \to \widetilde{\pi}_{1,*}(\widetilde{\pi}_2^*\widetilde{\CA}_{\bk,\bm,R})$ is an isomorphism
and $R^i\widetilde{\pi}_{1,*}(\widetilde{\pi}_2^*\widetilde{\CA}_{\bk,\bm,R}) = 0$ for all $i > 0$.
It follows also that if $\prod_\theta \theta(\mu)^{k_\theta+ 2 m_\theta}$ has trivial image in $R$ for all $\mu \in U \cap \CO_F^\times$, then $R^i{\pi}_{1,*}({\pi}_2^*{\CA}_{\bk,\bm,R}) = 0$ for all $i > 0$, and the base-change map
$({\pi}_{1,*}({\pi}_2^*{\CA}_{\bk,\bm,R}))_{R'}
\to {\pi}_{1,*}({\pi}_2^*{\CA}_{\bk,\bm,R'})$ is an isomorphism for all Noetherian $R$-algebras $R'$.
(Note however that the proof of Corollary~\ref{cor:flat1} does not carry over, as it would require analogous results for
$\widetilde{\CA}^{-1}_{\bk,\bm}\CK_{\widetilde{Y}_0(\gP)/\CO}$, for which the corresponding inequalities do not necessary hold.)

\end{remark}

We may also use the results above to construct integral models for Hilbert modular varieties of levels of
$U_0(\gP)$ and $U_1(\gP)$ which are finite and flat over the smooth models of level $U$.  Indeed consider
the Stein factorizations
$$Y_1(\gP)  \to Y'_1(\gP) \to Y \quad\mbox{and} \quad Y_0(\gP) \to Y'_0(\gP) \to Y,$$
where $Y'_1(\gP) := \SPEC (\varphi_*\CO_{Y_1(\gP)})$ and $Y'_0(\gP):= \SPEC( \pi_{1,*} \CO_{Y_0(\gP)})$.
Note that $Y_i(\gP)_K = Y'_i(\gP)_K$ for $i=1,2$ and that Corollaries~\ref{cor:flat1} and~\ref{cor:flat0}
immediately imply the following:

\begin{corollary}  \label{cor:models} The schemes $Y'_0(\gP)$ and $Y'_1(\gP)$ are finite and flat over $Y$; in particular
they are Cohen--Macaulay over $\CO$.
\end{corollary}

Note that the same conclusions apply to $Y''_0(\gP) := \SPEC( \pi_{2,*} \CO_{Y_0(\gP)})$; in fact $w_\gP$
induces an isomorphism $Y''_0(\gP) \stackrel{\sim}{\to} Y'_0(\gP)$ over $Y$.  Note also that the natural
map $Y'_1(\gP) \to Y$ factors through $Y'_0(\gP)$.  Furthermore if $U$, $U'$ are sufficiently small and 
$g \in \GL_2(\AA_{F,\f}^{(p)})$ is such that $g^{-1}U'g \subset U$, then (with the obvious notation) the morphisms
$\rho_g: Y_{U'_i(\gP)} \to Y_{U_i(\gP)}$ for $i=0,1$ induce morphisms $\rho_g':Y'_{U'_i(\gP)} \to Y'_{U_i(\gP)}$
satisfying the usual compatibilities.

Finally, we note the following consequence of Corollary~\ref{cor:models}, pointed out to us by G.~Pappas.
\begin{corollary}  \label{cor:normal}
The normalization of $Y$ in $Y_0(\gP)_K$ relative to $\pi_{1}$ (or to $\pi_{2}$) is flat over $Y$.
\end{corollary}
\begpf  First recall that $Y_0(\gP)$ is normal, since for example it is regular in codimension one and
Cohen--Macaulay.  Therefore it follows from \cite[Lem.~035L]{stacks} that $Y'_0(\gP)$ is normal.
Since $Y'_0(\gP)$ is also finite over $Y$, it is the normalization of $Y$ in $Y'_0(\gP)_K = Y_0(\gP)_K$
(relative to $\pi_{1}$),  and its flatness over $Y_U$ follows from Corollary~\ref{cor:models}.

The same argument applies with $\pi_1$ (resp.~$Y'_0(\gP)$) replaced by $\pi_2$ (resp.~$Y''_0(\gP)$).
\epf

\subsection{Hecke operators} \label{ss:hecke}
We now specialize to the case of $\gP = \gp$ and consider the morphisms
$$\pi_{1},\pi_{2}:Y_0(\gp)  \to Y,$$
$Y$ (resp.~$Y_0(\gp)$) for $Y_U$ (resp.~$Y_{U_0(\gp)}$) when $U$ is fixed. 

The pairing of Corollary~\ref{cor:flat0} thus defines an isomorphism
\begin{equation} \label{eqn:tr0}  \pi_{1,*}\CO_{Y_0(\gp)}  \stackrel{\sim}{\longrightarrow}
 \Shom_{\CO_{Y}}(\pi_{1,*}\CK_{Y_0(\gp)/\CO}, \,\CK_{Y/\CO}),\end{equation}
under which the image of the unit section is the trace morphism
\begin{equation}
\label{eqn:trace0} \Tr:  \pi_{1,*}\CK_{Y_0(\gp)/\CO}  \longrightarrow \CK_{Y/\CO}
\end{equation}
(defined as in (\ref{eqn:trace}) since $R^i\pi_{1,*} \CK_{Y_0(\gp)/\CO} =0 $ for all $i > 0$).
 On the other hand we have the morphism $\pi_1^* \CK_{Y/\CO} \to \CK_{Y_0(\gp)/\CO}$ defined in (\ref{eqn:integral}),
which by Proposition~\ref{prop:compatibility} extends the canonical isomorphism over $Y_0(\gp)_K$.   It follows that
the composite of its direct image with $\Tr$,
\begin{equation} \label{eqn:Ktrace} \pi_{1,*}\pi_1^*\CK_{Y/\CO}   \longrightarrow \pi_{1,*} \CK_{Y_0(\gp)/\CO}  \longrightarrow \CK_{Y/\CO}, \end{equation}
extends the trace morphism over $Y_K$ (associated to the finite flat morphism $\pi_{1,K}$).    Thus if $\CF$
is any locally free sheaf over $Y$, we obtain an extension of the trace over $Y_K$ to a morphism
$$\tr_{\CF}:  \pi_{1,*}\pi_1^*\CF \longrightarrow \CF$$
by tensoring (\ref{eqn:Ktrace}) with $\CK_{Y/\CO}^{-1} \otimes \CF$ and applying the projection formula.

Recall that we have the Kodaira--Spencer isomorphism $\CK_{Y/\CO} \cong \delta^{-1}\omega^{\otimes 2}$ on $Y$, and that
Theorem~\ref{thm:KSIwa} gives $\CK_{Y_0(\gp)/\CO} \cong \pi_1^*\omega \otimes \pi_2^*(\delta^{-1}\omega)$.  
Recall also that the universal isogeny on $\widetilde{Y}_0(\gp)$ induces an isomorphism 
$\widetilde{\pi}_2^*\widetilde{\delta} \stackrel{\sim}{\to} p^{f_\gp} \widetilde{\pi}_1^*\widetilde{\delta}$ which descends to $Y_0(\gp)$,
so composing with multiplication by $p^{-f_\gp}$ yields an isomorphism $\pi_2^*\delta \stackrel{\sim}{\to} \pi_1^*\delta$, or
equivalently $\pi_1^*\delta^{-1} \stackrel{\sim}{\to} \pi_2^*\delta^{-1}$.
We may therefore view~(\ref{eqn:tr0}) as giving rise to an isomorphism
\begin{equation}  \label{eqn:KStrace}\begin{array}{rcl}
\pi_{1,*}\CO_{Y_0(\gp)} & \stackrel{\sim}{\longrightarrow}  &
 \Shom_{\CO_{Y}}(\pi_{1,*}(\pi_1^*\omega\otimes_{\CO_{Y_0(\gp)}} \pi_2^*(\delta^{-1}\omega) ), \,\delta^{-1}\omega^{\otimes 2} ) \\
 &\stackrel{\sim}{\longrightarrow} &  \Shom_{\CO_{Y}}(\pi_{1,*}(\pi_1^*(\delta^{-1}\omega) \otimes_{\CO_{Y_0(\gp)}} \pi_2^*\omega ), \,\delta^{-1}\omega^{\otimes 2} ) \\
& \stackrel{\sim}{\longrightarrow} &  \Shom_{\CO_{Y}}(\delta^{-1}\omega \otimes_{\CO_{Y}}\pi_{1,*} \pi_2^*\omega ,\, \delta^{-1}\omega \otimes_{\CO_Y} \omega)\\
& \stackrel{\sim}{\longrightarrow} &  \Shom_{\CO_{Y}}( \pi_{1,*} \pi_2^*\omega ,\, \omega ), \end{array}
 \end{equation}
and we call the image of the unit section the {\em saving trace}
\begin{equation}  \label{eqn:stO} \st: \, \pi_{1,*}\pi_2^*\omega  \longrightarrow \omega.\end{equation}
It follows from Proposition~\ref{prop:compatibility} and the definition of the saving trace that the
diagram
\begin{equation} \label{eqn:compatibility}
\xymatrix{ \pi_{1,*} \pi_2^*\omega  \ar[r] \ar[d]_{\st}  & \pi_{1,*}\pi_1^* \omega \ar[d]^{\tr_\omega} \\
\omega  \ar[r]^{p^{f_\gp}} & \omega } \end{equation}
commutes, where the top arrow is the direct image of the descent of the morphism over $\widetilde{Y}_{U_0(\gp)}$
induced by the universal isogeny.

More generally consider the morphisms $\pi_{1},\pi_{2}:  Y_0(\gp)_R \to Y_R$
for any (Noetherian) $\CO$-algebra $R$ (omitting the subscript $R$ from the notation for the degeneracy maps).
It follows from Theorem~\ref{thm:KSIwa} and Corollary~\ref{cor:bc} (for $\CE = \CK_{Y_0(\gp)/\CO}$ and $i=0$)
that $(\pi_{1,*}\pi_2^*\omega)_R \stackrel{\sim}{\to} \pi_{1,*}\pi_2^*(\omega_R)$
so the base-change of (\ref{eqn:stO}) defines a saving trace over $R$:
$$\st_R: \, \pi_{1,*}\pi_2^*(\omega_R)  \longrightarrow \omega_R.$$
Note also that Theorem~\ref{thm:KSIwa} and Corollary~\ref{cor:bc} imply that $R^i\pi_{1,*}(\pi_2^*\omega_R) \cong 
\omega_R^{-1}R^i\pi_{1,*}\CK_{Y_0(\gp)_R/R} = 0$ for all $i > 0$.
More generally, it follows that 
if $\CF$ is any locally free sheaf on $Y_R$, then
$R^i\pi_{1,*}(\pi_1^*\CF \otimes_{\CO_{Y_0(\gp)_R}} \pi_2^*\omega_R)
\cong \CF \otimes_{\CO_{Y_R}} R^i\pi_{1,*}\pi_2^*\omega_R = 0$ for all $i > 0$, and hence that
\begin{equation} \label{eqn:degen}
 H^i(Y_0(\gp)_R,\,\pi_1^*\CF \otimes_{\CO_{Y_0(\gp)_R}} \pi_2^*\omega_R) \cong 
H^i(Y_R,\,\CF \otimes_{\CO_{Y_R}} \pi_{1,*}\pi_2^*\omega_R).\end{equation}

For $(\bk,\bm) \in \ZZ^\Sigma$, we let $\widetilde{\CA}_{\bk,\bm}$ denote the line bundle $\bigotimes_\theta(\CL_{\theta}^{\otimes k_\theta}\otimes\CN_\theta^{\otimes m_\theta})$ on $\widetilde{Y}$.  Recall that if $\mu^{\bk+2\bm}: = \prod_\theta \theta(\mu)^{k_\theta+ 2 m_\theta}$
has trivial image in $R$ for all $\mu \in U \cap \CO_F^\times$, then
$\widetilde{\CA}_{\bk,\bm,R}$ descends to
a line bundle on $Y_R$ which we denote $\CA_{\bk,\bm,R}$.
We will now explain how to use the saving trace to construct a Hecke operator $T_\gp$ on $H^i(Y_R,\CA_{\bk,\bm,R})$
for suitable $\bk,\bm \in \ZZ^\Sigma$ (and in particular whenever
$k_\theta \ge1$ and $m_\theta \ge 0$ for all $\theta$).  Our main interest is in the case $i=0$ and $F \neq \QQ$,
in which case 
$$M_{\bk,\bm}(U;R) :=  H^0(Y_R,\CA_{\bk,\bm,R})$$
is the space of Hilbert modular forms of weight $(\bk,\bm)$ and level $U$
with coefficients in $R$.

First recall that the universal isogeny $\psi:A_1 \to A_2$ over $\widetilde{Y}_0(\gp)$ induces morphisms
$\psi_\theta^*: \widetilde{\pi}_2^*\CP_{\theta} \to \widetilde{\pi}_1^*\CP_{\theta}$, and in turn $\widetilde{\pi}_2^*\CL_{\theta} \to \widetilde{\pi}_1^*\CL_{\theta}$ and
$\wedge^2\psi_\theta^*: \widetilde{\pi}_2^*\CN_{\theta} \to \widetilde{\pi}_1^*\CN_{\theta}$
such that 
$$\theta(\varpi_\gp) \widetilde{\pi}_1^*\CL_\theta \subset \psi_\theta^*(\widetilde{\pi}_2^*\CL_\theta) \subset \widetilde{\pi}_1^*\CL_\theta
\quad\mbox{and}\quad \wedge^2\psi_\theta^*(\widetilde{\pi}_2^*\CN_{\theta}) = \theta(\varpi_\gp)\widetilde{\pi}_1^*\CN_{\theta}.$$
Since $v_\CO(\theta(\varpi_\gp))$ is independent of $\theta \in \Sigma_\gp$ and trivial for $\theta \not\in \Sigma_\gp$, it follows that if
\begin{equation}
\label{eqn:optint}
\sum_{\theta \in \Sigma_\gp}  \min\{m_\theta,m_\theta+k_\theta-1\} \ge 0,
\end{equation}
then $\psi_\theta^*$ induces a morphism $\widetilde{\pi}_2^*\widetilde{\CA}_{\bk - \bf{1},\bm} \to \widetilde{\pi}_1^*\widetilde{\CA}_{\bk - \bone,\bm}$.
Furthermore its base-change to $\widetilde{Y}_0(\gp)_R$ descends to a morphism
\begin{equation}
\label{eqn:psi}  {\pi}_2^*{\CA}_{\bk - \bone,\bm,R} \to {\pi}_1^*{\CA}_{\bk - \bone,\bm,R}
\end{equation}
over $Y_R$.

For $(\bk,\bm)$ satisfying (\ref{eqn:optint}) we define the endomorphism $T_\gp$
of $H^i(Y_R,\CA_{\bk,\bm,R})$ as the composite
\begin{equation}\label{eqn:Tp}  \begin{array}{rcl}
H^i(Y_R,\,\CA_{\bk,\bm,R}) & \stackrel{\pi_2^*}{\longrightarrow} & H^i(Y_0(\gp)_R,\,\pi_2^*\CA_{\bk -\bone,\bm,R} \otimes_{\CO_{Y_0(\gp)_R}} \pi_2^*\omega_R) \\
& {\longrightarrow} & H^i(Y_0(\gp)_R,\, \pi_1^*\CA_{\bk-\bone,\bm,R}  \otimes_{\CO_{Y_0(\gp)_R}} \pi_2^*\omega_R) \\
& \stackrel{\sim}{\longrightarrow} & H^i(Y_R,\, \CA_{\bk-\bone,\bm,R} \otimes_{\CO_{Y_R}} \pi_{1,*}\pi_2^*\omega_R) \\
& \stackrel{1\otimes \st_R}{\longrightarrow}& H^i(Y_R,\,\CA_{\bk,\bm,R}),\end{array}\end{equation}
where the second arrow is induced by (\ref{eqn:psi}) and the third is (\ref{eqn:degen}).

It follows from the construction that the operator $T_\gp$ is compatible with base-change, in the sense
that if $R \to R'$ is any homomorphism of $\CO$-algebras, then the resulting diagram
$$ \xymatrix{H^i(Y_R,\CA_{\bk,\bm,R}) \otimes_R R' \ar[r]\ar[d]_{T_\gp \otimes 1}  & H^i(Y_{R'},\CA_{\bk,\bm,R'}) \ar[d]^{T_\gp} \\
H^i(Y_R,\CA_{\bk,\bm,R}) \otimes_R R' \ar[r]  & H^i(Y_{R'},\CA_{\bk,\bm,R'})}$$
commutes.  Furthermore it is straightforward to check that $T_\gp$ is compatible with the action of
$g \in \GL_2(\AA_{F,\f}^{(p)})$, in the sense that if $U$ and $U'$ are sufficiently small and $g^{-1}Ug \subset U'$,
then the resulting diagram
$$ \xymatrix{H^i(Y'_R,\CA'_{\bk,\bm,R}) \ar[r]\ar[d]_{T_\gp}  & H^i(Y_R,\rho_g^*\CA'_{\bk,\bm,R})  \ar[r]^-{\sim} &
H^i(Y_{R},\CA_{\bk,\bm,R}) \ar[d]^{T_\gp} \\
H^i(Y'_R,\CA'_{\bk,\bm,R})  \ar[r]  & H^i(Y_R,\rho_g^*\CA'_{\bk,\bm,R})  \ar[r]^-{\sim} &
H^i(Y_{R},\CA_{\bk,\bm,R})}$$
commutes, where $\rho_g: Y_R \to Y'_R := Y_{U',R}$ and $\rho_g^*\CA'_{\bk,\bm,R} \stackrel{\sim}{\to} \CA_{\bk,\bm,R}$
are defined in \S\ref{sec:bundles}.  Finally the commutativity of (\ref{eqn:compatibility}) implies that $T_\gp$ coincides
with the classical Hecke operator so denoted on the space  $H^0(Y_K,\CA_{\bk,\bm,K})$
of Hilbert modular forms of level $U$ and weight $(\bk,\bm)$ over $K$.

We record the result as follows:
\begin{theorem} \label{thm:hecke}  For sufficiently small open compact subgroups $U \subset \GL_2(\AA_{F,\f})$ containing
$\GL_2(\CO_{F,p})$, Noetherian $\CO$-algebras $R$, and $\bk, \bm \in \ZZ^\Sigma$ such that $\mu^{\bk + 2\bm}$ has
trivial image in $R$ for all $\mu \in U \cap \CO_F^\times$ and
$$\sum_{\theta \in \Sigma_\gp}  \min\{m_\theta,m_\theta+k_\theta-1\} \ge 0,$$
the operators $T_\gp$ defined on $H^i(Y_{U,R}, \CA_{\bk,\bm,R})$ 
by (\ref{eqn:Tp}) are compatible with base-change and the action of $\GL_2(\AA_{F,\f}^{(p)})$,
and coincide with the classical Hecke operator $T_\gp$ if $R$ is a $K$-algebra and $i=0$.
\end{theorem}

\begin{remark}  \label{rmk:optimality}
The inequality (\ref{eqn:optint}) is needed to ensure integrality, but 
results for more general $\bk$, $\bm$ follow from twisting arguments.  Indeed for fixed $\bk$
and varying $\bm$, the modules $H^i(Y_{U,R}, \CA_{\bk,\bm,R})$ are isomorphic for
sufficiently small $U$.  The isomorphisms are not canonical, but one can keep track of
its twisting effect on the action of $\GL_2(\AA_{F,\f}^{(p)})$ (see for example \cite[Lemma~4.6.1]{DS1}
and \cite[Prop.~3.2.2]{theta}), and the same considerations apply to its effect on $T_\gp$.
In particular if $R$ has characteristic zero (so $\bk + 2\bm$ is parallel), then we recover
the same integrality conditions as in \cite{FP} if $p$ is unramified in $F$.

On the other hand if $R$ has finite characteristic, then there is no restriction on $\bk$ and $\bm$,
and the critical case for $i=0$ becomes $\bm = \bO$.  Note that the hypothesis (\ref{eqn:optint})
is then equivalent to $k_\theta \ge 1$ for all $\theta \in \Sigma_\gp$, which is known to hold in the
cases of primary interest thanks to the main result of~\cite{DK2} and \cite[Prop.~1.13]{DDW} (see also \cite[Thm.~D]{theta}).

Note also that if $\bm = \bO$, then the inequality (\ref{eqn:optint}) becomes an equality.  By contrast
if the inequality is strict one finds that $T_\gp = 0$ if $R$ is a $k$-algebra, and more generally
that $T_\gp$ is nilpotent if $p^nR = 0$.
\end{remark}

\begin{remark} \label{rmk:adjoint}   Interchanging the roles of $\pi_1$ and $\pi_2$, one can similarly define a saving trace
$\st': \pi_{2,*}\pi_1^*\omega \to \omega$ (without using the isomorphism $\pi_2^*\delta \cong \pi_1^*\delta$).  A construction
analogous to the one above then defines an operator $T_\gp'$ on $H^i(Y_R,\CA_{\bk,\bm,R})$ whenever
$\sum_{\theta \in \Sigma_\gp}  \max\{m_\theta,m_\theta+k_\theta-1\} \le 0$,
and the analogue of the commutativity of (\ref{eqn:compatibility}) is then that of the diagram
\begin{equation} \label{eqn:compatibility2}
\xymatrix{ \pi_{2,*} \pi_2^*\omega  \ar[rr] \ar[dr]_{\tr'_\omega}  && \pi_{2,*}\pi_1^* \omega \ar[dl]^{\st'} \\
& \omega, & } \end{equation}
where $\tr'_\omega$ extends the trace relative to $\pi_{2,K}$.

Note that if $\bk = \bone$ and $\bm = \bO$
(or more generally  $k_\theta = 1$ for all $\theta \in \Sigma_\gp$ and $\sum_{\theta \in \Sigma_\gp} m_\theta = 0$), then
both $T_\gp$ and $T_\gp'$ are defined.   Using the commutativity of (\ref{eqn:compatibility}) and (\ref{eqn:compatibility2}),
it is straightforward to check that $T_\gp$ is the composite of $T_\gp'$ with the
automorphism of $H^i(Y_R,\omega_R)$ induced by the map $\sigma_\gp: Y \stackrel{\sim}{\to}{Y}$ obtained by descent
from $\underline{A} \mapsto \underline{A} \otimes_{\CO_F} \gp^{-1}$, together with the isomorphism 
$\sigma_\gp^*{\omega} = \Nm_{F/\Q}(\gp) \otimes_{\Z} {\omega} \stackrel{\sim}{\to} \omega$ defined by $p^{-f_\gp}$.
\end{remark}

\begin{remark} \label{rmk:compactifications}
The main results of this paper, in particular Theorem~3.2.1 and Theorem~5.3.1 (and its corollaries), are extended in \cite[\S5]{new} to toroidal compactifications. Furthermore the effect of the saving trace on $q$-expansions is given by Proposition~5.3.1 of \cite{new}, leading to that of $T_\gp$ in Proposition~6.8.1.  In particular this implies the commutativity of the Hecke operators $T_\gp$ (whenever defined) on $H^0(Y_R,\CA_{\bk,\bm,R})$ for varying $\gp$ in $S_p$.

By contrast, checking this directly from their construction leads to a formidable diagram whose commutativity ultimately
seems to require analogues of some of the results in
\S\ref{sec:vanishing} for direct images under the projections $Y_0(\gp\gp') \to Y_0(\gp')$.
While these might follow from arguments along similar lines to those for direct images under $Y_0(\gp) \to Y$ or from their compatibility with more general base-changes than we considered (more precisely, with respect to $Y_0(\gp)$ instead of $\CO$), we have not carried these out.  Consequently, we have not shown the commutativity of the operators $T_\gp$ on $H^i(Y_R,\CA_{\bk,\bm,R})$ for $i > 0$.
\end{remark}

\bibliographystyle{amsalpha} 

\bibliography{MRrefs_KS}

 \end{document}